\documentclass[a4paper,11pt,reqno]{amsart}
\usepackage[T1]{fontenc}
\usepackage{lmodern}
\usepackage{amsmath,amsfonts,amsthm,amssymb,mathtools,xcolor}
\usepackage{mathrsfs}
\usepackage[a4paper,margin=2.5cm,marginparwidth=2.1cm,marginparsep=2mm]{geometry}
\usepackage[hidelinks]{hyperref}
\usepackage{csquotes}
\usepackage{graphicx}
\usepackage{tikz}
\usepackage{calc}
\IfFileExists{mathabx.sty}{\usepackage{mathabx}}{}
\IfFileExists{ulem.sty}{\usepackage{ulem}}{}
\IfFileExists{stmaryrd.sty}{\usepackage{stmaryrd}}{}

\numberwithin{equation}{section}

\makeatletter
\renewcommand{\part}{\@startsection{part}{0}%
  \z@{1.5\linespacing\@plus\linespacing}{.8\linespacing}%
  {\normalfont\Large\bfseries\centering}}
\makeatother

\newcommand{\epsp}{\varepsilon'}
\newcommand{\Rret}{\mathcal R}
\newcommand{\rr}{\varrho}

\newtheorem{theorem}{Theorem}[section]
\newtheorem{proposition}[theorem]{Proposition}
\newtheorem{lemma}[theorem]{Lemma}
\newtheorem{corollary}[theorem]{Corollary}
\theoremstyle{definition}
\newtheorem{definition}[theorem]{Definition}
\theoremstyle{remark}
\newtheorem{remark}[theorem]{Remark}
\newcommand{\E}{\mathbb E}
\newcommand{\Ep}{E}
\newcommand{\Pp}{\mathbb P}
\newcommand{\R}{\mathbb R}
\newcommand{\N}{\mathbb N}
\newcommand{\Hh}{\mathcal H}
\newcommand{\D}{\mathcal D}
\newcommand{\cE}{\mathcal E}
\newcommand{\ip}[2]{\left\langle #1,#2\right\rangle}
\newcommand{\norm}[1]{\left\lVert #1\right\rVert}
\newcommand{\dd}{\,\mathrm d}

\newcommand{\ToLaw}{\xrightarrow{\,\mathrm{law}\,}}
\newcommand{\ToProb}{\xrightarrow{\,\Pp\,}}
\newcommand{\one}{\mathbf 1}
\newcommand{\CL}{\mathscr L}
\newcommand{\Icross}{\mathfrak I^{\times}}
\newcommand{\heat}{\mathsf p}

\title[Construction and approximation of fractional PAM]{Construction and discrete approximation of the
renormalized fractional parabolic Anderson model}
\author{Aur\'elien Deya}
\address{Universit\'e de Lorraine, CNRS, Institut \'Elie Cartan de Lorraine (IECL, UMR~7502), F-54000 Nancy, France}
\email{aurelien.deya@univ-lorraine.fr}
\date{13 November 2026}

\begin{document}
\begin{abstract}

For a directed polymer in a discrete fractional Gaussian environment,
we identify the scaling limit of its partition field with the
renormalized Stratonovich solution of the one-dimensional parabolic
Anderson model with fractional noise and flat initial condition.
The underlying walk is simple and symmetric, and the temporal and
spatial Hurst indices satisfy $H_0>\frac12$, $\frac12\le H_1<1$ and
$2H_0+H_1<2$. After an explicit deterministic normalization, the
interpolated partition field converges in law in
$C_{\rm loc}([0,\infty)\times\R)$, with convergence of all finite
moments of its supremum norm on compacts.

\smallskip

We construct the limiting field directly by smoothing the noise in
the equation and subtracting the corresponding stationary counterterm.
The classical solutions converge in $L^p(\Omega;C(\mathcal K))$ for
every $1\le p<\infty$ and every compact space-time rectangle $\mathcal K$,
including those that meet the initial line $t=0$. For
$2H_0+H_1>\frac53$, the limit agrees with the solution that Chen,
Ouyang, Tindel and the author constructed using regularity structures.
The construction also covers the new range $\frac32<2H_0+H_1\le\frac53$.
The proofs rely on Feynman-Kac representations and on
estimates and invariance principles for weighted self- and mutual
intersections of random walks and Brownian motions, together with
uniform increment estimates and Gaussian hypercontractivity.

\end{abstract}
\maketitle


\part{Introduction, main results and simulations}\label{chap:introduction-results}

\section{Introduction}\label{sec:article-introduction}

\subsection{Directed polymers and correlated disorder}

Directed polymers in random environments describe the competition
between the entropy of a fluctuating path and its preference for the
favorable regions of a disordered medium \cite{Comets}. Through
Feynman-Kac representations, their partition functions connect
discrete models of random paths to the multiplicative stochastic heat
equation, or parabolic Anderson model \cite{BC}.

When the disorder is independent in time and space, this connection
is by now well understood. Alberts, Khanin and Quastel \cite{AKQ}
identified the intermediate-disorder scaling
$\beta_N=\beta N^{-1/4}$ under which suitably normalized partition
functions converge to the It\^o stochastic heat equation driven by
space-time white noise, and the associated continuum polymer measure
is constructed in \cite{AKQ-CDRP}. Caravenna, Sun and Zygouras
developed a polynomial chaos approach to such scaling limits
\cite{CSZ-JEMS}, and later constructed the critical two-dimensional
stochastic heat flow \cite{CSZ-SHF}. Heavy-tailed disorder leads to
yet other limits \cite{DZ,BL}, among them a continuum polymer in a
L\'evy stable environment \cite{BL}.

Independence is, however, an idealization: the geometry of a
disordered material can produce correlations that persist across many
layers. If the directed polymer is viewed as an elastic line, time
is the longitudinal coordinate and the position of the walk is the
transverse one. A single extended defect may then influence the line
over a substantial longitudinal distance.

Vortex lines in type-II superconductors provide a
physical example. A magnetic field penetrates such a material along
flux lines, each of which behaves as an elastic line directed along
the field; the wandering of a single line in a disordered sample is a
problem of directed polymers in random media \cite{LNG-vortex}.
Irradiation by heavy ions creates columnar defects, that
is, long and thin tracks of damaged material, and such a track
attracts a flux line along its whole length, rather than at an
isolated point. In the description of Nelson and Vinokur
\cite{NV-columnar}, the direction of the field plays the role of an
imaginary time, and columns parallel to the field form a disorder that
does not depend on this time variable; when the columns are tilted at
random angles, the disorder is correlated over long ranges in both
space and time \cite{LNG-vortex}.
What appears as temporal memory in the
polymer description may thus originate in a fixed spatial environment.

\begin{figure}[tbp]
\centering
\begin{tikzpicture}[x=1cm,y=1cm,font=\small,
   grain/.style={fill=black,fill opacity=0.11,draw=black!60,draw opacity=1,line width=0.25pt},
   marked/.style={draw=blue!55!black,line width=1.2pt},
   section/.style={black!70,densely dotted,line width=0.6pt}]
 \begin{scope}[yshift=3.7cm]
  \begin{scope}
   \clip (0,0) rectangle (12,4);
   \path[grain] (-5.3,3.7) rectangle (2.6,4);
   \path[grain] (-1.4,0.9) rectangle (0.4,1.3);
   \path[grain] (-0.8,0.5) rectangle (0.6,1.1);
   \path[grain] (-0.7,0.8) rectangle (1,1.2);
   \path[grain] (-0.6,2.9) rectangle (0.8,3.3);
   \path[grain] (-0.5,0.7) rectangle (0.2,1);
   \path[grain] (-0.2,1.4) rectangle (2.4,1.8);
   \path[grain] (0.2,0.3) rectangle (1.2,0.7);
   \path[grain] (0.5,1.3) rectangle (1.2,1.9);
   \path[grain] (0.7,0.2) rectangle (2,0.7);
   \path[grain] (0.8,2.9) rectangle (1.5,3.2);
   \path[grain] (1.4,2) rectangle (2.1,2.5);
   \path[grain] (1.4,-11) rectangle (2.3,7);
   \path[grain] (1.5,1.4) rectangle (2.1,2.2);
   \path[grain] (1.8,2.4) rectangle (2.5,2.8);
   \path[grain] (2,0.9) rectangle (4.1,1.4);
   \path[grain] (2.2,2.5) rectangle (5.7,3);
   \path[grain] (2.3,2.3) rectangle (7.2,3.3);
   \path[grain] (2.6,3.8) rectangle (3.2,4.3);
   \path[grain] (2.7,1.2) rectangle (3.3,1.6);
   \path[grain] (2.7,-0.1) rectangle (4.3,0.2);
   \path[grain] (2.7,1) rectangle (5.5,3.4);
   \path[grain] (3,-0.1) rectangle (3.9,1.1);
   \path[grain] (3.2,0.5) rectangle (3.8,1.9);
   \path[grain] (3.4,3.5) rectangle (6,4.1);
   \path[grain] (4,-0.1) rectangle (4.6,0.3);
   \path[grain] (4.1,2.9) rectangle (4.8,3.2);
   \path[grain] (4.2,0.5) rectangle (5,0.9);
   \path[grain] (4.2,-0.2) rectangle (13,0.2);
   \path[grain] (4.5,3.1) rectangle (5.1,3.4);
   \path[grain] (4.8,2) rectangle (5.5,2.9);
   \path[grain] (4.9,-0.2) rectangle (6,0.2);
   \path[grain] (5,-0.2) rectangle (5.7,0.4);
   \path[grain] (5.1,1) rectangle (8.1,1.7);
   \path[grain] (5.5,-0.3) rectangle (6.4,1.9);
   \path[grain] (5.6,2.5) rectangle (8.9,3);
   \path[grain] (5.7,3.1) rectangle (6.4,3.5);
   \path[grain] (5.9,1.3) rectangle (9.1,1.7);
   \path[grain] (6,1.7) rectangle (6.9,2.8);
   \path[grain] (6.1,2.6) rectangle (6.9,2.9);
   \path[grain] (6.5,0.8) rectangle (7.1,1.4);
   \path[grain] (6.9,0.1) rectangle (7.9,0.5);
   \path[grain] (7.1,1.5) rectangle (8,2);
   \path[grain] (7.5,1.2) rectangle (8.2,1.7);
   \path[grain] (7.6,2.6) rectangle (8.4,3);
   \path[grain] (7.8,0.2) rectangle (8.5,0.6);
   \path[grain] (8,1.3) rectangle (8.8,1.8);
   \path[grain] (8.1,2.8) rectangle (8.9,3.9);
   \path[grain] (8.2,0.8) rectangle (8.9,2.4);
   \path[grain] (8.4,0.2) rectangle (11.2,0.6);
   \path[grain] (8.7,3.3) rectangle (10.2,3.6);
   \path[grain] (8.8,3.4) rectangle (9.6,4.3);
   \path[grain] (8.8,2.7) rectangle (9.9,3.2);
   \path[grain] (8.9,1.6) rectangle (9.6,2.6);
   \path[grain] (9.1,1.1) rectangle (10.4,2);
   \path[grain] (9.4,0.7) rectangle (10.5,1.2);
   \path[grain] (9.6,0.7) rectangle (10.6,1.6);
   \path[grain] (9.6,3.8) rectangle (17.1,4.5);
   \path[grain] (9.8,3.7) rectangle (10.6,4.2);
   \path[grain] (9.9,-6) rectangle (10.5,0.1);
   \path[grain] (10.1,3.1) rectangle (11.7,3.4);
   \path[grain] (10.5,3.5) rectangle (15.6,3.8);
   \path[grain] (10.6,3.8) rectangle (11.2,5.1);
   \path[grain] (10.8,0.5) rectangle (11.4,1.1);
   \path[grain] (11.7,2.5) rectangle (12.7,2.8);
   \path[grain] (11.7,2.3) rectangle (15.5,3.4);
  \end{scope}
  \draw[line width=0.5pt] (0,0) rectangle (12,4);
  \draw[marked] (2.3,2.3) rectangle (7.2,3.3);
  \draw[black!75,dashed,line width=0.5pt] (0,2.65)--(12,2.65);
  \draw[line width=0.5pt] (0,2.65)--(-0.08,2.65) node[left] {$x$};
  \draw[red!65!black,line width=0.6pt,line join=round]
   (0,1.3)--(0.05,1.39)--(0.1,1.3)--(0.15,1.22)--(0.2,1.3)
   --(0.25,1.22)--(0.3,1.13)--(0.35,1.22)--(0.4,1.13)--(0.45,1.22)
   --(0.5,1.13)--(0.55,1.05)--(0.6,0.96)--(0.65,0.88)--(0.7,0.96)
   --(0.75,0.88)--(0.8,0.96)--(0.85,1.05)--(0.9,1.13)--(0.95,1.05)
   --(1,0.96)--(1.05,0.88)--(1.1,0.79)--(1.15,0.88)--(1.2,0.79)
   --(1.25,0.88)--(1.3,0.79)--(1.35,0.88)--(1.4,0.79)--(1.45,0.71)
   --(1.5,0.79)--(1.55,0.88)--(1.6,0.96)--(1.65,0.88)--(1.7,0.79)
   --(1.75,0.88)--(1.8,0.79)--(1.85,0.71)--(1.9,0.62)--(1.95,0.54)
   --(2,0.62)--(2.05,0.71)--(2.1,0.79)--(2.15,0.88)--(2.2,0.96)
   --(2.25,1.05)--(2.3,1.13)--(2.35,1.05)--(2.4,1.13)--(2.45,1.22)
   --(2.5,1.3)--(2.55,1.22)--(2.6,1.13)--(2.65,1.05)--(2.7,1.13)
   --(2.75,1.05)--(2.8,1.13)--(2.85,1.05)--(2.9,1.13)--(2.95,1.05)
   --(3,0.96)--(3.05,0.88)--(3.1,0.96)--(3.15,1.05)--(3.2,0.96)
   --(3.25,0.88)--(3.3,0.96)--(3.35,0.88)--(3.4,0.79)--(3.45,0.71)
   --(3.5,0.62)--(3.55,0.54)--(3.6,0.62)--(3.65,0.54)--(3.7,0.45)
   --(3.75,0.54)--(3.8,0.62)--(3.85,0.54)--(3.9,0.62)--(3.95,0.71)
   --(4,0.79)--(4.05,0.71)--(4.1,0.62)--(4.15,0.54)--(4.2,0.62)
   --(4.25,0.71)--(4.3,0.79)--(4.35,0.71)--(4.4,0.79)--(4.45,0.88)
   --(4.5,0.79)--(4.55,0.88)--(4.6,0.96)--(4.65,0.88)--(4.7,0.79)
   --(4.75,0.71)--(4.8,0.79)--(4.85,0.88)--(4.9,0.79)--(4.95,0.71)
   --(5,0.79)--(5.05,0.88)--(5.1,0.96)--(5.15,1.05)--(5.2,0.96)
   --(5.25,0.88)--(5.3,0.96)--(5.35,0.88)--(5.4,0.96)--(5.45,1.05)
   --(5.5,1.13)--(5.55,1.05)--(5.6,0.96)--(5.65,1.05)--(5.7,0.96)
   --(5.75,1.05)--(5.8,1.13)--(5.85,1.05)--(5.9,1.13)--(5.95,1.22)
   --(6,1.3)--(6.05,1.39)--(6.1,1.47)--(6.15,1.55)--(6.2,1.47)
   --(6.25,1.55)--(6.3,1.64)--(6.35,1.72)--(6.4,1.64)--(6.45,1.55)
   --(6.5,1.64)--(6.55,1.55)--(6.6,1.47)--(6.65,1.39)--(6.7,1.47)
   --(6.75,1.55)--(6.8,1.64)--(6.85,1.55)--(6.9,1.64)--(6.95,1.55)
   --(7,1.47)--(7.05,1.39)--(7.1,1.3)--(7.15,1.22)--(7.2,1.3)
   --(7.25,1.39)--(7.3,1.3)--(7.35,1.22)--(7.4,1.13)--(7.45,1.05)
   --(7.5,1.13)--(7.55,1.05)--(7.6,0.96)--(7.65,1.05)--(7.7,0.96)
   --(7.75,1.05)--(7.8,1.13)--(7.85,1.05)--(7.9,0.96)--(7.95,1.05)
   --(8,0.96)--(8.05,0.88)--(8.1,0.96)--(8.15,1.05)--(8.2,0.96)
   --(8.25,0.88)--(8.3,0.96)--(8.35,1.05)--(8.4,1.13)--(8.45,1.22)
   --(8.5,1.13)--(8.55,1.22)--(8.6,1.3)--(8.65,1.22)--(8.7,1.3)
   --(8.75,1.22)--(8.8,1.13)--(8.85,1.22)--(8.9,1.3)--(8.95,1.22)
   --(9,1.13)--(9.05,1.22)--(9.1,1.13)--(9.15,1.22)--(9.2,1.13)
   --(9.25,1.22)--(9.3,1.3)--(9.35,1.22)--(9.4,1.3)--(9.45,1.39)
   --(9.5,1.47)--(9.55,1.55)--(9.6,1.64)--(9.65,1.55)--(9.7,1.64)
   --(9.75,1.72)--(9.8,1.64)--(9.85,1.72)--(9.9,1.81)--(9.95,1.89)
   --(10,1.98)--(10.05,1.89)--(10.1,1.81)--(10.15,1.72)--(10.2,1.81)
   --(10.25,1.89)--(10.3,1.98)--(10.35,2.06)--(10.4,1.98)
   --(10.45,2.06)--(10.5,2.15)--(10.55,2.23)--(10.6,2.15)
   --(10.65,2.23)--(10.7,2.32)--(10.75,2.23)--(10.8,2.15)
   --(10.85,2.23)--(10.9,2.15)--(10.95,2.06)--(11,2.15)--(11.05,2.06)
   --(11.1,2.15)--(11.15,2.06)--(11.2,2.15)--(11.25,2.23)--(11.3,2.32)
   --(11.35,2.4)--(11.4,2.49)--(11.45,2.57)--(11.5,2.66)--(11.55,2.57)
   --(11.6,2.66)--(11.65,2.57)--(11.7,2.66)--(11.75,2.57)--(11.8,2.66)
   --(11.85,2.57)--(11.9,2.49)--(11.95,2.57)--(12,2.49);
  \filldraw[fill=white,draw=black,line width=0.5pt] (3.45,2.65) circle[radius=1.7pt];
  \filldraw[fill=white,draw=black,line width=0.5pt] (6.25,3.05) circle[radius=1.7pt];
  \node[fill=white,fill opacity=0.85,text opacity=1,inner sep=1.5pt,anchor=south,outer sep=3pt] at (3.45,2.65) {$(t,x)$};
  \node[fill=white,fill opacity=0.85,text opacity=1,inner sep=1.5pt,anchor=north,outer sep=3pt] at (6.25,3.05) {$(t+r,x+z)$};
  \fill[blue!55!black] (2.3,2.3) circle[radius=1.6pt];
  \node[fill=white,fill opacity=0.85,text opacity=1,inner sep=1.5pt,anchor=north east,outer sep=2pt,text=blue!55!black] at (2.3,2.3) {$(s_i,y_i)$};
  \draw[blue!55!black,line width=0.5pt,<->] (2.3,3.46)--(7.2,3.46);
  \node[fill=white,fill opacity=0.85,text opacity=1,inner sep=1.5pt,text=blue!55!black] at (4.75,3.46) {$L_i$};
  \draw[blue!55!black,line width=0.5pt,<->] (7.36,2.3)--(7.36,3.3);
  \node[fill=white,fill opacity=0.85,text opacity=1,inner sep=1.5pt,text=blue!55!black,anchor=west,outer sep=1pt] at (7.36,2.8) {$M_i$};
  \node[rotate=90,anchor=south] at (-0.75,2) {\footnotesize transverse position};
  \node[anchor=east] at (-1.05,3.8) {(a)};
 \end{scope}
 \draw[section] (3.45,0)--(3.45,6.35);
 \begin{scope}[y=0.42cm]
  \draw[->,line width=0.5pt] (0,0)--(12.35,0);
  \draw[->,line width=0.5pt] (0,0)--(0,6.9) node[above] {$n(\cdot,x)$};
  \foreach \k in {0,...,6}{
    \draw[line width=0.5pt] (0,\k)--(-0.08,\k) node[left] {\footnotesize$\k$};
  }
  \draw[black!60,dashed,line width=0.6pt] (0,1.44)--(12,1.44)
       node[right,black] {$\lambda\,\mathbb E[L]\,\mathbb E[M]$};
  \draw[line width=0.9pt,line join=miter]
   (0,0) -- (1.4,0) -- (1.4,1) -- (1.8,1) -- (1.8,2) -- (2.2,2)
   -- (2.2,3) -- (2.3,3) -- (2.3,3) -- (2.5,3) -- (2.5,2) -- (2.7,2)
   -- (2.7,3) -- (4.8,3) -- (4.8,4) -- (5.5,4) -- (5.5,2) -- (5.6,2)
   -- (5.6,3) -- (5.7,3) -- (5.7,2) -- (6,2) -- (6,3) -- (6.1,3)
   -- (6.1,4) -- (6.9,4) -- (6.9,2) -- (7.2,2) -- (7.2,1) -- (7.6,1)
   -- (7.6,2) -- (8.4,2) -- (8.4,1) -- (8.9,1) -- (8.9,0) -- (11.7,0)
   -- (11.7,2) -- (12,2);
  \draw[<->,line width=0.7pt,blue!55!black] (3.45,1.44)--(3.45,3);
  \node[above,fill=white,inner sep=1.5pt,outer sep=2pt,text=blue!55!black] at (3.45,3) {$V(t,x)$};
  \node[below] at (3.45,0) {$t$};
  \node[anchor=north] at (6,-1.15) {\footnotesize longitudinal position};
  \node[anchor=east] at (-1.05,6.9) {(b)};
 \end{scope}
\end{tikzpicture}
\caption{Rectangular defects in the plane of the longitudinal and
transverse coordinates, for $U_i=1$.
(a)~Each defect is a rectangle $[s_i,s_i+L_i)\times[y_i,y_i+M_i)$; the
grey level at a point increases with the number of defects covering
it, and the broken line is a polymer path. The outlined defect covers
both $(t,x)$ and $(t+r,x+z)$: only such defects correlate the
potential at these two points.
(b)~The number $n(\cdot,x)$ of defects along the dashed line of~(a),
and its mean $\lambda\,\mathbb E[L]\,\mathbb E[M]$ (dashed); $V(t,x)$
is their difference. The lengths and widths are drawn from the laws
$\mathbb P(L>\ell)=(\ell_0/\ell)^{3/2}$, $\ell\ge\ell_0$, and
$\mathbb P(M>m)=(m_0/m)^{8/5}$, $m\ge m_0$, which correspond to
$H_0=\frac34$ and $H_1=\frac7{10}$.}
\label{fig:defects}
\end{figure}
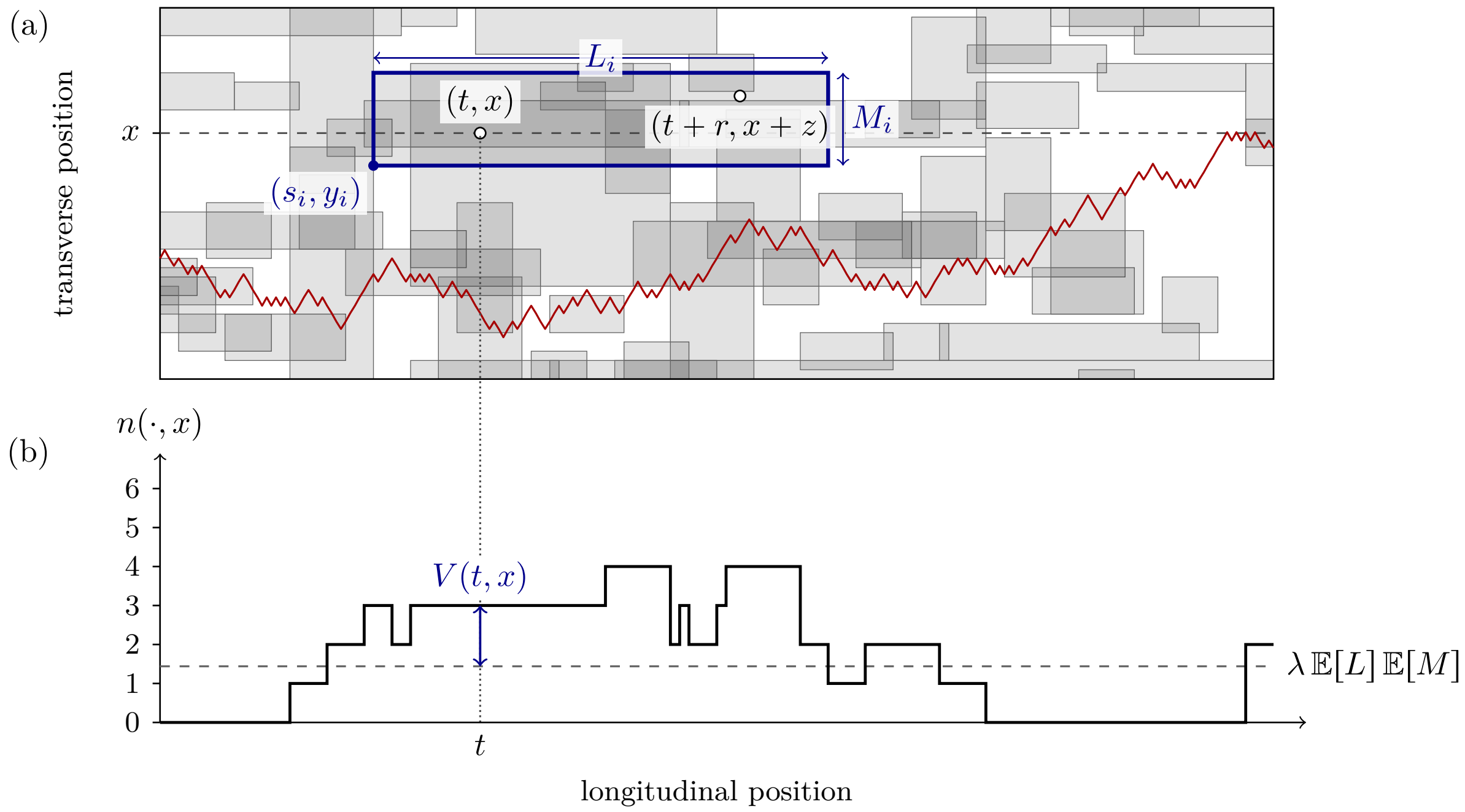

A simple mechanism shows how such memory can become algebraic, in time
and in space. Suppose that the defects are rectangles in the plane of
the longitudinal and transverse coordinates, as in
Figure~\ref{fig:defects}: a thin rectangle represents a column, and a
wider one a defect with a transverse extent, such as the twin
boundaries of \cite{NV-columnar}, which are planar defects aligned
with the field. The $i$th defect occupies the rectangle
\[
 D_i=[s_i,s_i+L_i)\times[y_i,y_i+M_i),
\]
of length $L_i$ and width $M_i$, and contributes the amount $U_i$ to
the potential at each of its points. The points $(s_i,y_i)$ form a
Poisson process of intensity $\lambda$ on $\R^2$, and the triples
$(L_i,M_i,U_i)$ are independent copies of $(L,M,U)$, independent of
these points, where $L$, $M$ and $U$ are independent, $L$ and $M$ are
positive with finite mean, and $U$ is bounded. The number $n(t,x)$ of
defects covering a point $(t,x)$ is then a Poisson variable with mean
$\lambda\,\mathbb E[L]\,\mathbb E[M]$, and the centered potential at
this point is
\[
 V(t,x)=\sum_{i:\,(t,x)\in D_i}U_i
        -\lambda\,\mathbb E[U]\,\mathbb E[L]\,\mathbb E[M].
\]
The field $V$ plays the role of the environment of the polymer. When
$U=1$, it records the fluctuations of the number of defects:
$V(t,x)=n(t,x)-\lambda\,\mathbb E[L]\,\mathbb E[M]$.

The values of $V$ at two points $(t,x)$ and $(t+r,x+z)$ are correlated
only through the defects that cover both of them. Indeed, up to
additive constants, these values are $A+C$ and $B+C$, where $C$ is the
sum of the $U_i$ over the defects covering both points, and $A$ and
$B$ are the sums over the defects covering only the first and only the
second, respectively. These three sums involve the points
$(s_i,y_i,L_i,M_i,U_i)$ of a Poisson process in disjoint regions; they
are therefore independent, and the covariance of the two values is the
variance of $C$. The $i$th defect contributes to $C$ exactly when
$(s_i,y_i)$ lies in a rectangle of sides $(L_i-|r|)_+$ and
$(M_i-|z|)_+$, so that the number of such defects is a Poisson
variable with mean
$\lambda\,\mathbb E[(L-|r|)_+]\,\mathbb E[(M-|z|)_+]$. Hence
\[
 \operatorname{Cov}\bigl(V(t,x),V(t+r,x+z)\bigr)
 =\operatorname{Var}(C)
 =\lambda\,\mathbb E[U^2]\,\mathbb E[(L-|r|)_+]\,
        \mathbb E[(M-|z|)_+],
\]
which is the product of a temporal and a spatial factor. At a fixed
transverse location $x$, the process $t\mapsto V(t,x)$ is a centered
infinite-source Poisson process of rate $\lambda\,\mathbb E[M]$, and
the formula at $z=0$ is its stationary covariance (see
\cite[Proposition~2.1]{FRS-Poisson}).

If $\mathbb P(L>\ell)\sim c\ell^{-(3-2H_0)}$ as $\ell\to\infty$, with
$\frac12<H_0<1$, the lengths have finite mean and
\[
 \mathbb E[(L-r)_+]=\int_r^\infty\mathbb P(L>\ell)\,\dd\ell
 \sim\frac{c}{2-2H_0}\,r^{2H_0-2},
       \qquad r\to\infty.
\]
In the same way, a tail of exponent $3-2H_1$ for the widths, with
$\frac12<H_1<1$, gives a spatial factor of order $|z|^{2H_1-2}$,
whereas defects of bounded width leave distant transverse locations
independent; this short-range case corresponds to $H_1=\frac12$ below.
Power-law distributions of the lengths and widths thus reproduce, in
each variable, the covariance decay of fractional Gaussian noise, and
lend the Hurst indices a geometric meaning: the larger $H_0$, the more
frequent the long defects.

The Gaussian description is itself natural at high density,
that is, when many defects overlap at a typical point. Although
$V$ is not Gaussian, the finite-dimensional distributions of
$\lambda^{-1/2}V$ converge, as $\lambda\to\infty$, to those of the
centered Gaussian field with covariance
$\mathbb E[U^2]\,\mathbb E[(L-|r|)_+]\,\mathbb E[(M-|z|)_+]$, since
the joint cumulants of order $k\ge3$ of $\lambda^{-1/2}V$ are of order
$\lambda^{1-k/2}$ (see \cite[(2.5)]{FRS-Poisson}).
Moreover, for the power-law lengths
considered above, the variance of the cumulative process
$\int_0^tV(s,x)\,\dd s$ at a fixed
location $x$ grows like $t^{2H_0}$ \cite{FRS-Poisson},
as it does for a
fractional Brownian motion of Hurst index $H_0$. Limit theorems of
this type, in which fractional Brownian motion arises from a dense
superposition of heavy-tailed durations, are classical in the
modelling of network traffic \cite{MRRS}; their
analogue for random fields generated by Poisson families of sets with
heavy-tailed volumes is established in \cite{KLNS}.

More broadly, quenched disorder with
power-law correlations is a classical theme in the physics of
disordered systems \cite{WH}. For directed polymers, Lacoin
\cite{Lacoin-correlation} showed that power-law spatial correlations of
a Gaussian environment affect the transition between weak and strong
disorder.

\subsection{The model and the main result}

These considerations lead us to a Gaussian field with a product
covariance as an effective model of long-range disorder. On the
lattice, we take a simple symmetric random walk $S$ on $\mathbb Z$ and
an independent Gaussian environment $\omega_h$ whose covariances in
time and space are those of fractional Brownian increments. The Hurst
indices $h=(H_0,H_1)$, which govern the longitudinal and the
transverse correlations respectively, are assumed to satisfy
\[
 H_0>\frac12,\qquad \frac12\le H_1<1,\qquad 2H_0+H_1<2.
\]
At $H_1=\frac12$ the environment is independent between sites, and
the limiting noise is white in space. With flat initial condition,
the formal continuum equation reads
\begin{equation}\label{intro:pam}
 \partial_t U=\frac12\Delta U+U\dot W_h-\infty\,U,
 \qquad U_0\equiv1.
\end{equation}
Its precise meaning, that is, the covariance of $\dot W_h$ and the
regularization behind the symbol $-\infty\,U$, is given in
Section~\ref{subsec:fractional-main}.

The lattice counterpart of $U$ is a renormalized partition function of
the polymer. Write $a_h=2H_0+H_1-1\in(\frac12,1)$, and let
$\beta_{N,h}=N^{-a_h/2}$ be the coupling strength; at the excluded
endpoint $H_0=H_1=\frac12$, this exponent would be the
intermediate-disorder exponent $\frac14$ of \cite{AKQ}. For $N\ge1$,
$t\ge0$ and $x\in\R$, set $N_t=\lfloor Nt\rfloor$ and
$\widetilde N_x=\lfloor x\sqrt N\rfloor$, and let
$\omega_h(n,j)$ be the value of the environment at time layer $n$
and site $j$. Its covariance is
\[
 \E[\omega_h(n,j)\omega_h(m,k)]=\gamma_{H_0}(n-m)\gamma_{H_1}(j-k),
\]
where $\gamma_H$ denotes the covariance of the unit increments of a
fractional Brownian motion of Hurst index $H$ (see
\eqref{fr:discrete-covariance}-\eqref{fr:fgn}). A path of $N_t$ steps started from the site
$\widetilde N_x$ has the scaled energy
\begin{equation}\label{intro:energy}
 X^N_{t,x}(S)=\beta_{N,h}\sum_{n=1}^{N_t}
       \omega_{h}(N_t-n+1,\widetilde{N}_x+S_n),
\end{equation}
and the partition field is defined by
\begin{equation}\label{intro:partition}
 \mathsf Z^N_t(x)
 =e^{-c_{\rm disc}(h)N_t\beta_{N,h}^2}\,
  E_S\Bigl[e^{X^N_{t,x}(S)}\Bigr],
\end{equation}
where $E_S$ denotes expectation over the walk, the environment being
fixed, and $c_{\rm disc}(h)$ is a deterministic constant, given in
\eqref{fr:counterterm-discrete}. This is the field
\eqref{fr:terminal-field} of Section~\ref{subsec:fractional-main},
where these conventions are discussed.

The exponential prefactor in \eqref{intro:partition}
cannot be dispensed with: without it, the partition function fails to
be tight as $N\to\infty$, at every point $(t,x)$ with $t>0$
(Proposition~\ref{prop:discrete-renormalization-needed}; see also
Remark~\ref{rem:discrete-mean} for the divergence of its mean).
The question is therefore how the polymer should be normalized.
In \eqref{intro:partition}, we retain
its ordinary Boltzmann weights $e^{X^N_{t,x}(S)}$,
and apply the same deterministic
correction, namely the factor
$e^{-c_{\rm disc}(h)N_t\beta_{N,h}^2}$, to all paths of a given length. Such a correction preserves
the relative weights of the paths in each realization of the
environment. It therefore cancels in the polymer Gibbs
measure, under which a path $S$ has probability
proportional to $e^{X^N_{t,x}(S)}$, and
in every ratio $\mathsf Z^N_t(x)/\mathsf Z^N_t(y)$ of partition
functions at a common time: it changes the normalization of the
partition function, and not the polymer itself.

What this choice retains is best seen through the energy of a single
path. Given the path, $X^N_{t,x}(S)$ is a centered
Gaussian variable, with variance (see \eqref{eq:GammaN})
\begin{equation}\label{intro:Gamma}
 \Gamma^N_t(S)=\beta_{N,h}^2\sum_{n,m=1}^{N_t}
     \gamma_{H_0}(n-m)\gamma_{H_1}(S_n-S_m).
\end{equation}
For an environment independent in time, only
the terms $n=m$ would remain in this sum, and this variance would take the same
value $N_t\beta_{N,h}^2$ for every path. Because the environment is correlated in time, the conditional
variance of this energy depends on the path, through the covariances
between pairs of its visits. The deterministic correction removes
only the divergent part of the mean of this self-interaction (Proposition~\ref{prop:counterterm}). Its
fluctuations about the mean survive under diffusive scaling, in the
form of a centered Brownian self-intersection: a contribution that
Wick ordering of each path weight would erase (see
\eqref{intro:wick-factorization} below).

To state the result, let
$\operatorname{Interp}_N$ denote bilinear interpolation in space and time (see \eqref{eq:bilinear-interpolation}). Theorem~\ref{thm:fractional-main} asserts that
\begin{equation}\label{intro:main-limit}
 \operatorname{Interp}_N\mathsf Z^N\ToLaw U
 \qquad\text{in }C_{\rm loc}([0,\infty)\times\R),
\end{equation}
where $U$ is the renormalized Stratonovich solution of
\eqref{intro:pam}.
The convergence thus holds in law for the entire field, the initial
time $t=0$ included, in the topology of uniform convergence on
compact space-time sets. Every finite moment of the supremum norm on
a compact space-time rectangle converges as well.

\subsection{The continuum field and its counterterm}

The limit $U$ is constructed here independently of the lattice,
directly on the Gaussian space of the noise. Let $U^\varepsilon$ solve
the classical equation with Gaussian-mollified potential
$\dot W_h^\varepsilon-c_h^\varepsilon$. Theorem~\ref{thm:continuous-summary} shows that, for
every $T,L>0$ and $1\le p<\infty$,
\[
 U^\varepsilon\longrightarrow U
 \qquad\text{in }L^p(\Omega;C([0,T]\times[-L,L])).
\]
The limit is a nonnegative
continuous field with moments of every finite order.

The renormalized equation \eqref{intro:pam} was previously studied by
Chen, Ouyang, Tindel and the author \cite{CDOT-KRP}. The solution is
constructed there in weighted spaces by means of regularity structures
\cite{Hairer-RS}, and, in dimension one, this construction covers the
range $\frac53<2H_0+H_1<2$.
The indices considered in the present paper satisfy
$\frac32<2H_0+H_1<2$. Our Brownian construction thus also covers the
rougher regime $\frac32<2H_0+H_1\le\frac53$.

The two counterterms, discrete and continuous, obey one and the same
stationary prescription. In both constructions we average the
conditional half-variance over paths and take its long-time growth
rate: per step for the random walk, per unit time for Brownian motion
(see \eqref{eq:stationary-counterterm-rates}). In this sense, the
lattice model provides a microscopic interpretation of the counterterm
in the equation.

Interface growth offers a familiar counterpart of this interpretation.
Through the Cole-Hopf correspondence, $\log\mathsf Z^N$ plays the
role of a height function, and the logarithm of a solution of
\eqref{intro:pam} formally solves a Kardar-Parisi-Zhang equation
driven by $\dot W_h$. Burgers and KPZ equations with noise correlated
in space and time were introduced by Medina, Hwa, Kardar and Zhang
\cite{MHKZ}, with directed polymers among their motivations, and the
effect of temporal correlations on the growth exponents is still under
investigation \cite{KS-KPZ,SC-KPZ}. Seen in this light, the counterterm
is a deterministic drift of the height, that is, a mean growth
velocity. Such a velocity is not universal, since its value
$c_{\rm disc}(h)$ depends on the lattice, whereas the limiting field
does not. We use this picture only as a heuristic guide, since strict
positivity of $U$ is not proved in this paper.

The counterterm is equally indispensable in the
continuum, up to the threshold $2H_0+H_1=2$. When
$2H_0+H_1>2$, Hu, Nualart and Song \cite{HNS} obtain a Feynman-Kac
formula without a counterterm. In the range considered here, by
contrast, Proposition~\ref{prop:renormalization-needed} shows that a
divergent renormalization is necessary even for convergence at a
single point of positive time.

\subsection{Related work}

The case of spatial white noise brings the role of the counterterm
into sharp focus. Rang, Song and Wang \cite{RSW} study polymers in a
Gaussian environment with independent sites and power-law temporal
correlations. For aperiodic walks in the domain of attraction of a
stable law of index $\rho$, their Theorem~1.1 establishes convergence in
law to the Stratonovich solution at $(1,x_0)$ for each fixed $x_0$,
under the condition $H_0-\frac1{2\rho}>\frac12$; in the diffusive
case $\rho=2$, this means $H_0>\frac34$. They also obtain a Skorokhod
limit by Wick ordering each path weight.

Our result is complementary. For the simple symmetric walk, which has
period two, \eqref{intro:main-limit} yields convergence of the entire
interpolated partition field to the renormalized Stratonovich solution
throughout the interval $\frac12<H_0<\frac34$, with a single
correction outside the path expectation.

The form of the correlations also affects the limit. Rang \cite{Rang}
obtains a limit driven by spatially fractional noise for a polymer
with algebraically decaying spatial correlations. For a polymer among
a Poisson field of independent walks, Shen, Song, Sun and Xu
\cite{SSSX} obtain in dimension $1+1$ a Stratonovich limit with
heat-kernel covariance and no divergent counterterm.

Finally, deterministic counterterms also arise in lattice models
evolving in continuous time. Chouk, Gairing and Perkowski \cite{CGP}
treat a small static potential in two dimensions, and Martin and
Perkowski \cite{MP-lattice} prove weak universality results on
Bravais lattices. For the dynamic potential generated by an exclusion
process, Erhard and Hairer \cite{EH} and Erhard, Hairer and Xu
\cite{EHX} obtain a renormalized equation driven by an
Ornstein-Uhlenbeck field in dimensions three and two, respectively.
In the present model, the polymer is discrete in both variables, and
the correlations determine how much weight each pair of visits
contributes to its self-interaction.

\subsection{Strategy of the proof}

The effect of temporal correlations on the normalization becomes
explicit on the lattice. Recall the scaled path energy
$X^N_{t,x}(S)$ from \eqref{intro:energy}, and its variance
$\Gamma^N_t(S)$ conditional on the walk from \eqref{intro:Gamma}.
Wick ordering subtracts half of $\Gamma^N_t(S)$ separately for each
path. The centered self-interaction is
$Y^N_t(S)=\Gamma^N_t(S)-E_S[\Gamma^N_t(S)]$. With
$R^N_t=\frac12E_S[\Gamma^N_t(S)]-c_{\rm disc}(h)N_t\beta_{N,h}^2$,
the exact identity
\begin{equation}\label{intro:wick-factorization}
 \mathsf Z^N_t(x)=e^{R^N_t}E_S\!\left[
      e^{X^N_{t,x}(S)-\frac12\Gamma^N_t(S)}
      e^{\frac12Y^N_t(S)}\right]
\end{equation}
separates the Wick exponential from the fluctuations retained by the
deterministic renormalization.

The identification of the limit also requires us to match the finite
parts of the two normalizations: the lattice remainder $R^N_t$ and
its analogue for the mollified equation converge to the same finite
$R_t$ (Propositions~\ref{prop:counterterm}
and~\ref{prop:R-eps-limit}).

To compare the two fields, we realize a copy of the lattice
environment on the Gaussian space of the continuum noise, and denote
the associated partition field by $\widehat{\mathsf Z}^N$. Its
interpolation has the same law as $\operatorname{Interp}_N\mathsf Z^N$
(Lemma~\ref{lem:environment-law}). On this common space we prove that
\[
 \operatorname{Interp}_N\widehat{\mathsf Z}^N\longrightarrow U
 \qquad\text{in }L^p(\Omega;C([0,T]\times[-L,L]))
\]
for every finite $p\ge1$ and $T,L>0$ (Theorem~\ref{thm:coupled-L2}).

The proof begins at a single point, where Gaussian duality reduces
the comparison to functionals of one or two paths. For the continuum
field, testing against a Gaussian exponential gives its
$\mathcal S$-transform. Up to the factor $e^{R_t}$, it is an average
over a Brownian path involving its centered self-intersection and a
linear occupation functional. The second moment introduces a second
independent path and the weighted mutual intersection of the pair
(see \eqref{eq:S-U-theo}-\eqref{eq:continuous-second-exact-theo}).
Brownian intersection functionals of this kind enter the moment
estimates of \cite{CDOT-moments} through a comparison of Stratonovich
and Skorokhod Feynman-Kac formulas; they also appear in the moment
representations of Gu and Xu \cite{GuXu} for the two-dimensional
static model.

On the discrete side, the identities of
Proposition~\ref{prop:discrete-identities} involve the corresponding
functionals of the walk. We prove their joint convergence to the
Brownian functionals, together with exponential moment
bounds that allow passage to the limit in the Feynman-Kac weights. For the centered self-intersection we use
the block decomposition introduced by Varadhan \cite{Varadhan} and
developed by Le Gall for random walks \cite{LeGall-RW}, together with
invariance principles for the walk and its local times
\cite{Revesz-1981}. Since Gaussian exponentials form a total family,
convergence of the $\mathcal S$-transforms together with the moment
bounds gives weak $L^2$ convergence at a point, and the two-path
calculation gives convergence of the $L^2$ norms.

It then remains to pass from a point to the entire field. By Nelson's
hypercontractivity \cite[Chapter~V]{Janson-Gaussian},
higher moments of space-time increments are bounded by another two-path
estimate (Propositions~\ref{prop:continuous-functional}
and~\ref{prop:discrete-increments}), and a
Kolmogorov-Chentsov argument strengthens pointwise
convergence to convergence in $L^p(\Omega;C(\mathcal K))$ on each compact
space-time rectangle $\mathcal K$ (Lemma~\ref{lem:compact-upgrade}).

\subsection{Organization of the paper}

Section~\ref{sec:main-results} gives the precise statements and
discusses their scope, including the questions of rates and
universality. Section~\ref{sec:numerics} illustrates the effect of
renormalization on terminal polymer profiles at $H_1=\frac12$.

The proofs begin with the Brownian self-intersection in
Section~\ref{sec:Brownian-construction} and the construction of the
continuum field in Section~\ref{sec:continuous-identification}.
On the lattice, Section~\ref{sec:discrete} derives the Gaussian
identities and the deterministic remainder.
Sections~\ref{subsec:legall-complet} and~\ref{sec:two-path} establish
the one- and two-path invariance principles, and
Section~\ref{sec:short} concludes the proof of
\eqref{intro:main-limit}.

Appendices~\ref{app:temporal}-\ref{app:compact-upgrade} collect the
kernel and moment bounds, the Feynman-Kac formula for unbounded
potentials, and the criterion for convergence on compact sets.

\section{Setting and main results}\label{sec:main-results}

\subsection{The models and the main result}
\label{subsec:fractional-main}

We begin with the two models, continuous and discrete, and with
the renormalization conventions that enter the main theorem. Throughout
the paper, the Hurst indices $h=(H_0,H_1)$ are assumed to satisfy
\begin{equation}\label{fr:range}
 H_0>\frac12,\qquad \frac12\le H_1<1,\qquad 2H_0+H_1<2.
\end{equation}
In particular, $H_0<\frac34$, since
$H_1\ge\frac12$.

All objects are constructed for each fixed $h$, with no uniformity
in $h$ claimed, and $h$ is dropped from the notation ($U$,
$c_{\rm disc}$, $\beta_N$, $\mathsf Z^N$, \ldots) when no confusion can
arise; we write $\mathsf Z^N_t(x;\omega)$ when the dependence on the
environment array has to be displayed.

The scaling of the problem is governed by a single exponent,
\begin{equation}\label{fr:exponents}
 a_{h}=2H_0+H_1-1\in\Bigl(\tfrac12,1\Bigr).
\end{equation}
The roles of the restrictions in
\eqref{fr:range} are discussed in Remark~\ref{fr:range-remark}.

We work in spatial dimension $d=1$, and the same family of covariance
kernels serves in time and in space: for $H\in[\frac12,1)$, we set
\begin{equation}\label{fr:spatial-kernel}
 K_H(z)=c_H|z|^{2H-2},\quad c_H=H(2H-1)\quad(\tfrac12<H<1),
 \qquad K_{1/2}=\delta_0.
\end{equation}
We then denote by $W_{h}$ the isonormal Gaussian process associated with the
centered stationary noise $\dot W_{h}$ whose covariance is
\begin{equation}\label{fr:noise-covariance}
 \E[\dot W_{h}(t,x)\dot W_{h}(s,y)]
 =K_{H_0}(t-s)K_{H_1}(x-y).
\end{equation}

With this noise, the continuous model is, formally, the parabolic Anderson equation
\[
\partial_tU=\frac12\Delta U+U\dot W_h-\infty\,U,
\]
with $U_0\equiv1$. The regularization that gives a meaning to this
equation, and the notion of renormalized Stratonovich solution, are
specified below.

We now turn to the discrete side. The polymer paths are those of the simple symmetric random walk
\begin{equation}\label{eq:simple-walk}
 S_0=0,\qquad S_n=X_1+\cdots+X_n,\qquad
 \Pp(X_i=1)=\Pp(X_i=-1)=\frac12,
\end{equation}
where $(X_i)_{i\ge1}$ are independent and identically distributed, and
which is taken independent of all Gaussian noises. This is the
canonical walk appearing in the elementary formulation of Donsker's
invariance principle. It is centered, has variance $1$, and has period
$2$. Writing $p_n(j)=\Pp(S_n=j)$, the return probabilities satisfy
\begin{equation}\label{eq:LLT}
 p_{2r+1}(0)=0,\qquad
 p_{2r}(0)=2^{-2r}\binom{2r}{r}
 \sim\frac1{\sqrt{\pi r}}
 =\sqrt{\frac2\pi}\,(2r)^{-\frac12}.
\end{equation}

This walk evolves in a centered Gaussian environment $\omega_{h}=(\omega_{h}(n,j))_{n\ge1,\ j\in\mathbb Z}$, independent of
$S$, with covariance
\begin{equation}\label{fr:discrete-covariance}
 \E[\omega_{h}(n,j)\omega_{h}(m,k)]
 =\gamma_{H_0}(n-m)\gamma_{H_1}(j-k),
\end{equation}
where, for $H\in[\frac12,1)$ and $k\in\mathbb Z$,
\begin{equation}\label{fr:fgn}
 \gamma_H(k)=\frac12\bigl(|k+1|^{2H}+|k-1|^{2H}-2|k|^{2H}\bigr).
\end{equation}

In particular, $\gamma_{H_0}(0)=\gamma_{H_1}(0)=1$, while $\gamma_{H}(r)\sim c_{H}|r|^{2H-2}$ as $|r|\to\infty$, when $H>\frac12$.

\smallskip 

With this notation, the coupling strength and the lattice counterterm
are defined by
\begin{align}
 \beta_{N,h}&=N^{-a_{h}/2},\label{fr:scaling}\\
 c_{\rm disc}(h)
 &=\frac12+\sum_{r=1}^\infty\gamma_{H_0}(r)q_{H_1}(r),
 \qquad q_{H_1}(r)=E_S[\gamma_{H_1}(S_r)],\label{fr:counterterm-discrete}
\end{align}
where $E_S$ averages over the walk with the environment fixed.

The summand in \eqref{fr:counterterm-discrete} is $O(r^{a_h-2})$
as $r\to\infty$. For $H_1>\frac12$, it is asymptotic to a positive
multiple of $r^{a_h-2}$ (see the proof of Proposition~\ref{prop:counterterm}).
For $H_1=\frac12$, $q_{1/2}(r)=p_r(0)$: by \eqref{eq:LLT} the odd
summands vanish, while along even $r$ the summand is asymptotic to
$c_{H_0}\sqrt{2/\pi}\,r^{a_h-2}$. Thus in both cases the series
converges exactly when $a_h<1$.

\smallskip\noindent\textit{The endpoint $H_1=\frac12$.}
Each object is defined for a general pair $h=(H_0,H_1)$, and
Theorems~\ref{thm:continuous-summary} and~\ref{thm:coupled-L2} include
the endpoint $H_1=\frac12$. There, the range \eqref{fr:range} reduces to
\begin{equation}\label{eq:H-range}
  \frac12<H_0<\frac34,
\end{equation}
the scaling exponent becomes $a_h=2H_0-\frac12$, and
$\gamma_{1/2}(k)=\one_{\{k=0\}}$. The coupling strength and the lattice
counterterm become $\beta_N=N^{-(H_0-\frac14)}$ and
\begin{equation}\label{eq:clat}
 c_{\rm disc}=\frac12+\sum_{r=1}^\infty\gamma_{H_0}(r)p_r(0),
\end{equation}
finite for $H_0<\frac34$. The remaining specializations at
$H_1=\frac12$ are recorded where the corresponding objects are defined.

To define the partition field, we fix a finite time horizon $t\ge0$ and a point $x\in\R$, and write
\begin{equation}\label{eq:Nt}
 N_t=\lfloor Nt\rfloor,\qquad t_N=\frac{N_t}{N},
\end{equation}
\begin{equation}\label{eq:xN}
 \widetilde{N}_x=\lfloor x\sqrt N\rfloor,\qquad
 x_N=\frac{\widetilde{N}_x}{\sqrt N}.
\end{equation}

The partition field is then given by
\begin{equation}\label{fr:terminal-field}
 \mathsf Z^N_t(x)
 =e^{-c_{\rm disc}(h)N_t\beta_{N,h}^2}
 E_S\left[\exp\left(\beta_{N,h}\sum_{n=1}^{N_t}
                    \omega_{h}(N_t-n+1,\widetilde{N}_x+S_n)\right)\right].
\end{equation}
The number of steps is thus $N_t$, while the coupling strength
$\beta_{N,h}$ keeps the scale $N$ of the lattice. The time layers are
read backwards from the terminal layer $N_t$; for $t=0$ the sum is
empty and $\mathsf Z^N_0(x)=1$.

The field $\mathsf Z^N$ is piecewise constant in $(t,x)$, and it is its interpolation that will be compared with the continuous field.
For a field $F$ on the grid
$N^{-1}\mathbb Z_{\ge0}\times N^{-1/2}\mathbb Z$, and for
$(t,x)\in[0,\infty)\times\R$, we use the rescaled coordinates $t_N,x_N$ from
\eqref{eq:Nt}-\eqref{eq:xN}, set
$\theta=N(t-t_N)$ and $\zeta=\sqrt N(x-x_N)$, and define
the bilinear interpolation
\begin{equation}\label{eq:bilinear-interpolation}
 \begin{split}
 (\operatorname{Interp}_N F)(t,x)
 &=\sum_{i,k\in\{0,1\}}\lambda_i(\theta)\lambda_k(\zeta)
       F\left(t_N+\frac{i}{N},x_N+\frac{k}{\sqrt N}\right),\\
 &\hspace{15mm}\lambda_0(u)=1-u,\qquad \lambda_1(u)=u.
 \end{split}
\end{equation}

This gives a continuous field agreeing with $F$ on the grid, and we
equip $C_{\rm loc}([0,\infty)\times\R)$ with the topology of uniform
convergence on compact sets. The interpolation of the field
\eqref{fr:terminal-field} is denoted by $\operatorname{Interp}_N\mathsf Z^N$.

We finally describe the regularization of the continuous equation. We write
\begin{equation}\label{eq:heat-kernel}
 \heat_\sigma(y)=(2\pi\sigma)^{-\frac12}e^{-\frac{y^2}{2\sigma}},
 \qquad \sigma>0,\ y\in\R,
\end{equation}
for the centered Gaussian density with variance $\sigma$.

We use the product Gaussian mollifier $\rho(t,x)=\heat_1(t)\heat_1(x)$ of \cite[Definition~3.4]{CDOT-KRP}, with parabolic rescaling
\begin{equation}\label{eq:parabolic-mollifier}
 \rho_\varepsilon(t,x)
 :=\varepsilon^{-3}\rho\left(\frac{t}{\varepsilon^2},\frac{x}{\varepsilon}\right)
 =\heat_{\varepsilon^4}(t)\heat_{\varepsilon^2}(x),\qquad \varepsilon>0.
\end{equation}
The smoothed noise is $\dot W_h^\varepsilon=\rho_\varepsilon*\dot W_h$
on $\R^2$, and the associated counterterm is defined by
\begin{equation}\label{fr:counterterm-continuous}
 \begin{split}
 c_h^\varepsilon&:=\varepsilon^{2a_h-2}c_h^1,\qquad \varepsilon>0,\\
 c_h^1&:=\int_0^\infty E_B\left[
   \E[\dot W_{h}^1(s,B_s)\dot W_{h}^1(0,0)]\right]\dd s.
 \end{split}
\end{equation}
Here $B$ is a standard Brownian motion independent of the noise.
The constant $c_h^1$ is finite and strictly positive; see
Proposition~\ref{prop:R-eps-limit}. Thus $c_h^\varepsilon$ tends to
$+\infty$ as $\varepsilon\downarrow0$, since $a_h<1$.

It remains to specify the class in which the regularized equation is
solved. A field $V$ is said to have \textit{at most exponential spatial growth on compact
time intervals} if, for every $T>0$, almost surely there are finite
constants $a_T,C_T\ge0$, which may depend on the realization, such that
\begin{equation}\label{eq:growth-class}
 \sup_{0\le t\le T}|V_t(x)|\le C_T e^{a_T|x|},
 \qquad x\in\R.
\end{equation}
For each $\varepsilon>0$, let $U^\varepsilon$ be the unique classical
solution in this class of
\begin{equation}\label{fr:regularized-equation}
 \partial_t U^\varepsilon
 =\tfrac12\Delta U^\varepsilon+(\dot W_h^\varepsilon-c_h^\varepsilon)U^\varepsilon,
 \qquad U^\varepsilon(0,\cdot)\equiv1.
\end{equation}
Lemma~\ref{lem:regularized-FK} proves existence and uniqueness and gives
the Feynman-Kac formula for $U^\varepsilon$.

We can now state our main result.

\begin{theorem}[Space-time fractional construction and approximation]
\label{thm:fractional-main}
Assume \eqref{fr:range}. There exists a nonnegative continuous random
field $U$, defined on the Gaussian space of $W_{h}$, with
$U_0(x)=1$, such that for every finite $T,L>0$, every compact
$\mathcal K\subset[0,T]\times[-L,L]$, and every $1\le p<\infty$,

\begin{equation}\label{fr:continuous-functional}
 U^\varepsilon\longrightarrow U
 \qquad\text{in }L^p(\Omega;C(\mathcal K)),
\end{equation}
\begin{equation}\label{fr:original-functional}
 \operatorname{Interp}_N\mathsf Z^N\ToLaw U
 \quad\text{in }C_{\rm loc}([0,\infty)\times\R),
\end{equation}
and
\begin{equation}\label{fr:compact-moments}
 \E\|\operatorname{Interp}_N\mathsf Z^N\|_{C(\mathcal K)}^p
 \longrightarrow\E\|U\|_{C(\mathcal K)}^p<\infty.
\end{equation}

If $2H_0+H_1>\frac53$, $U$
coincides with the solution of \cite[Definition~3.4]{CDOT-KRP}.
\end{theorem}

\begin{definition}[Solution obtained by the fixed regularization]
\label{def:continuous-solution}
The field $U$ of Theorem~\ref{thm:fractional-main} is called the
\textit{renormalized Stratonovich solution obtained by the stationary
Gaussian regularization with counterterm
\eqref{fr:counterterm-continuous}}. This means the limit of the
classical equations \eqref{fr:regularized-equation} in
$L^p(\Omega;C(\mathcal K))$ for every finite $p\ge1$ and every compact
space-time rectangle $\mathcal K$. The terminology refers to this
specified regularization and normalization.
\end{definition}

\begin{remark}
Both counterterms follow the same stationary prescription: they are
the asymptotic mean conditional half-variance per unit time, namely
\begin{equation}\label{eq:stationary-counterterm-rates}
\begin{aligned}
 c_h^\varepsilon
 &=\lim_{T\to\infty}\frac1{2T}E_B\!\left[
    \E\!\left[\left(\int_0^T\dot W_h^\varepsilon(s,B_s)\dd s\right)^2
             \Bigm|B\right]\right],\qquad\varepsilon>0,\\
 c_{\rm disc}(h)
 &=\lim_{m\to\infty}\frac1{2m}E_S\!\left[
    \E\!\left[\left(\sum_{n=1}^m\omega_h(n,S_n)\right)^2
             \Bigm|S\right]\right].
\end{aligned}
\end{equation}
The Gaussian variance is averaged over the path; the discrete formula
uses unit coupling, so that the rate at coupling $\beta_N$ is
$c_{\rm disc}(h)\beta_N^2$ per step. Stationarity of the covariance
and of path increments yields these identities from the integral and
series above. 
This prescription therefore differs from subtracting the full
conditional half-variance separately for each path.
\end{remark}

Theorem~\ref{thm:fractional-main} is made explicit by
Theorems~\ref{thm:continuous-summary} and~\ref{thm:coupled-L2} below,
which construct the continuous field and prove the discrete convergence
on a common Gaussian space, and by Lemma~\ref{lem:environment-law},
which transfers the latter to the original environment.

\

\subsection{The continuous model}
\label{sec:continuous}

We now describe the continuous field more precisely. This requires
three ingredients: the Gaussian space of the noise, the Brownian
functionals that appear in the limit, and the regularization and
renormalization conventions for the equation.

\subsubsection{The Gaussian space and the \texorpdfstring{$\mathcal S$}{S}-transform}

Recall that the Gaussian density $\heat_t$ of \eqref{eq:heat-kernel}
is the heat kernel of $\frac12\Delta$.
Let $W=\{W(f):f\in\Hh\}$ be the isonormal process associated with
\eqref{fr:noise-covariance}. The space $\Hh=\mathcal H_{h}$ is the
completion of $\D=C_c^\infty(\R\times\R)$ for the covariance inner
product
\begin{equation}\label{eq:H-inner}
 \ip{\varphi}{\psi}_{\Hh}
 =\int_{\R^4}
 \varphi(s,x)\psi(t,y)K_{H_0}(t-s)K_{H_1}(x-y)
                                       \dd x\dd y\dd s\dd t.
\end{equation}
At $H_1=\frac12$ the spatial kernel is $\delta_0$ and the double
spatial integral collapses to a single one; this is the form used in
Sections~\ref{sec:Brownian-construction}-\ref{sec:short}.
The noise is stationary on the full time axis; the SPDE is solved for
$t\ge0$. In particular, temporal mollification below is a convolution on
$\R$, without truncation at $0$.
The continuous path averages and the coupled discrete partition functions
are considered in
$L^2(\Omega,\sigma(W),\Pp)$. The notation $\E$ refers to Gaussian expectation;
$E_B$ and $E_S$, or simply $E$ when no ambiguity is possible, refer to
path expectation. All auxiliary paths are independent of $W$.

For $f\in\Hh$, define the Gaussian exponential
\begin{equation}\label{eq:stoch-exp}
 \cE(f)=\exp\left(W(f)-\frac12\norm f_{\Hh}^2\right).
\end{equation}
For $f,g\in\Hh$, the Gaussian moment-generating formula gives
\begin{equation}\label{eq:Gaussian-pairing}
 \E[\cE(f)]=1,\qquad
 \E[\cE(f)\cE(g)]=e^{\ip{f}{g}_{\Hh}}.
\end{equation}
In particular, $\cE(f)\in L^2(\Omega)$.

For $X\in L^2(\Omega,\sigma(W),\Pp)$,
its $\mathcal S$-transform is defined by
\begin{equation}\label{eq:S-transform}
 \mathcal S_X(\varphi)=\E\left[X\cE(\varphi)\right],\qquad \varphi\in\Hh.
\end{equation}

For background on Gaussian Hilbert spaces and Gaussian exponentials,
see \cite{Janson-Gaussian}; the $\mathcal S$-transform is treated in
\cite[Chapter~XVI]{Janson-Gaussian}. 

The $\mathcal S$-transform of the solution involves a linear functional
of the Brownian path, which we now introduce. For $\varphi\in\D$, put
\begin{equation}\label{eq:G-phi}
 G_\varphi(s,x)
 =\int_{\R^2}K_{H_0}(t-s)K_{H_1}(x-y)
                                      \varphi(t,y)\dd t\dd y.
\end{equation}
This function is bounded and globally uniformly continuous on $\R^2$.
Indeed, both kernels are locally integrable and bounded away from
their singularities, so that convolution with a smooth compactly
supported function, and the same argument applied to its derivatives,
give uniform bounds for $G_\varphi$ and its first derivatives. At
$H_1=\frac12$ the spatial integration disappears and
$G_\varphi(s,x)=\int_\R K_{H_0}(t-s)\varphi(t,x)\dd t$. For
$t>0$, $x\in\R$ and every continuous path $B$ on $[0,t]$, define
\begin{equation}\label{eq:ell-B-general}
 \ell^B_{t,x}(\varphi)=\int_0^tG_\varphi(s,x+B_{t-s})\dd s.
\end{equation}

These linear functionals are bounded in absolute value by
$t\|G_\varphi\|_\infty$. Unlike the path functionals $Y_t$ and $\alpha_t$ defined below, $\ell^B_{t,x}$ is not
homogeneous in $t$, the test function $\varphi$ being fixed.

We finally fix the constants that enter the deterministic remainder. Let $Z$ be a standard normal random
variable and set
\begin{equation}\label{eq:AH}
 \kappa_{H_1}=E[K_{H_1}(Z)]
 =\frac{H_1(2H_1-1)2^{H_1-1}\Gamma(H_1-\frac12)}{\sqrt\pi},
 \qquad A_{h}=c_{H_0}\kappa_{H_1},
\end{equation}
and, for every $t\ge0$,
\begin{equation}\label{eq:R-definition}
 R_t=-A_{h}\,t^{a_{h}}\left(\frac1{a_{h}}+\frac1{1-a_{h}}\right).
\end{equation}
Both are finite precisely under \eqref{fr:range}, the second because
$0<a_{h}<1$. The expression for $\kappa_{H_1}$ has limit
$\frac1{\sqrt{2\pi}}$ as $H_1\downarrow\frac12$; accordingly we write
$A_{H_0}=c_{H_0}/\sqrt{2\pi}$ for the value of $A_{h}$ at
$H_1=\frac12$, where $a_{h}=2H_0-\frac12$ and \eqref{eq:R-definition} reads
$R_t=-A_{H_0}t^{2H_0-\frac12}
\bigl((2H_0-\tfrac12)^{-1}+(\tfrac32-2H_0)^{-1}\bigr)$.
Propositions~\ref{prop:R-eps-limit} and~\ref{prop:counterterm} identify
$R_t$ as the common finite remainder of the continuous and discrete
renormalizations.

\subsubsection{Spatial kernels and mutual intersections}

For the regularized kernels, we use the convention
\begin{equation}\label{fr:kernels}
 K_H^{[v]}:=K_H*\heat_v\quad(v>0,\ H\in[\tfrac12,1)),
 \qquad K_H^{[0]}:=K_H.
\end{equation}
This notation applies to both the temporal and spatial kernels.
The covariance kernel of two spatial regularizations of sizes
$\varepsilon,\varepsilon'$ is
$K_{H_1}^{[\varepsilon^2+\varepsilon'^2]}$.
At $H_1=\frac12$, one has $K_{1/2}=\delta_0$,
$K_{1/2}^{[v]}=\heat_v$ for $v>0$, and
$\gamma_{1/2}=\one_{\{0\}}$.

The second moment of the solution involves a pair of paths. For two independent standard Brownian motions $B,\widetilde B$
started at zero and an arbitrary horizon $t>0$, we define their mutual intersection functional by
\begin{equation}\label{fr:alpha-continuous}
 \alpha_t(B,\widetilde B)
 =\int_{[0,t]^2}K_{H_0}(r-s)
                    K_{H_1}(B_r-\widetilde B_s)\dd r\dd s.
\end{equation}

At $H_1=\frac12$, where $K_{1/2}=\delta_0$, this definition is read
through the mutual occupation measure of the two paths; see
\eqref{eq:alpha-endpoint} in Section~\ref{sec:continuous-identification}.
Lemma~\ref{tp:Brownian-cutoff} proves finiteness and
integrability of $\alpha_t(B,\widetilde B)$ for this independent pair, and
Proposition~\ref{prop:continuous-mutual} gives the higher moment and
exponential bounds, together with the convergence of the regularized mutual
intersections.

\subsubsection{\texorpdfstring{The centered self-intersection}{The centered self-intersection}}\label{subsub:pre-cons}

We next introduce the centered Brownian self-intersection process
$(Y_t(B))_{t\ge0}$, which enters the renormalized Feynman-Kac
representation of the solution. We define it by decomposing the time
triangle into dyadic rectangles and centering the interaction on each
rectangle before summing. The convergence and integrability of this
construction are established in Section~\ref{sec:Brownian-construction}.

\begin{figure}[tbp]
\centering
\begingroup
\definecolor{dyadic1}{HTML}{440154}
\definecolor{dyadic2}{HTML}{482878}
\definecolor{dyadic3}{HTML}{3E4989}
\definecolor{dyadic4}{HTML}{31688E}
\definecolor{dyadic5}{HTML}{26828E}
\definecolor{dyadic6}{HTML}{1F9E89}
\definecolor{dyadic7}{HTML}{35B779}
\definecolor{dyadic8}{HTML}{6ECE58}
\definecolor{dyadic9}{HTML}{B5DE2B}
\definecolor{dyadic10}{HTML}{FDE725}
\begin{minipage}[b]{0.56\linewidth}
\centering
\begin{tikzpicture}[x=7cm,y=7cm,font=\small]
 \fill[black!3] (0,0)--(0,1)--(1,1)--cycle;
 \foreach \M in {1,...,10}{
   \pgfmathtruncatemacro{\nblocks}{2^(\M-1)}
   \foreach \k in {1,...,\nblocks}{
     \pgfmathsetmacro{\blockleft}{(2*\k-2)*2^(-\M)}
     \pgfmathsetmacro{\blockmid}{(2*\k-1)*2^(-\M)}
     \pgfmathsetmacro{\blockright}{2*\k*2^(-\M)}
     \fill[dyadic\M]
       (\blockleft,\blockmid) rectangle (\blockmid,\blockright);
   }
 }
 \draw[black!55,densely dashed,line width=0.4pt] (0,0)--(1,1);
 \draw[->] (0,0)--(1.045,0) node[right] {$s$};
 \draw[->] (0,0)--(0,1.045) node[above] {$t$};
 \foreach \pos/\ticktext in {0/0,0.5/{\frac12},1/1}{
   \draw (\pos,0)--(\pos,-0.012)
      node[below] {$\ticktext$};
 }
 \foreach \pos/\ticktext in {0.5/{\frac12},1/1}{
   \draw (0,\pos)--(-0.012,\pos)
      node[left] {$\ticktext$};
 }
 \node[text=white] at (0.25,0.75) {$M=1$};
 \node[text=white,font=\footnotesize] at (0.125,0.375) {$M=2$};
 \node[text=white,font=\footnotesize] at (0.625,0.875) {$M=2$};
 \node[text=black!55] at (0.68,0.25) {$s>t$};
 \node[rotate=45,fill=white,inner sep=1.5pt]
      at (0.72,0.72) {$t=s$};
 \draw[red!75!black,line width=0.8pt]
      (0,0) rectangle (0.03125,0.03125);
 \draw[->,red!75!black] (0.18,0.07)--(0.04,0.02);
 \node[anchor=west,text=red!75!black,font=\footnotesize]
      at (0.18,0.07) {zoom};
\end{tikzpicture}
\end{minipage}\hfill
\begin{minipage}[b]{0.40\linewidth}
\centering
\begin{tikzpicture}[x=4.6cm,y=4.6cm,font=\small]
 \fill[black!3] (0,0)--(0,1)--(1,1)--cycle;
 \foreach \M in {6,...,10}{
   \pgfmathtruncatemacro{\nblocks}{2^(\M-6)}
   \foreach \k in {1,...,\nblocks}{
     \pgfmathsetmacro{\blockleft}{(2*\k-2)*2^(5-\M)}
     \pgfmathsetmacro{\blockmid}{(2*\k-1)*2^(5-\M)}
     \pgfmathsetmacro{\blockright}{2*\k*2^(5-\M)}
     \fill[dyadic\M]
       (\blockleft,\blockmid) rectangle (\blockmid,\blockright);
   }
 }
 \draw[black!55,densely dashed,line width=0.4pt] (0,0)--(1,1);
 \draw[->] (0,0)--(1.05,0) node[right] {$s$};
 \draw[->] (0,0)--(0,1.05) node[above] {$t$};
 \foreach \pos/\ticktext in {0/0,0.5/{\frac1{64}},1/{\frac1{32}}}{
   \draw (\pos,0)--(\pos,-0.018) node[below] {$\ticktext$};
 }
 \foreach \pos/\ticktext in {0.5/{\frac1{64}},1/{\frac1{32}}}{
   \draw (0,\pos)--(-0.018,\pos) node[left] {$\ticktext$};
 }
 \node at (0.25,0.75) {$M=6$};
 \node[font=\footnotesize] at (0.125,0.375) {$M=7$};
 \node[font=\footnotesize] at (0.625,0.875) {$M=7$};
\end{tikzpicture}
\end{minipage}
\par\medskip
{\small Level $M$:\quad
\begin{tikzpicture}[x=1.05cm,y=0.8cm,baseline=0.02cm,font=\footnotesize]
 \foreach \M in {1,...,10}{
   \pgfmathsetmacro{\legendx}{\M-1}
   \fill[dyadic\M] (\legendx,0) rectangle ++(0.28,0.20);
   \node[anchor=west,inner sep=2pt]
      at (\legendx+0.30,0.10) {$\M$};
 }
\end{tikzpicture}}
\endgroup
\caption{The dyadic rectangles $I_{M,k}\times J_{M,k}$,
with $1\le M\le10$ and $1\le k\le2^{M-1}$, in the $(s,t)$-plane.
Colors identify the level $M$. The right panel enlarges
$[0,\frac1{32}]^2$, where levels $6$ through $10$ can be distinguished.
At any finite level, smaller triangles along $s=t$ remain uncovered;
the full decomposition is obtained by allowing all levels.}
\label{fig:dyadic-rectangles}
\end{figure}
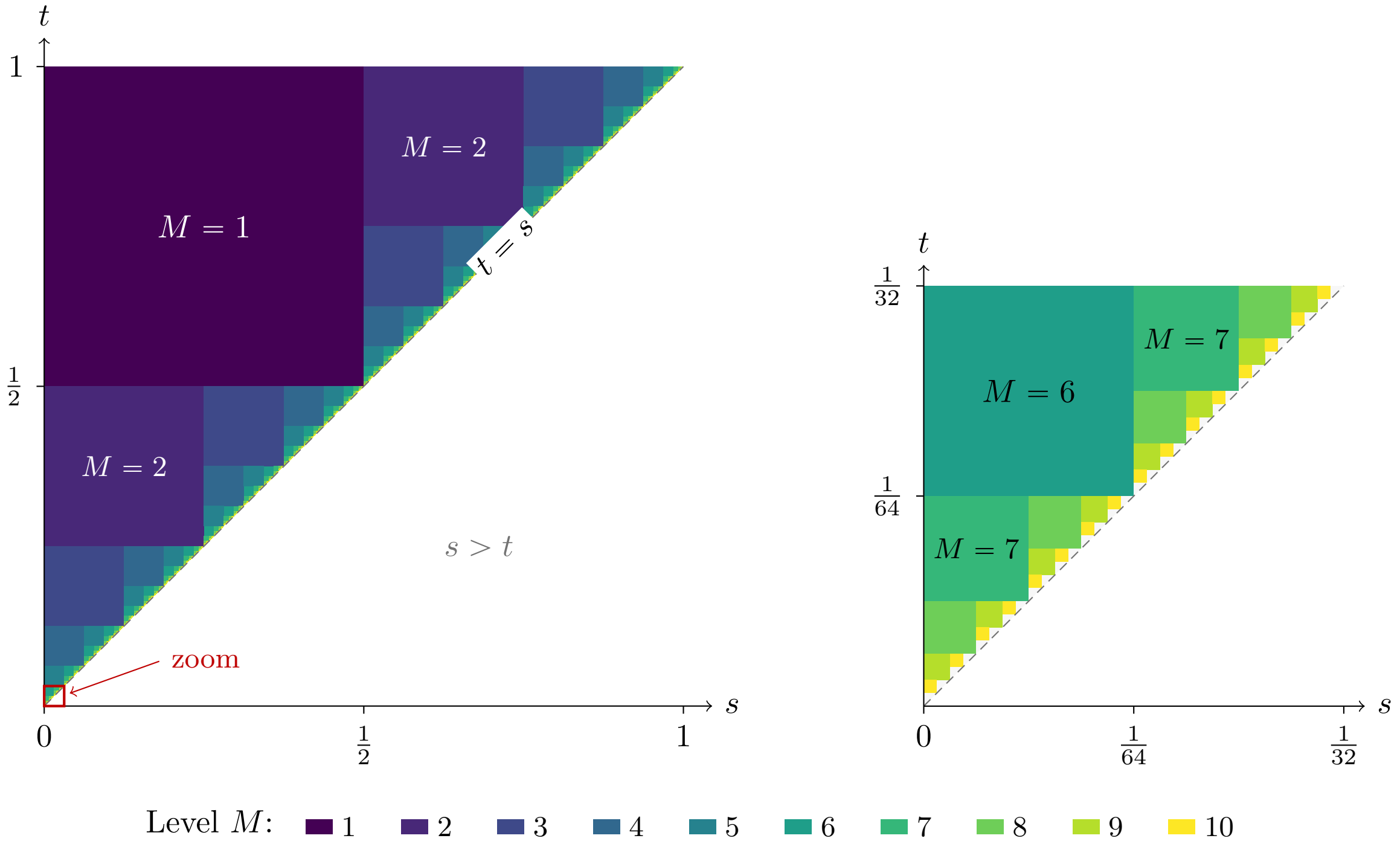

For $M\ge1$ and $1\le k\le2^{M-1}$, set
\begin{align}
& I_{M,k}=\big[\frac{2k-2}{2^{M}},\frac{2k-1}{2^{M}}\big), \qquad  J_{M,k}=\big[\frac{2k-1}{2^{M}},\frac{2k}{2^{M}}\big).\label{eq:dyadic-JI}
\end{align}
These rectangles form the disjoint partition
\begin{equation}\label{eq:dyadic-triangle-partition}
 \{(s,t):0\le s<t<1\}
 =\bigsqcup_{M\ge1}\bigsqcup_{k=1}^{2^{M-1}}
                          (I_{M,k}\times J_{M,k}).
\end{equation}
Indeed, for $s<t$, take the first dyadic subdivision at which the
two points lie in different halves of the same parent interval.
The half-open convention resolves boundary points and makes this
choice unique. This dyadic decomposition goes back to Varadhan's
construction of the renormalized self-intersection local time of planar
Brownian motion \cite{Varadhan}; see also \cite{LeGall-StFlour}.
Figure~\ref{fig:dyadic-rectangles} displays the first ten levels.

Writing $I_{M,k,t}=tI_{M,k}$ and $J_{M,k,t}=tJ_{M,k}$
for the rectangles scaled to $[0,t]$, define the Brownian blocks by
\begin{align}
 C_{M,k,t}(B)
 &:=\int_{I_{M,k,t}\times J_{M,k,t}}
            K_{H_0}(s-r)K_{H_1}(B_s-B_r)\dd r\dd s,
                                              \label{fr:Brownian-block}\\
 A_{M,k,t}&:=C_{M,k,t}-E_B[C_{M,k,t}].
 \label{eq:block-continuous}
\end{align}
Every block is finite, and its expectation is the explicit integral
\eqref{eq:Brownian-block-mean} of $(s-r)^{a_{h}-2}$ against the
constant $A_{h}$ of \eqref{eq:AH}, by
Lemma~\ref{lem:Brownian-blocks}. 
At $H_1=\frac12$, the blocks are read through the self-intersection
measure of $B$; see \eqref{eq:block-endpoint} below.

\begin{proposition}[Centered self-intersection]
\label{prop:Y}
Under \eqref{fr:range}, the series
\begin{equation}\label{fr:Y-definition}
 Y_t(B):=2\sum_{M\ge1}\sum_{k=1}^{2^{M-1}}A_{M,k,t}
\end{equation}
converges in every finite $L^p$ and defines a centered path
functional. Moreover, there exists a constant
$C_{h}>0$ such that, for $t\le T$ and every $\lambda\in\R$,
\begin{equation*}
 E[e^{\lambda Y_t(B)}]\le e^{C_{h}\lambda^2T^{2a_{h}}}.
\end{equation*}
\end{proposition}

The proof is given in Section~\ref{sec:Brownian-construction}; see Theorem~\ref{thm:Brownian-self-intersection}.

\subsubsection{Regularized noise and Stratonovich solution}\label{subsub:regul-noise-strato}

For fixed $h$, write $\dot W^\varepsilon=\dot W_h^\varepsilon$
for the noise smoothed with \eqref{eq:parabolic-mollifier}. Its temporal and spatial covariance factors are respectively
$K_{H_0}^{[2\varepsilon^4]}$ and $K_{H_1}^{[2\varepsilon^2]}$. Thus
\begin{equation}\label{eq:smoothed-noise-covariance}
 \E\left[\dot W^\varepsilon(s,x)\dot W^\varepsilon(t,y)\right]
 =K_{H_0}^{[2\varepsilon^4]}(t-s)K_{H_1}^{[2\varepsilon^2]}(x-y).
\end{equation}

The counterterm is fixed by the scaling definition
\eqref{fr:counterterm-continuous}. Averaging the covariance at unit
regularization along $B$ gives the explicit constant
\begin{equation}\label{eq:counterterm-unit}
 c_h^1=\kappa_{H_1}\int_0^\infty
       K_{H_0}^{[2]}(v)(v+2)^{H_1-1}\dd v\in(0,\infty),
\end{equation}
where $\kappa_{H_1}$ is defined in \eqref{eq:AH}.
Proposition~\ref{prop:R-eps-limit} verifies finiteness and identifies
the corresponding time-domain formula at every scale.

We can now state the construction of the continuous field, together
with the two identities that characterize it.

\begin{theorem}[Construction and characterization of the continuous field]
\label{thm:continuous-summary}
Assume \eqref{fr:range}. For each $\varepsilon>0$, let
$U^\varepsilon$ be the unique classical solution of
\eqref{fr:regularized-equation} in the growth class \eqref{eq:growth-class}.
The following conclusions hold on the given Gaussian space.
\begin{enumerate}
\item There is a continuous nonnegative field $U$, measurable
with respect to $W$, such that, at every fixed $(t,x)$,
\[
 U_t^\varepsilon(x)\longrightarrow U_t(x)
 \quad\text{in }L^p(\Omega),\qquad 1\le p<\infty,
 \quad\varepsilon\downarrow0.
\]
For every finite $T>0$ and
$1\le p<\infty$,
\[
 \sup_{\substack{\varepsilon>0,\ 0\le t\le T\\x\in\R}}
       \E[|U_t^\varepsilon(x)|^p]
 +\sup_{\substack{0\le t\le T\\x\in\R}}\E[|U_t(x)|^p]
 <\infty.
\]
Moreover, $\E[U_t(x)]>0$ for $t>0$.
\item For every finite $T,L>0$, set $\mathcal K_{T,L}=[0,T]\times[-L,L]$. Then, for $1\le p<\infty$,
\[
 \E\left[\|U^\varepsilon-U\|_{C(\mathcal K_{T,L})}^p\right]
 \longrightarrow0,
 \qquad \E\left[\sup_{|x|\le L}|U_t(x)-1|^p\right]
       \longrightarrow0\quad(t\downarrow0).
\]

\item At every point $(t,x)$ with $t>0$, for every $\varphi\in\D$, the
$\mathcal S$-transform satisfies
\begin{equation}\label{eq:S-U-theo}
 \mathcal S_{U_{t}(x)}(\varphi)
 =e^{R_{t}} E_B\left[\exp\left(\frac12Y_{t}(B)
                 +\ell^B_{t,x}(\varphi)\right)\right].
\end{equation}
For two independent standard Brownian motions $B,\widetilde B$ on $[0,t]$, the
second moment is
\begin{equation}\label{eq:continuous-second-exact-theo}
 \E[U_{t}(x)^2]
 =e^{2R_{t}}E_{B,\widetilde B}\left[
   \exp\left(\frac12\bigl(Y_{t}(B)+Y_{t}(\widetilde B)\bigr)
                    +\alpha_{t}(B,\widetilde B)\right)\right].
\end{equation}
Here $R_t$ is \eqref{eq:R-definition}, $\ell^B_{t,x}$ is
\eqref{eq:ell-B-general}, $Y_t(B)$ is \eqref{fr:Y-definition} and
$\alpha_t(B,\widetilde B)$ is \eqref{fr:alpha-continuous}.

\item If $2H_0+H_1>\frac53$, this field agrees at each fixed
$(t,x)$ almost surely with the solution constructed in \cite{CDOT-KRP}.
\end{enumerate}
\end{theorem}

\begin{remark}[Expected extensions of the pathwise construction]
\label{rem:pathwise-extensions}
The extension of the continuous construction to
$\frac32<2H_0+H_1\le\frac53$ was already expected from the perspective
of regularity structures. As explained in
\cite[Remark~3.6]{CDOT-KRP}, the deterministic solution theory extends
to noise regularities $\alpha>-\frac32$ upon adjoining components of
order three to the enhanced noise. In dimension one, this corresponds
to the threshold $2H_0+H_1>\frac32$; implementing this approach on the
whole space requires the corresponding stochastic construction in
suitable weighted spaces.

The nonlinear framework of \cite{De16}, together with the extensions
discussed there, also suggests that a local renormalized construction
should remain possible in this range for equations with nonlinear
multiplicative terms $F(x,u)\dot W$.
\end{remark}

Theorem~\ref{thm:continuous-summary} is proved in
Section~\ref{sec:continuous-identification}, throughout the range
\eqref{fr:range}.
Proposition~\ref{prop:CDOT-comparison} establishes the comparison with
the solution of \cite{CDOT-KRP} under the additional condition
$2H_0+H_1>\frac53$.

\

Theorem~\ref{thm:continuous-summary} also allows one to prove,
\textit{a posteriori}, that the counterterm cannot be dispensed with
when $a_{h}<1$. 

\begin{proposition}[Necessity of a renormalization]
\label{prop:renormalization-needed}
Assume \eqref{fr:range}. For $\varepsilon>0$, let $u^\varepsilon$ be the
unique classical solution of
\begin{equation}\label{eq:unrenormalized-PDE}
 \partial_tu^\varepsilon=\frac12\Delta u^\varepsilon+\dot W^\varepsilon u^\varepsilon,
 \qquad u^\varepsilon_0\equiv1,
\end{equation}
with at most exponential spatial growth on compact time intervals, in
the sense of \eqref{eq:growth-class}. Then, for every $t>0$, $x\in\R$
and $K>0$,
\[
 \limsup_{\varepsilon\downarrow0}\Pp\bigl(u^\varepsilon_t(x)\le K\bigr)
 \le\Pp\bigl(U_t(x)=0\bigr)<1,
\]
where $U$ is the field of Theorem~\ref{thm:continuous-summary}. In
particular, for no sequence $\varepsilon_n\downarrow0$ is the family
$(u^{\varepsilon_n}_t(x))_{n\ge1}$ tight, and none of these sequences
converges in law.
\end{proposition}

The proof is given in Section~\ref{subsec:renormalization-necessity}. Recall that when $a_{h}>1$, by contrast, no renormalization is
needed \cite{HNS}.

\subsection{Convergence of the discrete model}\label{sec:common-noise}

We now realize the discrete environment on the Gaussian space of
Section~\ref{sec:continuous}, so that the discrete and continuous objects
can be compared in $L^2(\Omega)$. The comparison is made at an
arbitrary point $(t,x)\in[0,\infty)\times\R$, using the notation
$N_t$ and $\widetilde{N}_x$ from \eqref{eq:Nt} and \eqref{eq:xN}.
We use the same isonormal process $W=W_{h}$ as for $U$. Introduce
the cells
\[
 \mathfrak{I}_{t,n}^N=\left(t-\frac nN,t-\frac{n-1}{N}\right],\quad 1\le n\le N_t,
 \qquad \mathfrak{J}_j^N=\left[\frac{j}{\sqrt N},\frac{j+1}{\sqrt N}\right),
\]
and set
\begin{equation}\label{eq:cell-functions}
 e_{t,n,j}^N=N^{H_0+\frac{H_1}2}\one_{\mathfrak{I}_{t,n}^N\times \mathfrak{J}_j^N}.
\end{equation}

The coupled environment is then defined by
\begin{equation}\label{eq:common-noise}
 \omega_t^N(n,j)=W(e_{t,n,j}^N).
\end{equation}

The terminal-time array $\omega_t^N$ realizes the covariance
\eqref{fr:discrete-covariance} exactly. 

\begin{proposition}[Exact isonormal coupling]\label{prop:coupling}
For all $1\le n,m\le N_t$ and $j,\ell\in\mathbb Z$,
\begin{equation}\label{eq:cell-covariance}
 \ip{e_{t,n,j}^N}{e_{t,m,\ell}^N}_{\Hh}
 =\gamma_{H_0}(n-m)\gamma_{H_1}(j-\ell).
\end{equation}
In particular, the Gaussian family \eqref{eq:common-noise} has the law of the environment
\eqref{fr:discrete-covariance}. Moreover, if
\begin{equation}\label{eq:hNS-explicit}
 h^{S,N}_{t,x}=\beta_N\sum_{n=1}^{N_t} e_{t,n,\widetilde{N}_x+S_n}^N
 =\sqrt N\sum_{n=1}^{N_t}\one_{\mathfrak{I}_{t,n}^N\times \mathfrak{J}_{\widetilde{N}_x+S_n}^N},
\end{equation}
then $W(h^{S,N}_{t,x})=\beta_N\sum_{n=1}^{N_t}\omega_t^N(n,\widetilde{N}_x+S_n)$.
\end{proposition}

The proof is given in Section~\ref{subsec:coupling-identification}.

Define the partition function associated with this realization by
\begin{equation}\label{eq:ZN-coupled}
 \widehat{\mathsf Z}^N_t(x)
 :=e^{-c_{\rm disc}N_t\beta_N^2}
 E_S\left[
   \exp\left(\beta_N\sum_{n=1}^{N_t}\omega_t^N(n,\widetilde{N}_x+S_n)\right)\right]
 {}=E_S\left[e^{W(h^{S,N}_{t,x})-c_{\rm disc}N_t\beta_N^2}\right].
\end{equation}

We note that $\widehat{\mathsf Z}^N_0(x)=1$.

\begin{lemma}[Equality in law of the environments and interpolated fields]
\label{lem:environment-law}
For each $N\ge1$ and each finite $t\ge0$,
\begin{equation}\label{eq:environment-law}
 \bigl(\omega_t^N(n,j)\bigr)_{\substack{1\le n\le N_t\\j\in\mathbb Z}}
 \overset{\rm law}=
 \bigl(\omega(n,j)\bigr)_{\substack{1\le n\le N_t\\j\in\mathbb Z}}.
\end{equation}
Consequently, for every $x\in\R$, the two variables
$\widehat{\mathsf Z}^N_t(x)$ and $\mathsf Z^N_t(x)$ have the same law. Moreover, for each $N\ge1$, the interpolations \eqref{eq:bilinear-interpolation}
of the two fields satisfy
\begin{equation}\label{eq:interpolated-field-law}
 \operatorname{Interp}_N\mathsf Z^N\overset{\rm law}=\operatorname{Interp}_N\widehat{\mathsf Z}^N
 \quad\text{in }C_{\rm loc}([0,\infty)\times\R).
\end{equation}

\end{lemma}

The proof is given in Section~\ref{subsec:coupling-identification}.

The interpolated field $\operatorname{Interp}_N\widehat{\mathsf Z}^N$
is defined on the same Gaussian space as $U$ and equals $1$ at time zero. The two fields can therefore be compared in a strong sense, and not
only in law.

\begin{theorem}[Convergence in the common Gaussian space]
\label{thm:coupled-L2}
Under \eqref{fr:range}, use the counterterm $c_{\rm disc}(h)$ of
\eqref{fr:counterterm-discrete} and the realization
\eqref{eq:common-noise}.
Let $U$ be the field of Theorem~\ref{thm:continuous-summary}. Then,
for every finite $t\ge0$ and every $x\in\R$,
\begin{equation}\label{eq:main-L2}
 \widehat{\mathsf Z}^N_t(x)\longrightarrow U_{t}(x)
 \quad\text{in }L^p(\Omega,\sigma(W),\Pp),\qquad 1\le p<\infty.
\end{equation}

Moreover, for every finite $T,L>0$ and $1\le p<\infty$,
\begin{equation}\label{eq:main-functional}
 \E\left[
   \|\operatorname{Interp}_N\widehat{\mathsf Z}^N-U\|_{C([0,T]\times[-L,L])}^p
 \right]\longrightarrow0.
\end{equation}
Thus the interpolated fields converge in probability in
$C_{\rm loc}([0,\infty)\times\R)$.

\end{theorem}

The proof of Theorem~\ref{thm:coupled-L2} occupies Sections~\ref{sec:discrete} to~\ref{sec:short}; it is completed in
Section~\ref{subsec:proof-coupled-L2}.

Theorem~\ref{thm:coupled-L2} also allows one to prove,
\textit{a posteriori}, that the lattice counterterm cannot be dispensed
with either. The statement is the exact counterpart of
Proposition~\ref{prop:renormalization-needed}.

\begin{proposition}[Necessity of the lattice counterterm]
\label{prop:discrete-renormalization-needed}
Assume \eqref{fr:range}. For $N\ge1$, $t\ge0$ and $x\in\R$, let
\begin{equation}\label{eq:unrenormalized-partition}
 \mathsf z^N_t(x)
 =E_S\left[\exp\left(\beta_{N,h}\sum_{n=1}^{N_t}
                    \omega_{h}(N_t-n+1,\widetilde{N}_x+S_n)\right)\right]
 =e^{c_{\rm disc}(h)N_t\beta_{N,h}^2}\,\mathsf Z^N_t(x)
\end{equation}
be the partition function \eqref{fr:terminal-field} without its
counterterm. Then, for every $t>0$, $x\in\R$ and $K>0$,
\[
 \limsup_{N\to\infty}\Pp\bigl(\mathsf z^N_t(x)\le K\bigr)
 \le\Pp\bigl(U_t(x)=0\bigr)<1,
\]
where $U$ is the field of Theorem~\ref{thm:continuous-summary}. In
particular, for no sequence $N_k\to\infty$ is the family
$(\mathsf z^{N_k}_t(x))_{k\ge1}$ tight, and none of these sequences
converges in law.
\end{proposition}

The proof is given in Section~\ref{subsec:discrete-necessity}.

\subsection{Comments and open questions}\label{sec:comments}

We conclude this section with some comments on the scope of
Theorems~\ref{thm:fractional-main},
\ref{thm:continuous-summary} and~\ref{thm:coupled-L2}, and on questions
that are not addressed in this article.

\begin{remark}[Rate of convergence]\label{rk:rate}
The proof of Theorem~\ref{thm:coupled-L2} is qualitative. Strong convergence in $L^2(\Omega)$ is obtained by combining
convergence of the second moments with that of the mixed moment.
The latter follows from convergence of the $\mathcal S$-transforms
by approximating $U_t(x)$ in $L^2$ with linear combinations of
Gaussian exponentials (Lemma~\ref{lem:total}), an argument that
provides no rate.
A rate of convergence would require a quantitative
control of the mixed moment
$\E[\widehat{\mathsf Z}^N_1(0)U_1(0)]$, and hence a quantitative
version of the joint invariance principle
(Sections~\ref{subsec:legall-complet} and~\ref{sec:two-path}); the
same applies at any other point.

Although no rate is established, one may expect the convergence to
slow down as $a_h\downarrow\frac12$, that is, as $(H_0,H_1)$ approaches
$(\frac12,\frac12)$. This heuristic is consistent with the noisier
profiles observed at smaller $H_0$ and fixed $N$ in
Section~\ref{sec:numerics}, where $H_1=\frac12$ and
$a_h-\frac12=2H_0-1$.
\end{remark}

\begin{remark}[Characterization and choice of regularization]\label{rk:solution}
The identity \eqref{eq:S-U-theo} characterizes
$U_t(x)$ uniquely in $L^2(\Omega,\sigma(W),\Pp)$: the Gaussian
exponentials with arguments in $\D$ have dense linear span in this
space (Lemma~\ref{lem:total}). Any two continuous fields on the same
Gaussian space satisfying this identity therefore agree outside one
null set, simultaneously at all space-time points.

The right-hand side of \eqref{eq:S-U-theo} contains no mollifier: it thus provides a criterion
for identifying limits obtained from other regularizations. Applying
that criterion requires convergence of their centered path
functionals and agreement of their deterministic remainders; the
characterization alone does not establish these facts. The matching
of remainders is proved here for the lattice and Gaussian
regularizations in Propositions~\ref{prop:R-eps-limit}
and~\ref{prop:counterterm}.

Well-posedness in an intrinsically defined solution space and
independence of the mollifier are distinct questions. When $2H_0+H_1>\frac53$, $U$ is identified with the
regularity-structures solution of \cite{CDOT-KRP}
(Theorem~\ref{thm:continuous-summary}(4)). We do not address these two questions
in the remaining range $\frac32<2H_0+H_1\le\frac53$.
\end{remark}

\begin{remark}[Parameter restrictions and excluded boundaries]
\label{fr:range-remark}\label{rk:endpoints}
The proof uses $2H_0-1>0$ to remove the temporal cutoff in mutual
intersections, $a_h-\frac12>0$ to sum the centered dyadic levels,
and $a_h<1$ for the stationary counterterm and the unequal-horizon
self-intersection estimates. The first two inequalities follow from
$H_0>\frac12$ and $H_1\ge\frac12$. The separate treatment of the two situations $H_1>\frac12$ and $H_1=\frac12$ is a
feature of the subsequent estimates, and this endpoint is included in the results.

At the excluded boundary $H_0=\frac12$, the temporal cutoff estimate
breaks down; the dyadic summability condition also fails at the corner
$h=(\frac12,\frac12)$. On the critical line $a_h=1$, the
function $f_h^{[\varepsilon]}$ of \eqref{eq:f-eps} behaves as
$A_h v^{-1}$ for large $v$. The stationary counterterm therefore
has a logarithmic divergence at large time separations, so that this
normalization cannot be used there without modification.
\end{remark}

\begin{remark}[Universality]\label{rk:universality}
The discrete model is fixed at the microscopic level by three choices:
the simple symmetric random walk, the Gaussian law of the environment,
and its exact covariance \eqref{fr:discrete-covariance}. 

The Gaussian character is used in an essential way. The exact
factorization \eqref{intro:wick-factorization} and the coupling of
Section~\ref{sec:common-noise}, which realizes the environment on the
Gaussian space of the limiting noise, both rely on it. For a
non-Gaussian environment with the same covariance, the conditional
log-Laplace transform of the path energy is no longer one half of its
conditional variance, and higher cumulants enter. In the white-noise
case, such contributions are negligible in the intermediate-disorder
regime under moment assumptions \cite{AKQ}, and the polynomial chaos
method of \cite{CSZ-JEMS} gives a general framework for replacing one
disorder by another; whether a comparable statement holds with
long-range temporal correlations, where the renormalization diverges,
is not addressed here. 

Changing the walk affects the result more
explicitly. The lattice counterterm \eqref{fr:counterterm-discrete}
depends on the law of the walk through $q_{H_1}(r)$.
The proofs also use invariance principles
specific to the simple random walk, jointly with its local times
\cite{Revesz-1981}. For another centered walk with finite variance, one
expects the same limit after the corresponding changes of the diffusion
constant and of the counterterm, but we have not verified the joint
invariance principles in that generality. 

Finally, replacing $\gamma_{H_0}$ and $\gamma_{H_1}$ by other covariances with the same
large-scale behavior would change the finite part of
$c_{\rm disc}(h)$, and would require the matching of remainders
carried out in Proposition~\ref{prop:counterterm} to be redone.
\end{remark}

\begin{remark}[Further directions]\label{rk:further-directions}
Extending the construction and the discrete approximation to
$H_1<\frac12$ would require different spatial estimates. In that
range, the spatial covariance is a distribution rather than a
nonnegative locally integrable kernel or a Dirac mass, so the
positive-kernel occupation arguments used here do not apply directly.

Theorems~\ref{thm:fractional-main} and~\ref{thm:continuous-summary}
already give moments of every finite order, including moments of
supremum norms on compacts. Questions about their sharp growth as
the order or the time horizon tends to infinity go beyond those
finiteness and convergence statements and are not treated here.
\end{remark}

\section{Numerical illustration of the discrete approximation}
\label{sec:numerics}

This section illustrates Theorem~\ref{thm:fractional-main} on the
lattice model itself, at $H_1=\frac12$: it displays renormalized
terminal profiles for three values of $H_0$, and compares them with
the unrenormalized profiles in the same environment. Each figure shows
a single realization of the environment. The computations are meant as
an illustration of the objects involved, and not as a numerical study
of the convergence: no rate and no statistical conclusion is drawn
from them.

\smallskip

We fix $N=1024$, $H_1=1/2$, $t=1$ and $x\in[-3,3]$, with flat
initial condition. 

\subsection{Sampling and exact walk averaging}
\label{subsec:numerical-method}

For each temporal index $H_0$, we simulate a centered
Gaussian array
$\bar\omega^{N,H_0}=(\bar\omega^{N,H_0}(i,j))_{i\ge1,\,j\in\mathbb Z}$
with covariance
\[
 \E[\bar\omega^{N,H_0}(i,j)\bar\omega^{N,H_0}(i',j')]
   =\gamma_{H_0}(i-i')\one_{\{j=j'\}}.
\]
This is the environment \eqref{fr:discrete-covariance}
at $h=(H_0,\frac12)$, indexed by forward time; it has the law of the
forward array \eqref{eq:forward-cell-array} used in
Section~\ref{subsec:coupling-identification}.

We keep $a_h=2H_0-1/2$ and
$\beta_N=N^{-(H_0-1/4)}$, and suppress $H_0$ in the notation below.

\begin{lemma}[Walk-averaging recursion]\label{lem:numerical-recursion}
For a fixed environment and $c\in\R$, define the path energy and
partition function by
\begin{equation}\label{eq:numerical-path-energy}
 \mathcal H_m^N(j,S):=\beta_N\sum_{k=1}^m
       \bar\omega^N(m+1-k,j+S_k),\qquad
 Z_m^{N,c}(j):=E_S\!\left[e^{\mathcal H_m^N(j,S)-cm\beta_N^2}\right].
\end{equation}
Then, for $m\ge0$ and $j\in\mathbb Z$,
\begin{equation}\label{eq:numerical-recursion}
 \begin{split}
 Z_0^{N,c}(j)&=1,\\
 Z_{m+1}^{N,c}(j)&=\frac12\sum_{\sigma=\pm1}
  e^{\beta_N\bar\omega^N(m+1,j+\sigma)-c\beta_N^2}
  Z_m^{N,c}(j+\sigma).
 \end{split}
\end{equation}
\end{lemma}
\begin{proof}
Condition on $S_1=\sigma$. Since the subsequent increments
are independent of $S_1$, the process $S'_r=S_{r+1}-\sigma$,
$r\ge0$, has the same law as the simple symmetric walk $S$
started at zero. Moreover,
\[
 \mathcal H_{m+1}^N(j,S)
 =\beta_N\bar\omega^N(m+1,j+\sigma)+\mathcal H_m^N(j+\sigma,S').
\]
Each value of $\sigma$ has probability $1/2$, which proves the
recursion. The initial value follows from the empty sum. The
environment is held fixed throughout, so no temporal independence
of the noise is used.
\end{proof}

To obtain the terminal grid $|j|\le J=3\sqrt N=96$, start
from $[-J-N,J+N]\cap\mathbb Z$ and remove one site at each end after
every step. This retains the full dependence cone without an
artificial boundary condition. 

The noise in the recursion is generated jointly in time. Set $L=2N$
and form the positive definite matrices
\[
 C_{H_0}(r,s)=\gamma_{H_0}\!\left(\min\{|r-s|,L-|r-s|\}\right),
 \qquad 0\le r,s<L.
\]

The choice $L=2N$ preserves all time differences in the first $N$
coordinates: for $0\le r,s<N$,
\[
 |r-s|\le N-1<L/2,\qquad
 C_{H_0}(r,s)=\gamma_{H_0}(r-s).
\]
With $L=N$, the circulant construction would instead identify large
time differences with small ones; for example, its $(0,N-1)$ entry
would be $\gamma_{H_0}(1)$ rather than $\gamma_{H_0}(N-1)$.
One could also sample directly from the $N\times N$ covariance
matrix $(\gamma_{H_0}(r-s))_{0\le r,s<N}$. The enlargement is used
here to obtain a circulant matrix, whose square-root multiplication
can be evaluated by fast Fourier transforms \cite{Dietrich-Newsam}.

Draw independent vectors $G^{(j)}\sim\mathcal N(0,I_L)$, one per
spatial site, and put
\begin{equation}\label{eq:numerical-noise-sampling}
 \bar\omega^{N,H_0}(i,j)=\bigl(C_{H_0}^{1/2}G^{(j)}\bigr)_{i-1},
 \qquad 1\le i\le N.
\end{equation}
Indeed, $\operatorname{Cov}(C_{H_0}^{1/2}G^{(j)})=C_{H_0}$ and,
since $|i-i'|<L/2$, its first $N$ coordinates have exactly the
required covariance. The same $G^{(j)}$ are used for every $H_0$;
each index still has its own covariance matrix.

\smallskip

\subsection{Renormalized profiles}
\label{subsec:numerical-counterterm}

For each of $H_0=11/16,5/8,9/16$, we use the full stationary
counterterm \eqref{fr:counterterm-discrete}:
\begin{equation}\label{fr:counterterm-discrete-bis}
 c_{\rm disc}(H_0,1/2)=\frac12+
       \sum_{k=1}^{\infty}\gamma_{H_0}(2k)2^{-2k}\binom{2k}{k}.
\end{equation}
All three series converge. The terminal interpolation with
$c=c_{\rm disc}(h)$ in \eqref{eq:numerical-recursion} is denoted by
$F_N^{\rm ren}(1,\cdot)$. It is the discrete approximation of $U_1$
specified in Theorem~\ref{thm:coupled-L2}.

We evaluate the full series through its Laplace integral
\begin{equation}\label{eq:numerical-counterterm-integral}
 c_{\rm disc}(H_0,1/2)=\frac12+
 \frac{H_0(2H_0-1)}{\Gamma(2-2H_0)}
 \int_0^\infty
 \frac{e^{-u}(1-e^{-u})^2}
 {u^{1+2H_0}\sqrt{1-e^{-2u}}(1+\sqrt{1-e^{-2u}})}\dd u.
\end{equation}

This follows from the Laplace representation of $\gamma_{H_0}$ and
$\sum_{k\ge0}p_{2k}(0)e^{-2ku}=(1-e^{-2u})^{-1/2}$.
The integrand is of order $u^{-a_h}$ at zero; the substitution
$u=v^{1/(1-a_h)}$ removes this singularity for quadrature.

The reason for this choice is the slow decay of the omitted tail
under direct truncation. Set
\[
 c^{(M)}(h):=\frac12+
       \sum_{k=1}^{M}\gamma_{H_0}(2k)2^{-2k}\binom{2k}{k}.
\]
By \eqref{eq:LLT} and
$\gamma_{H_0}(2k)\sim H_0(2H_0-1)(2k)^{2H_0-2}$,
\begin{equation}\label{eq:numerical-counterterm-tail}
 c_{\rm disc}(h)-c^{(M)}(h)
 \sim\frac{H_0(2H_0-1)2^{2H_0-2}}
            {\sqrt\pi(1-a_h)}\,M^{a_h-1}
 \qquad(M\to\infty).
\end{equation}
For example, at $H_0=11/16$ the tail decays only as $M^{-1/8}$.
Even $M=10^6$ leaves an error of approximately $0.13418$ in the
counterterm. At $N=1024$, using this truncated value in
\eqref{eq:numerical-recursion} multiplies the terminal profile by
\[
 \exp\!\left\{N^{1-a_h}
              \bigl(c_{\rm disc}(h)-c^{(M)}(h)\bigr)\right\}
 \simeq1.376.
\]

The integral representation includes this entire time tail.

\begin{figure}[htbp]
\centering
 \includegraphics[width=\textwidth]{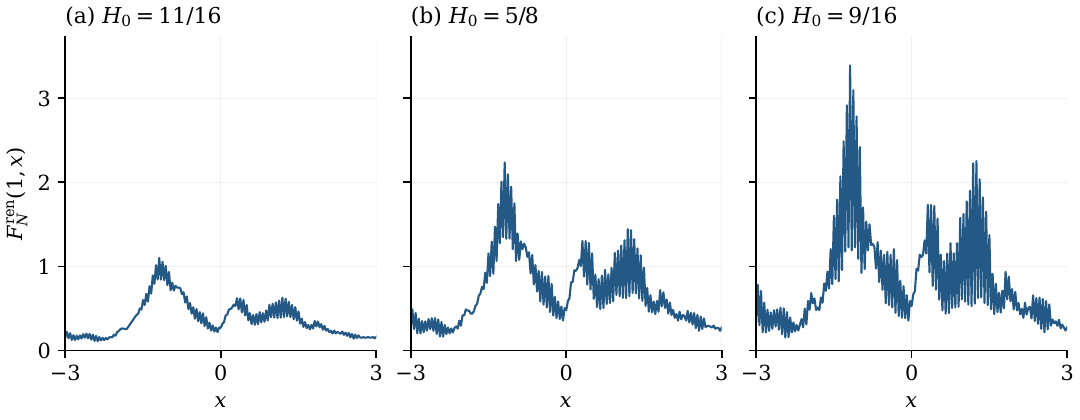}
 \caption{Renormalized profiles for $H_0=11/16,5/8,9/16$, with
 $N=1024$, $H_1=1/2$ and $t=1$, on a common vertical scale.}
 \label{fig:numerical-profiles}
\end{figure}

\subsection{A realization (Figure \ref{fig:numerical-profiles})}

The profiles remain positive by Lemma~\ref{lem:numerical-recursion}.
In this realization, oscillations between spatial parity classes
become more visible as $H_0$ decreases. 

These oscillations are a lattice effect: at step $k$, a walk started at $j$ occupies a site of
the parity of $j+k$, so that the two parity classes of terminal sites
read the environment on two disjoint sublattices.
We nevertheless keep both classes, because the interpolation
$\operatorname{Interp}_N$ of Theorem~\ref{thm:fractional-main} uses all
grid points: a plot restricted to one class would be smoother, but it
would no longer display the field to which the theorem applies.
Since $\operatorname{Interp}_N\mathsf Z^N$ converges in law in
$C_{\rm loc}([0,\infty)\times\R)$ to a continuous field
(Theorem~\ref{thm:fractional-main}), the largest difference between
neighboring grid values over a compact set tends to zero in
probability, and these oscillations disappear in the limit.

The noisier appearance of the profiles at smaller $H_0$
is consistent with the slower convergence expected as
$H_0\downarrow\frac12$.

\smallskip

\subsection{The effect of renormalization in one environment}
\label{subsec:numerical-effect}

For each of the three indices $H_0=11/16,5/8,9/16$, the two
recursions in the same environment satisfy
\begin{equation}\label{eq:numerical-raw-factor}
 F_N^{\rm raw}(1,x)=e^{d_N(H_0)}F_N^{\rm ren}(1,x),
 \qquad d_N(H_0):=c_{\rm disc}(h)N^{1-a_h}.
\end{equation}
Indeed, $Z_m^{N,0}=e^{cm\beta_N^2}Z_m^{N,c}$, and linear
interpolation preserves this identity at the terminal time.
Renormalization therefore changes the vertical scale, while
preserving the spatial shape of each profile.

\begin{figure}[htbp]
\centering
 \includegraphics[width=\textwidth]{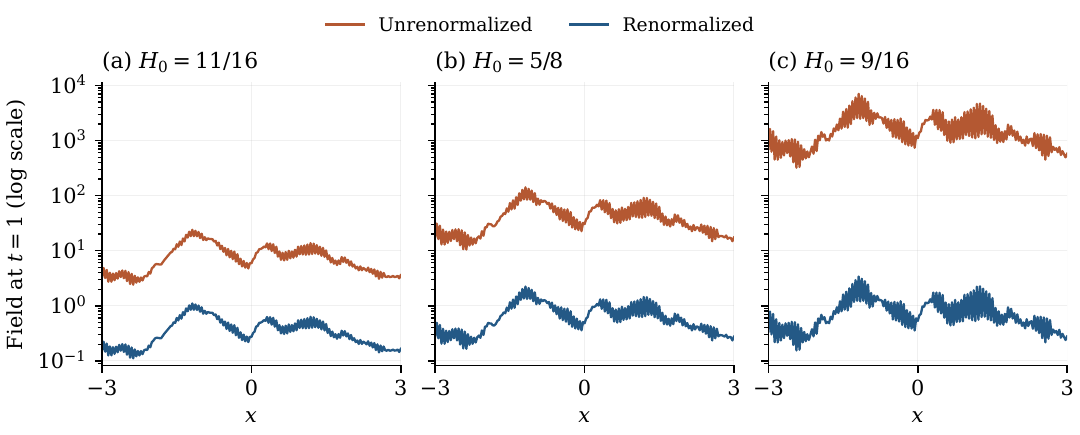}
 \caption{Unrenormalized and renormalized terminal profiles in
 the same environment, for $H_0=11/16,5/8,9/16$. The common
 logarithmic vertical scale keeps both curves visible. Their
 vertical separation is the deterministic quantity
 $d_N(H_0)/\log 10$ on the base-ten logarithmic scale.
 The blue profiles are those of Figure~\ref{fig:numerical-profiles}.}
 \label{fig:numerical-comparison}
\end{figure}

The numerical counterterms and renormalization factors are
\[
 \begin{array}{c|c|c|c}
 H_0 & c_{\rm disc}(h) & d_{1024}(H_0) & e^{d_{1024}(H_0)}\\ \hline
 11/16 & 1.29514 & 3.08038 & 21.7666\\
 5/8   & 0.733579 & 4.14975 & 63.4180\\
 9/16  & 0.568068 & 7.64299 & 2085.96
 \end{array}
\]

In particular, the factor removed
by renormalization exceeds $2000$ at $H_0=9/16$, despite a smaller
value of $c_{\rm disc}(h)$ than at the other two indices: the
exponent $1-a_h$ in \eqref{eq:numerical-raw-factor} is larger.


\clearpage
\part{Proofs}\label{chap:proofs}

\section{\texorpdfstring{Brownian self-intersections}{Brownian self-intersections}}
\label{sec:Brownian-construction}

The aim of this section is to construct the centered self-intersection
$Y_t(B)$ of Proposition~\ref{prop:Y} and to identify it as the limit of its
Gaussian regularizations. The main result is
Theorem~\ref{thm:Brownian-self-intersection}, which gives convergence in
every finite $L^p$, a quantitative bound on the dyadic remainder, and
uniform exponential estimates.

The construction uses the dyadic centering method of Varadhan and
Le Gall \cite{Varadhan,LeGall-StFlour}, applied to regularized
Brownian self-intersections with fractional covariance weights.

\

Throughout, we fix $h=(H_0,H_1)$ satisfying \eqref{fr:range}, and write
$a_h=2H_0+H_1-1\in(\frac12,1)$. For Brownian functionals, expectations
and $L^p$ norms refer to the Brownian probability space. Constants denoted by
$C_h$ may change from one occurrence to the next and depend only on
the fixed Hurst parameters, unless indicated otherwise.

Throughout the article, we fix the nondecreasing continuous cutoff
\begin{equation}\label{eq:cutoff}
 \chi(u)=\min\{1,\max\{0,u-1\}\},\qquad u\ge0,
\end{equation}
which vanishes on $[0,1]$ and equals one on $[2,\infty)$. It is used to
remove neighborhoods of the time diagonal.

For a standard Brownian motion $B$ started at zero, $t>0$ and
$\varepsilon>0$, let
\begin{equation}\label{eq:Gamma-eps}
 \Gamma_t^\varepsilon(B)
 =\int_{[0,t]^2}K_{H_0}^{[2\varepsilon^4]}(\tau-\sigma)
          K_{H_1}^{[2\varepsilon^2]}(B_\sigma-B_\tau)\dd\sigma\dd\tau,
\end{equation}
where $K_{H_0}^{[2\varepsilon^4]}$ and $K_{H_1}^{[2\varepsilon^2]}$ are the temporal
and spatial covariance factors of \eqref{eq:smoothed-noise-covariance}.
For fixed $\varepsilon>0$ these factors are bounded, so
$\Gamma_t^\varepsilon$ is integrable. Define
\begin{equation}\label{eq:Y-eps}
 Y_t^\varepsilon(B)
 =\Gamma_t^\varepsilon(B)-E_B[\Gamma_t^\varepsilon(B)].
\end{equation}
In the Feynman-Kac representation, $\Gamma_t^\varepsilon$ is the
conditional variance of the Gaussian noise integrated along the path.
Its centered fluctuation $Y_t^\varepsilon(B)$ remains after deterministic
renormalization; this connection is developed in
Section~\ref{sec:continuous-identification}.

One could in principle prove convergence of $Y_t^\varepsilon(B)$ directly,
but applying dominated convergence to \eqref{eq:Gamma-eps} is obstructed
by the singularity along the time diagonal. Indeed, the mean
unregularized interaction has density
$A_h|\tau-\sigma|^{a_h-2}$ off the diagonal, by
Lemma~\ref{lem:Brownian-mean-measure}; this density is not integrable
near the diagonal because $a_h<1$. Bounding the two terms in
\eqref{eq:Y-eps} separately therefore loses the cancellation due to
centering. The dyadic decomposition (leading to \eqref{eq:block-continuous}) isolates finite interactions
between adjacent time intervals. Centering these blocks before
summing, and using their independence at each level, yields both a
summable remainder and the uniform exponential bounds needed for the
Feynman-Kac representation.

The construction rests on one estimate. For $p\ge2$, a dyadic block $A_{M,k,t}$ (see \eqref{eq:block-continuous}) at
level $M$ has time scale $t2^{-M}$ and $L^p$ size at most
$C_h\sqrt p\,(t2^{-M})^{a_h}$. The $2^{M-1}$ centered blocks at that
level are independent. Their sub-Gaussian moment bounds therefore give
\[
 \left\|\sum_{k=1}^{2^{M-1}}A_{M,k,t}\right\|_p
 \le C_h\sqrt p\,2^{M/2}(t2^{-M})^{a_h}
 =C_h\sqrt p\,t^{a_h}2^{-M(a_h-\frac12)}.
\]
The last factor is summable precisely when $a_h>\frac12$. The same
estimate holds uniformly for the regularized blocks, so convergence
of each fixed block identifies the limit of the full series.

\

We first establish the occupation estimates and spatial approximation
results used by the blocks. We then prove their bounds and convergence,
and sum the centered levels.

\subsection{\texorpdfstring{Local times and the blocks at $H_1=\frac12$}{Local times and the blocks at H1=1/2}}
\label{subsec:Brownian-measures}

For $H_1>\frac12$, the blocks \eqref{fr:Brownian-block} are ordinary
integrals of a locally integrable kernel along the path. At
$H_1=\frac12$, where $K_{1/2}=\delta_0$, they have to be read through
the local times of $B$. The purpose of this subsection is to make this
reading precise. Along the way, it fixes the local-time notation used
in the occupation estimates of both regimes.

For a standard Brownian motion $B$ started at zero, choose its jointly
continuous occupation-density local time $L^B$, normalized by
\begin{equation}\label{temps-local}
 \int_0^t f(B_r)\dd r=\int_\R f(x)L_t^B(x)\dd x.
\end{equation}
On every finite interval $[0,T]$, let $\mu_x^B$ be the atomless
Stieltjes measure satisfying
\[
 \mu_x^B([0,t])=L_t^B(x),\qquad
 \mu_x^B([u,v))=L_v^B(x)-L_u^B(x).
\]
For two Brownian paths $B,\widetilde B$, define
\begin{equation}\label{eq:two-BM-occupation-measure}
 \mathfrak I_{B,\widetilde B}(A)
 :=\int_\R(\mu_x^B\otimes\mu_x^{\widetilde B})(A)\dd x,
 \qquad A\in\mathcal B([0,T]^2).
\end{equation}
Joint continuity of local times ensures measurability of the integrand.
Writing $L_I^B=L_v^B-L_u^B$ for $I=[u,v)$, one has
\begin{align*}
 \mathfrak I_{B,\widetilde B}(I\times J)
 &=\int_\R L_I^B(x)L_J^{\widetilde B}(x)\dd x,\\
 \mathfrak I_{B,\widetilde B}([0,T]^2)
 &\le T\|L_T^{\widetilde B}\|_\infty<\infty
 \quad\text{almost surely}.
\end{align*}
The last inequality follows from the occupation formula and the compact
spatial support of local times. For a single path, its self-intersection
measure is the special case
\begin{equation}\label{eq:BM-intersection-measure}
 \mathfrak I_B:=\mathfrak I_{B,B}.
\end{equation}
All these measures assign zero mass to horizontal and vertical lines.
They also assign zero mass to the time diagonal
$D=\{(s,t):s=t\}$: by Tonelli's theorem and atomlessness,
\[
 \mathfrak I_{B,\widetilde B}(D)
 =\int_\R\int_{[0,T]}
       \mu_x^{\widetilde B}(\{s\})\,\mu_x^B(\dd s)\dd x=0.
\]
The measure $\mathfrak I_{B,\widetilde B}$ plays, for a pair of paths,
the role played by the local time in \eqref{temps-local}. Indeed,
applying \eqref{temps-local} successively to the increments of each
path gives, for all nonnegative Borel functions $f$ on $[0,T]^2$ and
$g$ on $\R^2$,
\begin{equation}\label{eq:two-path-occupation}
 \int_{[0,T]^2}f(r,s)\,g(B_r,\widetilde B_s)\dd r\dd s
 =\int_{\R^2}g(x,y)
       \left(\int_{[0,T]^2}f\dd(\mu_x^B\otimes\mu_y^{\widetilde B})\right)
       \dd x\dd y.
\end{equation}
For bounded $g$, both sides define finite measures in $f$ that
coincide on rectangles $I\times J$; the general case follows by
monotone convergence. Taking formally $g(x,y)=\delta_0(x-y)$ concentrates the
right-hand side on the diagonal $x=y$ and produces
$\int_{[0,T]^2}f\dd\mathfrak I_{B,\widetilde B}$, by
\eqref{eq:two-BM-occupation-measure}. This is the meaning given below
to integrals of $\delta_0(B_r-\widetilde B_s)$, and
Lemma~\ref{lem:Brownian-spatial-regularization} shows that it is also
the limit obtained when $\delta_0$ is approximated by the heat kernels
$\heat_v$, $v\downarrow0$.

At $H_1=\frac12$, where
$K_{1/2}=\delta_0$
and $E_B[K_{1/2}(B_s-B_r)]=\heat_{s-r}(0)=(2\pi(s-r))^{-\frac12}$, the
definitions \eqref{fr:Brownian-block} and \eqref{eq:block-continuous}
now read
\begin{equation}\label{eq:block-endpoint}
\begin{split}
 C_{M,k,t}(B)&=\int_{I_{M,k,t}\times J_{M,k,t}}
       K_{H_0}(s-r)\,\mathfrak I_B(\dd r\dd s),\\
 E_B[C_{M,k,t}]&=A_{H_0}\int_{I_{M,k,t}\times J_{M,k,t}}
       (s-r)^{2H_0-\frac52}\dd r\dd s,
\end{split}
\end{equation}
where $\mathfrak I_B$ is the self-intersection measure defined in
\eqref{eq:BM-intersection-measure} and $A_{H_0}=c_{H_0}/\sqrt{2\pi}$,
in accordance with \eqref{eq:AH} and $a_{h}=2H_0-\frac12$.

\subsection{Occupation estimates and spatial regularization}
\label{subsec:Brownian-weighted-intersections}

Throughout this subsection the spatial kernel is $K_{H_1}$ and its
regularizations $K_{H_1}^{[v]}$, as in \eqref{fr:kernels}, with
$K_{H_1}^{[0]}=K_{H_1}$.

At $H_1=\frac12$, $K_{1/2}^{[v]}=\heat_v$ is a genuine function for
$v>0$, whereas $K_{1/2}^{[0]}=\delta_0$ is not; every integral against
it is read, once and for all, through the occupation densities of
Section~\ref{subsec:Brownian-measures}:
\begin{equation}\label{eq:endpoint-reading}
 \int_0^TK_{1/2}^{[0]}(B_r-z)\dd r:=L_T^B(z),
 \qquad
 \int_{[0,T]^2}f(r,s)K_{1/2}^{[0]}(B_r-\widetilde B_s)\dd r\dd s
 :=\int_{[0,T]^2}f\dd\mathfrak I_{B,\widetilde B}.
\end{equation}
The second definition is the reading of \eqref{eq:two-path-occupation}
at $g(x,y)=\delta_0(x-y)$ described in
Section~\ref{subsec:Brownian-measures}.
For $H_1>\frac12$, the corresponding occupation formulas, namely
\eqref{temps-local} and \eqref{eq:two-path-occupation} with
$g(x,y)=K_{H_1}(x-y)$, retain integration against the spatial kernel:
\begin{align*}
 \int_0^T K_{H_1}(B_r-z)\dd r
 &=\int_\R K_{H_1}(x-z)L_T^B(x)\dd x,\\
 \int_{[0,T]^2}f(r,s)K_{H_1}(B_r-\widetilde B_s)\dd r\dd s
 &=\int_{\R^2}K_{H_1}(x-y)
       \left(\int_{[0,T]^2}f\dd(\mu_x^B\otimes\mu_y^{\widetilde B})\right)
       \dd x\dd y.
\end{align*}
 We begin with the two properties of these kernels that every spatial estimate below uses.

\begin{lemma}[Bounds and translations of the spatial kernel]
\label{lem:spatial-kernel}
Uniformly in $v\ge0$ for $H_1>\frac12$ and in
$v>0$ for $H_1=\frac12$, for $x\ne0$, $R>0$ and $z\in\R$,
\begin{align}
 K_{H_1}^{[v]}(x)&\le C_{H_1}(|x|+\sqrt v)^{2H_1-2}\le C_{H_1}|x|^{2H_1-2},
                                             \label{eq:kernel-envelope}\\
 \int_{|x|\le R}K_{H_1}^{[v]}(x)\dd x&\le C_{H_1}R^{2H_1-1},
                                             \label{eq:kernel-mass}\\
 \|K_{H_1}^{[v]}(\cdot+z)-K_{H_1}^{[v]}\|_{L^1(\R)}
 &\le C_{H_1}|z|^{2H_1-1}.                      \label{eq:kernel-translation}
\end{align}
Moreover, for every $u>0$,
\begin{equation}\label{eq:kernel-at-zero}
 K_{H_1}^{[u]}(0)=\kappa_{H_1}u^{H_1-1},
\end{equation}
with $\kappa_{H_1}$ as in \eqref{eq:AH} for $H_1>\frac12$ and $\kappa_{1/2}=(2\pi)^{-1/2}$. For $H_1>\frac12$ one has
$K_{H_1}^{[v]}(x)\to K_{H_1}(x)$ for every $x\ne0$ as $v\downarrow0$;
at $H_1=\frac12$, the kernel $K_{1/2}^{[v]}=\heat_v$ converges weakly to
$K_{1/2}^{[0]}=\delta_0$.
\end{lemma}

\begin{proof}

For $H_1>\frac12$, all assertions are the case $H=H_1$ of
Lemma~\ref{lem:kernel-regularization} in Appendix~\ref{app:temporal},
the pointwise convergence even holding uniformly on $\{|x|\ge\eta\}$
for every $\eta>0$.

At $H_1=\frac12$ all four assertions are immediate for $v>0$:
$\heat_v(x)\le C(|x|+\sqrt v)^{-1}$ because
$\heat_1(u)\le C(1+|u|)^{-1}$, the mass bound holds because $\heat_v$
is a probability density, the translation bound holds because
$\|\heat_v(\cdot+z)-\heat_v\|_1\le2$, and \eqref{eq:kernel-at-zero} is
the value of the Gaussian density at the origin. The weak convergence
$\heat_v\to\delta_0$ is standard.
\end{proof}

\begin{lemma}[Brownian local-time bound]\label{lem:Brownian-local-time}
For a Brownian motion $B$ with variance parameter $1$, let
$L_t^B(x)$ be its occupation-density local time, normalized as in
Section~\ref{subsec:Brownian-measures}. There is a universal constant $C$ such that, uniformly in the
starting point, for every $p\ge2$ and $t\ge0$,
\begin{equation}\label{eq:Brownian-local-time-Lp}
 \sup_{x\in\R}\|L_t^B(x)\|_p\le C\sqrt{pt}.
\end{equation}
The same bound holds, with a modified constant, for a Brownian motion
of variance parameter $2$.
\end{lemma}

\begin{proof}
Let $E_z,\Pp_z$ denote expectation and probability for Brownian
motion started at $z$. At $z=0$, L\'evy's identity and the
reflection principle give
\[
 L_t^B(0)\overset{\rm law}=\max_{0\le s\le t}B_s
 \overset{\rm law}=|B_t|;
\]
see \cite[Theorems~6.10 and~2.21]{MP-Brownian}.
Given $p\ge2$, take $j=\lceil\frac p2\rceil$:
\[
 \|L_t^B(0)\|_p=\sqrt t\,\|B_1\|_p
 \le\sqrt t\,\|B_1\|_{2j}\le\sqrt{2jt}\le C\sqrt{pt}.
\]
For an arbitrary starting point $z$ and site $x$, let $T_x^B$
be the first hitting time of $x$. The strong Markov property \cite[Theorem~2.16]{MP-Brownian} and additivity of local time give
\[
 E_z[L_t^B(x)^p]
 =\int_{[0,t]}E_0[L_{t-u}^B(0)^p]\,\Pp_z(T_x^B\in\dd u)
 \le E_0[L_t^B(0)^p].
\]
Here local time at $x$ is zero up to its first visit, and the
inequality uses monotonicity in the remaining time. This proves
the uniform bound. For variance parameter $2$, write
$D_t=z+\sqrt2B_t$, where $B$ is standard and starts at zero.
The occupation-density normalization gives
\[
 L_t^D(x)=\frac1{\sqrt2}
                 L_t^B\left(\frac{x-z}{\sqrt2}\right),
\]
which proves the same bound uniformly in $z$ and $x$.
\end{proof}

The next estimate is the one that replaces, in the fractional regime,
the integrability of the spatial kernel.

\begin{lemma}[Occupation against the spatial kernel]
\label{lem:Brownian-occupation}
Let $B$ be a Brownian motion of variance
parameter $1$ or $2$, started at an arbitrary point. Uniformly in
$v\ge0$, in $z,z'\in\R$ and in the starting point, for every $p\ge2$
and $T>0$,
\begin{align}
 \left\|\int_0^TK_{H_1}^{[v]}(B_r-z)\dd r\right\|_p
 &\le C_{h}\,p^{1-H_1}T^{H_1},\label{eq:BM-occupation}\\
 \left\|\int_0^T\bigl[K_{H_1}^{[v]}(B_r-z)-K_{H_1}^{[v]}(B_r-z')\bigr]\dd r\right\|_p
 &\le C_{H_1}\sqrt p\,T^{\frac12}|z-z'|^{2H_1-1}.
                                    \label{eq:BM-occupation-increment}
\end{align}
Both bounds also hold conditionally on the position of $B$ at the
beginning of the time segment.
\end{lemma}

\begin{proof}
Consider first the case $H_1>\frac12$. By the occupation formula,
\[
 \int_0^TK_{H_1}^{[v]}(B_r-z)\dd r
 =\int_\R K_{H_1}^{[v]}(x-z)L_T^B(x)\dd x.
\]
Fix $R>0$ and split the spatial integral. On $|x-z|\le R$, Minkowski's
inequality, Lemma~\ref{lem:Brownian-local-time} and
\eqref{eq:kernel-mass} give an $L^p$ norm at most
\[
 \int_{|x-z|\le R}K_{H_1}^{[v]}(x-z)\|L^B_T(x)\|_p\dd x
 \le C\sqrt{pT}\int_{|y|\le R}K_{H_1}^{[v]}(y)\dd y
 \le C_h\sqrt{pT}\,R^{2H_1-1}.
\]
On $|x-z|>R$, the kernel is at most $C_{H_1}R^{2H_1-2}$ by
\eqref{eq:kernel-envelope}, and $\int_\R L^B_T(x)\dd x=T$, so that this
part is bounded pathwise by $C_{H_1}TR^{2H_1-2}$. Choosing
$R=\sqrt{\frac Tp}$ makes both terms equal to
$C_hp^{1-H_1}T^{H_1}$, which proves \eqref{eq:BM-occupation}.

For the increment, the same occupation formula, Minkowski's inequality
and \eqref{eq:kernel-translation} give
\begin{align*}
 \left\|\int_\R\bigl[K_{H_1}^{[v]}(x-z)-K_{H_1}^{[v]}(x-z')\bigr]L^B_T(x)\dd x\right\|_p&
 \le C\sqrt{pT}\,\|K_{H_1}^{[v]}(\cdot+z-z')-K_{H_1}^{[v]}\|_{L^1}\\
& \le C_{H_1}\sqrt{pT}\,|z-z'|^{2H_1-1},
\end{align*}
which is \eqref{eq:BM-occupation-increment}. Both estimates are
uniform in the starting point of $B$, by
Lemma~\ref{lem:Brownian-local-time}, hence hold conditionally on the
position at the beginning of the segment.

At $H_1=\frac12$ the argument is shorter. For $v=0$ the left-hand side
of \eqref{eq:BM-occupation} is $\|L_T^B(z)\|_p$, which
Lemma~\ref{lem:Brownian-local-time} bounds by $C\sqrt{pT}$, the value
of the right-hand side. For $v>0$,
$\int_0^T\heat_v(B_r-z)\dd r=\int_\R\heat_v(y)L_T^B(z+y)\dd y$, so
Minkowski's inequality and the same uniform bound give the same
estimate. Both are uniform in $z$ and in the starting point, hence
hold conditionally. Finally, \eqref{eq:BM-occupation-increment} follows from \eqref{eq:BM-occupation} at $z$ and at $z'$ and the triangle inequality, since $|z-z'|^{2H_1-1}=1$ for $z\ne z'$.
\end{proof}

Now, for two independent Brownian motions $B,\widetilde B$ started from
arbitrary real points, define, for $T>0$ and $v\ge0$,
\begin{equation}\label{eq:Brownian-weighted-functionals}
 \begin{split}
 \mathcal T_T^v(B,\widetilde B)
 &=\int_{[0,T]^2}K_{H_1}^{[v]}(B_r-\widetilde B_s)\dd r\dd s,\\
 \mathcal Q_T^v(B,\widetilde B)
 &=\int_{[0,T]^2}K_{H_0}(r+s)
                  K_{H_1}^{[v]}(B_r-\widetilde B_s)\dd r\dd s.
 \end{split}
\end{equation}
We omit the paths from the notation when they are fixed.

\begin{lemma}[Weighted Brownian intersections]
\label{lem:Brownian-weighted-block}
For two independent Brownian motions $B,\widetilde B$ started from
arbitrary real points, the following estimates hold uniformly in
$v\ge0$, for every $T>0$ and $p\ge2$:
\begin{align}
 \|\mathcal T_T^v\|_p
 &\le C_{h}\,p^{1-H_1}T^{1+H_1},\label{eq:BM-unweighted-block}\\
 \|\mathcal Q_T^v\|_p
 &\le C_{h}\,p^{1-H_1}T^{a_{h}}.\label{eq:BM-weighted-block}
\end{align}
In particular, both functionals are almost surely finite, including at
$v=0$.
\end{lemma}

\begin{proof}
Put $F_s=\int_0^TK_{H_1}^{[v]}(B_r-\widetilde B_s)\dd r$ and let $\nu_s$ be
the distribution of $\widetilde B_s$. The two paths are independent, so
conditioning on $\widetilde B_s=y$ and applying
\eqref{eq:BM-occupation} to $B$ with $z=y$ gives
\[
 E[|F_s|^p]
 =\int_\R E_B\left[
       \left|\int_0^TK_{H_1}^{[v]}(B_r-y)\dd r\right|^p
                    \right]\nu_s(\dd y)
 \le\bigl(C_hp^{1-H_1}T^{H_1}\bigr)^p,
\]
the bound in \eqref{eq:BM-occupation} being uniform in $y$.
Minkowski's integral inequality therefore gives
\[
 \left\|\int_0^T F_s\dd s\right\|_p
 \le\int_0^T\|F_s\|_p\dd s
 \le C_hp^{1-H_1}T^{1+H_1},
\]
which is \eqref{eq:BM-unweighted-block}.

\smallskip

For the weighted estimate, set $\delta_j=2^{-j}T$ and decompose the
square into
\[
 \{(r,s)\in[0,T]^2:T<r+s\le2T\}
 \quad\text{and}\quad
 \{(r,s)\in[0,T]^2:\tfrac12\delta_j<r+s\le\delta_j\},
 \quad j\ge0.
\]
On the $j$th annulus, $r,s\le\delta_j$ and
$K_{H_0}(r+s)\le c_{H_0}2^{2-2H_0}\delta_j^{2H_0-2}$.
Positivity of the kernel therefore bounds its contribution by
\[
 c_{H_0}2^{2-2H_0}\delta_j^{2H_0-2}
       \int_{[0,\delta_j]^2}K_{H_1}^{[v]}(B_r-\widetilde B_s)\dd r\dd s,
\]
whose $L^p$ norm is at most
$C_hp^{1-H_1}\delta_j^{2H_0-2}\delta_j^{1+H_1}
 =C_hp^{1-H_1}\delta_j^{a_{h}}$
by \eqref{eq:BM-unweighted-block} and $a_{h}=2H_0+H_1-1$. Minkowski's
inequality and the convergence of $\sum_{j\ge0}2^{-ja_h}$, which holds
because $a_{h}>0$, prove the desired bound on $\{0<r+s\le T\}$. On the
outer region $\{T<r+s\le2T\}$ the weight is at most $c_{H_0}T^{2H_0-2}$,
and a further application of \eqref{eq:BM-unweighted-block} concludes
the proof.
\end{proof}

\begin{lemma}[Removal of the spatial regularization]
\label{lem:Brownian-spatial-regularization}
Let $B,\widetilde B$ be independent Brownian motions started at arbitrary
real points, and let $T>0$. For every finite $p\ge2$,
\begin{equation}\label{eq:BM-spatial-regularization-convergence}
 \|\mathcal T_T^v-\mathcal T_T^0\|_p
 +\|\mathcal Q_T^v-\mathcal Q_T^0\|_p
 \longrightarrow0\qquad(v\downarrow0).
\end{equation}
More generally, for every bounded weight $f$ on $[0,T]^2$ which is
continuous off the origin,
\begin{equation}\label{eq:BM-bounded-weight-convergence}
 \int_{[0,T]^2}f(r,s)K_{H_1}^{[v]}(B_r-\widetilde B_s)\dd r\dd s
 \longrightarrow
 \int_{[0,T]^2}f(r,s)K_{H_1}(B_r-\widetilde B_s)\dd r\dd s
 \quad\text{in }L^p.
\end{equation}
At $H_1=\frac12$, the limiting integrals are read through
$\mathfrak I_{B,\widetilde B}$, as in \eqref{eq:endpoint-reading}.
\end{lemma}

\begin{proof}
Suppose first that $H_1>\frac12$.
For $(r,s)\ne(0,0)$ the variable $B_r-\widetilde B_s$ is Gaussian with
positive variance $r+s$, hence almost surely nonzero; the exceptional
set is Lebesgue-null in $[0,T]^2$. By
Lemma~\ref{lem:spatial-kernel}, $K_{H_1}^{[v]}\to K_{H_1}$ at every such
point, while $K_{H_1}^{[v]}\le C_{H_1}K_{H_1}/c_{H_1}$ everywhere by
\eqref{eq:kernel-envelope}. The majorant
$\int_{[0,T]^2}|f|\,K_{H_1}(B_r-\widetilde B_s)\dd r\dd s$ is integrable by
\eqref{eq:BM-unweighted-block}, so dominated convergence gives
\eqref{eq:BM-bounded-weight-convergence} almost surely and in $L^1$.
Since \eqref{eq:BM-unweighted-block} bounds the family uniformly in
$v$ at the order $2p$, the convergence is uniformly integrable at the
order $p$ and therefore holds in $L^p$. Taking $f\equiv1$ gives the
first half of \eqref{eq:BM-spatial-regularization-convergence}.

At $H_1=\frac12$ that domination is unavailable, $\heat_v$ converging
to zero at every point off the origin, and we argue instead with the
occupation densities themselves. Applying the occupation formula
successively in $r$ and in $s$ to the bounded integrand gives, for
$v>0$,
\[
 \int_{[0,T]^2}f(r,s)\heat_v(B_r-\widetilde B_s)\dd r\dd s
 =\int_{\R^2}\heat_v(x-y)\Phi(x,y)\dd x\dd y,
 \qquad
 \Phi(x,y)=\int_{[0,T]^2}f\dd(\mu_x^B\otimes\mu_y^{\widetilde B}),
\]
whereas the limit in \eqref{eq:BM-bounded-weight-convergence}, read
through \eqref{eq:endpoint-reading}, is $\int_\R\Phi(x,x)\dd x$ by the
definition \eqref{eq:two-BM-occupation-measure}. The measures $\mu_x^B$
depend weakly continuously on $x$: their distribution functions are
$t\mapsto L^B_t(x)$, which converge at every $t$ as $x$ varies, by
joint continuity of local times, and so do their total masses. Since
$f$ is bounded and continuous off the origin, and the limiting product
measure is atomless, $\Phi$ is therefore continuous; it is bounded by
$\|f\|_\infty L^B_T(x)L^{\widetilde B}_T(y)$ and supported in the
compact rectangle
\[
 B([0,T])\times\widetilde B([0,T])
 =
 \bigl[\min_{0\le r\le T}B_r,\max_{0\le r\le T}B_r\bigr]
 \times
 \bigl[\min_{0\le r\le T}\widetilde B_r,
       \max_{0\le r\le T}\widetilde B_r\bigr].
\]

 Writing
\[
 \int_{\R^2}\heat_v(x-y)\Phi(x,y)\dd x\dd y-\int_\R\Phi(x,x)\dd x
 =\int_\R\heat_v(z)\left(\int_\R
       \bigl[\Phi(x,x-z)-\Phi(x,x)\bigr]\dd x\right)\dd z,
\]
uniform continuity and compact support make the inner integral a
bounded function $F(z)$ with $F(z)\to0$ as $z\to0$. The bound is
uniform in $z$ because its $x$-support is contained in the compact set
$B([0,T])$. After the change of variables $z=\sqrt v\,u$,
the last expression becomes
$\int_\R\heat_1(u)F(\sqrt v\,u)\dd u$; dominated convergence now
gives almost sure convergence as $v\downarrow0$. The uniform bound
\eqref{eq:BM-unweighted-block}, at the order $2p$, upgrades it to $L^p$ as before. The rest of the proof
uses only positivity and \eqref{eq:BM-weighted-block}, and is therefore
common to the two regimes $H_1>\frac12$ and $H_1=\frac12$.

For the singular weight, recall the cutoff $\chi$ of \eqref{eq:cutoff}
and put
\[
 w_\eta(r,s)=K_{H_0}(r+s)\chi\left(\frac{r+s}{\eta}\right),
 \qquad w_\eta(0,0)=0,
\]
which is bounded for each $\eta>0$. Write
$\mathcal Q_T^{v,\eta}$ for the corresponding integral. Then
\eqref{eq:BM-bounded-weight-convergence} gives
\begin{equation}\label{eq:BM-weighted-cutoff-convergence}
 \mathcal Q_T^{v,\eta}\longrightarrow\mathcal Q_T^{0,\eta}
 \quad\text{in }L^p\qquad(v\downarrow0).
\end{equation}
For $0<\eta\le T/2$, the definition of $\chi$ gives, for $r,s\ge0$,
\[
 0\le1-\chi\left(\frac{r+s}{\eta}\right)
 \le\one_{\{r+s\le2\eta\}}
 \le\one_{[0,2\eta]^2}(r,s).
\]
With the convention \eqref{eq:endpoint-reading}, positivity therefore
gives the pathwise comparison
\begin{align*}
 0\le\mathcal Q_T^v-\mathcal Q_T^{v,\eta}
 &=\int_{[0,T]^2}K_{H_0}(r+s)
       \left[1-\chi\left(\frac{r+s}{\eta}\right)\right]
       K_{H_1}^{[v]}(B_r-\widetilde B_s)\dd r\dd s\\
 &\le\int_{[0,2\eta]^2}K_{H_0}(r+s)
       K_{H_1}^{[v]}(B_r-\widetilde B_s)\dd r\dd s
 =\mathcal Q_{2\eta}^v.
\end{align*}
The last functional uses the restrictions of the same two paths to
$[0,2\eta]$; at $H_1=\frac12$ and $v=0$, the comparison is between
integrals against the same positive measure $\mathfrak I_{B,\widetilde B}$.

Applying \eqref{eq:BM-weighted-block} on this shorter square yields
\begin{equation}\label{eq:BM-weighted-corner-uniform}
 \sup_{v\ge0}
 \|\mathcal Q_T^{v}-\mathcal Q_T^{v,\eta}\|_p
 \le C_hp^{1-H_1}\eta^{a_{h}}.
\end{equation}
Combining this with \eqref{eq:BM-weighted-cutoff-convergence} through
the triangle inequality gives
\[
 \limsup_{v\downarrow0}\|\mathcal Q_T^{v}-\mathcal Q_T^0\|_p
 \le 2C_hp^{1-H_1}\eta^{a_{h}},
\]
and letting $\eta\downarrow0$, which is legitimate because $a_{h}>0$,
proves the second half of
\eqref{eq:BM-spatial-regularization-convergence}.
\end{proof}

\begin{lemma}[Mean self- and mutual intersection measures]
\label{lem:Brownian-mean-measure}
Let $B,\widetilde B$ be independent standard Brownian motions started at
zero. For every $T>0$, every nonnegative Borel function $f$ on
$[0,T]^2$ and every $v\ge0$,
\begin{align}
 E\left[\int_{[0,T]^2}f(s,t)K_{H_1}^{[v]}(B_s-B_t)\dd s\dd t\right]
 &=\kappa_{H_1}\int_{[0,T]^2}f(s,t)(|t-s|+v)^{H_1-1}\dd s\dd t,
 \label{eq:BM-mean-occupation-measure}\\
 E\left[\int_{[0,T]^2}f(s,t)K_{H_1}^{[v]}(B_s-\widetilde B_t)\dd s\dd t\right]
 &=\kappa_{H_1}\int_{[0,T]^2}f(s,t)(s+t+v)^{H_1-1}\dd s\dd t.
 \label{tp:Brownian-mean-measure}
\end{align}
At $H_1=\frac12$ and $v=0$, the left-hand sides mean
$E[\int f\dd\mathfrak I_B]$ and
$E[\int f\dd\mathfrak I_{B,\widetilde B}]$, respectively.
Both identities hold in $[0,\infty]$. In particular,
\[
 E\left[\int_{[0,T]^2}K_{H_1}(B_s-B_t)\dd s\dd t\right]
 =\frac{2\kappa_{H_1}T^{1+H_1}}{H_1(H_1+1)}<\infty.
\]
\end{lemma}

\begin{proof}
Suppose first that $H_1>\frac12$, or that $H_1=\frac12$ and $v>0$.
For $s\ne t$ the increment $B_t-B_s$ is centered Gaussian with variance
$|t-s|$, so that, by \eqref{eq:kernel-at-zero} applied to the kernel
$K_{H_1}^{[v]}*\heat_{|t-s|}=K_{H_1}^{[v+|t-s|]}$,
\[
 E\left[K_{H_1}^{[v]}(B_s-B_t)\right]
 =K_{H_1}^{[v+|t-s|]}(0)
 =\kappa_{H_1}(|t-s|+v)^{H_1-1}.
\]
For the independent pair, $B_s-\widetilde B_t$ is centered
Gaussian with variance $s+t$, so the same computation gives
$E[K_{H_1}^{[v]}(B_s-\widetilde B_t)]
 =\kappa_{H_1}(s+t+v)^{H_1-1}$.
Tonelli's theorem for nonnegative integrands proves both mean
identities.

It remains to justify the case $H_1=\frac12$ and $v=0$ at the level
of occupation measures. First let $B,\widetilde B$ be two independent
standard Brownian motions started at zero. For a continuous function
$q$ on $[0,T]^2$, Lemma~\ref{lem:Brownian-spatial-regularization} gives
\[
 \int_{[0,T]^2}q(r,s)\heat_v(B_r-\widetilde B_s)\dd r\dd s
 \longrightarrow
 \int_{[0,T]^2}q\dd\mathfrak I_{B,\widetilde B}
 \quad\text{in }L^1.
\]
For $v>0$ the expectation on the left is
$\int_{[0,T]^2}q(r,s)\heat_{r+s+v}(0)\dd r\dd s$.
The bound
$\heat_{r+s+v}(0)\le(2\pi(r+s))^{-1/2}$ provides an integrable
majorant off the origin. Passing to the limit therefore gives
\[
 E\left[\int_{[0,T]^2}q\dd\mathfrak I_{B,\widetilde B}\right]
 =\int_{[0,T]^2}q(r,s)\heat_{r+s}(0)\dd r\dd s.
\]
Both sides define finite measures, the finiteness on the left
following from \eqref{eq:BM-unweighted-block} at $v=0$.
Continuous functions determine finite measures on $[0,T]^2$, so
\begin{equation}\label{eq:BM-independent-mean-endpoint}
 E[\mathfrak I_{B,\widetilde B}(\dd r,\dd s)]
 =\heat_{r+s}(0)\dd r\dd s
 =\frac{\dd r\dd s}{\sqrt{2\pi(r+s)}}.
\end{equation}
This proves \eqref{tp:Brownian-mean-measure} at $H_1=\frac12$ and
$v=0$, for every nonnegative Borel function.

To transfer it to a single path, recall the intervals of
\eqref{eq:dyadic-JI} and fix the scaled rectangle
$I_{M,k,T}\times J_{M,k,T}$, where $I_{M,k,T}=TI_{M,k}$ and
$J_{M,k,T}=TJ_{M,k}$. Put $\delta=T2^{-M}$ and
$\theta=T(2k-1)2^{-M}$, so that
$I_{M,k,T}=[\theta-\delta,\theta)$ and
$J_{M,k,T}=[\theta,\theta+\delta)$, and define
\[
 B^-_r=B_{\theta-r}-B_\theta,\qquad
 B^+_u=B_{\theta+u}-B_\theta,\qquad 0\le r,u\le\delta.
\]
These are independent standard Brownian motions, since they use
increments on disjoint time intervals. The changes of variables
$\sigma=\theta-r$ and $\tau=\theta+u$ in the occupation formula give,
for $z\in\R$ and $r,u\in[0,\delta]$,
\begin{align*}
 L_r^{B^-}(z)
 &=L_\theta^B(z+B_\theta)-L_{\theta-r}^B(z+B_\theta),\\
 L_u^{B^+}(z)
 &=L_{\theta+u}^B(z+B_\theta)-L_\theta^B(z+B_\theta).
\end{align*}
Indeed, integrating either identity against a nonnegative test
function of $z$ gives precisely the corresponding change of time
variable, and continuity of the local times identifies the densities
at every $z$. Thus the restrictions of $\mu_x^B$ to the two half
intervals are the images of $\mu_{x-B_\theta}^{B^-}$ and
$\mu_{x-B_\theta}^{B^+}$ under these time maps. Using
\eqref{eq:two-BM-occupation-measure} and then changing $x=z+B_\theta$
yields
\begin{align*}
 \int_{I_{M,k,T}\times J_{M,k,T}}f(\sigma,\tau)
                              \mathfrak I_B(\dd\sigma,\dd\tau)
 &=\int_\R\int_{[0,\delta]^2}f(\theta-r,\theta+u)
       \mu_{x-B_\theta}^{B^-}(\dd r)
       \mu_{x-B_\theta}^{B^+}(\dd u)\dd x\\
 &=\int_{[0,\delta]^2}f(\theta-r,\theta+u)
                              \mathfrak I_{B^-,B^+}(\dd r,\dd u).
\end{align*}
The choice of open or closed endpoints is immaterial because these
measures assign zero mass to horizontal and vertical lines.
Taking expectations and applying \eqref{eq:BM-independent-mean-endpoint}
identifies the mean on this rectangle as
$\heat_{\tau-\sigma}(0)\dd\sigma\dd\tau$.
The dyadic partition \eqref{eq:dyadic-triangle-partition}, rescaled to
$[0,T]$, together with symmetry and the zero mass of the diagonal,
then gives
\[
 E[\mathfrak I_B(\dd s,\dd t)]
 =\heat_{|t-s|}(0)\dd s\dd t.
\]
This proves \eqref{eq:BM-mean-occupation-measure} also at
$H_1=\frac12$ and $v=0$, for every nonnegative Borel $f$.

Finally, for $f\equiv1$ and $v=0$,
\[
 \int_{[0,T]^2}|t-s|^{H_1-1}\dd s\dd t
 =2\int_0^Tu^{H_1-1}(T-u)\dd u
 =\frac{2T^{1+H_1}}{H_1(H_1+1)},
\]
which is finite because $H_1>0$. At $H_1=\frac12$,
$\kappa_{1/2}=\frac1{\sqrt{2\pi}}$, and hence
$E[\mathfrak I_B([0,T]^2)]=\frac{8T^{3/2}}{3\sqrt{2\pi}}$.
\end{proof}

\subsection{Dyadic block estimates}
\label{subsec:Brownian-blocks}

Recall the intervals $I_{M,k},J_{M,k}$ of \eqref{eq:dyadic-JI},
their scaled versions $I_{M,k,t}=tI_{M,k}$,
$J_{M,k,t}=tJ_{M,k}$, and the unregularized blocks
$C_{M,k,t},A_{M,k,t}$ of
\eqref{fr:Brownian-block}-\eqref{eq:block-continuous}.
For their regularized counterparts, set
\begin{align*}
 C^{\varepsilon}_{M,k,t}
 &=\int_{I_{M,k,t}\times J_{M,k,t}}
       K_{H_0}^{[2\varepsilon^4]}(\tau-\sigma)
       K_{H_1}^{[2\varepsilon^2]}(B_\sigma-B_\tau)\dd\sigma\dd\tau,\\
 A^{\varepsilon}_{M,k,t}
 &=C^{\varepsilon}_{M,k,t}-E[C^{\varepsilon}_{M,k,t}].
\end{align*}

\begin{lemma}[Dyadic block estimates and convergence]
\label{lem:Brownian-blocks}
Let $t>0$, $\varepsilon>0$, $M\ge1$ and $1\le k\le2^{M-1}$. The variables
$C_{M,k,t}$ and $C^{\varepsilon}_{M,k,t}$ are integrable, and
\begin{equation}\label{eq:Brownian-block-mean}
 E[C_{M,k,t}]
 =A_h\int_{I_{M,k,t}\times J_{M,k,t}}
                 (\tau-\sigma)^{a_h-2}\dd\sigma\dd\tau.
\end{equation}
Uniformly in $t,M,k$, $\varepsilon>0$ and $p\ge2$,
\begin{equation}\label{eq:Brownian-block-Lp}
 \begin{split}
 &\|C_{M,k,t}\|_p+\|A_{M,k,t}\|_p
   +\|C^{\varepsilon}_{M,k,t}\|_p+\|A^{\varepsilon}_{M,k,t}\|_p\\
 &\hspace{25mm}\le C_h\sqrt p\,(t2^{-M})^{a_h}.
 \end{split}
\end{equation}
For fixed $t,M,k$ and every finite $p\ge2$,
\begin{equation}\label{eq:Brownian-block-reg-conv}
 \|C^{\varepsilon}_{M,k,t}-C_{M,k,t}\|_p
 +\|A^{\varepsilon}_{M,k,t}-A_{M,k,t}\|_p
 \longrightarrow0\qquad(\varepsilon\downarrow0).
\end{equation}
For fixed $M$ and $\varepsilon$, the pairs
$(C_{M,k,t},C^{\varepsilon}_{M,k,t})$ are independent in $k$,
and so are the corresponding centered pairs.
\end{lemma}

\begin{proof}
\textit{Step 1: moment bounds, means and independence.}
Put $\delta=t2^{-M}$, $\theta=t(2k-1)2^{-M}$, and define
\begin{equation}\label{eq:BM-block-half-paths}
 B^-_r=B_{\theta-r}-B_\theta,\qquad
 B^+_u=B_{\theta+u}-B_\theta,\qquad 0\le r,u\le\delta.
\end{equation}
These are independent standard Brownian motions. For
$\sigma=\theta-r$ and $\tau=\theta+u$, one has
$\tau-\sigma=r+u$ and $B_\sigma-B_\tau=B^-_r-B^+_u$.
At $H_1=\frac12$, the same change of variables holds for $\mathfrak I_B$, as shown in the proof of Lemma~\ref{lem:Brownian-mean-measure}. With $\mathcal Q_\delta^v$ defined in
\eqref{eq:Brownian-weighted-functionals}, this gives
\[
 C_{M,k,t}=\mathcal Q_\delta^0(B^-,B^+),\qquad
 0\le C^{\varepsilon}_{M,k,t}
      \le C_{H_0}\mathcal Q_\delta^{2\varepsilon^2}(B^-,B^+),
\]
where the second inequality follows from
$K_{H_0}^{[2\varepsilon^4]}(r+u)\le C_{H_0}K_{H_0}(r+u)$
in Lemma~\ref{lem:kernel-regularization}. 
Equation~\eqref{eq:BM-weighted-block} and
$p^{1-H_1}\le\sqrt p$ give the asserted bounds for the uncentered
blocks. Centering increases their norms by at most a factor of two,
which proves \eqref{eq:Brownian-block-Lp}.

Lemma~\ref{lem:Brownian-mean-measure} gives
\[
 E[C_{M,k,t}]
 =c_{H_0}\kappa_{H_1}\int_{I_{M,k,t}\times J_{M,k,t}}
       (\tau-\sigma)^{2H_0-2+H_1-1}\dd\sigma\dd\tau,
\]
which is \eqref{eq:Brownian-block-mean}, since
$2H_0+H_1-3=a_h-2$ and $A_h=c_{H_0}\kappa_{H_1}$.
Its finiteness also follows directly by writing $w=r+u$:
\[
 \int_{[0,\delta]^2}(r+u)^{a_h-2}\dd r\dd u
 =\int_0^{2\delta}\min\{w,2\delta-w\}w^{a_h-2}\dd w
 \le\frac{(2\delta)^{a_h}}{a_h}.
\]

For independence, let $\ell_{M,k}=t(2k-2)2^{-M}$ be the left endpoint of the
parent interval, and put
\[
 \mathcal G_{M,k}
 :=\sigma\bigl(B_{\ell_{M,k}+r}-B_{\ell_{M,k}}:
                                      0\le r\le t2^{1-M}\bigr).
\]
Each block, centered or uncentered and with or without regularization,
is $\mathcal G_{M,k}$-measurable: spatial translation cancels in
the path differences. For fixed $M$, these sigma-fields are independent
by independence of Brownian increments on disjoint parent intervals.
This proves the asserted independence in $k$.

\smallskip
\noindent\textit{Step 2: control of the singular corner.}
Fix $0<\eta<\frac\delta2$, recall the cutoff $\chi$ of
\eqref{eq:cutoff}, and define
\begin{align*}
 C^{(\eta)}_{M,k,t}
 &=\int_{I_{M,k,t}\times J_{M,k,t}}
       K_{H_0}(\tau-\sigma)
       \chi\left(\frac{\tau-\sigma}{\eta}\right)
                    K_{H_1}(B_\sigma-B_\tau)\dd\sigma\dd\tau,\\
 C^{\varepsilon;\eta}_{M,k,t}
 &=\int_{I_{M,k,t}\times J_{M,k,t}}
       K_{H_0}^{[2\varepsilon^4]}(\tau-\sigma)
       \chi\left(\frac{\tau-\sigma}{\eta}\right)
                    K_{H_1}^{[2\varepsilon^2]}(B_\sigma-B_\tau)
                                       \dd\sigma\dd\tau.
\end{align*}
In the half-path coordinates of \eqref{eq:BM-block-half-paths},
subtracting a cutoff block from the full block inserts the factor
$1-\chi((r+u)/\eta)$ into its integrand. By \eqref{eq:cutoff},
\[
 0\le1-\chi\left(\frac{r+u}{\eta}\right)
 \le\one_{\{r+u\le2\eta\}}
 \le\one_{[0,2\eta]^2}(r,u).
\]
Consequently, positivity gives
\begin{align*}
 0\le C_{M,k,t}-C^{(\eta)}_{M,k,t}
 &=\int_{[0,\delta]^2}K_{H_0}(r+u)
     \left[1-\chi\left(\frac{r+u}{\eta}\right)\right]
     K_{H_1}(B^-_r-B^+_u)\dd r\dd u\\
 &\le \int_{[0,2\eta]^2}K_{H_0}(r+u)
     K_{H_1}(B^-_r-B^+_u)\dd r\dd u
 =\mathcal Q_{2\eta}^0(B^-,B^+).
\end{align*}
For the regularized block, the same indicator bound and
$K_{H_0}^{[2\varepsilon^4]}(r+u)\le C_{H_0}K_{H_0}(r+u)$ give
\begin{align*}
 0\le C^{\varepsilon}_{M,k,t}-C^{\varepsilon;\eta}_{M,k,t}
 &\le C_{H_0}\int_{[0,2\eta]^2}K_{H_0}(r+u)
       K_{H_1}^{[2\varepsilon^2]}(B^-_r-B^+_u)\dd r\dd u\\
 &=C_{H_0}\mathcal Q_{2\eta}^{2\varepsilon^2}(B^-,B^+).
\end{align*}
Since $2\eta<\delta$, the shorter integrals use restrictions of the
same half-paths. At $H_1=\frac12$, the unregularized comparison is
understood through \eqref{eq:endpoint-reading}, using positivity of
the measure $\mathfrak I_{B^-,B^+}$.
Applying \eqref{eq:BM-weighted-block} at horizon $2\eta$ yields
\begin{equation}\label{eq:Brownian-corner}
 \|C_{M,k,t}-C^{(\eta)}_{M,k,t}\|_p
 +\sup_{\varepsilon>0}\|C^{\varepsilon}_{M,k,t}
                   -C^{\varepsilon;\eta}_{M,k,t}\|_p
 \le C_{h}\sqrt p\,\eta^{a_{h}}.
\end{equation}

\smallskip
\noindent\textit{Step 3: convergence at a fixed cutoff.}
For fixed $\eta$, set
\[
 d_{\varepsilon,\eta}
 :=\sup_{\eta\le w\le2\delta}
             |K_{H_0}^{[2\varepsilon^4]}(w)-K_{H_0}(w)|\stackrel{\varepsilon \to 0}{\longrightarrow}0
\]
by Lemma~\ref{lem:kernel-regularization}, and let
\[
 \overline C^{\varepsilon;\eta}_{M,k,t}
 :=\int_{I_{M,k,t}\times J_{M,k,t}}
     K_{H_0}(\tau-\sigma)
     \chi\left(\frac{\tau-\sigma}{\eta}\right)
                 K_{H_1}^{[2\varepsilon^2]}(B_\sigma-B_\tau)
                                       \dd\sigma\dd\tau.
\]
Then \eqref{eq:BM-unweighted-block} gives the quantitative estimate
\begin{align*}
 \|C^{\varepsilon;\eta}_{M,k,t}
                   -\overline C^{\varepsilon;\eta}_{M,k,t}\|_p
 &\le d_{\varepsilon,\eta}
           \|\mathcal T_\delta^{2\varepsilon^2}(B^-,B^+)\|_p\\
 &\le C_{h}\sqrt p\,\delta^{1+H_1}d_{\varepsilon,\eta}
\stackrel{\varepsilon \to 0}{\longrightarrow} 0.
\end{align*}
For the spatial limit, apply \eqref{eq:BM-bounded-weight-convergence}
to the bounded weight on $[0,\delta]^2$ given by
\[
 f_\eta(r,u)=K_{H_0}(r+u)\chi\bigl(\tfrac{r+u}\eta\bigr),
 \qquad f_\eta(0,0)=0.
\]
It yields $\overline C^{\varepsilon;\eta}_{M,k,t}\to C^{(\eta)}_{M,k,t}$
in $L^p$ as $\varepsilon\downarrow0$.

\smallskip
\noindent\textit{Step 4: removal of the cutoff and centering.}
The three approximation errors satisfy
\begin{align*}
 \|C^{\varepsilon}_{M,k,t}-C_{M,k,t}\|_p
 &\le C_{h}\sqrt p\,\eta^{a_{h}}
       +C_{h}\sqrt p\,\delta^{1+H_1}d_{\varepsilon,\eta}\\
 &\quad+\|\overline C^{\varepsilon;\eta}_{M,k,t}
                              -C^{(\eta)}_{M,k,t}\|_p.
\end{align*}
First let $\varepsilon\downarrow0$, then $\eta\downarrow0$, which is
legitimate because $a_{h}>0$. This proves
$C^{\varepsilon}_{M,k,t}\to C_{M,k,t}$ in $L^p$; centering gives
\[
 \|A^{\varepsilon}_{M,k,t}-A_{M,k,t}\|_p
 \le\|C^{\varepsilon}_{M,k,t}-C_{M,k,t}\|_p
          +|E[C^{\varepsilon}_{M,k,t}-C_{M,k,t}]|
 \le2\|C^{\varepsilon}_{M,k,t}-C_{M,k,t}\|_p\longrightarrow0.
\]
\end{proof}

The blocks also give a dyadic decomposition of $Y_t^\varepsilon(B)$. By \eqref{eq:Gamma-eps}, \eqref{eq:Y-eps} and
Lemma~\ref{lem:Brownian-mean-measure},
\begin{equation}\label{eq:Yeps-centered-integral}
 \begin{split}
 Y_t^\varepsilon(B)
 &=2\int_{0\le\sigma<\tau<t}K_{H_0}^{[2\varepsilon^4]}(\tau-\sigma)\\
 &\quad\times\left(
 K_{H_1}^{[2\varepsilon^2]}(B_\sigma-B_\tau)
 -\kappa_{H_1}(\tau-\sigma+2\varepsilon^2)^{H_1-1}\right)
 \dd\sigma\dd\tau.
 \end{split}
\end{equation}
Here symmetry accounts for the factor two, and the diagonal has
zero Lebesgue measure.
The partition \eqref{eq:dyadic-triangle-partition}, rescaled to
$[0,t]$, together with symmetry and Tonelli's theorem, gives
\[
 \Gamma_t^\varepsilon=2\sum_{M\ge1}\sum_k C^{\varepsilon}_{M,k,t},
 \qquad E[\Gamma_t^\varepsilon]
                 =2\sum_{M\ge1}\sum_k E[C^{\varepsilon}_{M,k,t}].
\]
For fixed $\varepsilon>0$, both nonnegative sums are finite, the first one in $L^1$,
since $K_{H_0}^{[2\varepsilon^4]}$ and $K_{H_1}^{[2\varepsilon^2]}$ are bounded and
the time square has finite area.
In particular, for $L\ge1$,
\begin{align*}
 \left\|Y_t^\varepsilon(B)
           -2\sum_{M\le L}\sum_k A^{\varepsilon}_{M,k,t}\right\|_1
 &\le2\left\|\Gamma_t^\varepsilon(B)
           -2\sum_{M\le L}\sum_k C^{\varepsilon}_{M,k,t}\right\|_1\\
 &=4\sum_{M>L}\sum_k E[C^{\varepsilon}_{M,k,t}]
 \longrightarrow0\qquad(L\to\infty).
\end{align*}
Indeed, the last expression is four times the tail of a nonnegative
series whose total sum is $E[\Gamma_t^\varepsilon]/2<\infty$; its tails
therefore vanish, with $\varepsilon>0$ held fixed.
Thus the centered decomposition is
\begin{equation}\label{eq:Yeps-dyadic}
 Y_t^\varepsilon(B)=2\sum_{M\ge1}\sum_k A^{\varepsilon}_{M,k,t},
 \qquad\text{initially in }L^1.
\end{equation}

\subsection{Construction of the centered self-intersection}
\label{subsec:Brownian-summation}

We use a quantitative form of standard sub-Gaussian estimates. The following statement is a slight variant of \cite[Proposition~2.5.2 and Proposition~2.6.1]{Vershynin}; its proof can be found in Appendix \ref{app:subgaussian}.

\begin{lemma}[Centered sub-Gaussian estimates]\label{lem:subgaussian}
Suppose $E[Z]=0$ and $\|Z\|_p\le c\sqrt p\,s$ for every $p\ge2$,
where $c>0$ and $s\ge0$. Then
\begin{equation}\label{eq:subgaussian-mgf}
 E[e^{\lambda Z}]\le\exp(8c^2\lambda^2s^2),
 \qquad\lambda\in\R.
\end{equation}
More generally, let $Z_1,\ldots,Z_n$ be independent random variables
such that, for every $i$ and $p\ge2$,
\[
 E[Z_i]=0,\qquad \|Z_i\|_p\le c\sqrt p\,s_i,
 \qquad s_i\ge0.
\]
Put $\Sigma_n=\sum_{i=1}^n Z_i$ and $V_n^2=\sum_{i=1}^n s_i^2$.
Then
\begin{align}
 E[e^{\lambda\Sigma_n}]
 &\le\exp(8c^2\lambda^2V_n^2),
       \qquad \lambda\in\R,\label{eq:subgaussian-sum-mgf}\\
 \|\Sigma_n\|_p
 &\le C c\sqrt p\,V_n,
       \qquad p\ge2,\label{eq:subgaussian-sum-Lp}
\end{align}
with an absolute constant $C$. If $V_n=0$, every $Z_i$ and their sum vanish almost surely.
\end{lemma}

\begin{corollary}[Positive sub-Gaussian variables]\label{tp:moment-to-exp}
If $X\ge0$ and $\|X\|_p\le c\sqrt p$ for every $p\ge2$, then, for some universal constant $C>0$,
\begin{equation}\label{tp:moment-to-exp-bound}
 E\left[e^{\lambda X}\right]\le e^{\sqrt2c\lambda+Cc^2\lambda^2},\qquad\lambda\ge0.
\end{equation}
\end{corollary}
\begin{proof}
We have $E\left[X\right]\le\sqrt2c$ and
$\|X-E\left[X\right]\|_p\le2c\sqrt p$. Apply Lemma~\ref{lem:subgaussian}
to $X-E\left[X\right]$, then multiply the resulting bound by $e^{\lambda E\left[X\right]}$.
\end{proof}

We now sum the dyadic levels. The following theorem proves
Proposition~\ref{prop:Y} and identifies the limit of the regularized
self-intersections.

\begin{theorem}[Construction and regularization of the centered self-intersection]
\label{thm:Brownian-self-intersection}
Let $t \geq 0$. The series defining $Y_t(B)$ in Proposition~\ref{prop:Y}
converges in every finite $L^p$, and $Y_t(B)$ is centered. For every
$p\ge2$ and every integer $L\ge0$,
\begin{equation}\label{eq:Brownian-Y-tail-Lp}
 \begin{split}
 &\left\|Y_t(B)-2\sum_{M\le L}\sum_k A_{M,k,t}\right\|_p\\
 &\quad+\sup_{\varepsilon>0}
   \left\|Y_t^\varepsilon(B)
              -2\sum_{M\le L}\sum_k A^{\varepsilon}_{M,k,t}\right\|_p
 \le C_{h}\sqrt p\,t^{a_{h}}2^{-L(a_{h}-\frac12)}.
 \end{split}
\end{equation}
Moreover,
\begin{equation}\label{eq:Yeps-dyadic-identification}
 Y_t^\varepsilon(B)\longrightarrow Y_t(B)=2\sum_{M\ge1}\sum_kA_{M,k,t}
 \quad\text{in every finite }L^p,\qquad \varepsilon\downarrow0,
\end{equation}
on the given Brownian probability space. Uniformly for $p\ge2$,
\begin{equation}\label{eq:Brownian-Y-subgaussian-Lp}
 \sup_{\varepsilon>0}\|Y_t^\varepsilon(B)\|_p+\|Y_t(B)\|_p
 \le C_{h}\sqrt p\,t^{a_{h}}.
\end{equation}
For every $\lambda\in\R$,
\begin{equation}\label{eq:Brownian-Y-exp-unif}
 \sup_{\varepsilon>0}E[e^{\lambda Y_t^\varepsilon(B)}]
 \le e^{C_{h}\lambda^2t^{2a_{h}}},
 \qquad E[e^{\lambda Y_t(B)}]\le e^{C_{h}\lambda^2t^{2a_{h}}}.
\end{equation}
All constants are uniform in $t$, which enters only through the
displayed powers.
\end{theorem}

\begin{proof}

For $t=0$ all blocks and $Y_0^\varepsilon(B)$ vanish, and every assertion is trivial; we therefore assume $t>0$.

\smallskip

\noindent
\textit{Step 1: estimate at a fixed level.}
At a fixed level $M$, Lemma~\ref{lem:Brownian-blocks} gives
independence of the centered blocks and their moment bounds.
Apply \eqref{eq:subgaussian-sum-Lp} with
$n=2^{M-1}$ and $s_k=(t2^{-M})^{a_{h}}$, both to $A_{M,k,t}$
and to $A^{\varepsilon}_{M,k,t}$. It follows that
\begin{align}
 &\left\|\sum_kA_{M,k,t}\right\|_p
       +\sup_{\varepsilon>0}
                    \left\|\sum_kA^{\varepsilon}_{M,k,t}\right\|_p\notag\\
 &\qquad\le C_{h}\sqrt p
              \left(2^{M-1}t^{2a_{h}}2^{-2Ma_{h}}\right)^{\frac12}
 \le C_{h}\sqrt p\,t^{a_{h}}2^{-M(a_{h}-\frac12)}.
 \label{eq:Brownian-level-Lp}
\end{align}
\smallskip
\noindent\textit{Step 2: summation of the levels.}
Write
\[
 \Sigma_{L,t}=2\sum_{M=1}^L\sum_kA_{M,k,t},\qquad
 \Sigma_{L,t}^\varepsilon
        =2\sum_{M=1}^L\sum_kA^{\varepsilon}_{M,k,t}.
\]
For $L'>L\ge0$, Minkowski's inequality and \eqref{eq:Brownian-level-Lp} give
\begin{align*}
 \|\Sigma_{L',t}-\Sigma_{L,t}\|_p
 +\sup_{\varepsilon>0}\|\Sigma_{L',t}^\varepsilon
                                      -\Sigma_{L,t}^\varepsilon\|_p
 &\le2C_{h}\sqrt p\,t^{a_{h}}\sum_{M=L+1}^{L'}2^{-M(a_{h}-\frac12)}\\
 &\le\frac{2C_{h}\sqrt p\,t^{a_{h}}}{1-2^{-(a_{h}-\frac12)}}\,2^{-(L+1)(a_{h}-\frac12)}
 \longrightarrow0\quad(L\to\infty),
\end{align*}
the series converging because $a_{h}>\frac12$.
Completeness of $L^p$ gives a limit $Y_t(B)$ for $\Sigma_{L,t}$ and, for
each fixed $\varepsilon$, a limit for $\Sigma_{L,t}^\varepsilon$.
These limits are consistent for different $p\ge2$, since
convergence in any such $L^p$ implies convergence in $L^1$.
Equation~\eqref{eq:Yeps-dyadic} already identifies the $L^1$ limit
of $\Sigma_{L,t}^\varepsilon$ as $Y_t^\varepsilon(B)$. Moreover,
\[
 |E[Y_t(B)]|
 =|E[Y_t(B)-\Sigma_{L,t}]|
 \le\|Y_t(B)-\Sigma_{L,t}\|_1\stackrel{L\to \infty}{\longrightarrow}0,
\]
so $Y_t(B)$ is centered. Letting $L'\to\infty$ in the partial-sum
estimate yields
\[
 \|Y_t(B)-\Sigma_{L,t}\|_p
 +\sup_{\varepsilon>0}\|Y_t^\varepsilon(B)-\Sigma_{L,t}^\varepsilon\|_p
 \le C_{h}\sqrt p\,t^{a_{h}}2^{-L(a_{h}-\frac12)},
\]
which is \eqref{eq:Brownian-Y-tail-Lp}.

\smallskip
\noindent\textit{Step 3: identification of the limit and exponential bounds.}
For fixed $L$, the sum is finite and
\[
 \|\Sigma_{L,t}^\varepsilon-\Sigma_{L,t}\|_p
 \le2\sum_{M=1}^L\sum_k
             \|A^{\varepsilon}_{M,k,t}-A_{M,k,t}\|_p
 \stackrel{\varepsilon \to 0}{\longrightarrow}0
\]
by \eqref{eq:Brownian-block-reg-conv}. Thus
\begin{align*}
 \|Y_t^\varepsilon(B)-Y_t(B)\|_p
 &\le\|Y_t^\varepsilon(B)-\Sigma_{L,t}^\varepsilon\|_p
        +\|\Sigma_{L,t}^\varepsilon-\Sigma_{L,t}\|_p
        +\|\Sigma_{L,t}-Y_t(B)\|_p,\\
 \limsup_{\varepsilon\downarrow0}\|Y_t^\varepsilon(B)-Y_t(B)\|_p
 &\le C_{h}\sqrt p\,t^{a_{h}}2^{-L(a_{h}-\frac12)}\longrightarrow0
                                      \quad(L\to\infty),
\end{align*}
which proves \eqref{eq:Yeps-dyadic-identification}.

Taking $L=0$ in \eqref{eq:Brownian-Y-tail-Lp} proves
\eqref{eq:Brownian-Y-subgaussian-Lp}.
All variables are centered, so Lemma~\ref{lem:subgaussian}
gives \eqref{eq:Brownian-Y-exp-unif}.
\end{proof}

\section{\texorpdfstring{Direct construction of the continuous field}{Direct construction of the continuous field}}
\label{sec:continuous-identification}

This section proves Theorem~\ref{thm:continuous-summary} by constructing
the limit of the regularized solutions on the given Gaussian space.
The argument has two main steps. First, the Feynman-Kac formula
(Lemma~\ref{lem:regularized-FK}) expresses the Gaussian moments of
the classical approximations in terms of a deterministic remainder,
the centered self-intersections of Section~\ref{sec:Brownian-construction},
and mutual intersections (Proposition~\ref{prop:regularized-moments}).
The joint limit of the mixed second moments makes the approximations
Cauchy in $L^2$ at each fixed point; uniform higher moments then give
convergence in every finite $L^p$
(Proposition~\ref{thm:continuous-construction}).
Proposition~\ref{prop:S-cont} records the resulting $\mathcal S$-transform
and second-moment identities, which will be compared with their
discrete counterparts in Sections~\ref{subsec:legall-complet}
and~\ref{sec:two-path}.

The second step is to estimate space-time increments. Bounds for
self-intersection increments (Lemma~\ref{lem:continuous-self-increments})
and mutual-intersection increments
(Lemma~\ref{lem:continuous-path-increments}), combined with Gaussian
hypercontractivity, yield uniform $L^q$ estimates for the field
(Proposition~\ref{prop:continuous-functional}).
A quantitative continuity criterion (Lemma~\ref{lem:compact-upgrade})
then gives a continuous field and convergence in
$L^p(\Omega;C(\mathcal K))$ on every compact rectangle
(Theorem~\ref{thm:continuous-field}). We then prove the necessity of the renormalization
(Proposition~\ref{prop:renormalization-needed}) and identify the field
of Definition~\ref{def:continuous-solution} with the solution of
\cite{CDOT-KRP} when $2H_0+H_1>\frac53$
(Proposition~\ref{prop:CDOT-comparison}).

\smallskip

\textit{For the whole section, we fix $h=(H_0,H_1)$ satisfying the
conditions in \eqref{fr:range}}.

\subsection{The regularized equation and its Feynman-Kac formula}

We use the following form of the Feynman-Kac formula, under linear
spatial growth bounds on the potential and its derivatives.

\begin{lemma}[Feynman-Kac formula under linear growth bounds]
\label{lem:FK-unbounded}
Let $T>0$ and let $q$ be a smooth function on
$[0,T]\times\R$ such that
\begin{equation}\label{eq:FK-linear-upper-bound}
 |q(s,x)|+|\partial_sq(s,x)|+|\partial_xq(s,x)|
       +|\partial_{xx}q(s,x)|\le A_T(1+|x|),
\end{equation}
for all $(s,x)\in[0,T]\times\R$ and some finite $A_T\ge0$.
Let $B$ be a standard Brownian motion
started at zero, and write $E_B$ for the expectation with respect to
$B$. Then the function
\begin{equation}\label{eq:bounded-potential-FK}
 u(t,x)=E_B\left[\exp\left(\int_0^t q(r,x+B_{t-r})\dd r\right)\right]
\end{equation}
is well defined on $[0,T] \times \R$, and it is the unique classical solution, in the growth class
\eqref{eq:growth-class}, of the problem
$$\partial_tu=\frac12\partial_{xx}u+qu, \qquad   u(0,\cdot)=1.$$
\end{lemma}

The proof is given in Appendix~\ref{app:FK-unbounded}: uniqueness
follows from Friedman's uniqueness theorem, based on the maximum
principle, while existence is obtained directly from the Brownian
representation.

For the Gaussian estimates below, we use the Fourier convention
\begin{equation}\label{eq:Fourier-convention}
 \widehat\varphi(\lambda,\xi)
 =\int_{\R^2}e^{-i(t\lambda+x\xi)}
                    \varphi(t,x)\dd t\dd x.
\end{equation}
With the convention \eqref{eq:Fourier-convention}, put
\begin{equation}\label{eq:spectral-measure}
 \mathfrak a_r=\frac{r\Gamma(2r)\sin(\pi r)}\pi,
 \qquad
 \mathfrak m_{h}(\dd\lambda,\dd\xi)
 =\mathfrak a_{H_0}\mathfrak a_{H_1}
       |\lambda|^{1-2H_0}|\xi|^{1-2H_1}\dd\lambda\dd\xi,
\end{equation}
which is the spectral measure of \eqref{fr:noise-covariance},
obtained by Fourier transforming $K_{H_0}(t)$ and
$K_{H_1}(x)$. At $H_1=\frac12$, one has $\mathfrak a_{1/2}=\frac1{2\pi}$ and
$|\xi|^{1-2H_1}=1$, so that $\mathfrak m_{h}$ reduces to
$$\mathfrak m_{h}=\frac{H_0\Gamma(2H_0)\sin(\pi H_0)}{2\pi^2}
|\lambda|^{1-2H_0}\dd\lambda\dd\xi,$$
the measure associated with
\eqref{fr:noise-covariance}.

\begin{lemma}[The classical regularized solution]\label{lem:regularized-FK}
Fix $\varepsilon>0$, let $B$ be a standard Brownian motion started
at zero and independent of $W$, and, for $t\ge0$ and $x\in\R$,
define the path function
\begin{equation}\label{eq:hB-tx}
 h^{B,\varepsilon}_{t,x}(r,y)
 =\int_0^t\rho_\varepsilon(s-r,x+B_{t-s}-y)\dd s.
\end{equation}

Then, almost surely, the process $U^\varepsilon$ defined for all
$t\ge0$ and $x\in\R$ by
\begin{equation}\label{eq:U-eps-tx}
 U_t^\varepsilon(x)
 =E_B\left[\exp\left(W(h^{B,\varepsilon}_{t,x})-t c_h^\varepsilon\right)\right],
 \qquad U_0^\varepsilon(x)=1
\end{equation}
is the unique classical solution, in the growth class
\eqref{eq:growth-class}, of the regularized problem \eqref{fr:regularized-equation}.
\end{lemma}

\begin{proof}
We first fix a $C^\infty$ version of $\dot W^\varepsilon$ on $\R^2$.
Such a version exists because the Gaussian mollifier gives finite
spectral moments of every order: Kolmogorov's continuity criterion
applied to the differentiated fields yields a smooth modification.
For every multi-index $\beta=(\beta_0,\beta_1)$, write
$D^\beta=\partial_s^{\beta_0}\partial_x^{\beta_1}$. The field
$D^\beta\dot W^\varepsilon$ is centered Gaussian and stationary, with
\[
 \sigma_{\varepsilon,\beta}^2
 :=\E\big[|D^\beta\dot W^\varepsilon(s,x)|^2\big]
 =\int_{\R^2}|\lambda|^{2\beta_0}|\xi|^{2\beta_1}
       e^{-\varepsilon^4\lambda^2-\varepsilon^2\xi^2}
       \,\mathfrak m_h(\dd\lambda,\dd\xi)<\infty.
\]

Set $V^\varepsilon(s,x)=\dot W^\varepsilon(s,x)-c_h^\varepsilon$.
We verify that, on a single event of probability one, for every
multi-index $\beta$ and every $T>0$ there is a finite random
constant $\Theta_{\varepsilon,\beta,T}$ such that
\begin{equation}\label{eq:regularized-potential-growth}
 |D^\beta V^\varepsilon(s,x)|
 \le\Theta_{\varepsilon,\beta,T}(1+|x|),
 \qquad (s,x)\in[0,T]\times\R.
\end{equation}
Fix $\beta$ and $T>0$. For $n=0,1,\ldots$, let
$I_n=[-(n+1),-n]\cup[n,n+1]$ and $Q_{T,n}=[0,T]\times I_n$.
The Sobolev embedding on $[0,T]\times[0,1]$, applied separately
to the two translated rectangles composing $Q_{T,n}$, gives a
constant $C_T$ independent of $n$ such that
\begin{equation}\label{eq:regularized-strip-Sobolev}
 \|D^\beta\dot W^\varepsilon\|_{L^\infty(Q_{T,n})}
 \le C_T\sum_{|\eta|\le2}
             \|D^{\beta+\eta}\dot W^\varepsilon\|_{L^2(Q_{T,n})}.
\end{equation}
Write $M_{\beta,T,n}=\|D^\beta\dot W^\varepsilon\|_{L^\infty(Q_{T,n})}$.
Squaring \eqref{eq:regularized-strip-Sobolev} and using
Cauchy-Schwarz, stationarity, and $|Q_{T,n}|=2T$, we obtain
\begin{align*}
 \E[M_{\beta,T,n}^2]
 &\le C_T\sum_{|\eta|\le2}\int_{Q_{T,n}}
       \E\big[|D^{\beta+\eta}\dot W^\varepsilon(s,x)|^2\big]\dd s\dd x\\
 &=2TC_T\sum_{|\eta|\le2}\sigma_{\varepsilon,\beta+\eta}^2
 =:C_{\varepsilon,\beta,T}<\infty,
\end{align*}
uniformly in $n$. Markov's inequality gives
\[
 \sum_{n\ge0}\Pp\big(M_{\beta,T,n}>1+n\big)
 \le C_{\varepsilon,\beta,T}\sum_{n\ge0}(1+n)^{-2}<\infty.
\]
By Borel-Cantelli, almost surely there is a finite integer
$n_0=n_0(\varepsilon,\beta,T)\geq 1$ such that
$M_{\beta,T,n}\le1+n$
for all $n\ge n_0$. Since $x\in I_n$ implies $n\le|x|$, this yields
\[
 |D^\beta\dot W^\varepsilon(s,x)|\le1+|x|,
 \qquad 0\le s\le T,\quad |x|\ge n_0.
\]
Consequently, \eqref{eq:regularized-potential-growth} holds with
\[
 \Theta_{\varepsilon,\beta,T}
 =|c_h^\varepsilon|
  +\|D^\beta\dot W^\varepsilon\|_{L^\infty([0,T]\times[-n_0,n_0])}+1<\infty.
\]
Let $\Omega_\varepsilon$ be the intersection of the smoothness event
and the probability-one events just obtained for all multi-indices
$\beta$ and all integer time horizons $T=m\ge1$. For arbitrary $T>0$, choose an
integer $m\ge T$ and restrict the estimate on $[0,m]$ to $[0,T]$.
Thus, on $\Omega_\varepsilon$,
\eqref{eq:regularized-potential-growth} holds for every $\beta$ and
every $T>0$.

Taking $\beta\in\{(0,0),(1,0),(0,1),(0,2)\}$ in
\eqref{eq:regularized-potential-growth} verifies
\eqref{eq:FK-linear-upper-bound} for $q=V^\varepsilon$ on every
finite time interval. The expectation in \eqref{eq:U-eps-tx} is precisely
\[
 E_B\left[\exp\left(\int_0^t
            V^\varepsilon(s,x+B_{t-s})\dd s\right)\right].
\]
Lemma~\ref{lem:FK-unbounded} therefore gives, on
$\Omega_\varepsilon$, finiteness, classical regularity, and
uniqueness in \eqref{eq:growth-class}, simultaneously for all
$t\ge0$ and $x\in\R$.
\end{proof}

We now relate the representation \eqref{eq:U-eps-tx} to the Brownian functionals of
Section~\ref{sec:Brownian-construction}. The covariance identity
\eqref{eq:smoothed-noise-covariance}, followed by reversal of both time
variables, gives
\[
 \|h^{B,\varepsilon}_{t,x}\|_{\Hh}^2=\Gamma_t^\varepsilon(B),
\]
with $\Gamma_t^\varepsilon$ defined by \eqref{eq:Gamma-eps}. The spatial
translation by $x$ cancels in the path differences. Since $B$ is
independent of $W$, conditionally on $B$ the variable
$W(h^{B,\varepsilon}_{t,x})$ is centered Gaussian with this variance.
Define the deterministic remainder
\begin{equation}\label{eq:R-eps}
 R^{\varepsilon}_t
 =\frac12E_B[\Gamma_t^\varepsilon(B)]-t\,c_h^\varepsilon.
\end{equation}

Together with \eqref{eq:Y-eps}, this gives the decomposition
\[
 W(h^{B,\varepsilon}_{t,x})-t\,c_h^\varepsilon
 =\left(W(h^{B,\varepsilon}_{t,x})
                 -\frac12\Gamma_t^\varepsilon(B)\right)
   +\frac12Y_t^\varepsilon(B)+R^{\varepsilon}_t.
\]
Thus \eqref{eq:stoch-exp} and \eqref{eq:U-eps-tx} yield the exact
factorization
\begin{equation}\label{eq:Gaussian-cont}
 U_t^\varepsilon(x)=e^{R^{\varepsilon}_t}
 E_B\left[e^{\frac12Y_t^\varepsilon(B)}
                    \cE(h^{B,\varepsilon}_{t,x})\right].
\end{equation}

This representation was already central to the moment analysis in
\cite[Section~4.2]{CDOT-moments}: it relates the regularized
Stratonovich and Skorohod Feynman-Kac formulas through the centered
Brownian self-intersection and a deterministic remainder.
It is also the continuous analogue of \eqref{intro:wick-factorization}.
Conditionally on $B$, the Gaussian exponential $\cE(h^{B,\varepsilon}_{t,x})$ has Gaussian expectation one. The deterministic counterterm is compared with the Brownian
mean of the conditional variance, leaving the random fluctuation
$Y_t^\varepsilon(B)$ in the exponent.

\subsection{The continuous deterministic remainder}
\label{subsec:continuous-counterterm}

We begin the analysis of \eqref{eq:Gaussian-cont} with the
convergence of its deterministic factor $e^{R^{\varepsilon}_t}$. By
\eqref{eq:R-eps}, $R^{\varepsilon}_t$ is the difference of
$\frac12E_B[\Gamma_t^\varepsilon(B)]$ and $t\,c_h^\varepsilon$, two
quantities that diverge as $\varepsilon\downarrow0$ (see
\eqref{eq:counterterm-time-domain}). Its convergence therefore rests on
an exact cancellation between them, and requires no Brownian estimate:
it follows from explicit computations with the covariance
\eqref{eq:smoothed-noise-covariance} and the scaling definition
\eqref{fr:counterterm-continuous}.

\begin{proposition}[Convergence of the remainder]
\label{prop:R-eps-limit}
Recall $\kappa_{H_1}$, $A_h$ and $R_t$ from
\eqref{eq:AH}-\eqref{eq:R-definition}. For $v>0$, define
\begin{equation}\label{eq:f-eps}
 f_h^{[\varepsilon]}(v)
 =\kappa_{H_1}\,K_{H_0}^{[2\varepsilon^4]}(v)
                    (v+2\varepsilon^2)^{H_1-1}.
\end{equation}

The counterterm defined by \eqref{fr:counterterm-continuous} is
finite and strictly positive, and admits the time-domain formula
\begin{equation}\label{eq:counterterm-time-domain}
 c_h^\varepsilon=\int_0^\infty f_h^{[\varepsilon]}(v)\dd v.
\end{equation}
Furthermore, for every $t\ge0$,
\begin{equation}\label{eq:R-eps-limit}
 R^{\varepsilon}_t\longrightarrow R_t,
 \qquad \sup_{\varepsilon>0}|R^{\varepsilon}_t|\le C_h t^{a_h}.
\end{equation}
At $H_1=\frac12$ one has $\kappa_{1/2}(v+2\varepsilon^2)^{-\frac12}
=\heat_{v+2\varepsilon^2}(0)$ and $2a_h-2=4H_0-3$.
\end{proposition}

\begin{proof}
\textit{Step 1: time-domain representation.}
Gaussian convolution in \eqref{eq:smoothed-noise-covariance} gives
\[
 f_h^{[\varepsilon]}(v)
 =E_B\left[\E\left[\dot W^\varepsilon(v,B_v)
                         \dot W^\varepsilon(0,0)\right]\right],
\]
since $E_B[K_{H_1}^{[2\varepsilon^2]}(B_v)]
=\kappa_{H_1}(v+2\varepsilon^2)^{H_1-1}$.
The envelope \eqref{eq:KH-envelope} of Lemma~\ref{lem:kernel-regularization}, applied with $H=H_0$, and $(v+2\varepsilon^2)^{H_1-1}\le v^{H_1-1}$ yield
\begin{equation}\label{eq:f-envelope}
 0\le f_h^{[\varepsilon]}(v)\le C_h v^{a_h-2},\qquad v>0.
\end{equation}
For fixed $\varepsilon$, $f_h^{[\varepsilon]}$ is bounded near zero and
strictly positive. Since $a_h<1$, it is integrable at infinity.
In particular, the unit-scale constant
$c_h^1=\int_0^\infty f_h^{[1]}(v)\dd v$ of
\eqref{fr:counterterm-continuous} is finite and strictly positive.

The substitution $v=\varepsilon^2z$ in the temporal convolution gives
\begin{align*}
 K_{H_0}^{[2\varepsilon^4]}(\varepsilon^2w)
 &=c_{H_0}\int_\R |\varepsilon^2w-v|^{2H_0-2}
                   \heat_{2\varepsilon^4}(v)\dd v\\
 &=c_{H_0}\varepsilon^{4H_0-4}\int_\R |w-z|^{2H_0-2}\heat_2(z)\dd z
 =\varepsilon^{4H_0-4}K_{H_0}^{[2]}(w).
\end{align*}
Thus the change of variables $v=\varepsilon^2w$ gives
\begin{align*}
 \int_0^\infty f_h^{[\varepsilon]}(v)\dd v
 &=\varepsilon^{2a_h-2}\kappa_{H_1}
       \int_0^\infty K_{H_0}^{[2]}(w)(w+2)^{H_1-1}\dd w\\
 &=\varepsilon^{2a_h-2}c_h^1=c_h^\varepsilon,
\end{align*}
by \eqref{eq:counterterm-unit} and the definition
\eqref{fr:counterterm-continuous}. This proves
\eqref{eq:counterterm-time-domain}.

\smallskip
\noindent\textit{Step 2: cancellation and limit of the remainder.}
The case $t=0$ is immediate; let $t>0$.
Brownian stationary increments,
\eqref{eq:smoothed-noise-covariance} and
Lemma~\ref{lem:Brownian-mean-measure} give
\[
 \frac12E_B\left[\Gamma_t^\varepsilon(B)\right]
 =\int_0^t(t-v)f_h^{[\varepsilon]}(v)\dd v.
\]
Combining this with \eqref{eq:counterterm-time-domain}
gives the exact cancellation
\begin{equation}\label{eq:R-eps-exact}
 R^{\varepsilon}_t
 =-\int_0^t v f_h^{[\varepsilon]}(v)\dd v
   -t\int_t^\infty f_h^{[\varepsilon]}(v)\dd v.
\end{equation}
For every $v>0$, $f_h^{[\varepsilon]}(v)\to A_h v^{a_h-2}$, by
Lemma~\ref{lem:kernel-regularization} and
$(v+2\varepsilon^2)^{H_1-1}\to v^{H_1-1}$.
The two integrands in \eqref{eq:R-eps-exact} are bounded respectively
by $C_h v^{a_h-1}$ on $(0,t)$ and $C_h v^{a_h-2}$ on $(t,\infty)$.
Both are integrable, the first because $a_h>0$ and the second because
$a_h<1$, and
\[
 \int_0^t v^{a_h-1}\dd v=\frac{t^{a_h}}{a_h},
 \qquad
 t\int_t^\infty v^{a_h-2}\dd v=\frac{t^{a_h}}{1-a_h}.
\]
Dominated convergence therefore proves
$R^{\varepsilon}_t\to-A_h t^{a_h}(\frac1{a_h}+\frac1{1-a_h})=R_t$,
which is \eqref{eq:R-eps-limit}, and the same majorants give its
uniform bound.
\end{proof}

\subsection{Mixed mutual intersections}

These functionals describe the cross covariances of Brownian path
integrals at two regularization scales. They will enter the mixed
second-moment identity \eqref{eq:continuous-mixed-second}, whose joint
limit makes $U_t^\varepsilon(x)$ a Cauchy family in $L^2(\Omega)$ in
Proposition~\ref{thm:continuous-construction}.

For $\varepsilon,{\epsp}>0$, Gaussian convolution in the covariance
inner product gives
\begin{equation}\label{eq:mixed-regularization}
 \E\left[\dot W^\varepsilon(s,x)\dot W^{\epsp}(t,y)\right]
 =K_{H_0}^{[\varepsilon^4+{\epsp}^4]}(t-s)
       K_{H_1}^{[\varepsilon^2+{\epsp}^2]}(x-y).
\end{equation}
For two independent Brownian motions $B,\widetilde B$ and a horizon
$t>0$, define
\begin{align}
 \alpha_t^{\varepsilon,{\epsp}}(B,\widetilde B)
 &=\ip{h^{B,\varepsilon}_{t,x}}{h^{\widetilde B,\epsp}_{t,x}}_{\Hh}\notag\\
 &=\int_{[0,t]^2}K_{H_0}^{[\varepsilon^4+{\epsp}^4]}(\tau-\sigma)
       K_{H_1}^{[\varepsilon^2+{\epsp}^2]}(B_\sigma-\widetilde B_\tau)
                                       \dd\sigma\dd\tau.
 \label{eq:alpha-mixed}
\end{align}
For the second equality, reverse both time variables in the definition
of the two path functions; the spatial point $x$ cancels. Write
\begin{equation}\label{eq:alpha-eps-again}
 \alpha_t^\varepsilon(B,\widetilde B)=\alpha_t^{\varepsilon,\varepsilon}(B,\widetilde B).
\end{equation}

\begin{lemma}
\label{lem:Brownian-lag}
Uniformly in $v\ge0$, for $p\ge2$ and $T>0$,
\begin{equation}\label{eq:Brownian-lag}
 \sup_{0\le u<T}
 \left\|\int_0^{T-u}
       K_{H_1}^{[v]}(B_{r+u}-\widetilde B_r)\dd r\right\|_p
 \le C_h p^{1-H_1}T^{H_1}.
\end{equation}
The same assertion holds with $B$ and $\widetilde B$ interchanged.
More generally, for $u,a\ge0$ and $\ell>0$, the integral over
$r\in[a,a+\ell]$ has conditional $L^p$ norm at most
$C_h p^{1-H_1}\ell^{H_1}$ given
\[
 \mathcal F_a^u
 =\sigma(B_r:0\le r\le a+u)\vee
       \sigma(\widetilde B_r:0\le r\le a).
\]
\end{lemma}

\begin{proof}
Conditionally on $\mathcal F_a^u$, the process
$(B_{a+u+r}-\widetilde B_{a+r})_{r\ge0}$ is a Brownian motion of
variance parameter $2$, started at $B_{a+u}-\widetilde B_a$.
Lemma~\ref{lem:Brownian-occupation} applies uniformly in this starting
point and gives, almost surely,
\begin{equation}\label{eq:Brownian-lag-conditional}
 E\left[\left|\int_a^{a+\ell}
       K_{H_1}^{[v]}(B_{r+u}-\widetilde B_r)\dd r\right|^p
       \Bigm|\mathcal F_a^u\right]^{\frac1p}
 \le C_h p^{1-H_1}\ell^{H_1}.
\end{equation}
Taking $a=0$, $\ell=T-u\le T$ and then expectations proves
\eqref{eq:Brownian-lag}. Exchanging the two paths gives the other
assertion.
\end{proof}

We now show that, as $\varepsilon,{\epsp}\downarrow0$, these
regularizations converge to the mutual intersection functional
\eqref{fr:alpha-continuous} of two independent Brownian motions,
\[
 \alpha_t(B,\widetilde B)
 =\int_{[0,t]^2}K_{H_0}(r-s)
                    K_{H_1}(B_r-\widetilde B_s)\dd r\dd s.
\]
At $H_1=\frac12$, where $K_{1/2}=\delta_0$, this definition is read
through the mutual occupation measure \eqref{eq:two-BM-occupation-measure},
in accordance with \eqref{eq:endpoint-reading}:
\begin{equation}\label{eq:alpha-endpoint}
\begin{split}
 \alpha_t(B,\widetilde B)
 &=\int_{[0,t]^2}K_{H_0}(r-s)\,
       \mathfrak I_{B,\widetilde B}(\dd r\dd s).
\end{split}
\end{equation}

As a preliminary to Proposition~\ref{prop:continuous-mutual}, we record
a cutoff estimate for the limiting functional $\alpha_t(B,\widetilde B)$; it will also
be used in the two-path invariance principle of
Section~\ref{sec:two-path}. For independent standard Brownian motions
$B,\widetilde B$ started at zero, a horizon $t>0$ and $0<\eta<t/2$, set
\begin{equation}\label{tp:alpha-cutoff}
 \alpha_t^{(\eta)}(B,\widetilde B)
 =\int_{[0,t]^2}K_{H_0}(\tau-\sigma)
       \chi\left(\frac{|\tau-\sigma|}{\eta}\right)
       K_{H_1}(B_\sigma-\widetilde B_\tau)\dd\sigma\dd\tau,
\end{equation}
where $\chi$ is the nondecreasing continuous cutoff of \eqref{eq:cutoff},
which vanishes on $[0,1]$ and equals one on $[2,\infty)$.

\begin{lemma}[Brownian mutual intersection and cutoff]\label{tp:Brownian-cutoff}
Let $B,\widetilde B$ be independent standard Brownian motions started
at zero.
For every horizon $t>0$, the functional $\alpha_t(B,\widetilde B)$ of
\eqref{fr:alpha-continuous} is finite almost surely and integrable,
with
\begin{equation}\label{tp:alpha-expectation}
 \Ep\left[\alpha_t(B,\widetilde B)\right]
 =\kappa_{H_1}\int_{[0,t]^2}K_{H_0}(\tau-\sigma)
                  (\sigma+\tau)^{H_1-1}\dd\sigma\dd\tau<\infty,
\end{equation}
and, for $0<\eta<\frac t2$,
\begin{equation}\label{tp:Brownian-cutoff-bound}
 0\le\Ep\left[\alpha_t(B,\widetilde B)-\alpha_t^{(\eta)}(B,\widetilde B)\right]
 \le C_{h}\,t^{H_1}\eta^{2H_0-1}.
\end{equation}
\end{lemma}

\begin{proof}
The mutual mean identity \eqref{tp:Brownian-mean-measure}, applied to
the nonnegative temporal weight, gives \eqref{tp:alpha-expectation}.
It applies also at $H_1=\frac12$ through the occupation-measure
interpretation in Lemma~\ref{lem:Brownian-mean-measure}.

To see that the mean is finite, substitute $u=|\tau-\sigma|$ and
$w=\sigma+\tau$. The two time triangles each have Jacobian
$\frac12$, and their common image is
$0\le u\le t$, $u\le w\le2t-u$. Consequently,
\begin{align*}
 \Ep[\alpha_t(B,\widetilde B)]
 &=c_{H_0}\kappa_{H_1}\int_0^{2t}w^{H_1-1}
       \left(\int_0^{\min(w,2t-w)}u^{2H_0-2}\dd u\right)\dd w\\
 &\le\frac{c_{H_0}\kappa_{H_1}}{2H_0-1}
             \int_0^{2t}w^{a_h-1}\dd w
 =\frac{c_{H_0}\kappa_{H_1}(2t)^{a_h}}
        {(2H_0-1)a_h}<\infty.
\end{align*}
Here the inequality enlarges the domain to $0\le u\le w\le2t$;
the integral is finite because $H_0>\frac12$ and $a_h>0$.

Since $\alpha_t(B,\widetilde B)\ge0$, a finite expectation also proves almost sure finiteness.

Finally, $\Ep\left[\alpha_t(B,\widetilde B)-\alpha_t^{(\eta)}(B,\widetilde B)\right]$ is bounded above by the same integral
restricted to $|\tau-\sigma|\le2\eta$. In the variables $(u,w)$, this integral is at most
\[
 C_{h}\int_0^{2\eta}u^{2H_0-2}
         \left(\int_u^{2t}w^{H_1-1}\dd w\right)\dd u
 \le\frac{C_{h}(2t)^{H_1}}{H_1}\int_0^{2\eta}u^{2H_0-2}\dd u
 \le C_{h}\,t^{H_1}\eta^{2H_0-1},
\]
which is \eqref{tp:Brownian-cutoff-bound}.
\end{proof}

The next statement proves the convergence and records the moment and
exponential bounds that $\alpha_t(B,\widetilde B)$ and its regularizations obey.

\begin{proposition}[Mixed mutual convergence and exponential bounds]
\label{prop:continuous-mutual}
Let $t>0$. For every $1\le p<\infty$,
\begin{equation}\label{eq:continuous-alpha-Lp}
 \alpha_t^{\varepsilon,{\epsp}}(B,\widetilde B)\longrightarrow\alpha_t(B,\widetilde B)
 \quad\text{in }L^p
 \quad\text{as }\varepsilon,{\epsp}\downarrow0,
\end{equation}
on the same pair $(B,\widetilde B)$. For $p\ge2$ and $\lambda\ge0$,
\begin{align}
 \sup_{\varepsilon,{\epsp}>0}\|\alpha_t^{\varepsilon,{\epsp}}(B,\widetilde B)\|_p
       +\|\alpha_t(B,\widetilde B)\|_p
   &\le C_{h}\sqrt p\,t^{a_{h}},\label{eq:continuous-alpha-bound}\\
 \sup_{\varepsilon,{\epsp}>0}
       E\left[e^{\lambda\alpha_t^{\varepsilon,{\epsp}}(B,\widetilde B)}\right]
   &\le e^{C_{h}(\lambda t^{a_{h}}+\lambda^2t^{2a_{h}})},\qquad
 E\left[e^{\lambda\alpha_t(B,\widetilde B)}\right]
   \le e^{C_{h}(\lambda t^{a_{h}}+\lambda^2t^{2a_{h}})}.
 \label{eq:cont-alpha-exp}
\end{align}
\end{proposition}

\begin{proof}
\textit{Step 1: uniform moment bounds.}
Split the time square into $\sigma<\tau$ and $\tau<\sigma$ and put
$u=|\tau-\sigma|$. Each resulting term in \eqref{eq:alpha-mixed} has
the form, possibly with $B$ and $\widetilde B$ interchanged,
\[
 \int_0^t K_{H_0}^{[\varepsilon^4+{\epsp}^4]}(u)
       \left(\int_0^{t-u}
         K_{H_1}^{[\varepsilon^2+{\epsp}^2]}(B_{r+u}-\widetilde B_r)
         \dd r\right)\dd u.
\]
Minkowski's inequality, Lemma~\ref{lem:Brownian-lag}, and the
envelope \eqref{eq:KH-envelope} of Lemma~\ref{lem:kernel-regularization},
applied with $H=H_0$ and $v=\varepsilon^4+{\epsp}^4$, bound its $L^p$
norm by
\begin{equation}\label{eq:alpha-lag-bound}
 C_{h}p^{1-H_1}t^{H_1}\int_0^tu^{2H_0-2}\dd u
 =C_{h}p^{1-H_1}t^{2H_0+H_1-1}
 =C_{h}p^{1-H_1}t^{a_{h}},
\end{equation}
which is finite because $H_0>\frac12$. Since $p^{1-H_1}\le\sqrt p$,
this proves the uniform regularized bound in
\eqref{eq:continuous-alpha-bound}. For $H_1>\frac12$, the same
computation with $v=0$ and the exact temporal kernel bounds $\alpha_t(B,\widetilde B)$
itself. At $H_1=\frac12$, the bound for $\alpha_t(B,\widetilde B)$ will follow from
the convergence argument below.

\smallskip
\noindent\textit{Step 2: convergence for $H_1>\frac12$.}
For
$(\sigma,\tau)\ne(0,0)$ the variable $B_\sigma-\widetilde B_\tau$ is
Gaussian with positive variance $\sigma+\tau$, hence almost surely
nonzero. By Fubini, almost surely the set of pairs where this
difference vanishes is Lebesgue-null in $[0,t]^2$. Outside this set
and the time diagonal, the two kernels converge:
$K_{H_0}^{[\varepsilon^4+{\epsp}^4]}(\tau-\sigma)\to K_{H_0}(\tau-\sigma)$
and
$K_{H_1}^{[\varepsilon^2+{\epsp}^2]}(B_\sigma-\widetilde B_\tau)
\to K_{H_1}(B_\sigma-\widetilde B_\tau)$, both by \eqref{eq:KH-convergence}
in Lemma~\ref{lem:kernel-regularization}, applied with $H=H_0$,
$v=\varepsilon^4+{\epsp}^4$ and with $H=H_1$, $v=\varepsilon^2+{\epsp}^2$.
By \eqref{eq:KH-envelope}, applied in the same way, the product is
dominated by
$C_h K_{H_0}(\tau-\sigma)K_{H_1}(B_\sigma-\widetilde B_\tau)$,
whose integral over $[0,t]^2$ is $C_h\alpha_t(B,\widetilde B)$ and is integrable by the
bound just proved. Dominated convergence therefore gives
\eqref{eq:continuous-alpha-Lp} almost surely and in $L^1$; the uniform
bound at the order $2p$ upgrades it to $L^p$, first for $p\ge2$ and
then for $1\le p<2$.

\smallskip
\noindent\textit{Step 3: convergence for $H_1=\frac12$.}
We first prove convergence in $L^1$ using a temporal
cutoff. Here $K_{H_1}^{[v]}=\heat_v$. For $0<\eta<\frac t2$, with
$\chi$ as in \eqref{eq:cutoff}, set
\[
 \alpha_t^{\varepsilon,{\epsp},(\eta)}(B,\widetilde B)
 =\int_{[0,t]^2}K_{H_0}^{[\varepsilon^4+{\epsp}^4]}(\tau-\sigma)
       \chi\left(\frac{|\tau-\sigma|}{\eta}\right)
       \heat_{\varepsilon^2+{\epsp}^2}(B_\sigma-\widetilde B_\tau)
                                       \dd\sigma\dd\tau,
\]
and recall $\alpha_t^{(\eta)}(B,\widetilde B)$ from \eqref{tp:alpha-cutoff}. With the
weight
\[
 w_\eta(\sigma,\tau)
 =K_{H_0}(\tau-\sigma)
     \chi\left(\frac{|\tau-\sigma|}{\eta}\right),
\]
extended by zero on the diagonal, the difference decomposes as
\[
 \alpha_t^{\varepsilon,{\epsp},(\eta)}(B,\widetilde B)-\alpha_t^{(\eta)}(B,\widetilde B)
 =D_1^{\varepsilon,{\epsp},(\eta)}+D_2^{\varepsilon,{\epsp},(\eta)},
\]
where
\begin{align*}
 D_1^{\varepsilon,{\epsp},(\eta)}
 &=\int_{[0,t]^2}\bigl(K_{H_0}^{[\varepsilon^4+{\epsp}^4]}
                   -K_{H_0}\bigr)(\tau-\sigma)
      \chi\left(\frac{|\tau-\sigma|}{\eta}\right)
      \heat_{\varepsilon^2+{\epsp}^2}(B_\sigma-\widetilde B_\tau)
                                       \dd\sigma\dd\tau,\\
 D_2^{\varepsilon,{\epsp},(\eta)}
 &=\int_{[0,t]^2}w_\eta(\sigma,\tau)
      \heat_{\varepsilon^2+{\epsp}^2}(B_\sigma-\widetilde B_\tau)
                                       \dd\sigma\dd\tau
   -\int_{[0,t]^2}w_\eta\dd\mathfrak I_{B,\widetilde B}.
\end{align*}
The cutoff vanishes for $|\tau-\sigma|\le\eta$, so that
\[
 |D_1^{\varepsilon,{\epsp},(\eta)}|
 \le d_{\varepsilon,{\epsp}}^{(\eta)}\,
     \mathcal T_t^{\varepsilon^2+{\epsp}^2}(B,\widetilde B),
 \qquad
 d_{\varepsilon,{\epsp}}^{(\eta)}
 =\sup_{|u|\ge\eta}\bigl|K_{H_0}^{[\varepsilon^4+{\epsp}^4]}(u)
                         -K_{H_0}(u)\bigr|,
\]
with $\mathcal T_t^v$ as in \eqref{eq:Brownian-weighted-functionals}.
Since $d_{\varepsilon,{\epsp}}^{(\eta)}\to0$ by \eqref{eq:KH-convergence},
applied with $H=H_0$ and $v=\varepsilon^4+{\epsp}^4$,
\eqref{eq:BM-unweighted-block} gives
$\|D_1^{\varepsilon,{\epsp},(\eta)}\|_p\le
C_hp^{1-H_1}t^{1+H_1}d_{\varepsilon,{\epsp}}^{(\eta)}\to0$. The weight $w_\eta$
is bounded by $K_{H_0}(\eta)$ and continuous on $[0,t]^2$, so that the
$H_1=\frac12$ case of \eqref{eq:BM-bounded-weight-convergence} gives
$\|D_2^{\varepsilon,{\epsp},(\eta)}\|_p\to0$. Thus,
for each fixed $\eta$, in every finite $L^p$,
\[
 \alpha_t^{\varepsilon,{\epsp},(\eta)}(B,\widetilde B)\longrightarrow\alpha_t^{(\eta)}(B,\widetilde B)
 \qquad(\varepsilon,{\epsp}\downarrow0).
\]
For the regularized integrals, the omitted strip has
$|\tau-\sigma|\le2\eta$. Restricting the $u$-integral in
\eqref{eq:alpha-lag-bound} to $[0,2\eta]$, the argument of Step~1
therefore gives, for $p\ge2$,
\[
 \sup_{\varepsilon,{\epsp}>0}
 \|\alpha_t^{\varepsilon,{\epsp}}(B,\widetilde B)
          -\alpha_t^{\varepsilon,{\epsp},(\eta)}(B,\widetilde B)\|_p
 \le C_h\sqrt p\,t^{1/2}\eta^{2H_0-1}.
\]
For the limiting intersection, defined through
$\mathfrak I_{B,\widetilde B}$, Lemma~\ref{tp:Brownian-cutoff} gives
\[
 \|\alpha_t(B,\widetilde B)-\alpha_t^{(\eta)}(B,\widetilde B)\|_1
 \le C_h t^{1/2}\eta^{2H_0-1}.
\]
Combining these two bounds with the fixed-cutoff convergence yields
\[
 \limsup_{\varepsilon,{\epsp}\downarrow0}
 \|\alpha_t^{\varepsilon,{\epsp}}(B,\widetilde B)-\alpha_t(B,\widetilde B)\|_1
 \le C_h t^{1/2}\eta^{2H_0-1}\longrightarrow0
 \qquad(\eta\downarrow0),
\]
because $H_0>\frac12$. Along a subsequence
$(\varepsilon_n,\epsp_n)$, the convergence also holds almost surely,
and Fatou's lemma gives, for every $p\ge2$,
\[
 \|\alpha_t(B,\widetilde B)\|_p
 \le\liminf_{n\to\infty}\|\alpha_t^{\varepsilon_n,\epsp_n}(B,\widetilde B)\|_p
 \le C_h\sqrt p\,t^{a_h},
\]
by the regularized bound of Step~1.
For each finite $p\ge2$, the uniform bound at order $2p$, together
with convergence in probability, gives convergence in $L^p$; for
$1\le p<2$, use the $L^2$ convergence. This proves both
\eqref{eq:continuous-alpha-Lp} and \eqref{eq:continuous-alpha-bound}
at $H_1=\frac12$.

\smallskip
\noindent\textit{Step 4: exponential moments.}
Apply the elementary moment-to-exponential estimate in
Corollary~\ref{tp:moment-to-exp} to the nonnegative variables just bounded,
after dividing by $t^{a_{h}}$. This proves \eqref{eq:cont-alpha-exp}.
\end{proof}

\begin{corollary}[The joint limit for two regularizations]
\label{lem:continuous-triple}
Let $t>0$. For every finite $p\ge1$,
\begin{equation}\label{eq:continuous-triple-Lp}
 \bigl(Y_t^\varepsilon(B),Y_t^{\epsp}(\widetilde B),
                    \alpha_t^{\varepsilon,{\epsp}}(B,\widetilde B)\bigr)
 \longrightarrow
 \bigl(Y_t(B),Y_t(\widetilde B),\alpha_t(B,\widetilde B)\bigr)
 \quad\text{in }L^p
\end{equation}
as $\varepsilon,{\epsp}\downarrow0$. For $b_1,b_2\in\R$ and $c\ge0$,
\begin{equation}\label{eq:continuous-mixed-exponential-limit}
 E\left[e^{b_1Y_t^\varepsilon(B)+b_2Y_t^{\epsp}(\widetilde B)
                       +c\alpha_t^{\varepsilon,{\epsp}}(B,\widetilde B)}\right]
 \longrightarrow
 E\left[e^{b_1Y_t(B)+b_2Y_t(\widetilde B)+c\alpha_t(B,\widetilde B)}\right].
\end{equation}
\end{corollary}

\begin{proof}
Theorem~\ref{thm:Brownian-self-intersection} and
Proposition~\ref{prop:continuous-mutual} give coordinatewise $L^p$
convergence on the same Brownian space, at the horizon $t$.

For any $r>1$, H\"older's inequality gives
\begin{align*}
 &E\left[\exp\left(r[b_1Y_t^\varepsilon(B)
                    +b_2Y_t^{\epsp}(\widetilde B)
                    +c\alpha_t^{\varepsilon,{\epsp}}(B,\widetilde B)]\right)\right]\\
 &\quad\le
 \left(E\left[e^{3rb_1Y_t^\varepsilon(B)}\right]\right)^{\frac13}
 \left(E\left[e^{3rb_2Y_t^{\epsp}(\widetilde B)}\right]\right)^{\frac13}
 \left(E\left[e^{3rc\alpha_t^{\varepsilon,{\epsp}}(B,\widetilde B)}\right]\right)^{\frac13}.
\end{align*}
All factors are uniformly bounded by \eqref{eq:Brownian-Y-exp-unif} and
\eqref{eq:cont-alpha-exp}. Thus the exponentials are uniformly
integrable. Their convergence in probability proves
\eqref{eq:continuous-mixed-exponential-limit}.
\end{proof}

\subsection{Pointwise construction}
\label{subsec:every-point}

We now combine the Feynman-Kac representation with the path estimates
to construct $U_t(x)$ at every fixed point. Besides the convergence
and uniform bound for $R^{\varepsilon}_t$ in \eqref{eq:R-eps-limit},
we shall repeatedly use the following consequences of
Theorem~\ref{thm:Brownian-self-intersection} and
Proposition~\ref{prop:continuous-mutual}. For $p\ge2$ and $t>0$,
\begin{equation}\label{eq:time-path-moments}
 \sup_{\varepsilon>0}\|Y_t^\varepsilon(B)\|_p
 +\sup_{\varepsilon,{\epsp}>0}\|\alpha_t^{\varepsilon,{\epsp}}(B,\widetilde B)\|_p
 \le C_{h}\sqrt p\,t^{a_{h}}.
\end{equation}

Their exponential versions are
\begin{align}
 \sup_{\varepsilon>0}E\left[e^{\lambda Y_t^\varepsilon(B)}\right]
    &\le e^{C_{h}\lambda^2t^{2a_{h}}},\qquad \lambda\in\R,\label{eq:time-path-exponentials-0}\\
 \sup_{\varepsilon,{\epsp}>0}E\left[e^{\lambda\alpha_t^{\varepsilon,{\epsp}}(B,\widetilde B)}\right]
    &\le e^{C_{h}(\lambda t^{a_{h}}+\lambda^2t^{2a_{h}})},
                   \qquad \lambda\ge0.
 \label{eq:time-path-exponentials}
\end{align}

\begin{proposition}[Moments before taking a solution limit]
\label{prop:regularized-moments}
Let $t>0$ and $x\in\R$. For every integer $m\ge1$ and independent standard
Brownian motions $B^1,\ldots,B^m$, started at zero and independent of $W$,
\begin{equation}\label{eq:continuous-m-moment}
 \E\left[\bigl(U_t^\varepsilon(x)\bigr)^m\right]
 =e^{mR^{\varepsilon}_t}E\left[\exp\left(
       \frac12\sum_{i=1}^mY_t^\varepsilon(B^i)
       +\sum_{1\le i<j\le m}
              \alpha_t^{\varepsilon}(B^i,B^j)\right)\right],
\end{equation}
with $\alpha_t^\varepsilon(B,\widetilde B)=\alpha_t^{\varepsilon,\varepsilon}(B,\widetilde B)$ as in
\eqref{eq:alpha-eps-again}. The right-hand side does not depend on $x$.
Consequently, for every finite $T>0$ and every finite $p\ge1$,
\begin{equation}\label{eq:continuous-uniform-moments}
 \sup_{\varepsilon>0}\
 \sup_{0\le t\le T}\ \sup_{x\in\R}
       \E\left[|U_t^\varepsilon(x)|^p\right]<\infty.
\end{equation}
\end{proposition}

\begin{proof}
Tonelli's theorem applies to the nonnegative Feynman-Kac integrands of
\eqref{eq:U-eps-tx}. Gaussian expectation then gives
\[
 \E\left[\exp\left(\sum_iW(h^{B^i,\varepsilon}_{t,x})
                             -mt\,c_h^{\varepsilon}\right)\right]
 =\exp\left(\frac12
      \Big\|\sum_i h^{B^i,\varepsilon}_{t,x}\Big\|_{\Hh}^2
                   -mt\,c_h^{\varepsilon}\right).
\]
By bilinearity, the identity $\|h^{B,\varepsilon}_{t,x}\|_{\Hh}^2
=\Gamma_t^\varepsilon(B)$ and \eqref{eq:alpha-mixed} with
${\epsp}=\varepsilon$,
\[
 \Big\|\sum_{i=1}^m h^{B^i,\varepsilon}_{t,x}\Big\|_{\Hh}^2
 =\sum_{i=1}^m\Gamma_t^\varepsilon(B^i)
  +2\sum_{1\le i<j\le m}\alpha_t^\varepsilon(B^i,B^j).
\]
By \eqref{eq:Y-eps} and \eqref{eq:R-eps},
$\frac12\Gamma_t^\varepsilon(B^i)-t\,c_h^\varepsilon
=\frac12Y_t^\varepsilon(B^i)+R^{\varepsilon}_t$ for each $i$, so that
\[
 \frac12\Big\|\sum_{i=1}^m h^{B^i,\varepsilon}_{t,x}\Big\|_{\Hh}^2
      -mt\,c_h^{\varepsilon}
 =mR^{\varepsilon}_t+\frac12\sum_{i=1}^mY_t^\varepsilon(B^i)
  +\sum_{1\le i<j\le m}\alpha_t^\varepsilon(B^i,B^j).
\]
Taking the expectation over the paths proves
\eqref{eq:continuous-m-moment}, initially with values in
$[0,\infty]$. The spatial point cancels in every path difference, so
the right-hand side does not involve $x$.

The path expectation in \eqref{eq:continuous-m-moment} is that of a
product of the factors $e^{\frac12Y_t^\varepsilon(B^i)}$ and
$e^{\alpha_t^\varepsilon(B^i,B^j)}$. H\"older's inequality, followed by
the uniform exponential bounds
\eqref{eq:time-path-exponentials-0}-\eqref{eq:time-path-exponentials},
bounds it by a constant depending only on $h$, $m$ and $T$, for every
$t\le T$. The
deterministic factor is bounded by \eqref{eq:R-eps-limit}.
Thus the right-hand side of \eqref{eq:continuous-m-moment} is finite and bounded uniformly in $\varepsilon$, in
$x$ and in $t\le T$. For arbitrary $p\ge1$, choose an integer $m\ge p$
and use monotonicity of probability-space $L^p$ norms.
\end{proof}

\begin{proposition}[Pointwise construction]
\label{thm:continuous-construction}\label{prop:L2-cont}
Let $t>0$ and $x\in\R$. There is a random variable $U_t(x)$, belonging
to every finite $L^p(\Omega,\sigma(W),\Pp)$, such that
\begin{equation}\label{eq:L2-cont}
 U_t^\varepsilon(x)\longrightarrow U_t(x)
 \quad\text{in }L^p(\Omega)
 \quad\text{as }\varepsilon\downarrow0,\qquad 1\le p<\infty.
\end{equation}
It is nonnegative, and its first two moments are
\begin{equation}\label{eq:field-first-second}
 \E\left[U_t(x)\right]=m_t:=e^{R_t}E\left[e^{\frac12Y_t(B)}\right]>0,
 \qquad
 \E\left[U_t(x)^2\right]=M_t,
\end{equation}
where, for two independent standard Brownian motions $B,\widetilde B$,
\begin{equation}\label{eq:Mt}
 M_t=e^{2R_t}E\left[\exp\left(
        \frac12\bigl(Y_t(B)+Y_t(\widetilde B)\bigr)
                       +\alpha_t(B,\widetilde B)\right)\right].
\end{equation}
The law of $U_t(x)$ does not depend on $x$. Set $U_0(x)=1$. Then, for every
finite $T>0$ and every finite $p\ge1$,
\begin{equation}\label{eq:continuous-limit-moments}
 \sup_{0\le t\le T}\ \sup_{x\in\R}
        \E\left[|U_t(x)|^p\right]<\infty.
\end{equation}
\end{proposition}

\begin{proof}
All the regularized solutions are defined from the same $W$. Fix
$(t,x)$. For $\varepsilon,{\epsp}>0$, Gaussian expectation and
Tonelli's theorem, exactly as in the previous proof with $m=2$ and two
different regularization parameters, give the mixed identity
\begin{equation}\label{eq:continuous-mixed-second}
 \E\left[U_t^\varepsilon(x)U_t^{\epsp}(x)\right]
 =e^{R^{\varepsilon}_t+R^{\epsp}_t}
 E\left[\exp\left(
   \frac12\bigl(Y_t^\varepsilon(B)+Y_t^{\epsp}(\widetilde B)\bigr)
      +\alpha_t^{\varepsilon,{\epsp}}(B,\widetilde B)\right)\right],
\end{equation}
whose right-hand side again does not involve $x$. By
Corollary~\ref{lem:continuous-triple} applied at the horizon $t$, and
by $R^{\varepsilon}_t\to R_t$ from \eqref{eq:R-eps-limit}, that
right-hand side converges, jointly as $\varepsilon,{\epsp}\downarrow0$, to the finite number $M_t$ of \eqref{eq:Mt}.
Consequently,
\begin{equation}\label{eq:continuous-Cauchy}
 \E\left[|U_t^\varepsilon(x)-U_t^{\epsp}(x)|^2\right]
 =\E\left[(U_t^\varepsilon(x))^2\right]
   +\E\left[(U_t^{\epsp}(x))^2\right]
   -2\E\left[U_t^\varepsilon(x)U_t^{\epsp}(x)\right]
 \longrightarrow M_t+M_t-2M_t=0.
\end{equation}
Completeness of $L^2(\Omega,\sigma(W),\Pp)$ gives a limit $U_t(x)$: one
may first take $\varepsilon_n=2^{-n}$, then use
\eqref{eq:continuous-Cauchy} to compare any sufficiently small
$\varepsilon$ with that sequence. This proves the full $L^2$ limit and
$\E[U_t(x)^2]=M_t$.

An almost surely convergent subsequence shows that $U_t(x)\ge0$ and,
by Fatou's lemma and \eqref{eq:continuous-uniform-moments}, that
$U_t(x)$ has every finite moment, with the uniform bound
\eqref{eq:continuous-limit-moments}. For $p>2$, choose $r>p$ and
$\theta\in(0,1)$ with $\frac1p=\frac\theta2+\frac{1-\theta}r$.
Interpolation gives
\[
 \|U_t^\varepsilon(x)-U_t(x)\|_p
 \le\|U_t^\varepsilon(x)-U_t(x)\|_2^\theta\,
       \|U_t^\varepsilon(x)-U_t(x)\|_r^{1-\theta}\longrightarrow0,
\]
since the second factor is bounded; for $1\le p\le2$, use the $L^2$
convergence directly. This is \eqref{eq:L2-cont}.

Taking $m=1$ in \eqref{eq:continuous-m-moment} gives the first
identity in \eqref{eq:field-first-second}: the left-hand side converges
by $L^1$ convergence, and the right-hand side converges by
\eqref{eq:R-eps-limit} and Corollary~\ref{lem:continuous-triple},
applied with $b_1=\frac12$ and $b_2=c=0$.

To prove spatial invariance in law,
use the stationarity of the Gaussian field $\dot W^\varepsilon$.
For every $b\in\R$, the shifted field
$(\dot W^\varepsilon(s,y+b))_{s,y}$ has the same law as
$(\dot W^\varepsilon(s,y))_{s,y}$. The Feynman-Kac representation
\eqref{eq:U-eps-tx} therefore gives
\[
 U_t^\varepsilon(x+b)\overset{\rm law}=U_t^\varepsilon(x).
\]
Passing to the limit in distribution at both points, using
\eqref{eq:L2-cont}, shows that the law of $U_t(x)$ does not depend on $x$.
\end{proof}

\subsection{The two path identities for the continuous field}
\label{subsec:two-continuous-identities}

 We now establish the two identities of Theorem~\ref{thm:continuous-summary}(3); their discrete counterparts are derived in Section~\ref{sec:discrete} and compared with them in Sections~\ref{subsec:legall-complet} and~\ref{sec:two-path}.

The second-moment formula is already part
of Proposition~\ref{thm:continuous-construction}. To establish the
$\mathcal S$-transform identity, we first control the regularization
of the linear occupation functional $\ell^B_{t,x}$ from
\eqref{eq:ell-B-general}, with $G_\varphi$ as in \eqref{eq:G-phi}.

\begin{lemma}[Regularized Brownian linear occupations]
\label{lem:continuous-linear-occupation}
For $\varphi\in\D$, $t>0$, $x\in\R$ and $\varepsilon>0$, define
\begin{equation}\label{eq:ell-B-regularized}
 \ell^{B,\varepsilon}_{t,x}(\varphi)
 :=\ip{h^{B,\varepsilon}_{t,x}}{\varphi}_{\Hh}
 =\int_0^{t}(\rho_\varepsilon*G_\varphi)(s,x+B_{t-s})\dd s.
\end{equation}
Then, uniformly over continuous paths $B$,
\begin{equation}\label{eq:ell-continuous-uniform-bound}
 |\ell^{B,\varepsilon}_{t,x}(\varphi)|\le t\|G_\varphi\|_\infty,
 \qquad |\ell^B_{t,x}(\varphi)|\le t\|G_\varphi\|_\infty,
\end{equation}
and
\begin{equation}\label{eq:ell-limit}
 \sup_B|\ell^{B,\varepsilon}_{t,x}(\varphi)-\ell^B_{t,x}(\varphi)|
 \le t\sqrt{\frac2\pi}\left(
       \varepsilon^2\|\partial_tG_\varphi\|_\infty
       +\varepsilon\|\partial_xG_\varphi\|_\infty\right)
 \longrightarrow0.
\end{equation}
\end{lemma}

\begin{proof}

Fubini in \eqref{eq:H-inner}, followed by \eqref{eq:hB-tx}, gives
\begin{align*}
 \ip{h^{B,\varepsilon}_{t,x}}{\varphi}_{\Hh}
 &=\int_{\R^2}h^{B,\varepsilon}_{t,x}(r,y)G_\varphi(r,y)\dd r\dd y\\
 &=\int_0^{t}\int_{\R^2}
     \rho_\varepsilon(s-r,x+B_{t-s}-y)G_\varphi(r,y)
                                    \dd r\dd y\dd s\\
 &=\int_0^{t}(\rho_\varepsilon*G_\varphi)(s,x+B_{t-s})\dd s.
\end{align*}

To justify the first use of Fubini, replace $\varphi$ by
$|\varphi|$. The corresponding covariance convolution
$G_{|\varphi|}$ is bounded, since $|\varphi|$ is bounded and compactly
supported. Moreover, $h^{B,\varepsilon}_{t,x}$ is nonnegative and has
Lebesgue integral $t$. Thus the absolute integral is at most
$t\|G_{|\varphi|}\|_\infty<\infty$. The subsequent interchange is
also legitimate because $G_\varphi$ is bounded and
$\rho_\varepsilon$ has mass one.
Since $\rho_\varepsilon$ is nonnegative,
$\|\rho_\varepsilon*G_\varphi\|_\infty
 \le\|G_\varphi\|_\infty$, which proves
\eqref{eq:ell-continuous-uniform-bound}.

The bounded derivatives of $G_\varphi$, established after
\eqref{eq:G-phi}, also give
\begin{align*}
 \|\rho_\varepsilon*G_\varphi-G_\varphi\|_\infty
 &\le\int_{\R^2}\rho_\varepsilon(r,y)
       \left(|r|\|\partial_tG_\varphi\|_\infty
             +|y|\|\partial_xG_\varphi\|_\infty\right)\dd r\dd y\\
 &=\sqrt{\frac2\pi}\left(
       \varepsilon^2\|\partial_tG_\varphi\|_\infty
       +\varepsilon\|\partial_xG_\varphi\|_\infty\right).
\end{align*}
Here we used
$\rho_\varepsilon(r,y)=\heat_{\varepsilon^4}(r)
                          \heat_{\varepsilon^2}(y)$ and
$\int_\R|z|\heat_v(z)\dd z=\sqrt{\frac{2v}\pi}$.
Integrating along the path over a time interval of length $t$
proves \eqref{eq:ell-limit}.
\end{proof}

\begin{proposition}[The continuous $\mathcal S$-transform and second moment]\label{prop:S-cont}
For every $\varphi\in\D$, $t>0$ and $x\in\R$,
\begin{equation}\label{eq:S-U}
 \mathcal S_{U_{t}(x)}(\varphi)
 =e^{R_{t}} E_B\left[\exp\left(\frac12 Y_{t}(B)
                    +\ell^B_{t,x}(\varphi)\right)\right].
\end{equation}
Moreover, for two independent standard Brownian motions,
\begin{equation}\label{eq:continuous-second-tx}
 \E\left[U_t(x)^2\right]
 =e^{2R_t}E_{B,\widetilde B}\left[\exp\left(
   \frac12\bigl(Y_t(B)+Y_t(\widetilde B)\bigr)
      +\alpha_t(B,\widetilde B)\right)\right].
\end{equation}
\end{proposition}

\begin{proof}
Fix $t>0$ and $x\in\R$. For $\varepsilon>0$, the factorization
\eqref{eq:Gaussian-cont}, Tonelli's theorem and
\eqref{eq:Gaussian-pairing} give
\begin{equation}\label{eq:S-U-eps}
 \mathcal S_{U^\varepsilon_{t}(x)}(\varphi)
 =e^{R^{\varepsilon}_t}E_B\left[\exp\left(\frac12 Y^\varepsilon_{t}(B)
                       +\ell^{B,\varepsilon}_{t,x}(\varphi)\right)\right].
\end{equation}
By Corollary~\ref{lem:continuous-triple} and \eqref{eq:ell-limit},
the exponent converges in probability to
$\frac12 Y_{t}(B)+\ell^B_{t,x}(\varphi)$.
Its exponential has uniformly bounded $r$th moments for every
$r>1$, by \eqref{eq:ell-continuous-uniform-bound} and
\eqref{eq:time-path-exponentials-0}.
Its expectation therefore converges; use also $R^{\varepsilon}_t\to R_{t}$.
On the left, \eqref{eq:L2-cont} and Cauchy-Schwarz give convergence
to $\mathcal S_{U_{t}(x)}(\varphi)$, since $\cE(\varphi)\in L^2$.
This proves \eqref{eq:S-U}.
Finally, \eqref{eq:continuous-second-tx} follows from
\eqref{eq:field-first-second} and the definition \eqref{eq:Mt} of $M_t$.
\end{proof}

\subsection{Functional convergence on space-time compacts}
\label{subsec:continuous-increments}

To pass from the pointwise limits to a continuous field, we shall prove
the following uniform estimate for $q\ge2$, $T>0$, $s,t\in[0,T]$ and $x,y\in\R$:
\[
 \sup_{\varepsilon>0}
 \|U_t^\varepsilon(x)-U_s^\varepsilon(y)\|_q
 \le C_{q,h,T}\bigl(|t-s|^{H_0-\frac12}+|x-y|^{\vartheta_h}\bigr),
\]
where the spatial exponent supplied by the occupation estimates is
\begin{equation}\label{eq:continuous-spatial-exponent}
 \vartheta_h=
 \begin{cases}
  H_1-\frac12,&H_1>\frac12,\\[2pt]
  \frac14,&H_1=\frac12.
 \end{cases}
\end{equation}
In both cases $\vartheta_h>0$.
Proposition~\ref{prop:continuous-functional} proves the estimate by
Gaussian hypercontractivity and Brownian path increments. We retain
this spatial exponent in the increment bounds, and choose a common
parabolic exponent only when applying the compactness criterion.

We first record a general tool, also used for the discrete
approximations: a quantitative form of the
Kolmogorov-Chentsov criterion; see \cite[Section~1.4]{Kunita}.

\begin{lemma}[From pointwise to compact-uniform convergence]
\label{lem:compact-upgrade}
Let $\mathcal K=[0,T]\times[-L,L]$, and let $(X_\iota)$ be continuous random
fields on one probability space, indexed by a sequence or a
regularization parameter tending to its limit. Suppose that, for
some $0<\eta\le1$ and every finite $q\ge2$,
\begin{align}
 \sup_\iota\|X_\iota(0,0)\|_q&<\infty,\notag\\
 \|X_\iota(t,x)-X_\iota(s,y)\|_q
 &\le C_{q,\mathcal K}\bigl(|t-s|^{\eta/2}+|x-y|^\eta\bigr)
       \quad((t,x),(s,y)\in\mathcal K),
 \label{eq:compact-increment-criterion}
\end{align}
where $C_{q,\mathcal K}$ does not depend on $\iota$.

Assume that $X_\iota(z)\to X(z)$ in every finite $L^p$ at each
fixed $z\in\mathcal K$. Then $X$ has a continuous modification, still denoted by $X$, and
\[
 \|X_\iota-X\|_{L^p(\Omega;C(\mathcal K))}\longrightarrow0
 \quad(1\le p<\infty).
\]
For each $p$ there is $\gamma>0$ such that
\begin{equation}\label{eq:compact-holder-control}
 \sup_\iota\E[\|X_\iota\|_{C^\gamma(\mathcal K)}^p]
       +\E[\|X\|_{C^\gamma(\mathcal K)}^p]<\infty,
\end{equation}
where $C^\gamma$ denotes the usual H\"older space.
\end{lemma}

The argument is relatively standard: we include in
Appendix~\ref{app:compact-upgrade} a short proof based on the
Garsia-Rodemich-Rumsey lemma.

At $H_1=\frac12$, the occupation increment estimate \eqref{eq:BM-occupation-increment} does not provide a positive power
of the spatial increment. The next lemma supplies that power through
the spatial regularity of Brownian local times. It will be used in
Lemma~\ref{lem:continuous-path-increments} below.

\begin{lemma}[Spatial increments of Brownian local times]
\label{lem:Brownian-local-time-increments}
Let $B$ be a Brownian motion of variance parameter $\sigma^2\in\{1,2\}$,
started at any deterministic point, and let $L^B$ be its jointly
continuous occupation-density local time, normalized as in
\eqref{temps-local}. For $\tau>0$, $q\ge2$ and $z,z'\in\R$,
\begin{equation}\label{eq:Brownian-local-time-increments}
 \|L^B_\tau(z)-L^B_\tau(z')\|_q
 \le C_q|z-z'|^{1/2}\tau^{1/4}.
\end{equation}
The constant is uniform in the starting point. Consequently, for every nonnegative $\kappa$ with $\|\kappa\|_1=1$,
\begin{equation}\label{eq:smoothed-local-time-increments}
 \left\|\int_0^\tau
       [\kappa(B_r+z)-\kappa(B_r+z')]\dd r\right\|_q
 \le C_q|z-z'|^{1/2}\tau^{1/4}.
\end{equation}
\end{lemma}

\begin{proof}
Put $d=|z-z'|$ and $I=[z\wedge z',z\vee z']$. Since
$\dd\langle B\rangle_r=\sigma^2\dd r$, the semimartingale local time of
$B$ at $a$ is $\sigma^2L^B_\tau(a)$, and Tanaka's formula \cite[Theorem~7.33]{MP-Brownian} reads
\[
 \sigma^2L^B_\tau(a)
 =|B_\tau-a|-|B_0-a|-\int_0^\tau\operatorname{sgn}(B_r-a)\dd B_r,
 \qquad a\in\R.
\]
Subtracting the formulas at $a=z$ and $a=z'$
gives
\begin{equation}\label{eq:Tanaka-difference}
 \sigma^2\bigl(L^B_\tau(z)-L^B_\tau(z')\bigr)
 =\bigl(|B_\tau-z|-|B_\tau-z'|\bigr)
  -\bigl(|B_0-z|-|B_0-z'|\bigr)-M_\tau,
\end{equation}
where
\[
 M_\tau=\int_0^\tau
   \bigl(\operatorname{sgn}(B_r-z)-\operatorname{sgn}(B_r-z')\bigr)\dd B_r.
\]
By the triangle inequality, each of the first two brackets in
\eqref{eq:Tanaka-difference} is bounded by $d$ in absolute value. The
integrand of $M$ is bounded by $2$ and vanishes when $B_r\notin I$, so
that, by the occupation formula \eqref{temps-local},
\[
 \langle M\rangle_\tau
 \le4\sigma^2\int_0^\tau\one_{\{B_r\in I\}}\dd r
 =4\sigma^2\int_IL^B_\tau(x)\dd x.
\]
Minkowski's inequality and Lemma~\ref{lem:Brownian-local-time}, applied
with $p=\max\{2,\frac q2\}$, bound the $L^{q/2}$ norm of the right-hand
side by $C_qd\sqrt\tau$, uniformly in the starting point. The
Burkholder-Davis-Gundy inequality therefore gives
$\|M_\tau\|_q\le C_q\|\langle M\rangle_\tau\|_{q/2}^{1/2}
\le C_qd^{1/2}\tau^{1/4}$, and \eqref{eq:Tanaka-difference}, with
$\sigma^2\ge1$, yields
\[
 \|L^B_\tau(z)-L^B_\tau(z')\|_q
 \le C_q(d+d^{1/2}\tau^{1/4}).
\]
For $d\le\sqrt{\tau}$ the first term is absorbed by the second.
For $d>\sqrt{\tau}$, the bound
$\|L^B_\tau(z)\|_q+\|L^B_\tau(z')\|_q\le C_q\sqrt{\tau}$
of Lemma~\ref{lem:Brownian-local-time} gives
\eqref{eq:Brownian-local-time-increments} directly.
Finally, the occupation formula rewrites the integral in
\eqref{eq:smoothed-local-time-increments} as
$\int\kappa(v)[L^B_\tau(v-z)-L^B_\tau(v-z')]\dd v$.
Minkowski's inequality proves the assertion.
\end{proof}

\begin{lemma}[Self-intersection and remainder increments]
\label{lem:continuous-self-increments}
For every $q\ge2$ and $s,t\ge0$, uniformly in $\varepsilon>0$,
\begin{equation}\label{eq:continuous-self-time-increment}
 \|Y_t^\varepsilon(B)-Y_s^\varepsilon(B)\|_q
       +|R^{\varepsilon}_t-R^{\varepsilon}_s|
 \le C_{q,h}|t-s|^{a_h}.
\end{equation}
\end{lemma}

\begin{proof}
Fix an arbitrary finite horizon $T>0$; the final bound will be
independent of $T$.
Throughout the proof, we abbreviate
$\mathsf K_0=K_{H_0}^{[2\varepsilon^4]}$ and
$\mathsf K_1=K_{H_1}^{[2\varepsilon^2]}$. Both kernels are even and
nonnegative, and \eqref{eq:KH-envelope}, applied with $H=H_0$, gives
\begin{equation}\label{eq:K0-envelope-proof}
 \mathsf K_0(u)\le C_h|u|^{2H_0-2}\quad(u\ne0),
 \qquad
 \int_{-T}^T\mathsf K_0(a-u)\dd u\le C_hT^{2H_0-1}\quad(|a|\le T),
\end{equation}
the second bound because $\int_{-2T}^{2T}|w|^{2H_0-2}\dd w
=C_{H_0}T^{2H_0-1}$.

Let $0\le s<t\le T$ and put $\delta=t-s$. Split $[0,t]^2$ into
$[0,s]^2$, $[s,t]^2$ and the two off-diagonal rectangles, which
contribute equally because the integrand of \eqref{eq:Gamma-eps} is
symmetric in $(\sigma,\tau)$:
\begin{equation}\label{eq:Gamma-splitting}
 \Gamma^\varepsilon_t(B)
 =\Gamma^\varepsilon_s(B)+\Gamma^\varepsilon_{[s,t]}(B)
  +2\,\mathcal C^\varepsilon_{s,\delta}(B),
\end{equation}
where
\begin{align*}
 \Gamma^\varepsilon_{[s,t]}(B)
 &=\int_{[s,t]^2}\mathsf K_0(\tau-\sigma)
       \mathsf K_1(B_\sigma-B_\tau)\dd\sigma\dd\tau,\\
 \mathcal C^\varepsilon_{s,\delta}(B)
 &=\int_0^s\int_s^t\mathsf K_0(\tau-\sigma)
       \mathsf K_1(B_\sigma-B_\tau)\dd\tau\dd\sigma.
\end{align*}
Centering \eqref{eq:Gamma-splitting} gives
\begin{equation}\label{eq:Y-splitting}
 Y^\varepsilon_t(B)-Y^\varepsilon_s(B)
 =\bigl(\Gamma^\varepsilon_{[s,t]}-E[\Gamma^\varepsilon_{[s,t]}]\bigr)
  +2\bigl(\mathcal C^\varepsilon_{s,\delta}
          -E[\mathcal C^\varepsilon_{s,\delta}]\bigr).
\end{equation}
The functional $\Gamma^\varepsilon_{[s,t]}(B)$ depends only on the
increments of $B$ on $[s,t]$, and $(B_{s+r}-B_s)_{r\ge0}$ is a standard
Brownian motion. Hence the first term of \eqref{eq:Y-splitting} has
the law of $Y^\varepsilon_\delta(B)$, and \eqref{eq:time-path-moments}
gives
\[
 \bigl\|\Gamma^\varepsilon_{[s,t]}-E[\Gamma^\varepsilon_{[s,t]}]\bigr\|_q
 =\|Y^\varepsilon_\delta(B)\|_q
 \le C_h\sqrt q\,\delta^{a_h}.
\]

For the cross term, set $B^-_r=B_{s-r}-B_s$ for $0\le r\le s$, and
$B^+_v=B_{s+v}-B_s$ for $v\ge0$. These are independent standard
Brownian motions, and we extend $B^-$ beyond time $s$ by an independent
Brownian motion. The changes of variables $\sigma=s-r$ and $\tau=s+v$
give $\tau-\sigma=r+v$ and $B_\sigma-B_\tau=B^-_r-B^+_v$. Since the
integrand is nonnegative,
\begin{equation}\label{eq:cross-term}
 \mathcal C^\varepsilon_{s,\delta}
 =\int_0^s\int_0^\delta\mathsf K_0(r+v)
       \mathsf K_1(B^-_r-B^+_v)\dd v\dd r
 \le\sum_{k\ge-1}\mathcal C^\varepsilon_k,
\end{equation}
where
\[
 \mathcal C^\varepsilon_k=\int_{I_k}\int_0^\delta\mathsf K_0(r+v)
       \mathsf K_1(B^-_r-B^+_v)\dd v\dd r,
\]
with $I_{-1}=[0,\delta]$ and $I_k=[2^k\delta,2^{k+1}\delta]$ for
$k\ge0$. On $I_{-1}$, \eqref{eq:K0-envelope-proof} gives
$\mathsf K_0(r+v)\le C_hK_{H_0}(r+v)$, so that, in the notation
\eqref{eq:Brownian-weighted-functionals} and by
\eqref{eq:BM-weighted-block},
\[
 \mathcal C^\varepsilon_{-1}\le C_h\,\mathcal Q^{2\varepsilon^2}_\delta(B^-,B^+),
 \qquad
 \|\mathcal C^\varepsilon_{-1}\|_q\le C_hq^{1-H_1}\delta^{a_h}
 \le C_h\sqrt q\,\delta^{a_h}.
\]
For $k\ge0$, put $R=2^k\delta$. On $I_k$ one has $r+v\ge R$, so that
\eqref{eq:K0-envelope-proof} gives
\[
 \mathcal C^\varepsilon_k
 \le C_hR^{2H_0-2}\int_0^\delta
     \left(\int_R^{2R}\mathsf K_1(B^-_r-B^+_v)\dd r\right)\dd v.
\]
Conditionally on $(B^-_r)_{r\le R}$ and on $B^+$, the process
$(B^-_{R+r})_{r\ge0}$ is a Brownian motion started at $B^-_R$. Since
\eqref{eq:BM-occupation} is uniform in the starting point and in the
target point, here $B^+_v$, the inner integral has $L^q$ norm at most
$C_hq^{1-H_1}R^{H_1}$. Minkowski's inequality in $v$ therefore gives
\[
 \|\mathcal C^\varepsilon_k\|_q
 \le C_h\sqrt q\,\delta R^{2H_0+H_1-2}
 =C_h\sqrt q\,\delta^{a_h}2^{-k(1-a_h)}.
\]
Since $a_h<1$, summing over $k\ge-1$ yields
$\|\mathcal C^\varepsilon_{s,\delta}\|_q\le C_h\sqrt q\,\delta^{a_h}$.
Centering costs at most a factor two, and \eqref{eq:Y-splitting}
gives
\[
 \|Y^\varepsilon_t(B)-Y^\varepsilon_s(B)\|_q\le C_h\sqrt q\,(t-s)^{a_h}.
\]

For the deterministic term, differentiating \eqref{eq:R-eps-exact} and
using \eqref{eq:f-envelope} give, for $u>0$,
\[
 \frac{\dd}{\dd u}R^{\varepsilon}_u
 =-\int_u^\infty f_h^{[\varepsilon]}(v)\dd v,
 \qquad
 \left|\frac{\dd}{\dd u}R^{\varepsilon}_u\right|
 \le C_h\int_u^\infty v^{a_h-2}\dd v
 =\frac{C_h}{1-a_h}\,u^{a_h-1}.
\]
Hence
\[
 |R^{\varepsilon}_t-R^{\varepsilon}_s|
 \le\frac{C_h}{a_h(1-a_h)}\bigl(t^{a_h}-s^{a_h}\bigr)
 \le C_h(t-s)^{a_h}.
\]
This completes the proof of \eqref{eq:continuous-self-time-increment}.
\end{proof}

For two independent standard Brownian motions $B,\widetilde B$ on $[0,T]$, started at zero, define
\begin{equation}\label{eq:continuous-two-point-intersection}
 A^\varepsilon_{t,s}(z)
 =\int_0^t\int_0^s
   K_{H_0}^{[2\varepsilon^4]}(t-s-r+v)
   K_{H_1}^{[2\varepsilon^2]}(z+B_r-\widetilde B_v)\dd v\dd r.
\end{equation}
Reversing the physical time variables in the path functions gives
\[
 A^\varepsilon_{t,s}(x-y)
   =\ip{h^{B,\varepsilon}_{t,x}}
            {h^{\widetilde B,\varepsilon}_{s,y}}_{\Hh}.
\]
In particular $A^\varepsilon_{t,t}(0)=\alpha_t^{\varepsilon,\varepsilon}(B,\widetilde B)$.

\begin{lemma}[Mutual-intersection moments and increments]
\label{lem:continuous-path-increments}
For $q\ge2$ and $s,s',t,t'\in[0,T]$, uniformly in $\varepsilon>0$,
\begin{align}
 \sup_{z\in\R}\|A^\varepsilon_{t,s}(z)\|_q
     &\le C_{h,T}\sqrt q, \label{eq:continuous-two-point-moment}\\
 \|A^\varepsilon_{t,s}(z)-A^\varepsilon_{t,s}(z')\|_q
     &\le C_{q,h,T}|z-z'|^{2\vartheta_h}, \label{eq:continuous-mutual-space}\\
 \sup_{z\in\R}
 \|A^\varepsilon_{t,s}(z)-A^\varepsilon_{t',s}(z)\|_q
     &\le C_{q,h,T}|t-t'|^{2H_0-1}, \label{eq:continuous-mutual-time}\\
 \sup_{z\in\R}
 \|A^\varepsilon_{t,s}(z)-A^\varepsilon_{t,s'}(z)\|_q
     &\le C_{q,h,T}|s-s'|^{2H_0-1}. \label{eq:continuous-mutual-time-bis}
\end{align}
Here $\vartheta_h$ is given by \eqref{eq:continuous-spatial-exponent}.
Moreover, for every $\lambda\ge0$,
\begin{equation}\label{eq:continuous-two-point-exponential}
 \sup_{\varepsilon>0}\,\sup_{s,t\in[0,T]}\,\sup_{z\in\R}
 E\Bigl[e^{\lambda A^\varepsilon_{t,s}(z)}\Bigr]
 \le e^{C_{h,T}(\lambda+\lambda^2)}.
\end{equation}
\end{lemma}

\begin{proof}
We again abbreviate
$\mathsf K_0=K_{H_0}^{[2\varepsilon^4]}$ and
$\mathsf K_1=K_{H_1}^{[2\varepsilon^2]}$, and use
\eqref{eq:K0-envelope-proof}.

\smallskip
\noindent\textit{Step 1: moments and spatial increments.}
Fix $s,t\in[0,T]$. In \eqref{eq:continuous-two-point-intersection},
use the variable $u=r-v$. For fixed $u$, the variable $v$ runs over
the interval
\[
 J_{u;t,s}=[\max(0,-u),\min(s,t-u)],
\]
which is empty unless $-s\le u\le t$. With $D^u_v=B_{v+u}-\widetilde B_v$,
\begin{equation}\label{eq:A-lag}
 A^\varepsilon_{t,s}(z)
 =\int_{-s}^{t}\mathsf K_0(t-s-u)\,G_u(z)\dd u,
 \qquad
 G_u(z)=\int_{J_{u;t,s}}\mathsf K_1(z+D^u_v)\dd v.
\end{equation}
Let $v_u$ be the left endpoint of $J_{u;t,s}$. Conditionally on
$(B_r)_{r\le v_u+u}$ and $(\widetilde B_v)_{v\le v_u}$, the process
$(D^u_{v_u+w})_{w\ge0}$ is a Brownian motion of variance parameter $2$
started at $D^u_{v_u}$. The bounds \eqref{eq:BM-occupation} and
\eqref{eq:BM-occupation-increment} are uniform in the starting point,
and $|J_{u;t,s}|\le T$. Hence, for every $u$ and $z,z'\in\R$,
\begin{align}
 \|G_u(z)\|_q&\le C_h\sqrt q\,T^{H_1},\label{eq:Gu-moment}\\
 \|G_u(z)-G_u(z')\|_q
 &\le C_{H_1}\sqrt q\,T^{\frac12}|z-z'|^{2H_1-1}
 \qquad\text{if }H_1>\tfrac12.\label{eq:Gu-increment}
\end{align}
For $H_1=\frac12$, the kernel $\mathsf K_1=\heat_{2\varepsilon^2}$ is
nonnegative with unit mass. Then
\eqref{eq:smoothed-local-time-increments} of
Lemma~\ref{lem:Brownian-local-time-increments}, applied to the same
conditional Brownian motion with $\kappa=\heat_{2\varepsilon^2}$,
replaces \eqref{eq:Gu-increment} by
\[
 \|G_u(z)-G_u(z')\|_q\le C_qT^{\frac14}|z-z'|^{\frac12}.
\]
By \eqref{eq:continuous-spatial-exponent}, the spatial powers in
the two regimes are respectively $2H_1-1=2\vartheta_h$ and
$\frac12=2\vartheta_h$. Thus, in both cases,
\begin{equation}\label{eq:Gu-increment-eta}
 \|G_u(z)-G_u(z')\|_q
 \le C_{q,h,T}|z-z'|^{2\vartheta_h}.
\end{equation}
Minkowski's inequality in \eqref{eq:A-lag}, together with \eqref{eq:Gu-moment} and
\eqref{eq:K0-envelope-proof}, gives
\[
 \|A^\varepsilon_{t,s}(z)\|_q
 \le\int_{-T}^T\mathsf K_0(t-s-u)\|G_u(z)\|_q\dd u
 \le C_h\sqrt q\,T^{2H_0-1+H_1}
 =C_h\sqrt q\,T^{a_h},
\]
which is \eqref{eq:continuous-two-point-moment}. The same argument,
with \eqref{eq:Gu-increment-eta} in place of \eqref{eq:Gu-moment},
proves \eqref{eq:continuous-mutual-space}.

\smallskip
\noindent\textit{Step 2: temporal increments and exponential moments.}
Let $0\le t<t'\le T$ and put $\delta=t'-t$. Splitting the integral in
$r$ at $t$ gives
$A^\varepsilon_{t',s}(z)-A^\varepsilon_{t,s}(z)=\Delta_1+\Delta_2$,
where
\begin{align*}
 \Delta_1&=\int_0^t\int_0^s
   \bigl[\mathsf K_0(t'-s-r+v)-\mathsf K_0(t-s-r+v)\bigr]
   \mathsf K_1(z+B_r-\widetilde B_v)\dd v\dd r,\\
 \Delta_2&=\int_t^{t'}\int_0^s\mathsf K_0(t'-s-r+v)
   \mathsf K_1(z+B_r-\widetilde B_v)\dd v\dd r.
\end{align*}
With the notation of Step~1,
\[
 \Delta_1=\int_{-s}^t
   \bigl[\mathsf K_0(t-s-u+\delta)-\mathsf K_0(t-s-u)\bigr]G_u(z)\dd u.
\]
The translation bound \eqref{eq:KH-translation}, applied with $H=H_0$
and regularization parameter $2\varepsilon^4$, gives
$\int_\R|\mathsf K_0(w+\delta)-\mathsf K_0(w)|\dd w
\le C_h\delta^{2H_0-1}$, so that Minkowski's inequality and
\eqref{eq:Gu-moment} yield
\[
 \|\Delta_1\|_q\le C_h\sqrt q\,T^{H_1}\delta^{2H_0-1}.
\]
Similarly,
\[
 \Delta_2=\int_{t-s}^{t'}\mathsf K_0(t'-s-u)\,G'_u(z)\dd u,
 \qquad
 G'_u(z)=\int_{J'_u}\mathsf K_1(z+D^u_v)\dd v,
\]
where $J'_u=\{v\in[0,s]:t\le v+u\le t'\}$ is an interval of length
at most $\delta$. Applied as in Step~1, \eqref{eq:BM-occupation} gives
$\|G'_u(z)\|_q\le C_h\sqrt q\,\delta^{H_1}$, and
\eqref{eq:K0-envelope-proof} then gives
\[
 \|\Delta_2\|_q\le C_h\sqrt q\,T^{2H_0-1}\delta^{H_1}.
\]
By \eqref{fr:range} and $H_1\ge\frac12$, one has
$2H_0-1<1-H_1\le H_1$. Since $\delta\le T$, it follows that
\[
 \delta^{H_1}=\delta^{2H_0-1}\delta^{H_1-2H_0+1}
 \le T^{H_1-2H_0+1}\delta^{2H_0-1}.
\]
The bounds on $\Delta_1$ and $\Delta_2$ are uniform in $z$, and prove
\eqref{eq:continuous-mutual-time}. Since $\mathsf K_0$ and
$\mathsf K_1$ are even, exchanging the two paths gives
\[
 A^\varepsilon_{t,s}(z)(B,\widetilde B)
 =A^\varepsilon_{s,t}(-z)(\widetilde B,B),
\]
so that the same bound holds for increments in the second horizon $s$, which proves \eqref{eq:continuous-mutual-time-bis}.

Finally, $A^\varepsilon_{t,s}(z)\ge0$ and, by
\eqref{eq:continuous-two-point-moment},
$\|A^\varepsilon_{t,s}(z)\|_q\le C_{h,T}\sqrt q$ for every $q\ge2$.
Corollary~\ref{tp:moment-to-exp}, applied with $c=C_{h,T}$, proves
\eqref{eq:continuous-two-point-exponential}.
\end{proof}

\begin{proposition}[Uniform space-time increment bounds]
\label{prop:continuous-functional}
For every $q\ge2$ and finite $T>0$,
\begin{equation}\label{eq:continuous-field-increments}
 \sup_{\varepsilon>0}
 \|U^\varepsilon_t(x)-U^\varepsilon_s(y)\|_q
 \le C_{q,h,T}\bigl(|t-s|^{H_0-\frac12}+|x-y|^{\vartheta_h}\bigr)
 \quad(s,t\in[0,T],\ x,y\in\R),
\end{equation}
with $\vartheta_h$ as in \eqref{eq:continuous-spatial-exponent}.
\end{proposition}

\begin{proof}

\smallskip

\noindent\textit{Step 1: reduction to a second moment.} We use two classical
properties of the Ornstein-Uhlenbeck semigroup $(P_\tau)_{\tau\ge0}$
on $L^2(\Omega,\sigma(W),\Pp)$, which acts as multiplication by
$e^{-n\tau}$ on the $n$th Wiener chaos. First, the $n$th chaos
component of the Gaussian exponential $\cE(h)$ of \eqref{eq:stoch-exp}
is the $n$th Wick power of $W(h)$ divided by $n!$
\cite[Chapter~III]{Janson-Gaussian}; it is homogeneous of degree $n$
in $h$. Hence
\begin{equation}\label{eq:OU-exponential}
 P_\tau\cE(h)=\cE(e^{-\tau}h),\qquad h\in\Hh,\ \tau\ge0.
\end{equation}
Second, by Nelson's hypercontractivity theorem
\cite[Chapter~V]{Janson-Gaussian}, for every
$G\in L^2(\Omega,\sigma(W),\Pp)$,
\begin{equation}\label{eq:Nelson}
 \|P_\tau G\|_q\le\|G\|_2
 \qquad\text{as soon as}\quad e^{-\tau}\le(q-1)^{-1/2}.
\end{equation}

\smallskip

For every $q\ge2$, define the amplified field
\[
 V^{\varepsilon,q}_t(x)
 =e^{R^{\varepsilon}_t}E_B\left[
       e^{\frac12Y^\varepsilon_t(B)}
       \cE\bigl(\sqrt{q-1}\,h^{B,\varepsilon}_{t,x}\bigr)\right].
\]
Since $P_\tau$ is a bounded linear operator on $L^2$, it commutes
with this path expectation, which is a Bochner integral in $L^2$. Taking $\tau_q=\frac12\log(q-1)$ in
\eqref{eq:OU-exponential} therefore gives
$P_{\tau_q}V^{\varepsilon,q}_t(x)=U^\varepsilon_t(x)$.
By \eqref{eq:Nelson},
\begin{equation}\label{eq:U-by-V}
 \|U^\varepsilon_t(x)-U^\varepsilon_s(y)\|_q
 =\bigl\|P_{\tau_q}\bigl(V^{\varepsilon,q}_t(x)-V^{\varepsilon,q}_s(y)\bigr)\bigr\|_q
 \le\|V^{\varepsilon,q}_t(x)-V^{\varepsilon,q}_s(y)\|_2.
\end{equation}

\smallskip
\noindent\textit{Step 2: the second moment.} Let $\widetilde B$ be an
independent copy of $B$. The integrands are nonnegative, so Tonelli's
theorem, \eqref{eq:Gaussian-pairing} and the identity
$\ip{h^{B,\varepsilon}_{t,x}}{h^{\widetilde B,\varepsilon}_{s,y}}_{\Hh}
=A^\varepsilon_{t,s}(x-y)$ following
\eqref{eq:continuous-two-point-intersection} give
\begin{equation}\label{eq:continuous-amplified-two-point}
 C^{\varepsilon}_q((t,x),(s,y))
 :=\E\bigl[V^{\varepsilon,q}_t(x)V^{\varepsilon,q}_s(y)\bigr]
 =E\left[e^{\Phi_q^\varepsilon((t,x),(s,y))}\right],
\end{equation}
where
\[
 \Phi_q^\varepsilon((t,x),(s,y))
 =R^{\varepsilon}_t+R^{\varepsilon}_s
  +\frac12Y^\varepsilon_t(B)+\frac12Y^\varepsilon_s(\widetilde B)
  +(q-1)A^\varepsilon_{t,s}(x-y).
\]
Only the Gaussian exponentials have been modified in $V^{\varepsilon,q}$:
this is why $q$ enters \eqref{eq:continuous-amplified-two-point} only
through the mutual-intersection term, whose exponential moments are
controlled by \eqref{eq:continuous-two-point-exponential}. Now
\begin{align*}
 \|V^{\varepsilon,q}_t(x)-V^{\varepsilon,q}_s(y)\|_2^2
 &=\bigl[C^\varepsilon_q((t,x),(t,x))-C^\varepsilon_q((t,x),(s,y))\bigr]\\
 &\quad+\bigl[C^\varepsilon_q((s,y),(s,y))-C^\varepsilon_q((t,x),(s,y))\bigr].
\end{align*}
We compare $C^\varepsilon_q((t,x),(t,x))$ with $C_q^\varepsilon((t,x),(s,y))$, using the
same Brownian paths in both expectations. The inequality
$|e^u-e^v|\le|u-v|(e^u+e^v)$ and the Cauchy-Schwarz inequality give
\begin{align*}
 &|C_q^\varepsilon((t,x),(t,x))-C_q^\varepsilon((t,x),(s,y))|\\
 &\le\bigl\|\Phi^\varepsilon_q((t,x),(t,x))-\Phi^\varepsilon_q((t,x),(s,y))\bigr\|_2\,
    \bigl\|e^{\Phi^\varepsilon_q((t,x),(t,x))}+e^{\Phi^\varepsilon_q((t,x),(s,y))}\bigr\|_2.
\end{align*}
By H\"older's inequality, \eqref{eq:R-eps-limit},
\eqref{eq:time-path-exponentials-0} and
\eqref{eq:continuous-two-point-exponential}, the second factor is
bounded by $C_{q,h,T}$, uniformly in $\varepsilon$, $(t,x)$ and $(s,y)$. In the
first factor,
\begin{align*}
 &|\Phi^\varepsilon_q((t,x),(t,x))-\Phi^\varepsilon_q((t,x),(s,y))|\\
 &\le|R^{\varepsilon}_t-R^{\varepsilon}_s|
   +\tfrac12|Y^\varepsilon_t(\widetilde B)-Y^\varepsilon_s(\widetilde B)|\\
 &\quad+(q-1)|A^\varepsilon_{t,t}(0)-A^\varepsilon_{t,s}(0)|
   +(q-1)|A^\varepsilon_{t,s}(0)-A^\varepsilon_{t,s}(x-y)|.
\end{align*}
By \eqref{eq:continuous-self-time-increment},
\eqref{eq:continuous-mutual-time-bis} and
\eqref{eq:continuous-mutual-space}, and since
$|t-s|^{a_h}\le T^{a_h-2H_0+1}|t-s|^{2H_0-1}$, this gives
\[
 |C^\varepsilon_q((t,x),(t,x))-C^\varepsilon_q((t,x),(s,y))|
 \le C_{q,h,T}\bigl(|t-s|^{2H_0-1}+|x-y|^{2\vartheta_h}\bigr).
\]
The same bound holds for $|C^\varepsilon_q((s,y),(s,y))-C^\varepsilon_q((t,x),(s,y))|$,
by symmetry, with \eqref{eq:continuous-mutual-time} in place of
\eqref{eq:continuous-mutual-time-bis}. Hence
$\|V^{\varepsilon,q}_t(x)-V^{\varepsilon,q}_s(y)\|_2
\le C_{q,h,T}(|t-s|^{H_0-\frac12}+|x-y|^{\vartheta_h})$, and \eqref{eq:U-by-V}
proves \eqref{eq:continuous-field-increments}.
\end{proof}

For the application of Lemma~\ref{lem:compact-upgrade}, choose the
common parabolic exponent
\begin{equation}\label{eq:increment-exponent-used}
 \eta=\min\{2H_0-1,\vartheta_h\}
 =\begin{cases}
  \eta_h=\min\{2H_0-1,H_1-\frac12\},&H_1>\frac12,\\[2pt]
  \eta_{H_0}=\min\{2H_0-1,\frac14\},&H_1=\frac12.
 \end{cases}
\end{equation}
This exponent is positive in both regimes and will also be used for
the discrete approximations.

\begin{theorem}[Construction of the continuous field and functional convergence]
\label{thm:continuous-field}
There is a continuous, nonnegative random field
$U=(U_t(x))_{t\ge0,x\in\R}$, measurable with respect to $W$, with
$U_0(x)=1$, whose value at each fixed $(t,x)$ is almost surely equal to the variable
constructed in Proposition~\ref{thm:continuous-construction}. For finite
$T,L>0$ and $p\ge1$, with
$\mathcal K_{T,L}=[0,T]\times[-L,L]$,
\begin{equation}\label{eq:field-functional-convergence}
 \E\left[\|U^\varepsilon-U\|_{C(\mathcal K_{T,L})}^p\right]
 \longrightarrow0
 \qquad(\varepsilon\downarrow0),
\end{equation}
so that $U$ and its regularizations belong to
$L^p(\Omega;C(\mathcal K_{T,L}))$.
The initial condition is attained uniformly on compacts:
\begin{equation}\label{eq:uniform-initial-continuous}
 \E\left[\sup_{|x|\le L}|U_t(x)-1|^p\right]\longrightarrow0
 \qquad(t\downarrow0).
\end{equation}
\end{theorem}

\begin{proof}

Fix $T,L>0$. Since $H_0-\frac12\ge\frac\eta2$ and
$\vartheta_h\ge\eta$, the increment estimate
\eqref{eq:continuous-field-increments}, restricted to
$\mathcal K_{T,L}$, implies \eqref{eq:compact-increment-criterion}
with the exponent \eqref{eq:increment-exponent-used} and a constant
depending on $T,L$. The regularized solutions are continuous,
$U_0^\varepsilon(0)=1$, and
Proposition~\ref{thm:continuous-construction} gives their pointwise
limits in every finite $L^p$. Lemma~\ref{lem:compact-upgrade}
therefore provides a continuous modification on $\mathcal K_{T,L}$
and \eqref{eq:field-functional-convergence}.

Continuous modifications obtained on the rectangles
$[0,k]\times[-k,k]$ agree almost surely on overlaps: they agree
on a countable dense set, and hence everywhere by continuity.
They can therefore be patched into one continuous, $W$-measurable
field on $[0,\infty)\times\R$. Its nonnegativity and its value $1$
on the initial line follow on a countable dense set from
Proposition~\ref{thm:continuous-construction}, and then everywhere by
continuity.

Finally, \eqref{eq:compact-holder-control} gives, for some $\gamma>0$,
$\|\sup_{|x|\le L}|U_t(x)-U_0(x)|\|_p\le C_{p,T,L}t^\gamma$, which
proves \eqref{eq:uniform-initial-continuous}.
\end{proof}

\subsection{Necessity of the renormalization}
\label{subsec:renormalization-necessity}

The continuous construction now proves the non-tightness assertion of
Proposition~\ref{prop:renormalization-needed}.
\begin{proof}[Proof of Proposition~\ref{prop:renormalization-needed}]
By \eqref{fr:counterterm-continuous} and
Proposition~\ref{prop:R-eps-limit},
$c_h^\varepsilon=\varepsilon^{2a_h-2}c_h^1$ with $0<c_h^1<\infty$;
since $a_h<1$, it tends
to $+\infty$ as $\varepsilon\downarrow0$. The proof of
Lemma~\ref{lem:regularized-FK}, in which the counterterm enters only as
a constant potential and may be replaced by $0$, shows that
\eqref{eq:unrenormalized-PDE} has a unique classical solution in the
growth class, given by
\[
 u^\varepsilon_t(x)
 =E_B\left[\exp\left(\int_0^t\dot W^\varepsilon(s,x+B_{t-s})\dd s\right)\right]
 =e^{tc_h^\varepsilon}U^\varepsilon_t(x),
\]
where $B$ is a standard Brownian motion independent of $W$, the second
identity following from \eqref{eq:U-eps-tx}. Fix $t>0$, $x\in\R$ and $K,\delta>0$. As soon
as $Ke^{-tc_h^\varepsilon}\le\delta$,
\[
 \Pp\bigl(u^\varepsilon_t(x)\le K\bigr)
 \le\Pp\bigl(U^\varepsilon_t(x)\le\delta\bigr).
\]
By Theorem~\ref{thm:continuous-summary}, $U^\varepsilon_t(x)\to U_t(x)$
in $L^1(\Omega)$, hence in law, and the Portmanteau theorem applied to
the closed set $(-\infty,\delta]$ gives
$\limsup_{\varepsilon\downarrow0}\Pp(u^\varepsilon_t(x)\le K)\le
\Pp(U_t(x)\le\delta)$. Letting $\delta\downarrow0$ and using
$U_t(x)\ge0$ yields the bound $\Pp(U_t(x)=0)$, which is smaller than
one because $\E[U_t(x)]>0$. Consequently,
$\liminf_{\varepsilon\downarrow0}\Pp(u^\varepsilon_t(x)>K)\ge
1-\Pp(U_t(x)=0)>0$ for every $K$, so that no sequence
$(u^{\varepsilon_n}_t(x))$ is tight. Since a sequence of real random
variables that converges in law is tight, none converges in law.
\end{proof}

\subsection{Comparison with the previously constructed solution}\label{subsec:compa-cdot}

The notion of solution is fixed in
Definition~\ref{def:continuous-solution}. We now compare this field
with the previously constructed solution.

For the comparison with \cite{CDOT-KRP}, we now identify the
spectral representation of the counterterm.

\begin{lemma}[Spectral representation of the counterterm]
\label{lem:counterterm-spectral}
Under \eqref{fr:range}, for every $\varepsilon>0$, the counterterm
\eqref{fr:counterterm-continuous} satisfies
\begin{equation}\label{eq:counterterm-spectral-match}
 c_h^\varepsilon
 =\int_{\R^2}
 \frac{|\widehat\rho(\varepsilon^2\lambda,\varepsilon\xi)|^2}
      {\frac12\xi^2+i\lambda}\,
          \mathfrak m_h(\dd\lambda,\dd\xi).
\end{equation}
Here
$|\widehat\rho(\varepsilon^2\lambda,\varepsilon\xi)|^2
=e^{-\varepsilon^4\lambda^2-\varepsilon^2\xi^2}$ and the Fourier
convention and spectral measure are those of
\eqref{eq:Fourier-convention}-\eqref{eq:spectral-measure}.
\end{lemma}

\begin{proof}
We first write the function $f_h^{[\varepsilon]}$ of \eqref{eq:f-eps}
as a Fourier integral. Let $\frac12<H<1$ and $v>0$. In the sense of
tempered distributions, the Fourier transform of $K_H$ is
$2\pi\mathfrak a_H|\lambda|^{1-2H}$; this is the computation behind
\eqref{eq:spectral-measure}. The Fourier transform of
$K_H^{[v]}=K_H*\heat_v$ is therefore the integrable function
$2\pi\mathfrak a_H|\lambda|^{1-2H}e^{-v\lambda^2/2}$. Since
$K_H^{[v]}$ is continuous, Fourier inversion holds at every point:
\begin{equation}\label{eq:KH-v-Fourier}
 K_H^{[v]}(u)
 =\mathfrak a_H\int_\R e^{iu\lambda}e^{-\frac v2\lambda^2}
          |\lambda|^{1-2H}\dd\lambda,
 \qquad u\in\R.
\end{equation}
This formula also holds for $H=\frac12$, where $K_{1/2}^{[v]}=\heat_v$
and $\mathfrak a_{1/2}=\frac1{2\pi}$. Apply \eqref{eq:KH-v-Fourier}
with $(H,v)=(H_0,2\varepsilon^4)$ and with $(H,v)=(H_1,2\varepsilon^2)$.
Since $e^{-\varepsilon^4\lambda^2-\varepsilon^2\xi^2}
=|\widehat\rho(\varepsilon^2\lambda,\varepsilon\xi)|^2$, this gives,
for all $u,y\in\R$,
\begin{equation}\label{eq:covariance-Fourier}
 K_{H_0}^{[2\varepsilon^4]}(u)K_{H_1}^{[2\varepsilon^2]}(y)
 =\int_{\R^2}e^{i(u\lambda+y\xi)}
      |\widehat\rho(\varepsilon^2\lambda,\varepsilon\xi)|^2
                     \mathfrak m_h(\dd\lambda,\dd\xi).
\end{equation}
By Step~1 of the proof of
Proposition~\ref{prop:R-eps-limit} and \eqref{eq:smoothed-noise-covariance},
\[
 f_h^{[\varepsilon]}(v)
 =E_B\Bigl[K_{H_0}^{[2\varepsilon^4]}(v)
           K_{H_1}^{[2\varepsilon^2]}(B_v)\Bigr],
\]
where $B$ is a standard Brownian motion. Insert
\eqref{eq:covariance-Fourier} with $u=v$ and $y=B_v$. The Brownian characteristic function
$E_B[e^{i\xi B_v}]=e^{-v\xi^2/2}$ and the invariance of $|\widehat\rho(\varepsilon^2\lambda,\varepsilon\xi)|^2\,\mathfrak m_h(\dd\lambda,\dd\xi)$ under $\lambda\mapsto-\lambda$ then give
\[
 f_h^{[\varepsilon]}(v)
 =\int_{\R^2}
      |\widehat\rho(\varepsilon^2\lambda,\varepsilon\xi)|^2
      e^{iv\lambda-\frac v2\xi^2}
                     \mathfrak m_h(\dd\lambda,\dd\xi)
 =\int_{\R^2}
      |\widehat\rho(\varepsilon^2\lambda,\varepsilon\xi)|^2
      e^{-v(\frac12\xi^2+i\lambda)}
                     \mathfrak m_h(\dd\lambda,\dd\xi).
\]

For $b>0$, the Gaussian factor and the
damping $e^{-bv}$ make Fubini's theorem applicable, so
\begin{equation}\label{eq:counterterm-Abel}
 \int_0^\infty e^{-bv}f_h^{[\varepsilon]}(v)\dd v
 =\int_{\R^2}
  \frac{|\widehat\rho(\varepsilon^2\lambda,\varepsilon\xi)|^2}
       {b+\frac12\xi^2+i\lambda}\,
                     \mathfrak m_h(\dd\lambda,\dd\xi).
\end{equation}

To pass to the limit as $b\downarrow0$, note that
\begin{align*}
 \int_{\R^2}
   \frac{e^{-\varepsilon^4\lambda^2-\varepsilon^2\xi^2}}
        {\sqrt{\xi^4/4+\lambda^2}}\,
           \mathfrak m_h(\dd\lambda,\dd\xi)
 &\le C_h\int_\R e^{-\varepsilon^2\xi^2}|\xi|^{1-2H_1}
       \left(\int_\R
          \frac{|\lambda|^{1-2H_0}}
               {\sqrt{\xi^4/4+\lambda^2}}\dd\lambda\right)\dd\xi\\
 &\le C_h\left(\int_0^1w^{1-2H_0}\dd w
                +\int_1^\infty w^{-2H_0}\dd w\right)
       \int_\R e^{-\varepsilon^2\xi^2}|\xi|^{1-2a_h}\dd\xi
 <\infty.
\end{align*}
The second inequality follows from the substitution $\lambda=\frac12\xi^2w$ in the inner integral.

Moreover, for $b>0$,
\[
 \left|\frac1{b+\frac12\xi^2+i\lambda}\right|
 \le\frac1{\sqrt{\xi^4/4+\lambda^2}}.
\]
Dominated convergence therefore applies to the right-hand side of
\eqref{eq:counterterm-Abel}. On its left-hand side, monotone
convergence and \eqref{eq:counterterm-time-domain} give
$\int_0^\infty e^{-bv}f_h^{[\varepsilon]}(v)\dd v\to c_h^\varepsilon$.
This proves \eqref{eq:counterterm-spectral-match}.
\end{proof}

\begin{proposition}[Agreement with the previously constructed solution]
\label{prop:CDOT-comparison}
If $2H_0+H_1>\frac53$, the solution above and the continuous version
of the solution constructed in \cite{CDOT-KRP}, for the same $W_h$
and the same counterterm, agree on one event of probability one
for all $t\ge0$ and $x\in\R$.
\end{proposition}

\begin{proof}
The stationary product Gaussian mollifier is that of
\cite[Definition~3.4]{CDOT-KRP}. By
\eqref{eq:counterterm-spectral-match}, our counterterm at
$\varepsilon=2^{-n}$ agrees with \cite[(3.12)]{CDOT-KRP}, also
recorded in \cite[(4.7) and (4.12)]{CDOT-moments}.

For this parameter range,
\cite[Theorems~2.11 and~3.3, Definition~3.4]{CDOT-KRP} give almost
sure convergence of the classical approximations in a weighted
supremum norm on each finite time interval, and hence uniformly on
space-time compacts. We choose the resulting continuous version
$\widetilde U$ of that solution. In particular, 
\begin{equation}\label{eq:continuous-CDOT-input}
 U^{2^{-n}}_t(x)\longrightarrow \widetilde U_t(x)
 \quad\text{almost surely at every fixed }(t,x).
\end{equation}
Theorem~\ref{thm:continuous-field} gives convergence of the same
sequence to $U_t(x)$ in probability. Uniqueness of a limit in
probability proves $U_t(x)=\widetilde U_t(x)$ almost surely.

Intersecting these probability-one events over a countable dense
set of space-time points, and then using continuity, gives equality
everywhere on each compact rectangle. A further countable
intersection over the rectangles $[0,k]\times[-k,k]$, $k\ge1$,
proves the simultaneous assertion.
\end{proof}

\section{\texorpdfstring{Discrete Gaussian identities and the lattice remainder}{Discrete Gaussian identities and the lattice remainder}}
\label{sec:discrete}

We now turn to the proof of Theorem~\ref{thm:coupled-L2}. The
strategy parallels the continuous construction: exact Gaussian
identities reduce the $\mathcal S$-transform and the second moment of
$\widehat{\mathsf Z}^N_t(x)$ to expectations of path functionals (this section),
whose limits are identified in Sections~\ref{subsec:legall-complet}
and~\ref{sec:two-path}. The deterministic remainder of the lattice
renormalization is computed in Section~\ref{subsec:discrete-counterterm}.

\textit{For the whole section, we fix $h=(H_0,H_1)$ satisfying the
conditions in \eqref{fr:range}}. The lattice spatial kernel is
$\gamma_{H_1}$ of \eqref{fr:fgn}, the coupling strength is
$\beta_N=\beta_{N,h}=N^{-a_{h}/2}$ of \eqref{fr:scaling}, and the
lattice counterterm is $c_{\rm disc}=c_{\rm disc}(h)$ of
\eqref{fr:counterterm-discrete}. For $H_1=\frac12$ one has
$\gamma_{H_1}=\one_{\{0\}}$, so that $\gamma_{H_1}(S_n-S_m)=\one_{\{S_n=S_m\}}$,
and $\beta_N^2=N^{\frac12-2H_0}$; every statement below is
to be read in that way at that value.

\subsection{The isonormal coupling and the field laws}
\label{subsec:coupling-identification}

We first justify the coupling and the equality in law stated in
Section~\ref{sec:common-noise}.
\begin{proof}[Proof of Proposition~\ref{prop:coupling}]
The covariance of fractional Brownian motion increments gives
\[
 \int_{\mathfrak{I}_{t,n}^N}\int_{\mathfrak{I}_{t,m}^N}K_{H_0}(r-s)\dd s\dd r
 =N^{-2H_0}\gamma_{H_0}(n-m),
\]
and the same identity in the spatial variable gives
\[
 \int_{\R^2}\one_{\mathfrak{J}_j^N}(x)\one_{\mathfrak{J}_{\ell}^N}(y)K_{H_1}(x-y)\dd x\dd y
 =N^{-H_1}\gamma_{H_1}(j-\ell).
\]
In both cases the increments are stationary, so the integral depends on
the two cells only through the difference between their indices, and not on the anchoring time
$t$. The squared prefactor in \eqref{eq:cell-functions} equals
$N^{2H_0+H_1}$ and exactly cancels these two factors.
This proves \eqref{eq:cell-covariance}. The last assertion follows from
the linearity of $W$ and from
$\beta_{N}N^{H_0+\frac{H_1}2}=N^{\frac12}$, a consequence of \eqref{fr:scaling} and \eqref{fr:exponents}. At $H_1=\frac12$ the spatial identity reads
$\int\one_{\mathfrak{J}_j^N}\one_{\mathfrak{J}_{\ell}^N}=N^{-\frac12}\one_{\{j=\ell\}}$.
\end{proof}

\begin{proof}[Proof of Lemma~\ref{lem:environment-law}]
Both arrays are centered Gaussian. Proposition~\ref{prop:coupling}
and \eqref{fr:discrete-covariance} identify all their covariances and therefore
all their finite-dimensional distributions. Since the index set is
countable, this proves \eqref{eq:environment-law}.

For a fixed $N$ and $t$, the simple-walk path average is a finite sum
over $2^{N_t}$ paths, hence a measurable function of finitely many
environment coordinates. Spatial translation preserves the Gaussian
array law, the spatial covariance depending on $j-\ell$ alone. The same
is true of the finite-horizon reversal $n\mapsto N_t-n+1$: it replaces the index difference $n-m$ by $m-n$, and $\gamma_{H_0}$ is even. This proves the
assertion about the two one-point laws.

For the joint field laws, we use a single array with forward time
indices:
\begin{equation}\label{eq:forward-cell-array}
 \bar\omega^N(k,j)
 :=W\left(N^{H_0+\frac{H_1}2}
       \one_{\left(\frac{k-1}N,\frac kN\right]\times \mathfrak{J}_j^N}\right),
 \qquad k\ge1,\quad j\in\mathbb Z.
\end{equation}
At every grid time $t=m/N$, its relation to the terminal-time array is
\begin{equation}\label{eq:forward-backward-cells}
 \omega_{m/N}^N(n,j)=\bar\omega^N(m-n+1,j),
 \qquad 1\le n\le m,\quad j\in\mathbb Z.
\end{equation}
The superscript $N$ records the mesh scale in both arrays. The covariance calculation in Proposition~\ref{prop:coupling} gives
\[
 \E[\bar\omega^N(k,j)\bar\omega^N(\ell,r)]
 =\gamma_{H_0}(k-\ell)\gamma_{H_1}(j-r).
\]
Thus the whole countable array $\bar\omega^N$ has the law of
$\omega$, simultaneously for all positive time indices and all
spatial indices. The relation \eqref{eq:forward-backward-cells} and the definitions of
the two partition fields give, almost surely and simultaneously for all
$m\in\mathbb Z_{\ge0}$ and $j\in\mathbb Z$,
\begin{equation}\label{eq:grid-field-identification}
 \widehat{\mathsf Z}^N_{m/N}\left(\frac j{\sqrt N}\right)
 =\mathsf Z^N_{m/N}\left(\frac j{\sqrt N};\bar\omega^N\right).
\end{equation}
For $m=0$ both sides equal one. The map from the environment array
to the grid partition field is measurable, since each coordinate is
a finite path sum. Bilinear interpolation is a measurable map from
the countable grid values into $C_{\rm loc}$: its restriction to each
compact rectangle depends continuously on finitely many grid values.
Applying these maps to $\bar\omega^N\overset{\rm law}=\omega$ proves
\eqref{eq:interpolated-field-law}. This uses one common array for all
horizons, rather than separate reversals of the environment at each
horizon.
\end{proof}

\subsection{Representation of $ \widehat{\mathsf Z}^N_t(x)$}\label{subsec:intro-yn-etc}

Recall from \eqref{eq:ZN-coupled} that the coupled partition field is
\[
 \widehat{\mathsf Z}^N_t(x)
 =e^{-c_{\rm disc}N_t\beta_N^2}
 E_S\left[
   \exp\left(\beta_N\sum_{n=1}^{N_t}\omega_t^N(n,\widetilde{N}_x+S_n)\right)\right]
 =E_S\left[e^{W(h^{S,N}_{t,x})-c_{\rm disc}N_t\beta_N^2}\right],
\]
where $\omega_t^N$ is the realization \eqref{eq:common-noise} of the
environment on the Gaussian space of $W$, and $h^{S,N}_{t,x}$ is the
cell function \eqref{eq:hNS-explicit}.

Throughout this section, and again in
Sections~\ref{subsec:general-occupation}, \ref{subsec:two-path-general}
and~\ref{sec:short}, we fix a space-time point $(t,x)$ with
$t>0$ and $x\in\R$, and we consider $N\ge\frac1t$, so that
$N_t\ge1$, with $N_t$ as in \eqref{eq:Nt}.

If $N_t=0$, all discrete sums below are empty, so that
$\widehat{\mathsf Z}^N_t(x)=1$ and all the associated path functionals
and remainders vanish. In particular, at $t=0$,
$\widehat{\mathsf Z}^N_0(x)=1=U_0(x)$.

For the cell function $h^{S,N}_{t,x}$ defined in \eqref{eq:hNS-explicit}, the conditional variance of $W(h^{S,N}_{t,x})$ is
\begin{equation}\label{eq:GammaN}
 \Gamma^N_t(S)=\norm{h^{S,N}_{t,x}}_{\Hh}^2
 =\beta_N^2\sum_{n,m=1}^{N_t}
 \gamma_{H_0}(n-m)\gamma_{H_1}(S_n-S_m),
\end{equation}
by \eqref{eq:cell-covariance}.
It does not depend on $x$, the spatial shift $\widetilde{N}_x$ cancelling in the
difference of the two walk positions. Set
\begin{align}
 Y^N_t&=\Gamma^N_t(S)-E_S\left[\Gamma^N_t(S)\right],\label{eq:YN}\\
 R^N_t&=\frac12E_S\left[\Gamma^N_t(S)\right]-c_{\rm disc}N_t\beta_N^2.\label{eq:RN}
\end{align}

The definitions of $\widehat{\mathsf Z}^N_t(x)$ and $\cE$ give
\begin{equation}\label{eq:Gaussian-disc}
 \widehat{\mathsf Z}^N_t(x)=e^{R^N_t}E_S\left[e^{\frac12 Y^N_t}\cE(h^{S,N}_{t,x})\right].
\end{equation}
This is the discrete analogue of \eqref{eq:Gaussian-cont}.

For $\varphi\in\D$, define
\begin{equation}\label{eq:ell-N}
 \ell^{N}_{t,x}(\varphi)=\ip{h^{S,N}_{t,x}}{\varphi}_{\Hh}.
\end{equation}
The dependence of $Y^N_t$ and $\ell^N_{t,x}$ on the reference walk
$S$ is left implicit.

Equations~\eqref{eq:H-inner} and \eqref{eq:ell-N} give
\[
 \ell^{N}_{t,x}(\varphi)
 =\int_{\R^2}h^{S,N}_{t,x}(s,y)G_\varphi(s,y)\dd s\dd y.
\]
With $t_N=N_t/N$ as in \eqref{eq:Nt}, the nonnegative cell function
$h^{S,N}_{t,x}$ has total mass
\[
 \int_{\R^2}h^{S,N}_{t,x}(s,y)\dd s\dd y
 =\sqrt N\sum_{n=1}^{N_t}\frac1N\frac1{\sqrt N}=t_N.
\]
Together with \eqref{eq:ell-B-general}, this gives the deterministic bounds
\begin{equation}\label{eq:ell-uniform-bound}
 |\ell^{N}_{t,x}(\varphi)|\le t_N\|G_\varphi\|_\infty\le t\|G_\varphi\|_\infty,
 \qquad |\ell^B_{t,x}(\varphi)|\le t\|G_\varphi\|_\infty.
\end{equation}

For two independent walks $S,\widetilde S$, both started at zero,
let $\widetilde Y^N_t$ denote the centered self-intersection of
$\widetilde S$. Define their mutual intersection by
\begin{equation}\label{eq:alphaN}
 \alpha^N_t(S,\widetilde S)
 =\ip{h^{S,N}_{t,x}}{h^{\widetilde S,N}_{t,x}}_{\Hh}
 =\beta_N^2\sum_{n,m=1}^{N_t}
 \gamma_{H_0}(n-m)\gamma_{H_1}(S_n-\widetilde S_m),
\end{equation}
which again does not depend on $x$.
Since $H_0>\frac12$ and $H_1\ge\frac12$, both $\gamma_{H_0}$ and $\gamma_{H_1}$ are
nonnegative, and hence $\alpha^N_t\ge0$.

\begin{proposition}[The two discrete identities]\label{prop:discrete-identities}
For every $N\ge1$, $\widehat{\mathsf Z}^N_t(x)\in L^2(\Omega,\sigma(W),\Pp)$. For $\varphi\in\D$,
\begin{align}
 \mathcal S_{\widehat{\mathsf Z}^N_t(x)}(\varphi)
 &=e^{R^N_t}E_S\left[
 \exp\left(\frac12 Y^N_t+\ell^{N}_{t,x}(\varphi)\right)\right],
 \label{eq:S-ZN}\\
 \E\left[\widehat{\mathsf Z}^N_t(x)^2\right]
 &=e^{2R^N_t}E_{S,\widetilde S}\left[
 \exp\left(\frac12\bigl(Y^N_t+\widetilde Y^N_t\bigr)
                         +\alpha^N_t(S,\widetilde S)\right)\right].
 \label{eq:second-moment-ZN}
\end{align}
\end{proposition}

\begin{proof}
For fixed $N$, the simple walk has finitely many paths of length $N_t$.
Thus \eqref{eq:Gaussian-disc} is a finite linear combination of Gaussian
exponentials, each in every finite $L^p(\Omega)$. The exchanges of Gaussian
and path expectations below are therefore justified. Taking
 the $\mathcal S$-transform of \eqref{eq:Gaussian-disc} and using
\eqref{eq:Gaussian-pairing} gives \eqref{eq:S-ZN}. Multiplying two copies
of \eqref{eq:Gaussian-disc}, averaging over independent walks, and using
\[
 \E\left[\cE(h^{S,N}_{t,x})\cE(h^{\widetilde S,N}_{t,x})\right]
 =\exp\left(\ip{h^{S,N}_{t,x}}{h^{\widetilde S,N}_{t,x}}_{\Hh}\right)
\]
gives \eqref{eq:second-moment-ZN}.
\end{proof}

\subsection{Exact time scaling}\label{subsec:discrete-scaling}

This subsection reduces a general horizon to the unit horizon. The
unit-horizon functionals, such as $Y^N_1,\alpha^N_1$ for the walk, and $Y_1(B),\alpha_1(B,\widetilde B)$
for Brownian motion, are studied in Sections~\ref{subsec:legall-complet}
and~\ref{sec:two-path}; the scaling identities below transfer the limits of $Y$ and $\alpha$ to a general horizon.

We first give the Brownian identities. Fix $t>0$.
With $\bar B_u=t^{-\frac12}B_{tu}$, $0\le u\le1$, direct changes of
variables give the pathwise identities
\begin{align}
 \Gamma_t^\varepsilon(B)
   &=t^{a_{h}}\Gamma_1^{\frac{\varepsilon}{\sqrt t}}(\bar B),&
 Y_t^\varepsilon(B)
   &=t^{a_{h}}Y_1^{\frac{\varepsilon}{\sqrt t}}(\bar B),\notag\\
 \alpha_t^{\varepsilon,{\epsp}}(B,\widetilde B)
   &=t^{a_{h}}\alpha_1^{\frac{\varepsilon}{\sqrt t},\frac{{\epsp}}{\sqrt t}}
                         (\bar B,\overline{\widetilde B}),&
 R^{\varepsilon}_t&=t^{a_{h}}R^{\frac{\varepsilon}{\sqrt t}}_1.
 \label{eq:time-scaling}
\end{align}
 The
last identity also uses
$c_h^{\frac{\varepsilon}{\sqrt t}}=t^{1-a_{h}}c_h^\varepsilon$, which
follows directly from \eqref{fr:counterterm-continuous}, together
with $a_h+(1-a_h)=1$. The identities hold for every positive
regularization parameter, including
$\frac{\varepsilon}{\sqrt t}>1$. Letting $\varepsilon,\varepsilon'\downarrow0$ in
\eqref{eq:time-scaling}, through
Theorem~\ref{thm:Brownian-self-intersection},
Proposition~\ref{prop:continuous-mutual} and \eqref{eq:R-eps-limit},
gives, almost surely for this fixed $t$, the corresponding
identities for the limiting functionals:
\begin{equation}\label{eq:Yt-alphat-Rt}
 Y_t(B)=t^{a_{h}}Y_1(\bar B),\qquad
 \alpha_t(B,\widetilde B)=t^{a_{h}}
              \alpha_1(\bar B,\overline{\widetilde B}),\qquad
 R_t=t^{a_{h}}R_1.
\end{equation}
Unlike $\ell^B_{t,x}$, these three functionals are therefore
homogeneous in the horizon, with the exponent $a_{h}$, which equals
$2H_0-\frac12$ at $H_1=\frac12$.

The discrete model has an exact scaling property of the same kind, with
the regularization parameter $\frac{\varepsilon}{\sqrt t}$ replaced by
the passage from $N$ to $N_t$, and it holds without any error term.

\begin{lemma}[Exact time scaling of the discrete functionals]
\label{lem:discrete-scaling}
Let $t>0$ and $N\ge\frac1t$. Recall $N_t,t_N$ from
\eqref{eq:Nt}; in this case,
\begin{equation}\label{eq:tN}
 0<t_N\le t,\qquad 0\le t-t_N<\frac1N.
\end{equation}
Set
\begin{equation}\label{eq:thetaN}
 \theta^N_t=t_N^{a_{h}}\in(0,t^{a_{h}}],\qquad
 \theta_t=t^{a_{h}},\qquad \theta^N_t\longrightarrow\theta_t.
\end{equation}

Then, pathwise,
\begin{equation}\label{eq:discrete-scaling}
 \Gamma^N_t(S)=\theta^N_t\Gamma^{N_t}_1(S),\quad
 Y^N_t=\theta^N_t Y^{N_t}_1,\quad
 \alpha^N_t(S,\widetilde S)=\theta^N_t\alpha^{N_t}_1(S,\widetilde S),
\end{equation}
and
\begin{equation}\label{eq:RN-scaling}
 R^N_t=\theta^N_t R^{N_t}_1.
\end{equation}

\end{lemma}

\begin{proof}
Since $\beta_N^2=N^{-a_{h}}$ and $\beta_{N_t}^2=N_t^{-a_{h}}$,
\[
 \beta_N^2
 =\left(\frac{N_t}N\right)^{a_{h}}\beta_{N_t}^2
 =\theta^N_t\beta_{N_t}^2.
\]
By \eqref{eq:GammaN} and \eqref{eq:alphaN}, the quantities
$\Gamma^N_t$ and $\alpha^N_t$ depend on $N$ only through the prefactor
$\beta_N^2$ and through the range $1\le n,m\le N_t$ of the summation
indices. Substituting the displayed identity therefore gives the first
and third identities in \eqref{eq:discrete-scaling}, and the second one
follows by centering. Taking path expectations in the first identity and
using $c_{\rm disc}N_t\beta_N^2=\theta^N_t\,c_{\rm disc}N_t\beta_{N_t}^2$
gives \eqref{eq:RN-scaling}. Finally $t_N\to t$, whence
$\theta^N_t\to\theta_t$.
\end{proof}

The linear occupation functional $\ell^{N}_{t,x}$ is not homogeneous
under this scaling, because the test function $\varphi$ is fixed; it is
treated directly in Section~\ref{subsec:general-occupation}.

\subsection{Convergence of \texorpdfstring{$R^N_1$}{the unit-time discrete remainder}}\label{subsec:discrete-counterterm}

We now compute the deterministic remainder $R^N_1$ introduced in
\eqref{eq:RN}. To display its unscaled covariance sum, set
\begin{equation}\label{eq:LambdaN}
 \Lambda_N(S)=\sum_{n,m=1}^N\gamma_{H_0}(n-m)\gamma_{H_1}(S_n-S_m),
\end{equation}
so that $\Gamma^N_1=\beta_N^2\Lambda_N$ by \eqref{eq:GammaN}, and recall $A_{h}$ from \eqref{eq:AH}.

For the simple symmetric walk, the transition probabilities
are explicitly given, for $r\ge0$ and $j\in\mathbb Z$, by
\begin{equation}\label{eq:LLT-j}
\begin{aligned}
 p_{2r+1}(j)
 &=\begin{cases}
  2^{-(2r+1)}\displaystyle\binom{2r+1}{r+(j+1)/2},
       & j\text{ odd and }|j|\le2r+1,\\
  0,   & \text{otherwise},
 \end{cases}\\[1ex]
 p_{2r}(j)
 &=\begin{cases}
  2^{-2r}\displaystyle\binom{2r}{r+j/2},
       & j\text{ even and }|j|\le2r,\\
  0,   & \text{otherwise}.
 \end{cases}
\end{aligned}
\end{equation}

\begin{proposition}[Limit of the discrete deterministic remainder]
\label{prop:counterterm}
For the simple symmetric walk \eqref{eq:simple-walk}, the series
defining $c_{\rm disc}(h)$ in \eqref{fr:counterterm-discrete}-\eqref{eq:clat}
converges, and the deterministic remainder $R^N_1$ defined in
\eqref{eq:RN} satisfies
\begin{equation}\label{eq:RN-limit}
 R^N_1\longrightarrow
 R_1=-A_{h}\left(\frac1{a_{h}}+\frac1{1-a_{h}}\right),
\end{equation}
where $R_1$ is the continuous remainder at time $1$ (see \eqref{eq:R-definition}).

\end{proposition}

\begin{proof}
Write $w_r=\gamma_{H_0}(r)q_{H_1}(r)$, with
$q_{H_1}(r)=E_S[\gamma_{H_1}(S_r)]$ as in
\eqref{fr:counterterm-discrete}.

\smallskip
\noindent\textit{Step 1: uniform bounds and convergence of the counterterm.}
The binomial probabilities in \eqref{eq:LLT-j} attain their maximum
at the central sites. Thus, for every integer $m\ge0$,
\[
 \sup_{j\in\mathbb Z}p_{2m}(j)=p_{2m}(0),\qquad
 \sup_{j\in\mathbb Z}p_{2m+1}(j)=p_{2m+1}(1)=p_{2m+2}(0).
\]
The last equality follows from
$2^{-(2m+1)}\binom{2m+1}{m+1}
=2^{-(2m+2)}\binom{2m+2}{m+1}$.
Therefore, \eqref{eq:LLT} gives the local limit upper bound
\begin{equation}\label{eq:LLT-sup}
 \sup_{j\in\mathbb Z}p_r(j)\le p_r(0)+p_{r+1}(0)\le Cr^{-\frac12}
 \qquad(r\ge1).
\end{equation}
The lattice bounds \eqref{fr:lattice-envelope} of
Lemma~\ref{fr:kernel-lemma} in Appendix~\ref{app:temporal} apply with
$H=H_1$, including $H_1=\frac12$. Split the spatial sum at
$|j|=\sqrt r$. On the inner part, use \eqref{eq:LLT-sup} and the
lattice mass bound; on the outer part, use the pointwise envelope
and $2H_1-2<0$. This gives, for every $r\geq 1$,
\begin{align*}
 q_{H_1}(r)
 &=\sum_{|j|\le\sqrt r}p_r(j)\gamma_{H_1}(j)
   +\sum_{|j|>\sqrt r}p_r(j)\gamma_{H_1}(j)\\
 &\le Cr^{-\frac12}\sum_{|j|\le\sqrt r}\gamma_{H_1}(j)
      +C_{H_1}r^{H_1-1}\sum_{|j|>\sqrt r}p_r(j)\\
 &\le C_{H_1}\bigl(r^{-\frac12}(\sqrt r)^{2H_1-1}
                           +r^{H_1-1}\bigr)
 \le C_{H_1}r^{H_1-1}.
\end{align*}
Here $\sqrt r\ge1$, as required by the mass bound, and the last
probability sum is at most one. Applying the same pointwise envelope
with $H=H_0$ yields, for $r\ge1$,
\begin{equation}\label{eq:w-envelope}
 0\le w_r=\gamma_{H_0}(r)q_{H_1}(r)
 \le C_h r^{2H_0-2}r^{H_1-1}
 =C_h r^{a_h-2}.
\end{equation}
Since $a_h<1$, the series $\sum_{r\ge1}w_r$ converges. This proves
the convergence of the counterterm in \eqref{fr:counterterm-discrete}.

\smallskip
\noindent\textit{Step 2: asymptotics and their common Ces\`aro limit.}
Suppose first that $H_1>\frac12$. Let $Z$ be standard normal and fix
$0<\delta<1<L$.

By \eqref{fr:lattice-asymptotic} in
Lemma~\ref{fr:kernel-lemma},
\[
 r^{1-H_1}\gamma_{H_1}(\sqrt r\,x)\stackrel{r\to \infty}{\longrightarrow} K_{H_1}(x)
 \quad\text{uniformly for }\delta\le|x|\le L.
\]
The central limit theorem for $S_r/\sqrt r$ therefore gives
\[
 r^{1-H_1}E_S\left[\gamma_{H_1}(S_r)
       \one_{\{\delta\sqrt r\le |S_r|\le L\sqrt r\}}\right]
 \stackrel{r\to \infty}{\longrightarrow}
 E\left[K_{H_1}(Z)\one_{\{\delta\le |Z|\le L\}}\right].
\]
The omitted neighborhood of zero is controlled by
\eqref{eq:LLT-sup} and \eqref{fr:lattice-envelope}:
\begin{align*}
 r^{1-H_1}E_S\left[\gamma_{H_1}(S_r)
                   \one_{\{|S_r|\le\delta\sqrt r\}}\right]
 &\le Cr^{\frac12-H_1}
                 \sum_{|j|\le\delta\sqrt r}\gamma_{H_1}(j)\\
 &\le C_{H_1}r^{\frac12-H_1}(1+\delta\sqrt r)^{2H_1-1}
 =C_{H_1}(\delta+r^{-\frac12})^{2H_1-1}.
\end{align*}
For the tail, the pointwise envelope gives
\[
 r^{1-H_1}E_S\left[\gamma_{H_1}(S_r)
                   \one_{\{|S_r|>L\sqrt r\}}\right]
 \le C_{H_1}L^{2H_1-2}.
\]
Letting $r\to\infty$, then $\delta\downarrow0$ and $L\to\infty$,
and using the integrability of $K_{H_1}(Z)$, proves
\[
 r^{1-H_1}q_{H_1}(r)\stackrel{r\to \infty}{\longrightarrow} E[K_{H_1}(Z)]=\kappa_{H_1},
\]
with $\kappa_{H_1}$ as in \eqref{eq:AH}.

For the temporal factor, Lemma~\ref{fr:kernel-lemma}, specifically
\eqref{fr:lattice-asymptotic} with $H=H_0$, gives
\[
 \gamma_{H_0}(r)\sim c_{H_0}r^{2H_0-2}.
\]
Combining the last two limits and $A_h=c_{H_0}\kappa_{H_1}$ proves
the asymptotic for $w_r$ when $H_1>\frac12$.
At $H_1=\frac12$, one has $q_{H_1}(r)=p_r(0)$ by
\eqref{fr:fgn}; the temporal asymptotic and the period-two return
probabilities in \eqref{eq:LLT} give the corresponding result.
Thus the two cases read
\begin{equation}\label{eq:ar-asympt}
\begin{aligned}
 w_r&\sim A_h r^{a_h-2}
       &&(H_1>\tfrac12),\\
 w_{2m+1}&=0,\qquad
 w_{2m}\sim 2A_{H_0}(2m)^{a_h-2}
       &&(H_1=\tfrac12),
\end{aligned}
\end{equation}
where $A_{H_0}=c_{H_0}/\sqrt{2\pi}$ is the value of $A_h$ at
$H_1=\frac12$.

Put $v_r=r^{2-a_h}w_r$. For $H_1>\frac12$, the convergence
$v_r\to A_h$ implies convergence of its Ces\`aro means to $A_h$.

For $H_1=\frac12$, one has $v_{2m+1}=0$ and
$v_{2m}\to2A_{H_0}$. For an integer $L\ge2$, set $M_L=\lfloor L/2\rfloor$.
Then
\[
 \frac1L\sum_{1\le r\le L}v_r
 =\frac1L\sum_{m=1}^{M_L}v_{2m}
 =\frac{M_L}{L}\left(\frac1{M_L}\sum_{m=1}^{M_L}v_{2m}\right)
 \longrightarrow\frac12\cdot2A_{H_0}=A_h,
\]
because $M_L/L\to\frac12$ and the arithmetic means of the convergent
sequence $(v_{2m})$ tend to $2A_{H_0}$. Thus, in both cases,
\begin{equation}\label{eq:w-cesaro}
 \frac1L\sum_{1\le r\le L}v_r
 =\frac1L\sum_{1\le r\le L}r^{2-a_h}w_r
 \longrightarrow A_h\qquad(L\to\infty).
\end{equation}

\smallskip
\noindent\textit{Step 3: an exact identity for the remainder.}
Independence and stationarity of the increments give
\[
 E_S[\Lambda_N]
 =N\gamma_{H_0}(0)\gamma_{H_1}(0)
       +2\sum_{r=1}^{N-1}(N-r)w_r.
\]
Recall from \eqref{fr:counterterm-discrete} that
$2c_{\rm disc}=1+2\sum_{r\ge1}w_r$. Since $\gamma_{H_0}(0)=\gamma_{H_1}(0)=1$, subtracting $2Nc_{\rm disc}$
from the preceding identity yields
\begin{equation}\label{eq:Lambda-remainder}
 E_S[\Lambda_N]-2Nc_{\rm disc}
 =-2\left(\sum_{r=1}^{N-1}rw_r
       +N\sum_{r=N}^{\infty}w_r\right).
\end{equation}
Using \eqref{eq:RN}, $\Gamma^N_1=\beta_N^2\Lambda_N$ and
$\beta_N^2=N^{-a_h}$, we obtain
\begin{align*}
 R^N_1
 &=\frac12\beta_N^2\bigl(E_S[\Lambda_N]-2Nc_{\rm disc}\bigr)\\
 &=-N^{-a_h}\sum_{r=1}^{N-1}rw_r
       -N^{1-a_h}\sum_{r=N}^{\infty}w_r.
\end{align*}

\smallskip

\noindent\textit{Step 4: convergence of the two weighted sums.}
By \eqref{eq:w-envelope} and \eqref{eq:w-cesaro}, the sequence
$(v_r)$ is bounded and has Ces\`aro limit $A_h$.

Fix $0<b<c<\infty$ and $f\in C^1([b,c])$. Define
\[
 D_0=0,\qquad D_n=\sum_{r=1}^n(v_r-A_h)=o(n).
\]
For all sufficiently large $N$, put
$m_N=\lceil bN\rceil\le n_N=\lfloor cN\rfloor$.
Summation by parts gives
\begin{align*}
 \frac1N\sum_{r=m_N}^{n_N}f(r/N)(v_r-A_h)
 &=\frac{f(n_N/N)D_{n_N}-f(m_N/N)D_{m_N-1}}N\\
 &\quad-\frac1N\sum_{r=m_N}^{n_N-1}D_r
           \bigl(f((r+1)/N)-f(r/N)\bigr).
\end{align*}
With both supremum norms taken over $[b,c]$, the mean value theorem
therefore yields
\[
 \left|\frac1N\sum_{r=m_N}^{n_N}f(r/N)(v_r-A_h)\right|
 \le\bigl(2\|f\|_\infty+(c-b)\|f'\|_\infty\bigr)
       \max_{m_N-1\le k\le n_N}\frac{|D_k|}{N}
 \longrightarrow0.
\]
Indeed, $D_k/k\to0$, while $m_N-1\to\infty$ and $k/N\le c$ on
this range. Adding the ordinary Riemann sum
$A_hN^{-1}\sum_{r=m_N}^{n_N}f(r/N)$ proves
\[
 \frac1N\sum_{bN\le r\le cN}f(r/N)v_r
 \stackrel{N\to\infty}{\longrightarrow} A_h\int_b^c f(u)\dd u.
\]

We apply this to the two sums
\begin{align*}
 N^{-a_h}\sum_{r=1}^{N-1}rw_r
 &=\frac1N\sum_{r=1}^{N-1}(r/N)^{a_h-1}v_r,\\
 N^{1-a_h}\sum_{r=N}^{\infty}w_r
 &=\frac1N\sum_{r=N}^{\infty}(r/N)^{a_h-2}v_r.
\end{align*}
For $0<\delta<1<L$, the omitted ends are bounded uniformly in $N$ by
\begin{align*}
 \frac1N\sum_{1\le r<\delta N}(r/N)^{a_h-1}v_r
 &\le C_hN^{-a_h}\sum_{1\le r<\delta N}r^{a_h-1}
 \le C_h\delta^{a_h},\\
 \frac1N\sum_{r>LN}(r/N)^{a_h-2}v_r
 &\le C_hN^{1-a_h}\sum_{r>LN}r^{a_h-2}
 \le C_hL^{a_h-1}.
\end{align*}
Since $0<a_h<1$, these bounds vanish as $\delta\downarrow0$ and
$L\to\infty$. The convergence on $[\delta,1]$ and $[1,L]$ thus gives
\[
 N^{-a_h}\sum_{r=1}^{N-1}rw_r\longrightarrow
 A_h\int_0^1u^{a_h-1}\dd u=\frac{A_h}{a_h},
\]
and
\[
 N^{1-a_h}\sum_{r=N}^{\infty}w_r\longrightarrow
 A_h\int_1^\infty u^{a_h-2}\dd u=\frac{A_h}{1-a_h}.
\]
Inserting these limits into the identity of Step~3 yields
\[
 R^N_1\longrightarrow
 -A_h\left(\frac1{a_h}+\frac1{1-a_h}\right)=R_1,
\]
where the last equality is \eqref{eq:R-definition} at $t=1$. This proves
\eqref{eq:RN-limit}.
\end{proof}

\begin{corollary}[The deterministic remainder at a general horizon]
\label{cor:RN-t}
For every $t>0$,
\begin{equation}\label{eq:RN-t-limit}
 R^N_t\longrightarrow t^{a_{h}}R_1=R_t.
\end{equation}
\end{corollary}

\begin{proof}
By \eqref{eq:RN-scaling}, $R^N_t=\theta^N_tR^{N_t}_1$. Since $N\mapsto N_t$ is nondecreasing and tends to infinity, Proposition~\ref{prop:counterterm} gives $R^{N_t}_1\to R_1$. Combined with $\theta^N_t\to\theta_t=t^{a_h}$ from \eqref{eq:thetaN}, this proves the convergence; the identity $t^{a_h}R_1=R_t$ is \eqref{eq:R-definition}.
\end{proof}

\section{\texorpdfstring{The one-path invariance principle}{The one-path invariance principle}}
\label{subsec:legall-complet}

This section identifies the limit of the $\mathcal S$-transform
\eqref{eq:S-ZN}, that is the limit of
\begin{equation}\label{gau-pai}
\mathcal S_{\widehat{\mathsf Z}^N_t(x)}(\varphi)=e^{R^N_t}E_S\left[
 \exp\left(\frac12 Y^N_t+\ell^{N}_{t,x}(\varphi)\right)\right].
\end{equation}

This limit requires joint convergence of the centered self-intersection and the linear occupation functionals,
with the limiting quantities evaluated on the same Brownian path,
and uniform integrability of the exponential. We first establish this
joint convergence at the unit horizon, together with the walk.

\textit{For the whole section, we fix $h=(H_0,H_1)$ satisfying the
conditions in \eqref{fr:range}}, and $S$ is the simple
symmetric walk \eqref{eq:simple-walk}. The Brownian blocks and the series $Y_1(B)$ built from them have been constructed in Section~\ref{sec:Brownian-construction}.
We now define their discrete counterparts and prove joint convergence
with the walk. A comparison of this type between renormalized
self-intersections of planar random walks and of planar Brownian motion
was carried out by Le Gall \cite{LeGall-RW}.

The discrete moment estimates use Lemma~\ref{lem:tanaka-local-time}.
Identification of the limit uses the joint invariance principle for
the walk and its local times in Lemma~\ref{lem:local-time-invariance}.

Constants may depend on $h$, but never on $N$, the dyadic level,
or the moment order.

\subsection{Exact discrete decomposition}

Use the intervals $I_{M,k},J_{M,k}$ introduced in
\eqref{eq:dyadic-JI} and the partition
\eqref{eq:dyadic-triangle-partition}.
Associate the index $i$ with $\tau^N_i=\frac{i-\frac12}{N}$.

Define the uncentered blocks by
\begin{equation}\label{eq:uncentered-block-discrete}
 C^{(N)}_{M,k}:=N^{-a_h}
 \sum_{\tau^N_i\in I_{M,k}}\sum_{\tau^N_j\in J_{M,k}}
 \gamma_{H_0}(j-i)\gamma_{H_1}(S_j-S_i).
\end{equation}
Their centered counterparts are
\begin{equation}\label{eq:block-discrete}
\begin{split}
 A^{(N)}_{M,k}
 &:=C^{(N)}_{M,k}-E[C^{(N)}_{M,k}]\\
 &=N^{-a_h}
 \sum_{\tau^N_i\in I_{M,k}}\sum_{\tau^N_j\in J_{M,k}}
 \gamma_{H_0}(j-i)\bigl(\gamma_{H_1}(S_j-S_i)-q_{H_1}(j-i)\bigr),
\end{split}
\end{equation}
where $q_{H_1}(r)=E_S[\gamma_{H_1}(S_r)]$. An empty sum is zero.

\begin{lemma}[Exact dyadic decomposition]\label{lem:dyadic-exact}
For every $N$,
\begin{equation}\label{eq:YN-dyadic}
 Y^N_1=2\sum_{M\ge1}\sum_{k=1}^{2^{M-1}}A^{(N)}_{M,k},
\end{equation}
and only finitely many levels contribute.
\end{lemma}

\begin{proof}
The diagonal term is zero after centering. By
\eqref{eq:dyadic-triangle-partition}, each pair $\tau^N_i<\tau^N_j$
belongs to exactly one rectangle $I_{M,k}\times J_{M,k}$.
Symmetry gives the factor $2$. If $2^{1-M}<\frac1N$, a parent
interval cannot contain two distinct grid points, so every block
at that level is zero.
\end{proof}

\subsection{The lattice kernel, local times and weighted intersections}

The bounds for $\gamma_{H_0},\gamma_{H_1}$ are collected in
Lemma~\ref{fr:kernel-lemma} of Appendix~\ref{app:temporal}, including
the case $H_1=\frac12$. We now establish the local-time and weighted
intersection estimates for the walk.

\begin{lemma}[Random-walk local-time bound]\label{lem:tanaka-local-time}
For a simple symmetric walk $S$ started from any integer, let
\[
 V_L^S(x)=\sum_{r=0}^{L-1}\one_{\{S_r=x\}},
 \qquad L\in\N,\ L\ge1.
\]
There is a universal constant $C$ such that, uniformly in the
starting point, for every $p\ge2$,
\begin{equation}\label{eq:tanaka-local-time}
 \sup_{x\in\mathbb Z}\|V_L^S(x)\|_p\le C\sqrt{pL}.
\end{equation}
\end{lemma}

\begin{proof}
Let $E_z,\Pp_z$ denote expectation and probability for the walk
started at $z$. Start the walk at zero and write
\[
 \Rret_n=\sum_{r=1}^n\one_{\{S_r=0\}},\qquad \Rret_0=0.
\]
The exact return-count law in
\cite[Section~3, display preceding (11)]{SS-local-time} reads
\begin{equation}\label{eq:return-count-law}
 q_{m,k}:=\Pp_0(\Rret_{2m}=k)=\Pp_0(\Rret_{2m+1}=k)
 =2^{-2m+k}\binom{2m-k}{m},\qquad 0\le k\le m.
\end{equation}
For $m\ge1$, the normalized central binomial coefficient satisfies
\[
 q_{m,0}=\prod_{j=1}^m\left(1-\frac1{2j}\right)
 \le\exp\left(-\frac12\sum_{j=1}^m\frac1j\right)
 \le(m+1)^{-\frac12}.
\]
Moreover, for $1\le k\le m$,
\begin{align*}
 \frac{q_{m,k}}{q_{m,0}}
 &=\prod_{\ell=0}^{k-1}\frac{2(m-\ell)}{2m-\ell}
 =\prod_{\ell=0}^{k-1}\left(1-\frac{\ell}{2m-\ell}\right)\\
 &\le\exp\left(-\sum_{\ell=0}^{k-1}\frac{\ell}{2m}\right)
 =\exp\left(-\frac{k(k-1)}{4m}\right).
\end{align*}
Since $k(k-1)\ge\frac12k^2$ for $k\ge2$, and since $q_{m,k}\le q_{m,0}\le m^{-1/2}$ for $k\in\{0,1\}$, enlarging the absolute constant gives
\[
 q_{m,k}\le C m^{-\frac12}\exp\left(-\frac{k^2}{8m}\right).
\]

Consequently, for every $u\ge0$ and $m\ge1$,
\begin{align}
 \Pp_0(\Rret_{2m}\ge u)
 &\le C m^{-\frac12}\sum_{k\ge\lceil u\rceil}
                      e^{-\frac{k^2}{8m}}\notag\\
 &\le C e^{-\frac{u^2}{16m}}
          m^{-\frac12}\sum_{k\ge0}e^{-\frac{k^2}{16m}}\notag\\
 &\le C e^{-\frac{u^2}{16m}}m^{-\frac12}
          \left(1+\int_0^\infty e^{-\frac{v^2}{16m}}\dd v\right)\notag\\
 &=C e^{-\frac{u^2}{16m}}
          \left(m^{-\frac12}+2\sqrt\pi\right)
 \le C e^{-\frac{u^2}{16m}}.
 \label{eq:return-count-tail}
\end{align}

For every integer $j\ge1$, integration of
\eqref{eq:return-count-tail}, with the substitution
$v=\frac{u^2}{16m}$, yields
\[
 E_0[\Rret_{2m}^{2j}]
 \le C\,2j\int_0^\infty u^{2j-1}e^{-\frac{u^2}{16m}}\dd u
 =C(16m)^j j!.
\]
Given $p\ge2$, take $j=\lceil\frac p2\rceil$ and use $j!\le j^j$:
\[
 \|\Rret_{2m}\|_p\le\|\Rret_{2m}\|_{2j}
 \le C\sqrt{jm}\le C\sqrt{pm}.
\]
Equation \eqref{eq:return-count-law} treats odd horizons as well, while
$\Rret_0=\Rret_1=0$. Since $V_L^S(0)=1+\Rret_{L-1}$, we obtain
$\|V_L^S(0)\|_p\le C\sqrt{pL}$.

For an arbitrary starting point $z$, let $T_x^S$ be the first
hitting time of $x$. Conditioning on that time and using
independence of subsequent increments gives
\[
 E_z[V_L^S(x)^p]
 =\sum_{j=0}^{L-1}\Pp_z(T_x^S=j)\,E_0[V_{L-j}^S(0)^p]
 \le E_0[V_L^S(0)^p].
\]
The inequality follows by monotonicity in the remaining time and
proves the asserted uniformity.
\end{proof}

The next estimate is the lattice counterpart of
Lemma~\ref{lem:Brownian-occupation}: it replaces, in the fractional
regime, the integrability of the spatial kernel, and it is again
obtained by splitting the occupation sum at a radius which is then
optimized.

\begin{lemma}[Lattice occupation against the spatial kernel]
\label{lem:RW-occupation}
Let $S$ be a simple symmetric walk started at any integer.
Uniformly in the starting point and in $z\in\mathbb Z$, for $p\ge2$
and every integer $L\ge1$,
\begin{equation}\label{eq:RW-occupation}
 \left\|\sum_{r=0}^{L-1}\gamma_{H_1}(S_r-z)\right\|_p
 \le C_{h}\,p^{1-H_1}L^{H_1}.
\end{equation}
The same bound holds on any segment of length $L$, conditionally
on the past at the beginning of that segment. 
\end{lemma}

\begin{proof}
If $H_1=\frac12$, then $\gamma_{1/2}(j)=\one_{\{j=0\}}$ for
$j\in\mathbb Z$, so Lemma~\ref{lem:tanaka-local-time} gives
\[
 \left\|\sum_{r=0}^{L-1}\gamma_{1/2}(S_r-z)\right\|_p
 =\|V_L^S(z)\|_p
 \le C\sqrt{pL}.
\]
Suppose now that $H_1>\frac12$. With $V_L^S$ as in that lemma,
the occupation identity gives
\begin{align*}
 \sum_{r=0}^{L-1}\gamma_{H_1}(S_r-z)
 &=\sum_{y\in\mathbb Z}\gamma_{H_1}(y-z)V_L^S(y)\\
 &=\sum_{|y-z|\le R}
                 \gamma_{H_1}(y-z)V_L^S(y)
   +\sum_{|y-z|>R}
                 \gamma_{H_1}(y-z)V_L^S(y)=: A_R+B_R.
\end{align*}
For the first term, Minkowski's inequality, the local-time bound
\eqref{eq:tanaka-local-time} and the lattice mass bound
\eqref{fr:lattice-envelope} yield
\begin{align*}
 \|A_R\|_p
 &\le\sum_{|y-z|\le R}\gamma_{H_1}(y-z)\|V_L^S(y)\|_p\\
 &\le C\sqrt{pL}\sum_{|j|\le R}\gamma_{H_1}(j)
 \le C_h\sqrt{pL}\,R^{2H_1-1}.
\end{align*}
For the second term, \eqref{fr:lattice-envelope} and $2H_1-2<0$
give $\gamma_{H_1}(y-z)\le C_hR^{2H_1-2}$ when $|y-z|>R$.
This implies the pathwise bound
\[
 0\le B_R
 \le C_hR^{2H_1-2}\sum_{y\in \mathbb Z}V_L^S(y)
 \le C_hLR^{2H_1-2}.
\]
Combining these estimates, we obtain, for every $R\ge1$,
\[
 \left\|\sum_{r=0}^{L-1}\gamma_{H_1}(S_r-z)\right\|_p
 \le C_h\left(\sqrt{pL}\,R^{2H_1-1}+LR^{2H_1-2}\right).
\]
If $L\ge p$, choose $R=\sqrt{L/p}\ge1$. Then
\begin{align*}
& \sqrt{pL}\,R^{2H_1-1}=p^{1-H_1}L^{H_1},\qquad 
 LR^{2H_1-2} =p^{1-H_1}L^{H_1},
\end{align*}
which proves \eqref{eq:RW-occupation} in this case.
If $L<p$, the pointwise envelope in \eqref{fr:lattice-envelope}
implies $\sup_{j\in\mathbb Z}\gamma_{H_1}(j)\le C_h$. Hence,
pathwise,
\[
 0\le\sum_{r=0}^{L-1}\gamma_{H_1}(S_r-z)
 \le C_hL
 =C_hL^{H_1}L^{1-H_1}
 \le C_hp^{1-H_1}L^{H_1}.
\]
This proves the bound for every integer $L\ge1$. All estimates are
uniform in the starting point and in $z$.

Finally, let $\mathcal F_n^S=\sigma(S_0,\ldots,S_n)$ for an integer
$n\ge0$, and write $E_a$ for expectation when $S_0=a$. The Markov
property and the uniform bound just proved, in either spatial regime,
give
\begin{align*}
 E\left[\left(\sum_{r=0}^{L-1}\gamma_{H_1}(S_{n+r}-z)\right)^p
          \biggm|\mathcal F_n^S\right]
 &=\left.
     E_a\left[\left(\sum_{r=0}^{L-1}\gamma_{H_1}(S_r-z)\right)^p\right]
     \right|_{a=S_n}\\
 &\le\left(C_hp^{1-H_1}L^{H_1}\right)^p
 \qquad\text{almost surely}.
\end{align*}
Taking the $p$th root proves the conditional version.
\end{proof}

The next lemma is the discrete counterpart of
Lemma~\ref{lem:Brownian-weighted-block}. The proof again combines the
occupation estimate with a decomposition into dyadic annuli in the
sum of the two time variables.

\begin{lemma}[Weighted random-walk intersections]\label{lem:weighted-block}
Let $S,\widetilde S$ be independent simple symmetric walks. For an integer $L\ge1$, put
\[
 T_L=\sum_{r,s=0}^{L-1}\gamma_{H_1}(S_r-\widetilde S_s),\qquad
 Q_L=\sum_{r,s=0}^{L-1}(1+r+s)^{2H_0-2}\gamma_{H_1}(S_r-\widetilde S_s).
\]
Uniformly in their starting points, for $p\ge2$,
\begin{equation}\label{eq:weighted-block-Lp}
 \|T_L\|_p\le C_{h}\,p^{1-H_1}L^{1+H_1},\qquad
 \|Q_L\|_p\le C_{h}\,p^{1-H_1}L^{a_{h}}.
\end{equation}
The second estimate also holds with the weight $\gamma_{H_0}(1+r+s)$. In particular,
since $p^{1-H_1}\le\sqrt p$, for every integer $m\ge1$,
\begin{equation}\label{eq:weighted-block-moments}
 E\left[Q_L^m\right]\le C_{h}^m(m!)^{\frac12}L^{a_{h}m}.
\end{equation}

\end{lemma}

\begin{proof}
Since $T_L=\sum_{s<L}\sum_{r<L}\gamma_{H_1}(S_r-\widetilde S_s)$ and the two walks
are independent, conditioning on $\widetilde S_s=z$ and applying
\eqref{eq:RW-occupation} to $S$ gives
\[
 \|{\textstyle\sum_{r<L}}\gamma_{H_1}(S_r-\widetilde S_s)\|_p^p
 =\sum_z\Pp(\widetilde S_s=z)
      E\left[\left|{\textstyle\sum_{r<L}}\gamma_{H_1}(S_r-z)\right|^p\right]
 \le \bigl(C_{h}p^{1-H_1}L^{H_1}\bigr)^p,
\]
the bound in \eqref{eq:RW-occupation} being uniform in $z$.
Minkowski's inequality over $s<L$ proves the first bound in
\eqref{eq:weighted-block-Lp}.

\smallskip

For $0\le j\le J_L:=\lfloor\log_2(2L-1)\rfloor$, define
\begin{align*}
 \mathcal A_{j,L}
 &=\{(r,s)\in\{0,\ldots,L-1\}^2:
                       2^j\le1+r+s<2^{j+1}\},\\
 Q_{L,j}
 &=\sum_{(r,s)\in\mathcal A_{j,L}}
       (1+r+s)^{2H_0-2}\gamma_{H_1}(S_r-\widetilde S_s).
\end{align*}
These sets partition all the pairs in $Q_L$, so
$Q_L=\sum_{j=0}^{J_L}Q_{L,j}$. Set
$\ell_{j,L}=\min\{L,2^{j+1}\}$.
For $(r,s)\in\mathcal A_{j,L}$, one has
$r,s<\ell_{j,L}$, and the negativity of $2H_0-2$ gives
$(1+r+s)^{2H_0-2}\le2^{j(2H_0-2)}$.
Since all summands are nonnegative, enlarging the summation domain
to the full square gives the pathwise comparison
\begin{equation}\label{eq:annulus-domination}
 0\le Q_{L,j}
 \le 2^{j(2H_0-2)}
       \sum_{r,s=0}^{\ell_{j,L}-1}\gamma_{H_1}(S_r-\widetilde S_s)
 =2^{j(2H_0-2)}T_{\ell_{j,L}}.
\end{equation}
Here $T_{\ell_{j,L}}$ is formed from the same two walks.
Apply the unweighted bound already proved, with time horizon
$\ell_{j,L}$:
\[
 \|Q_{L,j}\|_p
 \le C_{h}p^{1-H_1}2^{j(2H_0-2)}\ell_{j,L}^{1+H_1}
 \le C_{h}\,2^{1+H_1}p^{1-H_1}2^{j a_{h}}.
\]
Minkowski's inequality and the geometric sum now yield
\[
 \|Q_L\|_p
 \le C_{h}p^{1-H_1}\sum_{j=0}^{J_L}2^{j a_{h}}
 \le C_{h}p^{1-H_1}2^{J_La_{h}}
 \le C_{h}p^{1-H_1}L^{a_{h}},
\]
because $a_{h}>0$ and $2^{J_L}\le2L$.

For the weight $\gamma_{H_0}(1+r+s)$, the envelope
\eqref{fr:lattice-envelope} of Lemma~\ref{fr:kernel-lemma}, applied
with $H=H_0$, gives
\[
 0\le\gamma_{H_0}(1+r+s)\le C_{H_0}(2+r+s)^{2H_0-2}
 \le C_{H_0}(1+r+s)^{2H_0-2},
\]
so that the bound just proved for $Q_L$ applies.
For $m\ge2$, use $m^{\frac m2}\le e^{\frac m2}(m!)^{\frac12}$ to derive \eqref{eq:weighted-block-moments}; the case $m=1$ follows
from the $L^2$ bound.
\end{proof}

\subsection{Moments of the discrete blocks}

\begin{proposition}[Moments of a discrete block]\label{prop:block-moments}
Uniformly in $N,M,k$ and $p\ge2$,
\begin{equation}\label{eq:block-Lp}
 \|C^{(N)}_{M,k}\|_p+\|A^{(N)}_{M,k}\|_p
 \le C_h\sqrt p\,2^{-a_hM}.
\end{equation}
Hence, for every integer $m\ge1$,
\begin{equation}\label{eq:block-moment-strong}
 E\left[(C^{(N)}_{M,k})^m\right]
 +E\left[|A^{(N)}_{M,k}|^m\right]
 \le C_h^m(m!)^{\frac12}2^{-a_hMm}.
\end{equation}
\end{proposition}

\begin{proof}

Put $\delta=2^{-M}$ and $L=\lceil N\delta\rceil$, and denote the
sets of indices in the two halves of the block by
\[
 \mathcal I=\{i\in\{1,\ldots,N\}:\tau^{N}_i\in I_{M,k}\},
 \qquad
 \mathcal J=\{j\in\{1,\ldots,N\}:\tau^{N}_j\in J_{M,k}\}.
\]
If either set is empty, then $C^{(N)}_{M,k}=A^{(N)}_{M,k}=0$.
Otherwise, write
$n_-=\#\mathcal I$, $n_+=\#\mathcal J$ and
\[
 b=\max\mathcal I
   =\max\{i\in\{1,\ldots,N\}:\tau^{N}_i\in I_{M,k}\}.
\]
Since the two halves are adjacent intervals and the grid is increasing,
\[
 \mathcal I=\{b-n_-+1,\ldots,b\},\qquad
 \mathcal J=\{b+1,\ldots,b+n_+\},\qquad n_-,n_+\le L.
\]
Define the two walk segments
\[
 S^-_r=S_{b-r}-S_b\quad(0\le r<n_-),\qquad
 S^+_s=S_{b+1+s}-S_b\quad(0\le s<n_+).
\]
Here $S^-_0=0$ and $S^+_0=\varepsilon:=X_{b+1}\in\{-1,1\}$.
Using the increments of \eqref{eq:simple-walk}, we have
\[
 S^-_{r+1}-S^-_r=-X_{b-r},\qquad
 S^+_{s+1}-S^+_s=X_{b+2+s}
\]
whenever these increments are defined. These two families and
$X_{b+1}$ are independent. Thus, conditionally on $\varepsilon$,
$S^-$ and $S^+$ are independent simple symmetric walks started at
$0$ and $\varepsilon$, respectively. Extend each segment to index
$L-1$ by independent simple symmetric increments if necessary.

For $i=b-r$ and $j=b+1+s$, the temporal and spatial differences are
\[
 j-i=1+r+s,\qquad S_j-S_i=S^+_s-S^-_r.
\]

Consequently, the block \eqref{eq:uncentered-block-discrete} satisfies
\begin{align*}
 0\le C^{(N)}_{M,k}
 &=N^{-a_h}\sum_{i\in\mathcal I}\sum_{j\in\mathcal J}
          \gamma_{H_0}(j-i)\gamma_{H_1}(S_j-S_i)\\
 &=N^{-a_h}\sum_{r=0}^{n_--1}\sum_{s=0}^{n_+-1}
          \gamma_{H_0}(1+r+s)\gamma_{H_1}(S^-_r-S^+_s)\\
 &\le N^{-a_h}\sum_{r,s=0}^{L-1}
          \gamma_{H_0}(1+r+s)\gamma_{H_1}(S^-_r-S^+_s),
\end{align*}
where we used the evenness and nonnegativity of the kernels.
Apply Lemma~\ref{lem:weighted-block} conditionally on $\varepsilon$,
with $(S,\widetilde S)=(S^-,S^+)$ and weight
$\gamma_{H_0}(1+r+s)$. Uniformity in the starting points gives
\begin{align*}
 \bigl(E[(C^{(N)}_{M,k})^p\mid\varepsilon]\bigr)^{1/p}
 &\le C_hp^{1-H_1}N^{-a_h}L^{a_h},\\
 \|C^{(N)}_{M,k}\|_p
 &\le C_hp^{1-H_1}N^{-a_h}L^{a_h}.
\end{align*}
Stationarity of the walk increments yields
\begin{align*}
 E[C^{(N)}_{M,k}]
 &=N^{-a_h}\sum_{i\in\mathcal I}\sum_{j\in\mathcal J}
       \gamma_{H_0}(j-i)q_{H_1}(j-i),\\
 A^{(N)}_{M,k}&=C^{(N)}_{M,k}-E[C^{(N)}_{M,k}].
\end{align*}
A nonzero block requires $2N\delta\ge1$, hence
\[
 L=\lceil N\delta\rceil\le N\delta+1\le3N\delta.
\]
Since $H_1\ge\frac12$ and centering gives
$\|A^{(N)}_{M,k}\|_p\le2\|C^{(N)}_{M,k}\|_p$, the preceding bounds imply
\begin{align*}
 \|C^{(N)}_{M,k}\|_p+\|A^{(N)}_{M,k}\|_p
 &\le3\|C^{(N)}_{M,k}\|_p
 \le C_hp^{1-H_1}N^{-a_h}L^{a_h}\\
 &\le C_h\sqrt p\,N^{-a_h}(3N\delta)^{a_h}
 \le C_h\sqrt p\,\delta^{a_h}
 =C_h\sqrt p\,2^{-a_hM}.
\end{align*}
This proves \eqref{eq:block-Lp}. For $m\ge2$, use
$m^{m/2}\le e^{m/2}(m!)^{1/2}$ to obtain
\eqref{eq:block-moment-strong}; the case $m=1$ follows from the
$L^2$ bound.
\end{proof}

\subsection{Identification at fixed dyadic levels}

We first fix the notation for the rescaled walk and its local times.
For $N\ge1$, define the continuous, piecewise linear path $S^N$ by
\begin{equation}\label{eq:rescaled-walk-path}
 S^N_{(j+\theta)/N}
 =N^{-\frac12}\bigl((1-\theta)S_j+\theta S_{j+1}\bigr),
 \qquad 0\le j<N,\quad 0\le\theta\le1.
\end{equation}
Let $\mathcal L^N$ be the continuous field obtained by bilinear
interpolation on the grid $N^{-1}\mathbb Z\times N^{-1/2}\mathbb Z$
from the values
\begin{equation}\label{eq:rescaled-local-time}
 \mathcal L^N_{j/N}\left(\frac{\ell}{\sqrt N}\right)
 =N^{-\frac12}\sum_{i=1}^j\one_{\{S_i=\ell\}},
 \qquad 0\le j\le N,\quad\ell\in\mathbb Z.
\end{equation}
We also use the step field
\begin{equation}\label{eq:step-local-time}
 \mathcal L^{N,\mathrm{step}}_t(x)
 =N^{-\frac12}\sum_{i=1}^{N_t}
       \one_{\{S_i=\widetilde{N}_x\}},
 \qquad t\in[0,1],\quad x\in\R.
\end{equation}
Here $N_t$ and $\widetilde{N}_x$ are defined in \eqref{eq:Nt} and \eqref{eq:xN}, respectively.

The continuous field is used in the joint invariance principle below;
the step field gives exact identities for discrete occupation measures.

\begin{lemma}[Classical joint invariance of local times]\label{lem:local-time-invariance}

Let $S^N$ and $\mathcal L^N$ be defined by
\eqref{eq:rescaled-walk-path}-\eqref{eq:rescaled-local-time}.
Let $B$ be a standard Brownian motion started at zero, and let $L^B$
be its jointly continuous occupation-density local time, normalized
by \eqref{temps-local}. Put $\CL=C([0,1]\times\R)$ with the topology
of local uniform convergence. Then
\begin{equation}\label{eq:joint-local-time-invariance}
 (S^N,\mathcal L^N)\ToLaw(B,L^B)
 \qquad\text{in }C([0,1])\times\CL.
\end{equation}
Moreover, the step field of \eqref{eq:step-local-time} is
asymptotically equivalent to $\mathcal L^N$: for every $0<R<\infty$,
\begin{equation}\label{eq:local-time-step-equivalence}
 \sup_{\substack{0\le t\le1\\|x|\le R}}
 \bigl|\mathcal L^{N,\mathrm{step}}_t(x)-\mathcal L^N_t(x)\bigr|
 \ToProb0
 \qquad\text{as }N\to\infty.
\end{equation}

\end{lemma}

\begin{proof}

By the simultaneous strong approximation of the walk and its local
times in \cite[Chapter~10, Theorem~10.1 and Remark~3]{Revesz-book},
one can construct $S$ and a standard Brownian motion $W$ on the same
probability space so that, almost surely,
\[
 e_n:=|S_n-W_n|
 +\sup_{\ell\in\mathbb Z}
   \left|\sum_{i=1}^n\one_{\{S_i=\ell\}}-L_n^W(\ell)\right|
 =o(\sqrt n)
 \qquad\text{as }n\to\infty.
\]
The uniformity over $0\le j\le N$ needed below follows from this
almost-sure estimate. Indeed, on the event where it holds, fix
$\varepsilon>0$ and choose a finite $J\ge1$ such that
$e_j\le\varepsilon\sqrt j$ for every $j\ge J$. For $N\ge J$, one has
\[
 \frac1{\sqrt N}\max_{0\le j\le N}e_j
 \le\frac1{\sqrt N}\max_{0\le j<J}e_j+\varepsilon.
\]
The finite maximum on the right does not depend on $N$. Letting
$N\to\infty$ and then $\varepsilon\downarrow0$ gives
$\max_{0\le j\le N}e_j=o(\sqrt N)$ almost surely. In particular,
\[
 \max_{0\le j\le N}|S_j-W_j|
 +\max_{0\le j\le N}\sup_{\ell\in\mathbb Z}
   \left|\sum_{i=1}^j\one_{\{S_i=\ell\}}-L_j^W(\ell)\right|
 \le2\max_{0\le j\le N}e_j=o(\sqrt N).
\]
Set $W^N_t=N^{-1/2}W_{Nt}$. With this notation, one has
\[
 L_t^{W^N}(x)=N^{-1/2}L_{Nt}^W(\sqrt N x),
\]
and thus the total error at the grid points satisfies,
almost surely,
\begin{align*}
 D_N
 &:=\max_{0\le j\le N}|S^N_{j/N}-W^N_{j/N}|
   +\max_{0\le j\le N}\sup_{\ell\in\mathbb Z}
      \left|\mathcal L^N_{j/N}\left(\frac{\ell}{\sqrt N}\right)
             -L^{W^N}_{j/N}\left(\frac{\ell}{\sqrt N}\right)\right|\\
 &=\frac1{\sqrt N}\biggl(
      \max_{0\le j\le N}|S_j-W_j|
      +\max_{0\le j\le N}\sup_{\ell\in\mathbb Z}
        \left|\sum_{i=1}^j\one_{\{S_i=\ell\}}-L_j^W(\ell)\right|
      \biggr)\longrightarrow0.
\end{align*}
Set $\delta_N=N^{-1}+N^{-1/2}$. For $R>0$, put
$K_R=[0,1]\times[-R,R]$ and
$\|g\|_{K_R}=\sup_{(t,x)\in K_R}|g_t(x)|$. For
$f\in C([0,1])$ and $g\in\CL$, define the moduli of continuity
\begin{align*}
 \omega(f,\delta)
 &:=\sup_{\substack{s,t\in[0,1]\\|s-t|\le\delta}}
          |f(s)-f(t)|,\qquad 
 \omega_R(g,\delta)
 :=\sup_{\substack{(s,x),(t,y)\in K_R\\
                    |s-t|+|x-y|\le\delta}}
          |g_s(x)-g_t(y)|.
\end{align*}
Both tend to zero as $\delta\downarrow0$, by uniform continuity
on compact sets.

Recall the weights $\lambda_0(u)=1-u$ and $\lambda_1(u)=u$ from
\eqref{eq:bilinear-interpolation}. For $j\in\{0,\ldots,N-1\}$
and $\theta\in[0,1]$, linear interpolation gives
\begin{align*}
 S^N_{(j+\theta)/N}-W^N_{(j+\theta)/N}&\quad=\sum_{a\in\{0,1\}}\lambda_a(\theta)
           \bigl(S^N_{(j+a)/N}-W^N_{(j+a)/N}\bigr)\\
 &\qquad\quad+\sum_{a\in\{0,1\}}\lambda_a(\theta)
           \bigl(W^N_{(j+a)/N}-W^N_{(j+\theta)/N}\bigr).
\end{align*}
Using $\sum_{a\in\{0,1\}}\lambda_a(\theta)=1$ and
$|(j+a)/N-(j+\theta)/N|\le N^{-1}$, we obtain
\begin{equation}\label{eq:local-time-interpolation-errors-1}
 \|S^N-W^N\|_\infty
 \le D_N+\omega(W^N,N^{-1}).
\end{equation}
Similarly, for $\ell\in\mathbb Z$ and $\theta'\in[0,1]$,
bilinear interpolation gives
\begin{align*}
 &\mathcal L^N_{(j+\theta)/N}
       \left(\frac{\ell+\theta'}{\sqrt N}\right)
   -L^{W^N}_{(j+\theta)/N}
       \left(\frac{\ell+\theta'}{\sqrt N}\right)\\
 &\quad=\sum_{a,b\in\{0,1\}}\lambda_a(\theta)\lambda_b(\theta')
       \biggl[
         \mathcal L^N_{(j+a)/N}\left(\frac{\ell+b}{\sqrt N}\right)
         -L^{W^N}_{(j+a)/N}\left(\frac{\ell+b}{\sqrt N}\right)
       \biggr]\\
 &\qquad\quad+\sum_{a,b\in\{0,1\}}\lambda_a(\theta)\lambda_b(\theta')
       \biggl[
         L^{W^N}_{(j+a)/N}\left(\frac{\ell+b}{\sqrt N}\right)
         -L^{W^N}_{(j+\theta)/N}
              \left(\frac{\ell+\theta'}{\sqrt N}\right)
       \biggr].
\end{align*}
If the point $((j+\theta)/N,(\ell+\theta')/\sqrt N)$ belongs
to $K_R$, the four grid points in these sums lie in $K_{R+1}$
and are at distance at most $\delta_N$ from it in the metric
$|s-t|+|x-y|$. Since $\sum_{a,b\in\{0,1\}}\lambda_a(\theta)\lambda_b(\theta')=1$, this yields
\begin{equation}\label{eq:local-time-interpolation-errors-2}
 \|\mathcal L^N-L^{W^N}\|_{K_R}
 \le D_N+\omega_{R+1}(L^{W^N},\delta_N).
\end{equation}
Brownian scaling gives
\[
 (W^N,L^{W^N})\overset{\rm law}=(B,L^B),
\]
and thus, for every $\varepsilon>0$, 
\begin{align*}
 \Pp\bigl(\omega(W^N,N^{-1})>\varepsilon\bigr)
 &=\Pp\bigl(\omega(B,N^{-1})>\varepsilon\bigr)\longrightarrow0,\\
 \Pp\bigl(\omega_{R+1}(L^{W^N},\delta_N)>\varepsilon\bigr)
 &=\Pp\bigl(\omega_{R+1}(L^B,\delta_N)>\varepsilon\bigr)
   \longrightarrow0,
\end{align*}
since the moduli on the right tend to zero almost surely.
Together with $D_N\to0$, this proves that both errors in
\eqref{eq:local-time-interpolation-errors-1} and \eqref{eq:local-time-interpolation-errors-2} tend to zero in
probability, for every $R>0$.

To pass to the stated function spaces, use the metric
\[
 d_{\CL}(g,h)=\sum_{m=1}^{\infty}2^{-m}
                 \bigl(1\wedge\|g-h\|_{K_m}\bigr),
\]
which induces local uniform convergence. Each term for
$(g,h)=(\mathcal L^N,L^{W^N})$ tends to zero in probability,
and the tail after $m=M$ is at most $2^{-M}$. Hence
$d_{\CL}(\mathcal L^N,L^{W^N})\ToProb0$. Thus the distance
between $(S^N,\mathcal L^N)$ and $(W^N,L^{W^N})$ tends to zero
in probability in the product space. The latter pair has the
fixed law of $(B,L^B)$, so \eqref{eq:joint-local-time-invariance}
follows.

Finally, we prove the step-field equivalence on the same coupling.

With $t_N$ and $x_N$ as defined in \eqref{eq:Nt} and \eqref{eq:xN},
the step field satisfies the exact identity
\begin{equation}\label{eq:local-time-step-grid-identity}
 \mathcal L^{N,\mathrm{step}}_t(x)=\mathcal L^N_{t_N}(x_N).
\end{equation}
For $(t,x)\in K_R$, these definitions give
\[
 (t_N,x_N)\in K_{R+1},\qquad
 |t-t_N|+|x-x_N|\le N^{-1}+N^{-1/2}=\delta_N.
\]
Hence
\begin{align*}
 |\mathcal L^{N,\mathrm{step}}_t(x)-L^{W^N}_t(x)|
 &\le |\mathcal L^N_{t_N}(x_N)-L^{W^N}_{t_N}(x_N)|
      +|L^{W^N}_{t_N}(x_N)-L^{W^N}_t(x)|\\
 &\le D_N+\omega_{R+1}(L^{W^N},\delta_N).
\end{align*}

Consequently,
\[
 \|\mathcal L^{N,\mathrm{step}}-L^{W^N}\|_{K_R}
 \le D_N+\omega_{R+1}(L^{W^N},\delta_N).
\]
Together with \eqref{eq:local-time-interpolation-errors-2}, the
triangle inequality now gives
\begin{equation}\label{oub}
 \|\mathcal L^{N,\mathrm{step}}-\mathcal L^N\|_{K_R}
 \le2D_N+2\omega_{R+1}(L^{W^N},\delta_N)\ToProb0,
\end{equation}
which proves \eqref{eq:local-time-step-equivalence}, since its left-hand side is a functional of the walk alone, whose law is the same on the coupling.
\end{proof}

\begin{lemma}[Convergence at fixed levels]\label{lem:fixed-level}
For every $M_0<\infty$ and $\varphi_1,\ldots,\varphi_r\in\D$,
\begin{align}
 &\left(S^N,(A^{(N)}_{M,k})_{\substack{1\le M\le M_0\\1\le k\le2^{M-1}}},
       (\ell^{N}_{1,0}(\varphi_j))_{j\le r}\right)\notag\\
 &\hspace{1cm}\ToLaw
 \left(B,(A_{M,k,1})_{\substack{1\le M\le M_0\\1\le k\le2^{M-1}}},
       (\ell^B_{1,0}(\varphi_j))_{j\le r}\right).
 \label{eq:fixed-level-convergence}
\end{align}
\end{lemma}

\begin{proof}
We treat the two regimes separately: the spatial rescaling of $\gamma_{H_1}$ is a
kernel convergence for $H_1>\frac12$ and an occupation-measure
convergence for $H_1=\frac12$. The temporal rescaling, the corner
estimates and the passage to the centered blocks are common.

With the cutoff $\chi$ of \eqref{eq:cutoff}, define the truncated versions of the blocks
\eqref{eq:uncentered-block-discrete}-\eqref{eq:block-discrete} by
\begin{equation}\label{eq:discrete-cutoff-blocks}
\begin{split}
 C^{(N),\eta}_{M,k}
 &:=N^{-a_h}
    \sum_{\tau_i^N\in I_{M,k}}\sum_{\tau_j^N\in J_{M,k}}
      \gamma_{H_0}(j-i)\gamma_{H_1}(S_j-S_i)
      \chi\left(\frac{j-i}{\eta N}\right),\\
 A^{(N),\eta}_{M,k}
 &:=C^{(N),\eta}_{M,k}-E[C^{(N),\eta}_{M,k}].
\end{split}
\end{equation}

\smallskip\noindent\textit{\underline{The case $H_1=\frac12$}.}
Introduce the finite positive measures
\[
 \mathfrak I_N=N^{-\frac32}\sum_{i,j=1}^N
          \one_{\{S_i=S_j\}}\delta_{(\tau^N_i,\tau^N_j)}
\]
and the step occupation fields
\[
\mathcal L^N_I(x)=N^{-\frac12}\sum_{\tau^N_i\in I}
                       \one_{\{S_i=\widetilde{N}_x\}}.
\]
For intervals $I,J$, there is an exact identity
\begin{equation}\label{eq:discrete-intersection-measure}
 \mathfrak I_N(I\times J)
 =\int_\R\mathcal L^N_I(x)\mathcal L^N_J(x)\dd x.
\end{equation}

Since $\tau_i^N=(i-\frac12)/N$, one has
\begin{equation}\label{midpo}
 \sup_{\substack{t\in[0,1]\\x\in\R}}
 \Big|N^{-\frac12}\sum_{\substack{1\le i\le N\\\tau_i^N<t}}
       \one_{\{S_i=\widetilde{N}_x\}}
       -\mathcal L^{N,\mathrm{step}}_t(x)\Big|
 \le N^{-\frac12}.
\end{equation}
By Lemma~\ref{lem:local-time-invariance} and the Skorokhod representation
theorem, we may work on a coupling on which
\[
 \|S^N-B\|_\infty+\|\mathcal L^N-L^B\|_{K_R}\longrightarrow0
 \quad\text{a.s., for every integer }R\ge1.
\]
Here $K_R$, $\delta_N$ and $\omega_R$ are as in the preceding proof.
For $I=[u,v)\subset[0,1]$, write $L_I^B=L_v^B-L_u^B$. Taking increments in \eqref{midpo} and using the grid identity \eqref{eq:local-time-step-grid-identity}, which gives, as in the proof of \eqref{oub}, $\|\mathcal L^{N,\mathrm{step}}-L^B\|_{K_R}\le\|\mathcal L^N-L^B\|_{K_{R+1}}+\omega_{R+1}(L^B,\delta_N)$, yields
\begin{align*}
 \sup_{|x|\le R}|\mathcal L_I^N(x)-L_I^B(x)|
 &\le 2N^{-\frac12}
       +2\|\mathcal L^{N,\mathrm{step}}-L^B\|_{K_R}\\
 &\le 2N^{-\frac12}+2\|\mathcal L^N-L^B\|_{K_{R+1}}
       +2\omega_{R+1}(L^B,\delta_N)\\
 &=:\varepsilon_N(R)\longrightarrow0\quad\text{a.s.}
\end{align*}

This bound is uniform over $I$. Almost surely, choose a finite integer
$R>1+\|B\|_\infty+\sup_N\|S^N\|_\infty$; all these occupation fields
then vanish outside $[-R,R]$.
Using \eqref{eq:discrete-intersection-measure} and the Brownian
rectangle identity following \eqref{eq:two-BM-occupation-measure}
(with $\widetilde B=B$), we obtain, for half-open intervals $I,J$,
\begin{align}
 |\mathfrak I_N(I\times J)-\mathfrak I_B(I\times J)|
 &=\left|\int_{-R}^R
       (\mathcal L_I^N\mathcal L_J^N-L_I^B L_J^B)(x)\dd x\right|
       \notag\\
 &\le 2R\varepsilon_N(R)
       \bigl(2\|L_1^B\|_\infty+\varepsilon_N(R)\bigr)
       \longrightarrow0\quad\text{a.s.},
 \label{eq:fixed-level-rectangle-convergence}
\end{align}
where we used $0\le L_I^B,L_J^B\le L_1^B$.
Continuity of local times also gives, for every $t\in[0,1]$,
\[
 \mathfrak I_B(\{t\}\times[0,1])
 =\int_\R\mu_x^B(\{t\})L_1^B(x)\dd x
 =0=\mathfrak I_B([0,1]\times\{t\}).
\]
Taking $I=J=[0,1)$ in
\eqref{eq:fixed-level-rectangle-convergence} yields in particular
\begin{align*}
\mathfrak M_N:=\mathfrak I_N([0,1]^2)
 &=\int_\R\bigl(\mathcal L^N_{[0,1)}(x)\bigr)^2\dd x\\
 &\longrightarrow\int_\R\bigl(L_1^B(x)\bigr)^2\dd x
   =\mathfrak I_B([0,1]^2)=:\mathfrak M\quad\text{a.s.}
\end{align*}
Finally, for $f\in C([0,1]^2)$, let $I_a^q=[a/q,(a+1)/q)$,
$0\le a<q$, with the last interval closed at $1$, and set
\[
 f_q=\sum_{a,b=0}^{q-1}f(a/q,b/q)\one_{I_a^q\times I_b^q},
 \qquad\|f-f_q\|_\infty\longrightarrow0.
\]
For fixed $q$, rectangle convergence implies
$\int f_q\,\dd(\mathfrak I_N-\mathfrak I_B)\to0$ almost surely.
Consequently,
\begin{align*}
 \left|\int f\,\dd(\mathfrak I_N-\mathfrak I_B)\right|
 &\le(\mathfrak M_N+\mathfrak M)\|f-f_q\|_\infty
       +\left|\int f_q\,\dd(\mathfrak I_N-\mathfrak I_B)\right|,\\
 \limsup_{N\to\infty}
 \left|\int f\,\dd(\mathfrak I_N-\mathfrak I_B)\right|
 &\le2\mathfrak M\|f-f_q\|_\infty\xrightarrow[q\to\infty]{}0.
\end{align*}
Together with the path convergence on the same coupling, this proves
\begin{equation}\label{measu-conv}
 (S^N,\mathfrak I_N)\ToLaw(B,\mathfrak I_B)
 \quad\text{in }C([0,1])\times\mathcal M_+([0,1]^2),
\end{equation}
where $\mathcal M_+([0,1]^2)$ carries the weak topology of finite positive
measures.

For fixed $\eta>0$ and $N\eta\ge2$,
\eqref{fr:lattice-asymptotic} in Lemma~\ref{fr:kernel-lemma}, applied
with $H=H_0$, gives
\begin{equation}\label{eq:gamma-away-diagonal}
\begin{split}
 \rho_N(\eta)
 &:=\sup_{\substack{r\in\N\\r\ge\eta N}}
     \left|N^{2-2H_0}\gamma_{H_0}(r)
                  -K_{H_0}\left(\frac rN\right)\right|\\
 &\le C_{H_0}\eta^{2H_0-4}N^{-2}\longrightarrow0.
\end{split}
\end{equation}
Since $\tau_j^N-\tau_i^N=(j-i)/N$, $\gamma_{1/2}=\one_{\{0\}}$
on the lattice, and $N^{-a_h}=N^{-3/2}N^{2-2H_0}$ when
$H_1=\frac12$, the exact integral identity is
\begin{equation}\label{eq:discrete-cutoff-measure-identity}
 C^{(N),\eta}_{M,k}
 =\int_{I_{M,k}\times J_{M,k}}
       N^{2-2H_0}\gamma_{H_0}\bigl(N(t-s)\bigr)
       \chi\left(\frac{t-s}{\eta}\right)
       \mathfrak I_N(\dd s\dd t).
\end{equation}
Set
\[
 w_\eta(s,t)=K_{H_0}(t-s)
                  \chi\left(\frac{|t-s|}{\eta}\right)
 \quad(s\ne t),\qquad w_\eta(s,s)=0.
\]
This function is bounded and continuous on $[0,1]^2$, and the
Brownian cutoff block introduced in the proof of
Lemma~\ref{lem:Brownian-blocks}, whose defining formula makes sense for every $\eta>0$, satisfies
\[
 C^{(\eta)}_{M,k,1}
 =\int_{I_{M,k}\times J_{M,k}}
                w_\eta(s,t)\,\mathfrak I_B(\dd s\dd t).
\]
On the coupling used to prove \eqref{measu-conv}, the bound 
\eqref{eq:gamma-away-diagonal} and the convergence of the total
masses give
\begin{equation}\label{precis-1}
 \left|C^{(N),\eta}_{M,k}
       -\int_{I_{M,k}\times J_{M,k}}
                   w_\eta(s,t)\,\mathfrak I_N(\dd s\dd t)\right|
 \le\rho_N(\eta)\mathfrak M_N\longrightarrow0\quad\text{a.s.}
\end{equation}
Moreover, the rectangle boundaries have zero $\mathfrak I_B$-mass,
so weak convergence on this coupling (see \eqref{measu-conv}) yields
\begin{equation}\label{precis-2}
 \int_{I_{M,k}\times J_{M,k}}
          w_\eta(s,t)\,\mathfrak I_N(\dd s\dd t)
 \longrightarrow
 \int_{I_{M,k}\times J_{M,k}}
          w_\eta(s,t)\,\mathfrak I_B(\dd s\dd t)
 =C^{(\eta)}_{M,k,1}\quad\text{a.s.}
\end{equation}
Thus, writing
$\mathcal F_{M_0}=\{(M,k):1\le M\le M_0,\ 1\le k\le2^{M-1}\}$,
we obtain from \eqref{precis-1}, \eqref{precis-2} and $\|S^N-B\|_\infty\to0$ on this coupling,
\begin{equation}\label{eq:fixed-level-uncentered-cutoff-convergence}
 \left(S^N,(C^{(N),\eta}_{M,k})_{(M,k)\in\mathcal F_{M_0}}\right)
 \ToLaw
 \left(B,(C^{(\eta)}_{M,k,1})_{(M,k)\in\mathcal F_{M_0}}\right).
\end{equation}

\smallskip\noindent\textit{\underline{The case $H_1>\frac12$}.} Fix $\eta,\delta>0$ and a block $(M,k)$, and abbreviate
$I=I_{M,k}$, $J=J_{M,k}$. Insert a spatial cutoff into
\eqref{eq:discrete-cutoff-blocks}, setting
\begin{equation}\label{eq:double-cutoff-block}
 C^{(N),\eta,\delta}_{M,k}
 :=N^{-a_h}\sum_{\tau_i^N\in I}\sum_{\tau_j^N\in J}
 \gamma_{H_0}(j-i)\gamma_{H_1}(S_j-S_i)\chi\left(\frac{j-i}{\eta N}\right)
          \chi\left(\frac{|S_j-S_i|}{\delta\sqrt N}\right).
\end{equation}
With $w_\eta$ as above, put
\[
 g_\delta(x)=K_{H_1}(x)\chi\left(\frac{|x|}{\delta}\right)
 \quad(x\ne0),\qquad g_\delta(0)=0,
\]
and, for $v\in C([0,1])$, define
\[
 \Phi^{\eta,\delta}_{M,k}(v)
 :=\int_{I\times J}w_\eta(s,t)g_\delta(v(t)-v(s))\dd s\dd t,
 \qquad C^{(\eta,\delta)}_{M,k,1}:=\Phi^{\eta,\delta}_{M,k}(B).
\]
Both $w_\eta$ and $g_\delta$ are bounded and Lipschitz. In particular,
\begin{equation}\label{eq:cutoff-block-functional-continuity}
 |\Phi^{\eta,\delta}_{M,k}(v)-\Phi^{\eta,\delta}_{M,k}(\widetilde v)|
 \le 2|I||J|\,\|w_\eta\|_\infty
       \operatorname{Lip}(g_\delta)\|v-\widetilde v\|_\infty.
\end{equation}

We first compare the lattice kernels with their continuous limits.
By \eqref{fr:lattice-asymptotic}, for $\delta\sqrt N\ge2$,
\begin{equation}\label{eq:G-away-diagonal}
\begin{split}
 \sigma_N(\delta)
 &:=\sup_{\substack{z\in\mathbb Z\\|z|\ge\delta\sqrt N}}
      \left|N^{1-H_1}\gamma_{H_1}(z)
                  -K_{H_1}\left(\frac z{\sqrt N}\right)\right|\\
 &\le C_{H_1}\delta^{2H_1-4}N^{-1}\longrightarrow0.
\end{split}
\end{equation}
The temporal estimate \eqref{eq:gamma-away-diagonal} applies unchanged.
Define the midpoint sum
\[
 R_N^{\eta,\delta}
 :=\frac1{N^2}\sum_{\tau_i^N\in I}\sum_{\tau_j^N\in J}
    w_\eta(\tau_i^N,\tau_j^N)
    g_\delta\left(S^N_{j/N}-S^N_{i/N}\right).
\]
Using $S^N_{i/N}=S_i/\sqrt N$ and
$N^{-a_h}=N^{-2}N^{2-2H_0}N^{1-H_1}$, the two kernel estimates give
\begin{align*}
 |C^{(N),\eta,\delta}_{M,k}-R_N^{\eta,\delta}|
 &\le\rho_N(\eta)\bigl(\|g_\delta\|_\infty+\sigma_N(\delta)\bigr)
       +\|w_\eta\|_\infty\sigma_N(\delta)\\
 &\longrightarrow0.
\end{align*}
Indeed, each cutoff kills the corresponding singular region, and the
number of summands is at most $N^2$.

To control the quadrature error, let
$\Delta_i^N=((i-1)/N,i/N]$ and
\[
 I_N=\bigcup_{\tau_i^N\in I}\Delta_i^N,\qquad
 J_N=\bigcup_{\tau_j^N\in J}\Delta_j^N.
\]
Only the cells containing the endpoints of $I,J$ can differ, so,
with $|\cdot|$ denoting Lebesgue measure,
\[
 |I_N\mathbin{\triangle}I|+|J_N\mathbin{\triangle}J|\le\frac4N.
\]
For $s\in\Delta_i^N$ and $t\in\Delta_j^N$,
\[
 |s-\tau_i^N|+|t-\tau_j^N|\le\frac1N,\qquad
 |S^N_s-S^N_{i/N}|+|S^N_t-S^N_{j/N}|\le\frac2{\sqrt N}.
\]
Decomposing the integral in $\Phi^{\eta,\delta}_{M,k}(S^N)$ along the cells $\Delta_i^N\times \Delta_j^N$, and using the Lipschitz constant of
$w_\eta$ for the metric $|s-s'|+|t-t'|$, yields
\begin{align*}
 |R_N^{\eta,\delta}-\Phi^{\eta,\delta}_{M,k}(S^N)|
 &\le\frac4N\|w_\eta\|_\infty\|g_\delta\|_\infty
       +\frac1N\operatorname{Lip}(w_\eta)\|g_\delta\|_\infty+\frac2{\sqrt N}\|w_\eta\|_\infty
                             \operatorname{Lip}(g_\delta)\le C_{h,\eta,\delta}N^{-1/2}.
\end{align*}
Consequently,
\begin{equation}\label{eq:fixed-level-riemann-error}
 |C^{(N),\eta,\delta}_{M,k}-\Phi^{\eta,\delta}_{M,k}(S^N)|
 \stackrel{N\to \infty}{\longrightarrow}0.
\end{equation}

\smallskip

Since $\mathcal F_{M_0}$ is finite, Donsker's theorem and the continuous
mapping theorem, using \eqref{eq:cutoff-block-functional-continuity}, give
\[
 \left(S^N,(\Phi^{\eta,\delta}_{M,k}(S^N))_{(M,k)\in\mathcal F_{M_0}}\right)
 \ToLaw
 \left(B,(C^{(\eta,\delta)}_{M,k,1})_{(M,k)\in\mathcal F_{M_0}}\right).
\]
For fixed $\eta,\delta>0$, the error in
\eqref{eq:fixed-level-riemann-error} is bounded by a deterministic
quantity tending to zero, independently of the walk. Slutsky's theorem
therefore yields
\begin{equation}\label{eq:double-cutoff-convergence}
 \left(S^N,(C^{(N),\eta,\delta}_{M,k})_{(M,k)\in\mathcal F_{M_0}}\right)
 \ToLaw
 \left(B,(C^{(\eta,\delta)}_{M,k,1})_{(M,k)\in\mathcal F_{M_0}}\right).
\end{equation}

\smallskip\noindent\textit{\underline{Removing the spatial cutoff at fixed $\eta$ ($H_1>\frac12$)}.}
For every integer $r\ge\eta N$, \eqref{eq:LLT-sup} and
\eqref{fr:lattice-envelope} give
\begin{equation}\label{eq:spatial-cutoff-error}
\begin{split}
 &N^{1-H_1}E\left[\gamma_{H_1}(S_r)
                    \one_{\{|S_r|\le2\delta\sqrt N\}}\right]\\
 &\quad=N^{1-H_1}\sum_{|z|\le2\delta\sqrt N}
                       \gamma_{H_1}(z)p_r(z)\\
 &\quad\le C_{H_1}N^{1-H_1}r^{-1/2}
                         (1+2\delta\sqrt N)^{2H_1-1}\\
 &\quad\le C_{H_1,\eta}(\delta+N^{-1/2})^{2H_1-1}.
\end{split}
\end{equation}
Moreover, the temporal envelope in \eqref{fr:lattice-envelope} implies
\[
 0\le N^{2-2H_0}\gamma_{H_0}(r)
               \chi\left(\frac r{N\eta}\right)
 \le C_{H_0}\eta^{2H_0-2}\qquad(r\ge1).
\]
Since $0\le1-\chi(|z|/\delta)\le\one_{\{|z|\le2\delta\}}$,
stationarity of the walk increments and
\eqref{eq:spatial-cutoff-error} therefore yield
\begin{align*}
 0&\le E[C^{(N),\eta}_{M,k}-C^{(N),\eta,\delta}_{M,k}]\\
 &\le\frac{C_{H_0,\eta}}{N^2}
       \sum_{\substack{\tau_i^N\in I,\ \tau_j^N\in J\\j-i\ge\eta N}}
       N^{1-H_1}E\left[\gamma_{H_1}(S_{j-i})
            \one_{\{|S_{j-i}|\le2\delta\sqrt N\}}\right]\\
 &\le C_{h,\eta}(\delta+N^{-1/2})^{2H_1-1}.
\end{align*}
For the Brownian counterpart, whenever $t-s\ge\eta$,
\begin{align*}
 E\left[K_{H_1}(B_t-B_s)\one_{\{|B_t-B_s|\le2\delta\}}\right]
 &=c_{H_1}\int_{-2\delta}^{2\delta}
              |z|^{2H_1-2}\heat_{t-s}(z)\dd z\\
 &\le\frac{2c_{H_1}(2\delta)^{2H_1-1}}
             {(2H_1-1)\sqrt{2\pi\eta}}.
\end{align*}
Thus
\[
 0\le E[C^{(\eta)}_{M,k,1}-C^{(\eta,\delta)}_{M,k,1}]
 \le C_{h,\eta}\delta^{2H_1-1}.
\]

Let $f$ be bounded and $1$-Lipschitz on
$C([0,1])\times\mathbb R^{\mathcal F_{M_0}}$, for the sum of the
uniform norm on paths and the $\ell^1$ norm on the family of blocks.
The triangle inequality and \eqref{eq:double-cutoff-convergence} give
\begin{align*}
 &\limsup_{N\to\infty}
   \Bigl|E\bigl[f(S^N,(C^{(N),\eta}_{M,k})_{(M,k)\in\mathcal F_{M_0}})\bigr]\\
 &\hspace{3cm}{}-E\bigl[f(B,(C^{(\eta)}_{M,k,1})_{(M,k)\in\mathcal F_{M_0}})\bigr]\Bigr|\\
 &\quad\le\sum_{(M,k)\in\mathcal F_{M_0}}
       \limsup_{N\to\infty}E[C^{(N),\eta}_{M,k}-C^{(N),\eta,\delta}_{M,k}]\\
 &\qquad{}+\sum_{(M,k)\in\mathcal F_{M_0}}
       E[C^{(\eta)}_{M,k,1}-C^{(\eta,\delta)}_{M,k,1}]\\
 &\quad\le C_{h,\eta,M_0}\delta^{2H_1-1}
       \longrightarrow0\qquad(\delta\downarrow0).
\end{align*}
Since $2H_1-1>0$, we have proved, for every fixed $\eta>0$,
\[
 \left(S^N,(C^{(N),\eta}_{M,k})_{(M,k)\in\mathcal F_{M_0}}\right)
 \ToLaw
 \left(B,(C^{(\eta)}_{M,k,1})_{(M,k)\in\mathcal F_{M_0}}\right).
\]
Thus \eqref{eq:fixed-level-uncentered-cutoff-convergence} holds also
for $H_1>\frac12$.

\smallskip\noindent\textit{\underline{Removal of the temporal cutoff}.}
We now treat both spatial regimes together. Fix
$0<\eta<2^{-M_0-1}$.

For a fixed $(M,k)\in\mathcal F_{M_0}$, use the index sets
$\mathcal I,\mathcal J$ introduced at the beginning of the proof of
Proposition~\ref{prop:block-moments}. If either set is empty, then
\[
 C^{(N)}_{M,k}=C^{(N),\eta}_{M,k}=0.
\]
We may therefore assume that both sets are nonempty. Let
$\theta=(2k-1)2^{-M}$ be the common endpoint of the intervals
$I_{M,k}$ and $J_{M,k}$. For $i\in\mathcal I$ and $j\in\mathcal J$,
the terms affected by the cutoff satisfy
\begin{gather*}
 0\le1-\chi\left(\frac{j-i}{N\eta}\right)
       \le\one_{\{1\le j-i<2\eta N\}},\\
 j-i<2\eta N\quad\Longrightarrow\quad
 \theta-2\eta<\tau_i^N<\theta\le\tau_j^N<\theta+2\eta.
\end{gather*}
If $2\eta N\le1$, there is no such pair and
$C^{(N)}_{M,k}=C^{(N),\eta}_{M,k}$. Otherwise set
\[
 L_{\eta,N}=\lceil2\eta N\rceil\le4\eta N.
\]
We also use the index $b=\max\mathcal I$ and the walk segments
$S^-,S^+$ defined at the beginning of the proof of
Proposition~\ref{prop:block-moments}. For $i=b-r\in\mathcal I$ and
$j=b+1+s\in\mathcal J$,
\[
 j-i=1+r+s<2\eta N
 \quad\Longrightarrow\quad0\le r,s<L_{\eta,N}.
\]
Extend these walk segments by independent increments when necessary.
Positivity gives the pathwise comparison
\[
 0\le C^{(N)}_{M,k}-C^{(N),\eta}_{M,k}
 \le N^{-a_h}\sum_{r,s=0}^{L_{\eta,N}-1}
        \gamma_{H_0}(1+r+s)\gamma_{H_1}(S^-_r-S^+_s).
\]
Conditionally on $S^+_0=X_{b+1}$, the two walks are independent, so Lemma~\ref{lem:weighted-block} yields, for every $p\ge2$,
\begin{align}
 \|C^{(N)}_{M,k}-C^{(N),\eta}_{M,k}\|_p
 &\le C_hp^{1-H_1}N^{-a_h}L_{\eta,N}^{a_h}\notag\\
 &\le C_h\sqrt p\,(4\eta)^{a_h}
 \le C_h\sqrt p\,\eta^{a_h}.
 \label{eq:discrete-temporal-cutoff-error}
\end{align}
On the Brownian side, $2\eta<2^{-M}$ for every $M\le M_0$,
so \eqref{eq:Brownian-corner} at $t=1$ gives
\begin{equation}\label{eq:Brownian-temporal-cutoff-error}
 \|C_{M,k,1}-C^{(\eta)}_{M,k,1}\|_p
 \le C_h\sqrt p\,\eta^{a_h}.
\end{equation}

\smallskip\noindent\textit{\underline{Centering and (partial) conclusion}.}
We now pass to centered blocks for both spatial regimes. Recall that, since
$0\le\chi\le1$,
\[
 0\le C^{(N),\eta}_{M,k}\le C^{(N)}_{M,k},\qquad
 0\le C^{(\eta)}_{M,k,1}\le C_{M,k,1}.
\]
Proposition~\ref{prop:block-moments} and
\eqref{eq:Brownian-block-Lp} therefore give
\begin{equation}\label{manque-1}
 \sup_{N\ge1}\|C^{(N),\eta}_{M,k}\|_2
       +\|C^{(\eta)}_{M,k,1}\|_2
 \le C_h2^{-a_hM},
\end{equation}
uniformly in $\eta>0$. This bound allows us to pass from
convergence in law to convergence of the means. Indeed, fix $\eta>0$.
For every $R>0$, the function $z\mapsto(z\vee0)\wedge R$ is bounded
and continuous, so \eqref{eq:fixed-level-uncentered-cutoff-convergence}
implies
\[
 E[C^{(N),\eta}_{M,k}\wedge R]
 \longrightarrow E[C^{(\eta)}_{M,k,1}\wedge R].
\]
Since $0\le z-(z\wedge R)\le z^2/R$ for $z\ge0$, the triangle
inequality and \eqref{manque-1} give
\begin{align*}
 &\limsup_{N\to\infty}
   \left|E[C^{(N),\eta}_{M,k}]-E[C^{(\eta)}_{M,k,1}]\right|\\
 &\quad\le\frac1R\left(
       \sup_{N\ge1}E[(C^{(N),\eta}_{M,k})^2]
       +E[(C^{(\eta)}_{M,k,1})^2]\right)
 \le\frac{C_h2^{-2a_hM}}R.
\end{align*}
Letting $R\to\infty$ proves, for every fixed $\eta>0$,
\[
 E[C^{(N),\eta}_{M,k}]\longrightarrow E[C^{(\eta)}_{M,k,1}].
\]
With $A^{(\eta)}_{M,k,1}
=C^{(\eta)}_{M,k,1}-E[C^{(\eta)}_{M,k,1}]$, subtraction of these
deterministic means and Slutsky's theorem yield
\begin{equation}\label{eq:fixed-level-centered-cutoff-convergence}
 \left(S^N,(A^{(N),\eta}_{M,k})_{(M,k)\in\mathcal F_{M_0}}\right)
 \ToLaw
 \left(B,(A^{(\eta)}_{M,k,1})_{(M,k)\in\mathcal F_{M_0}}\right).
\end{equation}

Centering costs at most a factor of two: for any $X\in L^p$,
\[
 \|X-E[X]\|_p\le\|X\|_p+|E[X]|\le2\|X\|_p.
\]
Applying this to the two differences of uncentered blocks and
using \eqref{eq:discrete-temporal-cutoff-error}-\eqref{eq:Brownian-temporal-cutoff-error}
gives
\begin{equation}\label{eq:fixed-block-corner}
 \sup_N\|A^{(N)}_{M,k}-A^{(N),\eta}_{M,k}\|_p
 +\|A_{M,k,1}-A^{(\eta)}_{M,k,1}\|_p
 \le C_h\sqrt p\,\eta^{a_h}.
\end{equation}
For a bounded $1$-Lipschitz function $f$ as above, the same comparison,
using \eqref{eq:fixed-level-centered-cutoff-convergence} and then
\eqref{eq:fixed-block-corner} with $p=2$, yields
\begin{align*}
 &\limsup_{N\to\infty}
   \Bigl|E\bigl[f(S^N,(A^{(N)}_{M,k})_{(M,k)\in\mathcal F_{M_0}})\bigr]\\
 &\hspace{3cm}{}-E\bigl[f(B,(A_{M,k,1})_{(M,k)\in\mathcal F_{M_0}})\bigr]\Bigr|\\
 &\quad\le\sum_{(M,k)\in\mathcal F_{M_0}}\Bigl(
       \sup_N\|A^{(N)}_{M,k}-A^{(N),\eta}_{M,k}\|_2
       +\|A_{M,k,1}-A^{(\eta)}_{M,k,1}\|_2\Bigr)\\
 &\quad\le C_{h,M_0}\eta^{a_h}\longrightarrow0
                  \qquad(\eta\downarrow0).
\end{align*}
Hence
\begin{equation}\label{eq:fixed-level-block-convergence}
 \left(S^N,(A^{(N)}_{M,k})_{(M,k)\in\mathcal F_{M_0}}\right)
 \ToLaw
 \left(B,(A_{M,k,1})_{(M,k)\in\mathcal F_{M_0}}\right).
\end{equation}

\

\noindent\textit{\underline{About the linear occupation functionals}.}
For $\varphi\in\D$, define
\[
 F_\varphi(v)=\int_0^1G_\varphi(1-u,v(u))\dd u,
 \qquad v\in C([0,1]).
\]
Recall from \eqref{eq:G-phi} that $G_\varphi$ is bounded and uniformly
continuous. Denote its modulus of continuity by
\[
 \omega_\varphi(\delta)
 :=\sup_{|s-s'|+|y-y'|\le\delta}
       |G_\varphi(s,y)-G_\varphi(s',y')|,
 \qquad \omega_\varphi(\delta)\longrightarrow0
       \quad(\delta\downarrow0).
\]
Then
\[
 |F_\varphi(v)-F_\varphi(\widetilde v)|
 \le\omega_\varphi(\|v-\widetilde v\|_\infty),
\]
so $F_\varphi$ is continuous for the uniform norm. The substitution
$s=1-u$ in \eqref{eq:ell-B-general} also gives
$F_\varphi(B)=\ell^B_{1,0}(\varphi)$.

To compare $F_\varphi(S^N)$ with the discrete linear functional,
\eqref{eq:H-inner}, \eqref{eq:ell-N} and \eqref{eq:hNS-explicit} give
\begin{equation}\label{eq:ell-cell-sum}
 \ell^{N}_{1,0}(\varphi)
 =\sqrt N\sum_{n=1}^N
   \int_{\mathfrak{I}_{1,n}^N}\int_{\mathfrak{J}_{S_n}^N}
          G_\varphi(s,y)\dd y\dd s.
\end{equation}
For $(s,y)\in\mathfrak{I}_{1,n}^N\times\mathfrak{J}_{S_n}^N$,
\[
 |s-(1-n/N)|\le N^{-1},\qquad
 |y-S_n/\sqrt N|\le N^{-1/2}.
\]
Each cell has Lebesgue measure $N^{-3/2}$, so
\[
 \left|\ell^{N}_{1,0}(\varphi)
       -\frac1N\sum_{n=1}^N
          G_\varphi\left(1-\frac nN,\frac{S_n}{\sqrt N}\right)\right|
 \le\omega_\varphi(N^{-1}+N^{-1/2}).
\]
For $u\in[(n-1)/N,n/N]$, the interpolation in
\eqref{eq:rescaled-walk-path} and $|S_n-S_{n-1}|=1$ imply
\[
 \left|S^N_u-\frac{S_n}{\sqrt N}\right|\le N^{-1/2},
 \qquad
 \left|(1-u)-\left(1-\frac nN\right)\right|\le N^{-1}.
\]
Integrating over these $N$ time intervals yields
\[
 \left|\frac1N\sum_{n=1}^N
          G_\varphi\left(1-\frac nN,\frac{S_n}{\sqrt N}\right)
       -F_\varphi(S^N)\right|
 \le\omega_\varphi(N^{-1}+N^{-1/2}).
\]
Adding the two errors gives
\begin{equation}\label{eq:ell-path-approximation}
 \left|\ell^{N}_{1,0}(\varphi)-F_\varphi(S^N)\right|
 \le2\omega_\varphi(N^{-1}+N^{-1/2})\longrightarrow0,
\end{equation}
uniformly over all walk trajectories.

\

\noindent\textit{\underline{Conclusion}.}
We now work in the product space
$C([0,1])\times\R^{\mathcal F_{M_0}}\times\R^r$, equipped with
the uniform norm on paths and the maximum norm on each
finite-dimensional factor. Since every $F_{\varphi_j}$ is continuous,
the map
\[
 \left(v,(a_{M,k})_{(M,k)\in\mathcal F_{M_0}}\right)
 \longmapsto
 \left(v,(a_{M,k})_{(M,k)\in\mathcal F_{M_0}},
             (F_{\varphi_j}(v))_{j\le r}\right)
\]
is continuous from $C([0,1])\times\R^{\mathcal F_{M_0}}$ into
this product space. Thus the continuous mapping theorem applied to
\eqref{eq:fixed-level-block-convergence} gives
\begin{align*}
 &\left(S^N,(A^{(N)}_{M,k})_{(M,k)\in\mathcal F_{M_0}},
              (F_{\varphi_j}(S^N))_{j\le r}\right)\\
 &\hspace{1cm}\ToLaw
 \left(B,(A_{M,k,1})_{(M,k)\in\mathcal F_{M_0}},
              (\ell^B_{1,0}(\varphi_j))_{j\le r}\right),
\end{align*}
where we have used $F_{\varphi_j}(B)=\ell^B_{1,0}(\varphi_j)$
for each $j$.

It remains to replace $F_{\varphi_j}(S^N)$ by
$\ell^{N}_{1,0}(\varphi_j)$ in the left-hand tuple.
Since $r$ is fixed, \eqref{eq:ell-path-approximation} implies
\[
 \max_{1\le j\le r}
 \left|\ell^{N}_{1,0}(\varphi_j)-F_{\varphi_j}(S^N)\right|
 \le 2\max_{1\le j\le r}
          \omega_{\varphi_j}(N^{-1}+N^{-1/2})
 \longrightarrow0.
\]
Slutsky's theorem therefore yields precisely
\eqref{eq:fixed-level-convergence}.
\end{proof}

\subsection{Summation of the discrete levels and the joint limit}

\begin{theorem}[Invariance principle for the centered self-intersection]
\label{thm:YN-invariance}\label{prop:one-path}
Let $B$ be a standard Brownian motion, and let $Y_1(B)$ be the Brownian functional constructed in
Theorem~\ref{thm:Brownian-self-intersection}.
For every $p\ge2$ and integer $M_0\ge0$,
\begin{equation}\label{eq:YN-tail-Lp}
 \sup_N\left\|Y^N_1-2\sum_{M\le M_0}\sum_k A^{(N)}_{M,k}\right\|_p
 \le C_{h}\sqrt p\,2^{-M_0(a_{h}-\frac12)}.
\end{equation}
For any finite family $\varphi_1,\ldots,\varphi_r\in\D$,
\begin{equation}\label{eq:path-YN-conv}
 (S^N,Y^N_1,\ell^{N}_{1,0}(\varphi_1),\ldots,\ell^{N}_{1,0}(\varphi_r))
 \ToLaw
 (B,Y_1(B),\ell^B_{1,0}(\varphi_1),\ldots,\ell^B_{1,0}(\varphi_r)).
\end{equation}
Furthermore, for $p\ge2$,
\begin{equation}\label{eq:YN-subgaussian-Lp}
 \sup_N\|Y^N_1\|_p\le C_{h}\sqrt p,
\end{equation}
and, for every $\lambda\in\R$,
\begin{equation}\label{eq:YN-exp-unif}
 \sup_NE[e^{\lambda Y^N_1}]\le e^{C_{h}\lambda^2}.
\end{equation}
\end{theorem}

\begin{proof}
At a fixed level $M$, the centered discrete blocks depend only on
increments inside disjoint parent intervals, so they are independent
in $k$.
Their moment bounds in Proposition~\ref{prop:block-moments} allow us
to apply \eqref{eq:subgaussian-sum-Lp} directly, with
$n=2^{M-1}$ and $s_k=2^{-a_hM}$. This gives, uniformly in $N$,
\begin{equation}\label{eq:level-subgaussian-Lp}
 \begin{split}
 \left\|\sum_{k=1}^{2^{M-1}} A^{(N)}_{M,k}\right\|_p
 &\le C_h\sqrt p\,
         \left(2^{M-1}2^{-2a_hM}\right)^{1/2}\\
 &\le C_h\sqrt p\,2^{-M(a_h-\frac12)},\qquad p\ge2.
 \end{split}
\end{equation}

By \eqref{eq:YN-dyadic}, summing over $M>M_0$ with Minkowski's
inequality proves \eqref{eq:YN-tail-Lp}, since $a_h>\frac12$.
Taking $M_0=0$ gives \eqref{eq:YN-subgaussian-Lp}, and centering
allows us to apply Lemma~\ref{lem:subgaussian} once more to obtain
\eqref{eq:YN-exp-unif}.

We now prove the joint convergence. For an integer $M_0\ge1$, put
\[
 \Sigma^{(N)}_{M_0}
    =2\sum_{M=1}^{M_0}\sum_{k=1}^{2^{M-1}}A^{(N)}_{M,k},
 \qquad
 \Sigma_{M_0}
    =2\sum_{M=1}^{M_0}\sum_{k=1}^{2^{M-1}}A_{M,k,1}.
\]
For each fixed $M_0$, Lemma~\ref{lem:fixed-level} and continuity
of the finite summation map give
\begin{align*}
 &\left(S^N,\Sigma^{(N)}_{M_0},
                (\ell^{N}_{1,0}(\varphi_j))_{j\le r}\right)\ToLaw
 \left(B,\Sigma_{M_0},
                (\ell^B_{1,0}(\varphi_j))_{j\le r}\right)
       \qquad(N\to\infty).
\end{align*}
On the other hand, \eqref{eq:YN-tail-Lp} and
\eqref{eq:Brownian-Y-tail-Lp}, with $p=2$ and $t=1$, yield
\[
 \sup_N\|Y^N_1-\Sigma^{(N)}_{M_0}\|_2
       +\|Y_1(B)-\Sigma_{M_0}\|_2
 \le C_h2^{-M_0(a_h-\frac12)}.
\]

Let $f$ be a bounded $1$-Lipschitz function on
$C([0,1])\times\R\times\R^r$, for the sum of the uniform norm,
the absolute value and the maximum norm on the respective factors.
Replacing only the self-intersection coordinate and using
$\|X\|_1\le\|X\|_2$, we obtain
\begin{align*}
 &\Bigl|E\bigl[f(S^N,Y^N_1,
                 (\ell^{N}_{1,0}(\varphi_j))_{j\le r})\bigr]-E\bigl[f(B,Y_1(B),
                 (\ell^B_{1,0}(\varphi_j))_{j\le r})\bigr]\Bigr|\\
 &\quad\le\|Y^N_1-\Sigma^{(N)}_{M_0}\|_2
                  +\|Y_1(B)-\Sigma_{M_0}\|_2\\
 &\qquad+\Bigl|E\bigl[f(S^N,\Sigma^{(N)}_{M_0},
                 (\ell^{N}_{1,0}(\varphi_j))_{j\le r})\bigr]-E\bigl[f(B,\Sigma_{M_0},
                 (\ell^B_{1,0}(\varphi_j))_{j\le r})\bigr]\Bigr|.
\end{align*}
For fixed $M_0$, the last absolute difference tends to zero as
$N\to\infty$ by the preceding joint convergence. The tail estimate
is uniform in $N$, so
\begin{align*}
 &\limsup_{N\to\infty}
   \Bigl|E\bigl[f(S^N,Y^N_1,
                 (\ell^{N}_{1,0}(\varphi_j))_{j\le r})\bigr]-E\bigl[f(B,Y_1(B),
                 (\ell^B_{1,0}(\varphi_j))_{j\le r})\bigr]\Bigr|\le C_h2^{-M_0(a_h-\frac12)}.
\end{align*}
We can now let $M_0\to\infty$: the right-hand side tends to zero
because $a_h>\frac12$. This proves
\eqref{eq:path-YN-conv}.
\end{proof}

Theorem~\ref{thm:YN-invariance}, including the uniform exponential bound \eqref{eq:YN-exp-unif}, together with Proposition~\ref{prop:counterterm}, suffices to pass to the limit
in the $\mathcal S$-transform \eqref{gau-pai} at $(t,x)=(1,0)$.

\subsection{\texorpdfstring{The $\mathcal S$-transform at a general point}{The S-transform at a general point}}
\label{subsec:general-occupation}

We now identify the limit of \eqref{gau-pai} at an arbitrary point
$(t,x)$, with $t>0$ and $x\in\R$ fixed.

Recall $t_N=N_t/N$ from \eqref{eq:Nt} and $x_N=\frac{\widetilde{N}_x}{\sqrt N}$ from \eqref{eq:xN}. We have
\begin{equation}\label{eq:xiN}
 0\le x-x_N<\frac1{\sqrt N},\qquad x_N\longrightarrow x.
\end{equation}

Recall that, for every integer $m\ge1$, \eqref{eq:rescaled-walk-path}
defines the continuous, piecewise linear path $S^m$ on $[0,1]$ by
\begin{equation}\label{eq:bar-SN}
 S^m_{(j+\theta)/m}
 =m^{-\frac12}\bigl((1-\theta)S_j+\theta S_{j+1}\bigr),
 \qquad 0\le j<m,\quad 0\le\theta\le1.
\end{equation}
In the rest of this subsection, we take $N$ large enough that
$N_t\ge1$ and apply this definition with $m=N_t$. Since $t_N=N_t/N$,
\begin{equation}\label{eq:two-scalings}
 \frac{S_j}{\sqrt N}=\sqrt{t_N}\,S^{N_t}_{j/N_t},
 \qquad 0\le j\le N_t.
\end{equation}

\begin{lemma}[Occupation functionals at a general point]
\label{lem:ell-general}
Fix $\varphi\in\D$ and, for $w\in C([0,1])$, set
\begin{align*}
 F^N_{t,x}(w)&=t_N\int_0^1
     G_\varphi\bigl(t-t_Nu,\, x_N+\sqrt{t_N}\,w_u\bigr)\dd u,\\
 F_{t,x}(w)&=t\int_0^1
     G_\varphi\bigl(t-tu,\,  x+\sqrt t\,w_u\bigr)\dd u.
\end{align*}
Then the following hold.
\begin{enumerate}
\renewcommand{\theenumi}{\roman{enumi}}
\renewcommand{\labelenumi}{(\theenumi)}
\item Uniformly over the paths of the walk,
\begin{equation}\label{eq:ell-cell-riemann}
 \bigl|\ell^{N}_{t,x}(\varphi)-F^N_{t,x}(S^{N_t})\bigr|
 \le\frac{2t}N\|\partial_tG_\varphi\|_\infty
   +\frac{2t}{\sqrt N}\|\partial_xG_\varphi\|_\infty
 \longrightarrow0.
\end{equation}
\item The map $F_{t,x}$ is Lipschitz on $C([0,1])$, and for every $L>0$,
\begin{equation}\label{eq:FN-to-F}
 \sup_{\|w\|_\infty\le L}|F^N_{t,x}(w)-F_{t,x}(w)|\longrightarrow0.
\end{equation}
\item If $\bar B$ is a standard Brownian motion on $[0,1]$ and
$B_r=\sqrt t\,\bar B_{r/t}$, $0\le r\le t$, then
$F_{t,x}(\bar B)=\ell^B_{t,x}(\varphi)$.
\end{enumerate}
\end{lemma}

\begin{proof}
(i) By \eqref{eq:hNS-explicit} and \eqref{eq:ell-N},
\[
 \ell^{N}_{t,x}(\varphi)
 =\sqrt N\sum_{n=1}^{N_t}
   \int_{\mathfrak{I}_{t,n}^N}\int_{\mathfrak{J}_{\widetilde{N}_x+S_n}^N}G_\varphi(s,y)\dd y\dd s.
\]
For $1\le n\le N_t$ and
$(s,y)\in\mathfrak{I}_{t,n}^N\times\mathfrak{J}_{\widetilde{N}_x+S_n}^N$,
the cell definitions give
\[
 \left|s-\left(t-\frac nN\right)\right|\le\frac1N,
 \qquad
 \left|y-\left(x_N+\frac{S_n}{\sqrt N}\right)\right|
       \le\frac1{\sqrt N}.
\]
The bounded derivatives of $G_\varphi$ therefore imply
\begin{align*}
 &\left|G_\varphi(s,y)
       -G_\varphi\left(t-\frac nN,x_N+\frac{S_n}{\sqrt N}\right)\right|\le\frac1N\|\partial_tG_\varphi\|_\infty
                +\frac1{\sqrt N}\|\partial_xG_\varphi\|_\infty.
\end{align*}
Since
\[
 \sqrt N\int_{\mathfrak{I}_{t,n}^N}
           \int_{\mathfrak{J}_{\widetilde{N}_x+S_n}^N}\dd y\dd s
 =\sqrt N\,\frac1N\,\frac1{\sqrt N}=\frac1N,
\]
we obtain the error bound for one cell:
\begin{align*}
 &\Biggl|\sqrt N\int_{\mathfrak{I}_{t,n}^N}
           \int_{\mathfrak{J}_{\widetilde{N}_x+S_n}^N}
                 G_\varphi(s,y)\dd y\dd s-\frac1N
       G_\varphi\left(t-\frac nN,x_N+\frac{S_n}{\sqrt N}\right)\Biggr|\\
 &\qquad\le\frac1N\left(
       \frac1N\|\partial_tG_\varphi\|_\infty
       +\frac1{\sqrt N}\|\partial_xG_\varphi\|_\infty\right).
\end{align*}
Summing over $n=1,\ldots,N_t$ and using $N_t/N=t_N$ yields
\[
 \Bigl|\ell^{N}_{t,x}(\varphi)
   -\frac1N\sum_{n=1}^{N_t}
      G_\varphi\Bigl(t-\frac nN,\, x_N+\frac{S_n}{\sqrt N}\Bigr)\Bigr|
 \le t_N\Bigl(\frac1N\|\partial_tG_\varphi\|_\infty
   +\frac1{\sqrt N}\|\partial_xG_\varphi\|_\infty\Bigr).
\]

Next, put $u_n=n/N_t$ for $0\le n\le N_t$. By
\eqref{eq:two-scalings},
\[
 t_Nu_n=\frac nN,\qquad
 \sqrt{t_N}\,S^{N_t}_{u_n}=\frac{S_n}{\sqrt N},\qquad
 t_N(u_n-u_{n-1})=\frac1N\quad(1\le n\le N_t).
\]
Thus the preceding discrete sum is the right-endpoint Riemann sum
for $F^N_{t,x}(S^{N_t})$. For $u\in[u_{n-1},u_n]$,
\eqref{eq:bar-SN} and $|S_n-S_{n-1}|=1$ give
\begin{align*}
& |t_Nu-t_Nu_n|\le\frac{t_N}{N_t}=\frac1N,\qquad
 \sqrt{t_N}\left|S^{N_t}_u-S^{N_t}_{u_n}\right|
=\frac{n-N_tu}{\sqrt N}|S_n-S_{n-1}|
 \le\frac1{\sqrt N}.
\end{align*}
Using these bounds and integrating over the $N_t$ intervals yields
\begin{align*}
 \left|\frac1N\sum_{n=1}^{N_t}
       G_\varphi\left(t-\frac nN,x_N+\frac{S_n}{\sqrt N}\right)
                 -F^N_{t,x}(S^{N_t})\right|
 &\le t_N\sum_{n=1}^{N_t}\int_{u_{n-1}}^{u_n}
       \left(\frac1N\|\partial_tG_\varphi\|_\infty
           +\frac1{\sqrt N}\|\partial_xG_\varphi\|_\infty\right)\dd u\\
 &=t_N\left(\frac1N\|\partial_tG_\varphi\|_\infty
           +\frac1{\sqrt N}\|\partial_xG_\varphi\|_\infty\right).
\end{align*}
Adding the cell error and the Riemann error, and using $t_N\le t$,
proves \eqref{eq:ell-cell-riemann}.

(ii) For $w,w'$,
$|F_{t,x}(w)-F_{t,x}(w')|\le t^{\frac32}\|\partial_xG_\varphi\|_\infty\|w-w'\|_\infty$,
so $F_{t,x}$ is Lipschitz. For $\|w\|_\infty\le L$, comparing the two
integrands gives
\[
 |F^N_{t,x}(w)-F_{t,x}(w)|
 \le|t_N-t|\,\|G_\varphi\|_\infty
  +t\Bigl(|t_N-t|\,\|\partial_tG_\varphi\|_\infty
   +\bigl(|x_N-x|+L|\sqrt{t_N}-\sqrt t|\bigr)
                  \|\partial_xG_\varphi\|_\infty\Bigr),
\]
which tends to zero by \eqref{eq:tN} and \eqref{eq:xiN}.

(iii) Substituting $s=t-tu$ gives
\[
F_{t,x}(\bar B)=\int_0^tG_\varphi(s,x+\sqrt t\,\bar B_{(t-s)/t})\dd s=\int_0^tG_\varphi(s,x+B_{t-s})\dd s,
\]
which is $\ell^B_{t,x}(\varphi)$ by \eqref{eq:ell-B-general}.
\end{proof}

We can now state the main conclusion of this section: the
$\mathcal S$-transforms of the discrete partition field converge to
those of the continuous solution. The joint invariance principle
identifies the limiting path functionals, and the exponential moment
bounds justify passage to the limit in the path expectation.

\begin{corollary}[Identification of the limiting $\mathcal S$-transform]
\label{prop:S-convergence}
Assume \eqref{fr:range}, and let $U$ be the continuous field of
Theorem~\ref{thm:continuous-summary}. Then, for every $\varphi\in\D$,
$t>0$ and $x\in\R$,
\begin{equation}\label{eq:S-convergence}
 \mathcal S_{\widehat{\mathsf Z}^N_t(x)}(\varphi)
 \longrightarrow\mathcal S_{U_t(x)}(\varphi)
 \qquad(N\to\infty).
\end{equation}
The limiting transform has the Brownian representation
\[
 \mathcal S_{U_t(x)}(\varphi)
 =e^{R_t}E_B\left[
       \exp\left(\frac12Y_t(B)+\ell^B_{t,x}(\varphi)\right)\right],
\]
where $B$ is a standard Brownian motion, as established in
Proposition~\ref{prop:S-cont}.
\end{corollary}

\begin{proof}

Let $\bar B$ be a standard Brownian motion on $[0,1]$, and set
\[
 B_r=\sqrt t\,\bar B_{r/t},\qquad 0\le r\le t.
\]
Then $B$ is a standard Brownian motion on $[0,t]$, and
Lemma~\ref{lem:ell-general}(iii) gives
$F_{t,x}(\bar B)=\ell^B_{t,x}(\varphi)$.
As $N\to\infty$ we have $N_t\to\infty$, so
Theorem~\ref{thm:YN-invariance}, applied along the indices $N_t$, yields
\[
 (S^{N_t},Y^{N_t}_1)\ToLaw(\bar B,Y_1(\bar B))
 \quad\text{in }C([0,1])\times\R.
\]
The map $(w,y)\mapsto(F_{t,x}(w),y)$ is continuous by
Lemma~\ref{lem:ell-general}(ii). The continuous mapping theorem
therefore gives
\[
 (F_{t,x}(S^{N_t}),Y^{N_t}_1)
 \ToLaw(\ell^B_{t,x}(\varphi),Y_1(\bar B)).
\]

We next replace $F_{t,x}(S^{N_t})$ by
$\ell^{N}_{t,x}(\varphi)$. For each $L>0$, on the event
$\{\|S^{N_t}\|_\infty\le L\}$, Lemma~\ref{lem:ell-general}(i) gives
\begin{align*}
\bigl|\ell^{N}_{t,x}(\varphi)-F_{t,x}(S^{N_t})\bigr|
 &\le
   \bigl|\ell^{N}_{t,x}(\varphi)-F^N_{t,x}(S^{N_t})\bigr|
   +\bigl|F^N_{t,x}(S^{N_t})-F_{t,x}(S^{N_t})\bigr|\\
 &\le\frac{2t}{N}\|\partial_tG_\varphi\|_\infty
       +\frac{2t}{\sqrt N}\|\partial_xG_\varphi\|_\infty
       +\sup_{\|w\|_\infty\le L}
                |F^N_{t,x}(w)-F_{t,x}(w)|.
\end{align*}
For fixed $L$, the last bound tends to zero by
\eqref{eq:FN-to-F}. Thus, for every $\varepsilon>0$ and $L>0$,
for all sufficiently large $N$,
\[
 \Pp\left(\bigl|\ell^{N}_{t,x}(\varphi)
                     -F_{t,x}(S^{N_t})\bigr|>\varepsilon\right)
 \le\Pp\left(\|S^{N_t}\|_\infty>L\right).
\]

The weak convergence of $S^{N_t}$ implies tightness, hence
\[
 \lim_{L\to\infty}\sup_{N:\,N_t\ge1}
       \Pp\left(\|S^{N_t}\|_\infty>L\right)=0.
\]
Taking first the limit superior as $N\to\infty$ and then
$L\to\infty$ proves
\[
 \ell^{N}_{t,x}(\varphi)-F_{t,x}(S^{N_t})\ToProb0.
\]
Slutsky's lemma now yields the required joint convergence:
\[
 (\ell^{N}_{t,x}(\varphi),Y^{N_t}_1)
 \ToLaw(\ell^B_{t,x}(\varphi),Y_1(\bar B))
 \quad\text{in }\R^2.
\]

\

By Lemma~\ref{lem:discrete-scaling}, the exponent in \eqref{eq:S-ZN} is
$\frac{\theta^N_t}2Y^{N_t}_1+\ell^{N}_{t,x}(\varphi)$, with
$\theta^N_t\to\theta_t$. Slutsky's lemma therefore yields
\[
 \frac{\theta^N_t}2Y^{N_t}_1+\ell^{N}_{t,x}(\varphi)
 \ToLaw
 \frac{\theta_t}2Y_1(\bar B)+\ell^B_{t,x}(\varphi)
 =\frac12Y_t(B)+\ell^B_{t,x}(\varphi),
\]
the last equality being \eqref{eq:Yt-alphat-Rt}.
Write $X_N=\exp(\frac{\theta^N_t}2Y^{N_t}_1+\ell^{N}_{t,x}(\varphi))$.
For any $r>1$, \eqref{eq:ell-uniform-bound}, $\theta^N_t\le t^{a_{h}}$ and
\eqref{eq:YN-exp-unif} imply
\[
 \sup_NE_S[X_N^r]
 \le e^{rt\|G_\varphi\|_\infty}
     \sup_NE_S\bigl[e^{\frac{r\theta^N_t}2Y^{N_t}_1}\bigr]
 \le\exp\left(rt\|G_\varphi\|_\infty
                 +\frac14C_{h}r^2t^{2a_{h}}\right)<\infty.
\]
In particular,
$(X_N)_{N\geq 1}$ is uniformly integrable. Weak convergence and this uniform integrability give
convergence of $E_S[X_N]$. Combining \eqref{eq:S-ZN},
Corollary~\ref{cor:RN-t} and \eqref{eq:S-U} proves
\eqref{eq:S-convergence}.
\end{proof}

\section{\texorpdfstring{The two-path invariance principle}{The two-path invariance principle}}
\label{sec:two-path}

The aim of this section is to prove convergence of the second moment:
\[
 \E\bigl[\widehat{\mathsf Z}^N_t(x)^2\bigr]
 \longrightarrow\E\bigl[U_t(x)^2\bigr],
 \qquad t>0,\quad x\in\R.
\]
The identity \eqref{eq:second-moment-ZN}, that is
\[
\E\left[\widehat{\mathsf Z}^N_t(x)^2\right]=e^{2R^N_t}E_{S,\widetilde S}\left[\exp\left(\frac12\bigl(Y^N_t+\widetilde Y^N_t\bigr)+\alpha^N_t(S,\widetilde S)\right)\right],
\]
expresses this moment in terms
of the centered self-intersections of two independent walks and their
mutual intersection $\alpha^N_t$. The exact scaling of Lemma~\ref{lem:discrete-scaling} reduces
the problem to the unit horizon. Passing to the limit requires joint
convergence of these three functionals and uniform integrability of
the exponential. The resulting invariance
principle is stated and proved in Theorem~\ref{tp:main}; the final
subsection deduces the second-moment convergence at a general point.

\textit{For the whole section, we fix $h=(H_0,H_1)$ satisfying the
conditions in \eqref{fr:range}}. The spatial kernel is
$\gamma_{H_1}$, as in Sections~\ref{sec:discrete}
and~\ref{subsec:legall-complet}; for $H_1=\frac12$ it is $\one_{\{0\}}$,
and the spatial sums below are then collision counts.

Let $S,\widetilde S$ be independent simple symmetric walks and let
$B,\widetilde B$ be independent standard Brownian motions, all started
at zero. The variables $\Gamma^N_1,Y^N_1,\alpha^N_1$ are those defined in
\eqref{eq:GammaN}, \eqref{eq:YN}, and \eqref{eq:alphaN}, and
$\beta_N^2=N^{-a_{h}}$. Path expectations are denoted by $E$,
whereas $\E$ denotes expectation with respect to the Gaussian environment.

Throughout, we use the rescaled paths $S^N,\widetilde S^N$ of
\eqref{eq:rescaled-walk-path} and the continuous local-time fields
$\mathcal L^N,\widetilde{\mathcal L}^{N}$ of
\eqref{eq:rescaled-local-time}. Their joint convergence in
Lemma~\ref{lem:local-time-invariance} is needed to identify the limit
of $\alpha^N_1$.

\subsection{\texorpdfstring{Moment and exponential bounds}{Moment and exponential bounds}}

We use the lattice estimates of Lemma~\ref{fr:kernel-lemma}
in Appendix~\ref{app:temporal}, with $H=H_0$ for the temporal
covariance and $H=H_1$ for the spatial one. In particular,
\eqref{fr:lattice-envelope} controls their sums, while
\eqref{fr:lattice-asymptotic} gives, for every $0<\eta<1$,
\begin{equation}\label{tp:gamma-away}
 \sup_{\eta N\le |r|\le N}
 \left|N^{2-2H_0}\gamma_{H_0}(r)
                  -K_{H_0}\left(\frac rN\right)\right|
 \longrightarrow0.
\end{equation}
Indeed, once $\eta N\ge2$, the rescaled remainder in that expansion
is bounded by $C_{H_0}\eta^{2H_0-4}N^{-2}$.

\begin{lemma}[Uniform local-time bound for a difference walk]
\label{tp:difference-local-time}
Let $D=(D_j)_{j\ge0}$ be a random walk on $\mathbb Z$, with
$D_0=d$ under $\Pp_d$ and independent, identically distributed
increments $\Delta D_j=D_j-D_{j-1}$ satisfying
\[
 \Pp_d(\Delta D_j=2)=\Pp_d(\Delta D_j=-2)=\frac14,\qquad
 \Pp_d(\Delta D_j=0)=\frac12.
\]
Thus $D$ stays in the parity class $d+2\mathbb Z$.
For every integer $L\ge1$ and every real $p\ge2$,
\begin{equation}\label{tp:difference-local-time-bound}
 \sup_{d,z\in\mathbb Z}
 \left\|\sum_{j=1}^{L}\one_{\{D_j=z\}}\right\|_{L^p(\Pp_d)}
 \le C\sqrt{pL}.
\end{equation}
\end{lemma}

\begin{proof}
Let $U$ be a simple symmetric walk started at $d$. The increments
$U_{2j}-U_{2j-2}$ are independent and have the law specified above:
each is the sum of two independent symmetric signs. Hence
\[
 (D_j)_{j\ge0}\stackrel{\mathrm{law}}=(U_{2j})_{j\ge0}.
\]
Therefore Lemma~\ref{lem:tanaka-local-time}, with horizon $2L+1$, gives
\begin{align*}
 &\left\|\sum_{j=1}^{L}\one_{\{D_j=z\}}\right\|_p
 =\left\|\sum_{j=1}^{L}\one_{\{U_{2j}=z\}}\right\|_p\le\left\|\sum_{k=0}^{2L}\one_{\{U_k=z\}}\right\|_p
 \le C\sqrt{p(2L+1)}
 \le C\sqrt{pL}.
\end{align*}
That lemma is uniform in both the starting point $d$ and the visited
site $z$, as required.
\end{proof}

The next estimate is the counterpart, for a difference walk, of
Lemma~\ref{lem:RW-occupation}.

\begin{corollary}[Occupation of a difference walk against the kernel]
\label{tp:difference-occupation}
Let $D$ be as in Lemma~\ref{tp:difference-local-time}. Uniformly in the
starting point and in $z\in\mathbb Z$, for $p\ge2$ and $L\ge1$,
\begin{equation}\label{tp:difference-occupation-bound}
 \left\|\sum_{j=1}^{L}\gamma_{H_1}(D_j-z)\right\|_p
 \le C_{h}\,p^{1-H_1}L^{H_1}.
\end{equation}
At $H_1=\frac12$ this is \eqref{tp:difference-local-time-bound}.
\end{corollary}

\begin{proof}
The case $H_1=\frac12$ follows from
\eqref{tp:difference-local-time-bound}, since
$\gamma_{1/2}(D_j-z)=\one_{\{D_j=z\}}$.
Suppose $H_1>\frac12$. The translated walk $D-z$ has the same
increment law and starts at $d-z$, so we may take $z=0$.
Put
\[
 V^D_L(y)=\sum_{j=1}^{L}\one_{\{D_j=y\}},
 \qquad
 \sum_{y\in\mathbb Z}V^D_L(y)=L,\qquad
 \sup_{d,y\in\mathbb Z}\|V^D_L(y)\|_{L^p(\Pp_d)}
 \le C\sqrt{pL}.
\]
The last bound is \eqref{tp:difference-local-time-bound}.
For any real $R\ge1$, split the occupation sum as
\[
 \sum_{j=1}^{L}\gamma_{H_1}(D_j)
 =\sum_{|y|\le R}\gamma_{H_1}(y)V^D_L(y)+\sum_{|y|>R}\gamma_{H_1}(y)V^D_L(y)=: A_R+B_R.
\]
Minkowski's inequality and \eqref{fr:lattice-envelope} give
\begin{align*}
 \|A_R\|_p
 &\le\sum_{|y|\le R}\gamma_{H_1}(y)\|V^D_L(y)\|_p\le C\sqrt{pL}\sum_{|y|\le R}\gamma_{H_1}(y)
 \le C_h\sqrt{pL}\,R^{2H_1-1}.
\end{align*}
On $|y|>R$, \eqref{fr:lattice-envelope} gives
$\gamma_{H_1}(y)\le C_hR^{2H_1-2}$, hence, pathwise,
\[
 0\le B_R
 \le C_hR^{2H_1-2}\sum_{|y|>R}V^D_L(y)
 \le C_hLR^{2H_1-2}.
\]
Consequently,
\[
 \left\|\sum_{j=1}^{L}\gamma_{H_1}(D_j)\right\|_p
 \le C_h\left(\sqrt{pL}\,R^{2H_1-1}+LR^{2H_1-2}\right).
\]
If $L\ge p$, choose $R=\sqrt{L/p}\ge1$. Then
\[
 \sqrt{pL}\,R^{2H_1-1}=p^{1-H_1}L^{H_1},
 \qquad
 LR^{2H_1-2}=p^{1-H_1}L^{H_1},
\]
which proves the required estimate. If $L<p$, use
$\sup_y\gamma_{H_1}(y)\le C_h$ to obtain directly
\[
 0\le\sum_{j=1}^{L}\gamma_{H_1}(D_j)
 \le C_hL
 =C_hL^{H_1}L^{1-H_1}
 \le C_hp^{1-H_1}L^{H_1}.
\]
All constants are independent of the starting point and of the
translation $z$.
\end{proof}

For $|r|<N$, define the weighted intersection count with index difference $r=n-m$ by
\begin{equation}\label{tp:lag-count}
 L_{N,r}
 =\sum_{\substack{1\le n,m\le N\\n-m=r}}
       \gamma_{H_1}(S_n-\widetilde S_m).
\end{equation}

\begin{lemma}[Weighted mutual-intersection estimates]\label{tp:lag-estimates}
Uniformly in $N\ge1$, $|r|<N$, and $p\ge2$,
\begin{equation}\label{tp:lag-bound}
 \norm{L_{N,r}}_p\le C_{h}\,p^{1-H_1}(N-|r|)^{H_1}
 \le C_{h}\,p^{1-H_1}N^{H_1}.
\end{equation}
Consequently, since $p^{1-H_1}\le\sqrt p$,
\begin{equation}\label{tp:alpha-Lp}
 \sup_N\norm{\alpha^N_1}_p\le C_{h}\,p^{1-H_1}\le C_{h}\sqrt p,
\end{equation}
and for every $\lambda\ge0$,
\begin{equation}\label{tp:alpha-exp}
 \sup_N\Ep\left[e^{\lambda\alpha^N_1}\right]
 \le e^{C_h(\lambda+\lambda^2)}.
\end{equation}
\end{lemma}

\begin{proof}
For $r\ge0$,
\[
 L_{N,r}=\sum_{j=1}^{N-r}\gamma_{H_1}(S_{j+r}-\widetilde S_j).
\]
Conditionally on $S_r$, the differences
$D_j=S_{j+r}-\widetilde S_j$ form a walk with the increment law in
Lemma~\ref{tp:difference-local-time}, started from $D_0=S_r$.
This follows because the future increments of $S$ are independent of $S_r$
and of the whole walk $\widetilde S$. The uniformity in the starting point
of \eqref{tp:difference-occupation-bound} gives \eqref{tp:lag-bound}.
Exchange the two walks for $r<0$.

By nonnegativity of $\gamma_{H_0}$ and of $\gamma_{H_1}$, and Minkowski's inequality,
\[
 \norm{\alpha^N_1}_p
 \le N^{-a_{h}}\sum_{|r|<N}\gamma_{H_0}(r)\norm{L_{N,r}}_p
 \le C_{h}p^{1-H_1}N^{H_1-a_{h}}
                   \sum_{|r|<N}\gamma_{H_0}(r).
\]
The mass bound \eqref{fr:lattice-envelope}, with $R=N$, bounds
the last sum by $C_{H_0}N^{2H_0-1}$. Since
$H_1-a_h=1-2H_0$, this proves \eqref{tp:alpha-Lp}.
Finally, $\alpha^N_1$ is nonnegative, so Corollary~\ref{tp:moment-to-exp}
applied to \eqref{tp:alpha-Lp} gives \eqref{tp:alpha-exp}.
\end{proof}

The corresponding Brownian moment and exponential estimates have
already been established in Proposition~\ref{prop:continuous-mutual}.

\subsection{The temporal cutoff and its error}

Recall the cutoff $\chi$ and the Brownian functional $\alpha_1^{(\eta)}(B,\widetilde B)$
from \eqref{tp:alpha-cutoff}. For $0<\eta<\frac12$, define its discrete
counterpart by
\begin{equation}\label{tp:alphaN-cutoff}
 \alpha^{N;\eta}_1
 =N^{-a_{h}}\sum_{n,m=1}^N\gamma_{H_0}(n-m)
       \chi\left(\frac{|n-m|}{\eta N}\right)
       \gamma_{H_1}(S_n-\widetilde S_m).
\end{equation}

\begin{proposition}[Uniform discrete cutoff error]\label{tp:discrete-cutoff}
For $p\ge2$, $N\ge1$, and $0<\eta<\frac12$,
\begin{equation}\label{tp:discrete-cutoff-bound}
 \norm{\alpha^N_1-\alpha^{N;\eta}_1}_p
 \le C_{h}\,p^{1-H_1}(\eta+N^{-1})^{2H_0-1}.
\end{equation}
In particular,
\begin{equation}\label{tp:cutoff-vanishing}
 \lim_{\eta\downarrow0}\limsup_{N\to\infty}
 \norm{\alpha^N_1-\alpha^{N;\eta}_1}_p=0.
\end{equation}
\end{proposition}

\begin{proof}
All summands are nonnegative, and $1-\chi(u)=0$ for $u\ge2$. Therefore,
\[
 0\le\alpha^N_1-\alpha^{N;\eta}_1
 \le N^{-a_{h}}\sum_{\substack{|r|<N\\|r|\le2\eta N}}
                       \gamma_{H_0}(r)L_{N,r}.
\]
Apply Minkowski's inequality, \eqref{tp:lag-bound}, and
\eqref{fr:lattice-envelope} with
$R=\max\{1,\lfloor2\eta N\rfloor\}\le2\eta N+1$ to get
\[
 \norm{\alpha^N_1-\alpha^{N;\eta}_1}_p
 \le C_{h}\,p^{1-H_1}N^{H_1-a_{h}}(2\eta N+1)^{2H_0-1}.
\]
Since $H_1-a_{h}=1-2H_0$ and
$N^{1-2H_0}(2\eta N+1)^{2H_0-1}=(2\eta+N^{-1})^{2H_0-1}$, this is
\eqref{tp:discrete-cutoff-bound}, after changing the constant.
The bound also covers $2\eta N<1$: in that case only the synchronous
diagonal can remain. Positivity of $2H_0-1$ proves \eqref{tp:cutoff-vanishing}.
\end{proof}

\subsection{\texorpdfstring{Joint convergence with the local-time fields}{Joint convergence with the local-time fields}}

To identify the limit of $\alpha^N_1$ jointly with the
self-intersections, we need the local times and the self-intersection of
each walk jointly with the walk itself. The limiting coordinates are
measurable, but not continuous, functions of the Brownian path; the
following elementary lemma allows us to combine two such joint limits.

\begin{lemma}[Joining limits that are functions of a common coordinate]
\label{tp:joining-lemma}
Let $\mathcal X,\mathcal Y,\mathcal Z$ be Polish spaces, and let $f:\mathcal X\to \mathcal Y$, $g:\mathcal X\to \mathcal Z$ be Borel maps. Let $X$ be an $\mathcal X$-valued random variable.
Suppose $(X_N,\xi_N,\zeta_N)$ takes values in $\mathcal X\times\mathcal Y\times\mathcal Z$ and
\[
 (X_N,\xi_N)\ToLaw(X,f(X)),\qquad
 (X_N,\zeta_N)\ToLaw(X,g(X)).
\]
Then
\[
 (X_N,\xi_N,\zeta_N)\ToLaw(X,f(X),g(X)).
\]
\end{lemma}

\begin{proof}
The assumed convergences imply tightness of each coordinate and hence of
the full triple. Consider a weakly convergent subsequence, with limit
$(X_*,\xi_*,\zeta_*)$. By continuity of the coordinate projections,
\[
 (X_*,\xi_*)\stackrel{\mathrm{law}}=(X,f(X)),\qquad
 (X_*,\zeta_*)\stackrel{\mathrm{law}}=(X,g(X)).
\]
Choose compatible metrics
$d_{\mathcal Y}$ and $d_{\mathcal Z}$, and set
\[
 \varphi(x,y)=\min\{1,d_{\mathcal Y}(y,f(x))\},\qquad
 \psi(x,z)=\min\{1,d_{\mathcal Z}(z,g(x))\}.
\]
Both functions are nonnegative, bounded and Borel. The two equalities
of laws give
\[
 \Ep[\varphi(X_*,\xi_*)]=\Ep[\varphi(X,f(X))]=0,\qquad
 \Ep[\psi(X_*,\zeta_*)]=\Ep[\psi(X,g(X))]=0.
\]
It follows that $\xi_*=f(X_*)$ and $\zeta_*=g(X_*)$ almost surely.

Now let $F\in C_b(\mathcal X\times\mathcal Y\times\mathcal Z)$.
Since $X_*$ has the law of $X$ and
$x\mapsto F(x,f(x),g(x))$ is a bounded Borel function,
\begin{align*}
 \Ep[F(X_*,\xi_*,\zeta_*)]
 &=\Ep[F(X_*,f(X_*),g(X_*))]\\
 &=\Ep[F(X,f(X),g(X))].
\end{align*}
Thus every subsequential limit has the prescribed law.

Tightness and uniqueness of subsequential limits give
convergence of the full sequence.
\end{proof}

\begin{lemma}[Borel versions of local times and intersection measures]
\label{tp:borel-local-time}
Let $\CL=C([0,1]\times\R)$ be endowed with local uniform convergence,
and equip $\mathcal M_+([0,1]^2)$ with the weak topology.
There exist Borel maps
\[
 L:C([0,1])\to\CL,\qquad
 (b,\widetilde b)\longmapsto\mathfrak I_{b,\widetilde b}
 \quad\text{from }C([0,1])^2\text{ to }\mathcal M_+([0,1]^2),
\]
whose values at Brownian paths agree almost surely with the local-time
field of \eqref{temps-local} and the intersection measure of
\eqref{eq:two-BM-occupation-measure}, respectively.
In particular, $b\mapsto\mathfrak I_b:=\mathfrak I_{b,b}$ is Borel.

\end{lemma}

\begin{proof}
The space $\CL$ is Polish; let $d_{\CL}$ be the complete metric
used in the proof of Lemma~\ref{lem:local-time-invariance}.
For every continuous path $b$ and every $m\ge1$, define
\begin{equation}\label{el-em}
 L^{(m)}(b)(t,x):=\int_0^t\heat_{1/m}(b_s-x)\dd s,
 \qquad (t,x)\in[0,1]\times\R.
\end{equation}
Each map $L^{(m)}:C([0,1])\to\CL$ is continuous, since
$\heat_{1/m}$ has a bounded derivative. By completeness, the set
of paths on which $(L^{(m)}(b))_{m\ge1}$ converges in $\CL$ is
\[
 \mathcal A
 =\bigcap_{k\ge1}\ \bigcup_{M\ge1}\ \bigcap_{m,n\ge M}
 \left\{b:d_{\CL}\bigl(L^{(m)}(b),L^{(n)}(b)\bigr)<\frac1k\right\},
\]
and is therefore Borel. Set
\[
 L(b):=
 \begin{cases}
  \displaystyle\lim_{m\to\infty}L^{(m)}(b),&b\in\mathcal A,\\
  0,&b\notin\mathcal A,
 \end{cases}
\]
where $0$ denotes the zero field in $\CL$.
This map is Borel: it is the pointwise limit, in $\CL$, of the
Borel maps $b\mapsto\one_{\mathcal A}(b)L^{(m)}(b)$.
For the occupation-density formula and the jointly continuous version
used here, see \cite[Theorems~6.18 and~6.19, p.~165]{MP-Brownian}.
In particular, for Brownian motion,
\[
 L^{(m)}(B)(t,x)
 =\int_\R\heat_{1/m}(y-x)L_t^B(y)\dd y.
\]
The joint continuity of $L^B$ and its compact spatial support on
$[0,1]$ imply convergence to $L^B$ in $\CL$ almost surely.
Thus $\Pp(B\in\mathcal A)=1$ and $L(B)=L^B$ almost surely.
The values assigned outside $\mathcal A$ provide a measurable extension;
no local time is required to exist for those paths.

We now construct the intersection measures for arbitrary continuous
paths. Write $L_t(b)(x)=L(b)(t,x)$. For every $b$, this field is
nonnegative, continuous, nondecreasing in $t$, and zero at $t=0$.
Moreover, it vanishes outside the compact set $b([0,1])$. Indeed,
this is immediate when $b\notin\mathcal A$; when $b\in\mathcal A$
and $d=\operatorname{dist}(x,b([0,1]))>0$,
\[
 0\le L^{(m)}(b)(t,x)
 \le\sqrt{\frac m{2\pi}}e^{-md^2/2}\longrightarrow0.
\]
Let $\mu_x^b$ be the finite atomless Stieltjes measure on $[0,1]$
with $\mu_x^b([0,t])=L_t(b)(x)$, and define
\begin{equation}\label{tp:path-intersection-measure}
 \mathfrak I_{b,\widetilde b}(A)
 :=\int_\R(\mu_x^b\otimes\mu_x^{\widetilde b})(A)\dd x,
 \qquad A\in\mathcal B([0,1]^2).
\end{equation}
The interval masses of $\mu_x^b$ are continuous in $x$, so these are
measurable kernels and the last formula defines a measure. Its total
mass is finite, since
\[
 \mathfrak I_{b,\widetilde b}([0,1]^2)
 =\int_\R L_1(b)(x)L_1(\widetilde b)(x)\dd x<\infty
\]
by continuity and compact spatial support.

For intervals $I=[u,v)$ and $J=[u',v')$, including the endpoint $1$
when needed, atomlessness gives
\begin{align*}
 \mathfrak I_{b,\widetilde b}(I\times J)
 =\lim_{R\to\infty}\int_{-R}^R
 &\bigl(L_v(b)(x)-L_u(b)(x)\bigr)\bigl(L_{v'}(\widetilde b)(x)-L_{u'}(\widetilde b)(x)\bigr)\dd x.
\end{align*}
For fixed $R$, the integral is a continuous function of the two fields
in $\CL$, and hence a Borel function of $(b,\widetilde b)$.
Thus all these rectangle masses are Borel. Approximate each
$f\in C([0,1]^2)$ uniformly by functions constant on finite dyadic
partitions. Finiteness of the masses then shows that
$(b,\widetilde b)\mapsto\int f\dd\mathfrak I_{b,\widetilde b}$ is Borel.
A countable uniformly dense family of such test functions, including
$1$, generates the Borel sigma-field for the weak topology; hence
$(b,\widetilde b)\mapsto\mathfrak I_{b,\widetilde b}$ is Borel.
The diagonal map $b\mapsto(b,b)$ gives the assertion for $\mathfrak I_b$.
For Brownian paths, \eqref{tp:path-intersection-measure} agrees with
\eqref{eq:two-BM-occupation-measure} because $L(B)=L^B$ almost surely.
\end{proof}

\begin{proposition}[The augmented pair of paths]\label{tp:augmented-pair}
Define
\begin{align*}
 \mathsf T_N
 &=(S^N,\widetilde S^N,\mathcal L^N,\widetilde{\mathcal L}^{N},
                         Y^N_1,\widetilde Y^N_1),\\
 \mathsf T
 &=(B,\widetilde B,L^B,L^{\widetilde B},Y_1(B),Y_1(\widetilde B)).
\end{align*}
Then
\begin{equation}\label{tp:augmented-convergence}
 \mathsf T_N\ToLaw\mathsf T
 \quad\text{in }C([0,1])^2\times\CL^2\times\R^2.
\end{equation}
Here $\CL$ carries the topology specified in Lemma~\ref{tp:borel-local-time}.

\end{proposition}

\begin{proof}

Let $L$ and $b\mapsto\mathfrak I_b$ be the Borel maps supplied by
Lemma~\ref{tp:borel-local-time}. 

Each block $C_{M,k,1}$ in \eqref{fr:Brownian-block} is a $[0,\infty]$-valued Borel path functional:
for $H_1=\frac12$, use \eqref{eq:block-endpoint} and measurability of
integration against a fixed nonnegative Borel function; for
$H_1>\frac12$, apply measurability of parameter integrals to
\eqref{fr:Brownian-block}. By Lemma~\ref{lem:Brownian-blocks}, these
blocks are finite almost surely, so replacing infinite values by zero
and centering gives real-valued Borel versions of $A_{M,k,1}$.

By Proposition~\ref{prop:Y}, a fixed subsequence of the partial sums
$2\sum_{M=1}^{M_0}\sum_{k=1}^{2^{M-1}}A_{M,k,1}$ in
\eqref{fr:Y-definition} converges almost surely to $Y_1(B)$.
Its limit, extended by zero on the Borel set where the
subsequence does not converge, defines a Borel map
$G:C([0,1])\to\R$ with $G(B)=Y_1(B)$ almost surely.

We can now apply Lemma~\ref{tp:joining-lemma} with
$X_N=S^N$, $\xi_N=\mathcal L^N$, $\zeta_N=Y^N_1$,
$f=L$ and $g=G$, using
\eqref{eq:joint-local-time-invariance} and the first two coordinates
of \eqref{eq:path-YN-conv}. We obtain
\[
 (S^N,\mathcal L^N,Y^N_1)\ToLaw(B,L^B,Y_1(B)).
\]

For each $N$, the corresponding triple constructed from $\widetilde S$
is independent and has the same law. Write $\mu_N,\mu$ for the laws
of the triples in the last display, on the Polish space
$\mathbb X=C([0,1])\times\CL\times\R$. Since $\mu_N\ToLaw\mu$, the standard product-measure convergence
theorem \cite[Theorem~2.8(ii)]{Billingsley-convergence} gives
\[
 \mu_N\otimes\mu_N\ToLaw\mu\otimes\mu
 \qquad\text{on }\mathbb X^2.
\]
Reordering the coordinates gives \eqref{tp:augmented-convergence}.
\end{proof}

\subsection{\texorpdfstring{Convergence of the mutual occupation measures}{Convergence of the mutual occupation measures}}

We first establish the joint convergence of the mutual occupation
measures needed for the case $H_1=\frac12$ in the proof of
Theorem~\ref{tp:main}.

For the case $H_1=\frac12$, abbreviate the Brownian occupation
measure of \eqref{eq:two-BM-occupation-measure} by
\begin{equation}\label{tp:Brownian-measure}
 \Icross:=\mathfrak I_{B,\widetilde B}
 \quad\text{on }[0,1]^2.
\end{equation}
With $\tau^N_i=\frac{i-\frac12}{N}$, its discrete counterpart is
\begin{equation}\label{tp:discrete-measure}
 \Icross_N
 =N^{-\frac32}\sum_{i,j=1}^N
       \one_{\{S_i=\widetilde S_j\}}\delta_{(\tau^N_i,\tau^N_j)}.
\end{equation}
Let $\mathcal M_+([0,1]^2)$ denote the space of finite positive measures with
the weak topology.

\begin{proposition}[Joint measure convergence]\label{tp:measure-convergence}
Let $H_1=\frac12$. With the notation of Proposition~\ref{tp:augmented-pair},
\begin{equation}\label{tp:joint-measure-limit}
 (\mathsf T_N,\Icross_N)\ToLaw(\mathsf T,\Icross)
\end{equation}
in the product with $\mathcal M_+([0,1]^2)$.
\end{proposition}

\begin{proof}

By \eqref{tp:augmented-convergence}, the Skorokhod representation theorem
provides copies of $\mathsf T_N,\mathsf T$ on one probability space such that
all their coordinates converge almost surely in the stated topologies.
For each $N$, recover the walk positions from the path coordinates by
$S_i=\sqrt N\,S^N_{i/N}$ and
$\widetilde S_j=\sqrt N\,\widetilde S^N_{j/N}$, and define $\Icross_N$
by \eqref{tp:discrete-measure}. This is a Borel reconstruction from
$(S^N,\widetilde S^N)$, so the enlarged tuple
$(\mathsf T_N,\Icross_N)$ retains its original joint law.
At the limit, define $\Icross$ from $(L^B,L^{\widetilde B})$ as in
\eqref{eq:two-BM-occupation-measure}.

Define the midpoint step field
\begin{equation}\label{tp:midpoint-field}
 V^N_t(x)
 =N^{-\frac12}\sum_{\tau^N_i<t}
       \one_{\{S_i=\widetilde{N}_x\}},
 \qquad t\in[0,1],\quad x\in\R.
\end{equation}
If $k_N(t)=\#\{i\in\{1,\ldots,N\}:\tau^N_i<t\}$, then
$|\frac{k_N(t)}{N}-t|\le\frac1{2N}$ and, exactly,
\[
 V^N_t(x)=\mathcal L^N_{k_N(t)/N}(x_N).
\]
For each fixed
$0<\rr<\infty$, local uniform convergence of $\mathcal L^N$ and
uniform continuity of $L^B$ give, with the norms and moduli from the
proof of Lemma~\ref{lem:local-time-invariance},
\begin{align}
 \sup_{\substack{0\le t\le1\\|x|\le\rr}}
 |V^N_t(x)-L^B_t(x)|
 &\le\|\mathcal L^N-L^B\|_{K_{\rr+1}}
   +\omega_{\rr+1}\left(L^B,\frac1{2N}+\frac1{\sqrt N}\right)
\longrightarrow0\quad\text{a.s.}
 \label{tp:step-field-limit}
\end{align}
The same argument applies to the second walk and its field $\widetilde V^N_t(x)=\widetilde{\mathcal L}^{N}_{k_N(t)/N}(x_N)$.

For $I=[u,v)\subset[0,1]$, put
\[
 V^N_I(x)=V^N_v(x)-V^N_u(x)
 =N^{-\frac12}\sum_{\tau^N_i\in I}
          \one_{\{S_i=\widetilde{N}_x\}}.
\]
The half-open convention includes $u$ and excludes $v$, including when a
midpoint is exactly an endpoint. Spatial cells have length $N^{-\frac12}$, so
for intervals $I,J$ we have the exact identity
\begin{equation}\label{tp:rectangle-identity}
 \Icross_N(I\times J)
 =\int_\R V^N_I(x)\widetilde V^N_J(x)\dd x.
\end{equation}

On this representation, $S^N\to B$ and
$\widetilde S^N\to\widetilde B$ uniformly almost surely.
Choose a finite integer
$\rr>2+\max\{\|B\|_\infty,\|\widetilde B\|_\infty\}$.
Since $S_i=\sqrt N\,S^N_{i/N}$ and
$\widetilde S_j=\sqrt N\,\widetilde S^N_{j/N}$, the unions of spatial cells
\[
 \bigcup_{i=1}^N
   \left[\frac{S_i}{\sqrt N},\frac{S_i+1}{\sqrt N}\right),
 \qquad
 \bigcup_{j=1}^N
   \left[\frac{\widetilde S_j}{\sqrt N},
         \frac{\widetilde S_j+1}{\sqrt N}\right)
\]
are contained in $[-\rr,\rr]$ for all sufficiently large $N$; the fields $V^N$ and $\widetilde V^N$
vanish outside these respective sets. The fields $L^B$ and
$L^{\widetilde B}$ vanish outside $B([0,1])$ and
$\widetilde B([0,1])$, respectively, which are also contained in
$[-\rr,\rr]$.

Since \eqref{tp:step-field-limit} holds simultaneously for every integer $\rr$,
uniform convergence of the two fields implies
\begin{equation}\label{tp:rectangle-limit}
 \Icross_N(I\times J)
 \longrightarrow
 \int_\R L^B_I(x)L^{\widetilde B}_J(x)\dd x
 =\Icross(I\times J)
 \quad\text{almost surely}.
\end{equation}
This also holds for the total masses, by using $I=J=[0,1)$.
Neither measure charges $s=1$ or $t=1$; for the Brownian measure,
this follows from the atomlessness of the Stieltjes measures associated
with the local times, as recalled after \eqref{eq:two-BM-occupation-measure}.

Finally, let $f$ be continuous on $[0,1]^2$. Approximate $f$
uniformly by a function constant on the rectangles of a finite time
partition. Equation \eqref{tp:rectangle-limit} gives convergence of its
integrals. The error is bounded by the uniform approximation error times
the total masses, which are bounded along the convergent sequence.
Letting the mesh of the partition tend to zero proves
$\int f\dd\Icross_N\to\int f\dd\Icross$ almost surely.

Consequently, $(\mathsf T_N,\Icross_N)\to(\mathsf T,\Icross)$
almost surely on this representation. Since each enlarged tuple has
its original joint law by the Borel reconstruction above, this proves
\eqref{tp:joint-measure-limit} for the original laws.
\end{proof}

\subsection{\texorpdfstring{The joint two-path limit}{The joint two-path limit}}

We now prove the joint invariance principle at the unit horizon.
The proof establishes convergence at a fixed temporal cutoff separately
in the two spatial regimes, and then removes the temporal cutoff.

\begin{theorem}[Joint two-path invariance]\label{tp:main}
Under the standing assumptions,
\begin{equation}\label{tp:joint-limit}
\begin{split}
 &\bigl(S^N,\widetilde S^N,Y^N_1,\widetilde Y^N_1,
                         \alpha^N_1(S,\widetilde S)\bigr)\\
 &\qquad\ToLaw
 \bigl(B,\widetilde B,Y_1(B),Y_1(\widetilde B),
                         \alpha_1(B,\widetilde B)\bigr)
\end{split}
\end{equation}
in $C([0,1])^2\times\R^3$, with the uniform topology on each path space.
\end{theorem}

\begin{proof}
We prove the stronger convergence
\begin{equation}\label{tp:full-augmented-limit}
 (\mathsf T_N,\alpha^N_1)\ToLaw(\mathsf T,\alpha_1(B,\widetilde B))
\end{equation}
which retains the two local-time fields.

Equip the product state space with a compatible product metric such that
the distance between $(u,a)$ and $(u,b)$ is at most $|a-b|$.
We first show that, for every fixed $0<\eta<\frac12$,
\begin{equation}\label{tp:joint-cutoff-limit}
 (\mathsf T_N,\alpha^{N;\eta}_1)\ToLaw(\mathsf T,\alpha_1^{(\eta)}(B,\widetilde B)),
\end{equation}
where the truncated mutual intersections are defined by
\eqref{tp:alphaN-cutoff} and \eqref{tp:alpha-cutoff}.

\textit{Step 1: fixed temporal cutoff, $H_1=\frac12$.}
Define
\[
 w_\eta(s,t)
 =K_{H_0}(t-s)
     \chi\left(\frac{|t-s|}{\eta}\right),
\]
extended by zero on the diagonal.
This is a bounded continuous function on $[0,1]^2$.
Since $\tau^N_i-\tau^N_j=\frac{i-j}{N}$ and, at $H_1=\frac12$,
$N^{-a_{h}}=N^{-\frac32}N^{2-2H_0}$, definitions
\eqref{tp:alphaN-cutoff} and \eqref{tp:discrete-measure} give
\begin{equation}\label{tp:weight-error}
 \left|\alpha^{N;\eta}_1-\int w_\eta\dd\Icross_N\right|
 \le\delta_{N,\eta}\,\Icross_N([0,1]^2),
\end{equation}
where
\[
 \delta_{N,\eta}
 =\sup_{\eta N\le|r|\le N}
 \left|N^{2-2H_0}\gamma_{H_0}(r)-K_{H_0}\left(\frac rN\right)\right|.
\]
By \eqref{tp:gamma-away}, $\delta_{N,\eta}\to0$.
Proposition~\ref{tp:measure-convergence} implies tightness of the total
masses, so the right-hand side of \eqref{tp:weight-error} tends to zero
in probability. The continuous mapping theorem applied to
\eqref{tp:joint-measure-limit} gives
\[
 \left(\mathsf T_N,\int w_\eta\dd\Icross_N\right)
 \ToLaw
 \left(\mathsf T,\int w_\eta\dd\Icross\right).
\]
The last integral is $\alpha_1^{(\eta)}(B,\widetilde B)$. Slutsky's theorem proves
\eqref{tp:joint-cutoff-limit}.

\smallskip\noindent\textit{Step 2: fixed temporal cutoff, $H_1>\frac12$.}
Here $\Icross_N$ is replaced by a rescaling of the kernel $\gamma_{H_1}$, and a
spatial cutoff is inserted alongside the temporal one: for $\delta>0$, define
\begin{align}
 \alpha^{N;\eta,\delta}_1
 &=N^{-a_h}\sum_{n,m=1}^N
     \gamma_{H_0}(n-m)\chi\left(\frac{|n-m|}{\eta N}\right)\gamma_{H_1}(S_n-\widetilde S_m)
     \chi\left(\frac{|S_n-\widetilde S_m|}{\delta\sqrt N}\right),
 \label{tp:double-cutoff-discrete}\\
 \alpha_1^{(\eta,\delta)}(B,\widetilde B)
 &=\int_{[0,1]^2}K_{H_0}(t-s)
     \chi\left(\frac{|t-s|}{\eta}\right)K_{H_1}(B_s-\widetilde B_t)
     \chi\left(\frac{|B_s-\widetilde B_t|}{\delta}\right)\dd s\dd t.
 \label{tp:double-cutoff-Brownian}
\end{align}

Retain the function $w_\eta$ defined above and put
\[
 g_\delta(z)=K_{H_1}(z)\chi\left(\frac{|z|}{\delta}\right)
 \quad(z\ne0),\qquad g_\delta(0)=0.
\]
Thus both truncated kernels are extended by zero at their singularities.
They are bounded and Lipschitz, so the map
\[
 \Phi_{\eta,\delta}(w,\widetilde w)
 :=\int_{[0,1]^2}w_\eta(s,t)
        g_\delta(w_s-\widetilde w_t)\dd s\dd t
\]
is bounded and continuous on $C([0,1])^2$, and
$\Phi_{\eta,\delta}(B,\widetilde B)=\alpha_1^{(\eta,\delta)}(B,\widetilde B)$.

First compare the discrete kernels with their continuous limits.
Define the right-endpoint sum
\[
 R_N^{\eta,\delta}
 :=\frac1{N^2}\sum_{n,m=1}^N
      w_\eta\left(\frac nN,\frac mN\right)
      g_\delta\left(\frac{S_n-\widetilde S_m}{\sqrt N}\right).
\]
Since $N^{-a_h}=N^{-2}N^{2-2H_0}N^{1-H_1}$, evenness of the kernels
and the uniform bounds
$\rho_N(\eta)$ and $\sigma_N(\delta)$ of
\eqref{eq:gamma-away-diagonal} and \eqref{eq:G-away-diagonal} give
\begin{align}
 |\alpha^{N;\eta,\delta}_1-R_N^{\eta,\delta}|
 &\le\rho_N(\eta)\bigl(\|g_\delta\|_\infty+\sigma_N(\delta)\bigr)
      +\|w_\eta\|_\infty\sigma_N(\delta)\longrightarrow0.
 \label{tp:double-cutoff-kernel-error}
\end{align}
Indeed, each cutoff vanishes in the region excluded from its kernel
estimate, and the sum has $N^2$ terms. The bounds are uniform over all
walk positions.

For the quadrature error, let $\Delta_n^N=((n-1)/N,n/N]$.
For $(s,t)\in\Delta_n^N\times\Delta_m^N$, linear interpolation gives
\[
 \left|s-\frac nN\right|+\left|t-\frac mN\right|\le\frac2N,
 \qquad
 \left|S^N_s-\frac{S_n}{\sqrt N}\right|
 +\left|\widetilde S^N_t-\frac{\widetilde S_m}{\sqrt N}\right|
 \le\frac2{\sqrt N}.
\]
Integrating over these cells, with the Lipschitz constant of $w_\eta$
taken for the distance $|s-s'|+|t-t'|$, yields
\begin{align}
 |R_N^{\eta,\delta}-\Phi_{\eta,\delta}(S^N,\widetilde S^N)|
 &\le\frac2N\operatorname{Lip}(w_\eta)\|g_\delta\|_\infty
      +\frac2{\sqrt N}\|w_\eta\|_\infty
                          \operatorname{Lip}(g_\delta)\longrightarrow0.
 \label{tp:double-cutoff-quadrature-error}
\end{align}
Proposition~\ref{tp:augmented-pair}, the continuous mapping theorem,
and these two deterministic error bounds therefore imply
\begin{equation}\label{tp:joint-double-cutoff-limit}
 (\mathsf T_N,\alpha^{N;\eta,\delta}_1)
 \ToLaw(\mathsf T,\alpha_1^{(\eta,\delta)}(B,\widetilde B)).
\end{equation}

We next remove the spatial cutoff, keeping $\eta>0$ fixed. Write
$D_{n,m}=S_n-\widetilde S_m$. Independence and symmetry give
$\Pp(D_{n,m}=z)=p_{n+m}(z)$. If $|n-m|>\eta N$, then
$n+m\ge\eta N$, and \eqref{eq:LLT-sup} and
\eqref{fr:lattice-envelope} imply
\begin{align}
N^{1-H_1}E\left[\gamma_{H_1}(D_{n,m})
            \one_{\{|D_{n,m}|\le2\delta\sqrt N\}}\right]
 &=N^{1-H_1}\sum_{|z|\le2\delta\sqrt N}
                  \gamma_{H_1}(z)p_{n+m}(z)
      \notag\\
 &\le C_{H_1}\eta^{-1/2}N^{1/2-H_1}
                         (1+2\delta\sqrt N)^{2H_1-1}
      \notag\\
 &\le C_{H_1,\eta}(\delta+N^{-1/2})^{2H_1-1}.
 \label{tp:spatial-cutoff-error}
\end{align}
Moreover, \eqref{fr:lattice-envelope} gives the uniform temporal bound
\[
 0\le N^{2-2H_0}\gamma_{H_0}(n-m)
          \chi\left(\frac{|n-m|}{\eta N}\right)
 \le C_{H_0,\eta} \mathbf{1}_{\{|n-m|>\eta N\}}.
\]
Since $0\le1-\chi(u)\le\one_{\{u\le2\}}$ for $u\ge0$,
\begin{align}
 0\le E[\alpha^{N;\eta}_1-\alpha^{N;\eta,\delta}_1]
 &\le\frac{C_{H_0,\eta}}{N^2}
       \sum_{\substack{1\le n,m\le N\\|n-m|>\eta N}}
       N^{1-H_1}E\left[\gamma_{H_1}(D_{n,m})
            \one_{\{|D_{n,m}|\le2\delta\sqrt N\}}\right]
       \notag\\
 &\le C_{h,\eta}(\delta+N^{-1/2})^{2H_1-1}.
 \label{tp:discrete-spatial-cutoff-L1}
\end{align}
For the Brownian error, $B_s-\widetilde B_t$ has density
$\heat_{s+t}$, and $s+t\ge|t-s|>\eta$ wherever $w_\eta(s,t)\ne0$.
Hence
\begin{align}
 0\le E[\alpha_1^{(\eta)}(B,\widetilde B)-\alpha_1^{(\eta,\delta)}(B,\widetilde B)]
 &\le\|w_\eta\|_\infty
       \int_{\substack{(s,t)\in[0,1]^2\\|t-s|>\eta}}
       \int_{|z|\le2\delta}K_{H_1}(z)\heat_{s+t}(z)\dd z\dd s\dd t
       \notag\\
 &\le\frac{\|w_\eta\|_\infty}{\sqrt{2\pi\eta}}
       \int_{|z|\le2\delta}K_{H_1}(z)\dd z
 \le C_{h,\eta}\delta^{2H_1-1}.
 \label{tp:Brownian-spatial-cutoff-L1}
\end{align}

Finally, for every bounded $1$-Lipschitz test function $f$, the triangle inequality
and \eqref{tp:joint-double-cutoff-limit} give
\begin{align*}
 &\limsup_{N\to\infty}
   \left|\Ep[f(\mathsf T_N,\alpha^{N;\eta}_1)]
          -\Ep[f(\mathsf T,\alpha_1^{(\eta)}(B,\widetilde B))]\right|\\
 &\quad\le\limsup_{N\to\infty}
            E[\alpha^{N;\eta}_1-\alpha^{N;\eta,\delta}_1]
          +E[\alpha_1^{(\eta)}(B,\widetilde B)-\alpha_1^{(\eta,\delta)}(B,\widetilde B)]\\
 &\quad\le C_{h,\eta}\delta^{2H_1-1}
       \xrightarrow[\delta\downarrow0]{}0.
\end{align*}
Since $2H_1-1>0$, this proves \eqref{tp:joint-cutoff-limit}, with
$N\to\infty$ first and $\delta\downarrow0$ second.

\smallskip\noindent\textit{Step 3: removal of the temporal cutoff.}
Let $f$ be bounded by $1$ and
$1$-Lipschitz for this metric. The triangle inequality gives
\begin{align*}
 &|\Ep\left[f(\mathsf T_N,\alpha^N_1)\right]-\Ep\left[f(\mathsf T,\alpha_1(B,\widetilde B))\right]|\\
 &\quad\le
 \Ep\left[|\alpha^N_1-\alpha^{N;\eta}_1|\right]
 +\left|\Ep\left[f(\mathsf T_N,\alpha^{N;\eta}_1)\right]-\Ep\left[f(\mathsf T,\alpha_1^{(\eta)}(B,\widetilde B))\right]\right|\\
 &\qquad+\Ep\left[|\alpha_1^{(\eta)}(B,\widetilde B)-\alpha_1(B,\widetilde B)|\right].
\end{align*}
By \eqref{tp:joint-cutoff-limit}, proved in Steps 1-2, the middle term
tends to zero for fixed $\eta$. By Proposition~\ref{tp:discrete-cutoff}, used with $p=2$,
and Lemma~\ref{tp:Brownian-cutoff} at the unit horizon,
\[
 \limsup_{N\to\infty}
 |\Ep\left[f(\mathsf T_N,\alpha^N_1)\right]-\Ep\left[f(\mathsf T,\alpha_1(B,\widetilde B))\right]|
 \le C_{h}\eta^{2H_0-1}.
\]
Letting $\eta\downarrow0$ proves convergence of all bounded Lipschitz test
functions, hence \eqref{tp:full-augmented-limit}. Projection gives
\eqref{tp:joint-limit}. The order of limits is $N\to\infty$ first and
$\eta\downarrow0$ second.
\end{proof}

\subsection{The second moment at a general point}
\label{subsec:two-path-general}

We now combine the joint convergence with the exponential bounds
to pass to the second moment. The following corollary allows the
coefficients to vary with $N$, as required by the exact time scaling
of Lemma~\ref{lem:discrete-scaling}. 

\begin{corollary}[Mixed limits with converging coefficients]
\label{tp:variable-coefficients}
Let $\lambda_{1,N}\to\lambda_1$ and $\lambda_{2,N}\to\lambda_2$ in $\R$,
and let $\vartheta_N\to\vartheta$ in $[0,\infty)$. Then, for any sequence of integers $(M_N)$ such that $M_N \to \infty$, one has
\begin{equation}\label{tp:variable-coefficient-limit}
\Ep\bigl[e^{\lambda_{1,N}Y^{M_N}_1+\lambda_{2,N}\widetilde Y^{M_N}_1+\vartheta_N\alpha^{M_N}_1}\bigr]
 \longrightarrow
 \Ep\bigl[e^{\lambda_1Y_1(B)+\lambda_2 Y_1(\widetilde B)+\vartheta\alpha_1(B,\widetilde B)}\bigr].
\end{equation}

\end{corollary}

\begin{proof}
Let $X_N$ denote the variable inside the first expectation, and let
$X$ denote its proposed limit. The joint convergence \eqref{tp:joint-limit}, which holds along the indices $M_N\to\infty$, Slutsky's lemma and the continuous mapping
theorem give $X_N\ToLaw X$.
Choose $\Lambda<\infty$ such that
$|\lambda_{1,N}|,|\lambda_{2,N}|,\vartheta_N\le\Lambda$ for all $N$.
For any $r>1$, H\"older's inequality, \eqref{eq:YN-exp-unif} and
\eqref{tp:alpha-exp} give
\begin{align}
 E[X_N^r]
 &\le
 \left(E[e^{3r\lambda_{1,N}Y^{M_N}_1}]\right)^{1/3}
 \left(E[e^{3r\lambda_{2,N}\widetilde Y^{M_N}_1}]\right)^{1/3}
 \left(E[e^{3r\vartheta_N\alpha^{M_N}_1}]\right)^{1/3}\notag\\
 &\le\exp\!\left(C_h[r\Lambda+r^2\Lambda^2]\right)<\infty,
 \label{tp:mixed-uniform-Lr}
\end{align}
uniformly in $N$. Thus $(X_N)$ is uniformly integrable, and weak
convergence yields $E[X_N]\to E[X]$. 
\end{proof}

\begin{corollary}[Convergence of the $L^2$ norms]
\label{prop:second-moment-convergence}
For every $t>0$ and $x\in\R$, the coupled partition functions
satisfy
\begin{equation}\label{eq:second-moment-convergence}
 \E[\widehat{\mathsf Z}^N_t(x)^2]
 \longrightarrow
 e^{2R_t}E_{B,\widetilde B}\left[
    e^{\frac12(Y_t(B)+Y_t(\widetilde B))+\alpha_t(B,\widetilde B)}\right]
 =\E[U_t(x)^2].
\end{equation}
Consequently,
\[
 \|\widehat{\mathsf Z}^N_t(x)\|_{L^2(\Omega)}
       \longrightarrow\|U_t(x)\|_{L^2(\Omega)},
 \qquad \sup_N\|\widehat{\mathsf Z}^N_t(x)\|_{L^2(\Omega)}<\infty.
\]
\end{corollary}

\begin{proof}
By Lemma~\ref{lem:discrete-scaling}, the path exponential in
\eqref{eq:second-moment-ZN} is
\[
 \exp\Bigl(\frac{\theta^N_t}2\bigl(Y^{N_t}_1+\widetilde Y^{N_t}_1\bigr)
          +\theta^N_t\alpha^{N_t}_1(S,\widetilde S)\Bigr).
\]
Apply Corollary~\ref{tp:variable-coefficients} with $M_N=N_t$ and coefficients
$(\lambda_{1,N},\lambda_{2,N},\vartheta_N)
=(\frac{\theta^N_t}2,\frac{\theta^N_t}2,\theta^N_t)$ converging to
$(\frac{\theta_t}2,\frac{\theta_t}2,\theta_t)$ by \eqref{eq:thetaN}:
\[
 E_{S,\widetilde S}\left[
   e^{\frac{\theta^N_t}2(Y^{N_t}_1+\widetilde Y^{N_t}_1)
                     +\theta^N_t\alpha^{N_t}_1(S,\widetilde S)}\right]
 \longrightarrow
 E_{B,\widetilde B}\left[
   e^{\frac{\theta_t}2(Y_1(B)+Y_1(\widetilde B))+\theta_t\alpha_1(B,\widetilde B)}\right].
\]
By \eqref{eq:Yt-alphat-Rt}, the right-hand side is the Brownian
expectation in \eqref{eq:second-moment-convergence}. Multiplication by
$e^{2R^N_t}\to e^{2R_t}$, given by Corollary~\ref{cor:RN-t}, and the
identification \eqref{eq:continuous-second-tx} prove
\eqref{eq:second-moment-convergence}. Taking square roots gives
the norm convergence. Each norm is finite by
Proposition~\ref{prop:discrete-identities}, so the convergent
sequence of norms is bounded.
\end{proof}

\section{\texorpdfstring{Conclusion: proof of Theorem~\ref{thm:coupled-L2} and Theorem~\ref{thm:fractional-main}}{Conclusion of the coupled convergence proof}}

\label{sec:short}

\textit{For the whole section, we fix $h=(H_0,H_1)$ satisfying the
conditions in \eqref{fr:range}}.

\smallskip

In the pointwise arguments we keep $t>0$ and $x\in\R$ fixed.
Corollaries~\ref{prop:S-convergence} and
\ref{prop:second-moment-convergence} already give convergence of
the $\mathcal S$-transforms and of the $L^2$ norms, at that point. We combine
these two convergences using totality of Gaussian exponentials, then
use uniform higher moments to obtain convergence in every finite $L^p$.

\subsection{Preliminaries}

\begin{lemma}[Totality of Gaussian exponentials]\label{lem:total}
The linear span of $\{\cE(\varphi):\varphi\in\D\}$ is dense in
$L^2(\Omega,\sigma(W),\Pp)$.
\end{lemma}

\begin{proof}
We use the classical totality \cite[Lemma~1.1.2, p.~5]{Nualart-Malliavin} of the Gaussian exponentials $\{\cE(\varphi):\varphi\in\Hh\}$. Multiplication by the
nonzero deterministic factors in \eqref{eq:stoch-exp} preserves their
linear span. Since $\D$ is dense in $\Hh$ and
$f\mapsto\cE(f)$ is continuous from $\Hh$ into $L^2(\Omega)$ by
\eqref{eq:Gaussian-pairing}, the same totality holds with arguments
restricted to $\D$.
\end{proof}

\begin{proposition}[Uniform moments of the partition functions]
\label{prop:discrete-uniform-moments}
For every $1\le p<\infty$ and every finite $T>0$,
\begin{equation}\label{eq:discrete-uniform-moments}
 \sup_N\E[\widehat{\mathsf Z}^N_t(x)^p]<\infty,
 \quad\text{uniformly in }t\in[0,T]\text{ and }x\in\R.
\end{equation}
More precisely, for every integer $m\ge1$ and independent copies
$S^1,\ldots,S^m$ of the simple walk,
\begin{equation}\label{eq:discrete-m-moment}
 \E[\widehat{\mathsf Z}^N_t(x)^m]
 =e^{mR^N_t}E\left[\exp\left(
       \frac12\sum_{i=1}^mY^N_t(S^i)
       +\sum_{1\le i<j\le m}\alpha^N_t(S^i,S^j)\right)\right].
\end{equation}
\end{proposition}

\begin{proof}
Multiply $m$ copies of \eqref{eq:Gaussian-disc} and take Gaussian expectation.
Conditionally on $S^1,\ldots,S^m$, put $h_i=h^{S^i,N}_{t,x}$.
These are deterministic elements of $\Hh$, and
$\sum_iW(h_i)=W(\sum_i h_i)$ is centered Gaussian with variance
$\|\sum_i h_i\|_{\Hh}^2$. Thus \eqref{eq:stoch-exp} gives
\begin{align*}
 \E\left[\prod_{i=1}^m\cE(h_i)\right]
 &=\exp\left(-\frac12\sum_{i=1}^m\|h_i\|_{\Hh}^2\right)
       \E\left[\exp\left(W\left(\sum_{i=1}^m h_i\right)\right)\right]\\
 &=\exp\left(\frac12\left\|\sum_{i=1}^m h_i\right\|_{\Hh}^2
                    -\frac12\sum_{i=1}^m\|h_i\|_{\Hh}^2\right)\\
 &=\exp\left(\sum_{1\le i<j\le m}\ip{h_i}{h_j}_{\Hh}\right).
\end{align*}
Each inner product is $\alpha^N_t(S^i,S^j)$ by
\eqref{eq:alphaN}. Substituting this identity into the product of the
path representations proves \eqref{eq:discrete-m-moment}.

If $N_t=0$, including $0<t<\frac1N$, the partition function equals
$1$. Otherwise, by Lemma~\ref{lem:discrete-scaling}, the path exponential equals
$\exp(\frac{\theta^N_t}2\sum_iY^{N_t}_1(S^i)
+\theta^N_t\sum_{i<j}\alpha^{N_t}_1(S^i,S^j))$, with $\theta^N_t\le T^{a_{h}}$.
Put $\mathsf m=m+\binom m2$. H\"older's inequality bounds its path expectation by
\[
 \prod_{i=1}^m
       \left(E[e^{\frac {\mathsf m\theta^N_t}2Y^{N_t}_1(S^i)}]\right)^{\frac1{\mathsf m}}
 \prod_{i<j}
       \left(E[e^{\mathsf m\theta^N_t\alpha^{N_t}_1(S^i,S^j)}]\right)
                                       ^{\frac1{\mathsf m}}.
\]
Each factor is uniformly bounded in $N$ and in
$\theta^N_t\in(0,T^{a_{h}}]$
by \eqref{eq:YN-exp-unif}
or \eqref{tp:alpha-exp}, with parameters bounded in terms of $\mathsf m,T,h$.
The deterministic factor is bounded uniformly in $t\in[0,T]$ by \eqref{eq:RN-scaling} and \eqref{eq:RN-limit}, since $|R^N_t|=\theta^N_t|R^{N_t}_1|\le T^{a_h}\sup_{N'}|R^{N'}_1|<\infty$.
This proves the uniform estimate for integer $m$. For general $p\ge1$,
choose an integer $m\ge p$ and use monotonicity of $L^p$ norms.
None of the bounds depends on $x$, and the case $t=0$ is
trivial.
\end{proof}

\subsection{Discrete increment estimates and interpolation}
\label{subsec:discrete-increments}

We retain $a_{h}$ from
\eqref{fr:exponents} and the increment exponent $\eta$ of
\eqref{eq:increment-exponent-used}.

The spatial increment of the mutual intersection rests on one of the next
two estimates, according to the regime $H_1>\frac12$ or $H_1=\frac12$. The first is the lattice
counterpart of \eqref{eq:BM-occupation-increment} and carries the
exponent $2H_1-1$; the second is its substitute for $H_1=\frac12$,
where that exponent vanishes, and carries the exponent $\frac12$ at
the price of a parity restriction.

\begin{lemma}[Kernel increments along a difference walk]
\label{lem:difference-kernel-increments}
Let $D$ have the increment law of
Lemma~\ref{tp:difference-local-time}, with any starting point. For
$L\ge1$, $q\ge2$ and $z,z'\in\mathbb Z$,
\begin{equation}\label{eq:difference-kernel-increments}
 \left\|\sum_{\ell=1}^{L}
    \bigl[\gamma_{H_1}(D_\ell-z)-\gamma_{H_1}(D_\ell-z')\bigr]\right\|_q
 \le C_{H_1}\sqrt q\,L^{\frac12}|z-z'|^{2H_1-1}.
\end{equation}
The bound is uniform in the starting point, and there is no
restriction on the parity of the displacement.
\end{lemma}

\begin{proof}
Write the sum as $\sum_y[\gamma_{H_1}(y-z)-\gamma_{H_1}(y-z')]V_L(y)$ with
$V_L(y)=\sum_{\ell\le L}\one_{\{D_\ell=y\}}$. Minkowski's inequality,
the bound $\sup_y\|V_L(y)\|_q\le C\sqrt{qL}$ of
Lemma~\ref{tp:difference-local-time}, and the lattice translation
bound \eqref{fr:lattice-translation} give the result.
\end{proof}

\begin{lemma}[Spatial increments of difference-walk local times]
\label{lem:difference-local-time-increments}
Let $D$ have the increment law of
Lemma~\ref{tp:difference-local-time}, with any starting point, and
put $V_L(z)=\sum_{\ell=1}^L\one_{\{D_\ell=z\}}$.
If $z-z'\in2\mathbb Z$, then, for $L\ge1$ and $q\ge2$,
\begin{equation}\label{eq:difference-local-time-increments}
 \|V_L(z)-V_L(z')\|_q
 \le C_q|z-z'|^{1/2}L^{1/4}.
\end{equation}
The bound is uniform in the starting point. 
\end{lemma}

\begin{proof}
If $z=z'$ there is nothing to prove. Otherwise put $\delta=|z-z'|\in2\N$, so that $\delta\ge2$, and $x_0=D_0$.
As in the proof of Lemma~\ref{tp:difference-local-time}, we may realize
$(D_j)_{j\geq 0}=(U_{2j})_{j\geq 0}$, where $U$ is a simple symmetric walk started at $x_0$.

Let $k_r$ denote the mass function of its displacement over $2r$ steps:
\[
 k_r(u):=\Pp_{x_0}(U_{2r}-U_0=u)=p_{2r}(u),\qquad u\in\mathbb Z.
\]
The explicit binomial formula \eqref{eq:LLT-j} gives, for every
integer $r\ge0$,
\[
 k_r(u)=
 \begin{cases}
  2^{-2r}\displaystyle\binom{2r}{r+u/2},
       &u\in2\mathbb Z\text{ and }|u|\le2r,\\[1ex]
  0,   &\text{otherwise}.
 \end{cases}
\]
Set
\[
 a_r=k_r(0)+k_r(\delta),\qquad b_r=k_r(0)-k_r(\delta).
\]
The binomial probabilities are symmetric and maximal at zero, so
$a_r,b_r\ge0$. Moreover, \eqref{eq:LLT-sup}, together with $k_0(0)=1$,
gives
\[
 \sum_{r=0}^L a_r\le C\sqrt L,\qquad
 \sum_{r=1}^L c_r\le C\sqrt L,
 \quad c_r:=k_r(z-x_0)+k_r(z'-x_0).
\]
The second bound is uniform in $x_0$.

To evaluate $\sum_{r\ge0}b_r$, we use an equivalent Fourier
representation of this explicit binomial law: it will turn the
summation over $r$ into a geometric series. Write
$U_{2r}-U_0=\varepsilon_1+\cdots+\varepsilon_{2r}$, where the
$\varepsilon_j$ are independent symmetric signs. Since
$E_{x_0}[e^{i\xi\varepsilon_j}]=\frac12(e^{i\xi}+e^{-i\xi})=\cos\xi$,
independence gives $E_{x_0}[e^{i\xi(U_{2r}-U_0)}]=(\cos\xi)^{2r}$. Fourier inversion on $\mathbb Z$ then yields, for every $u\in\mathbb Z$,
\[
 k_r(u)=\frac1{2\pi}\int_{-\pi}^{\pi}
                  e^{-iu\xi}(\cos\xi)^{2r}\dd\xi
       =\frac1{2\pi}\int_{-\pi}^{\pi}
                  \cos(u\xi)(\cos\xi)^{2r}\dd\xi,
\]
because the imaginary part is odd. Subtracting the formulas at
$u=0$ and $u=\delta$ gives
\[
 b_r=\frac1{2\pi}\int_{-\pi}^{\pi}
              (1-\cos(\delta\xi))(\cos\xi)^{2r}\dd\xi.
\]
Each integrand is nonnegative. Hence, summing first over
$0\le r\le M$ and then applying monotone convergence as
$M\to\infty$, we obtain
\begin{align*}
& \sum_{r=0}^\infty b_r=\frac1{2\pi}\int_{-\pi}^{\pi}
       \sum_{r=0}^\infty(1-\cos(\delta\xi))(\cos\xi)^{2r}\dd\xi=\frac1{2\pi}\int_{-\pi}^{\pi}
       \frac{1-\cos(\delta\xi)}{1-\cos^2\xi}\dd\xi.
\end{align*}
To evaluate the integral, write $\delta=2v$ with $v\ge1$ an integer.
For $\sin\xi\ne0$,
\[
 \frac{1-\cos(2v\xi)}{1-\cos^2\xi}
 =2\left|\frac{\sin(v\xi)}{\sin\xi}\right|^2,
 \qquad
 \frac{\sin(v\xi)}{\sin\xi}
 =e^{i(v-1)\xi}\frac{1-e^{-2iv\xi}}{1-e^{-2i\xi}}
 =\sum_{j=0}^{v-1}e^{i(v-1-2j)\xi}.
\]
Expanding the squared modulus and integrating term by term in this
finite sum therefore proves
\[
 \sum_{r=0}^\infty b_r
 =2\sum_{j,k=0}^{v-1}
       \frac1{2\pi}\int_{-\pi}^{\pi}e^{2i(k-j)\xi}\dd\xi
 =2\sum_{j=0}^{v-1}1=2v=\delta.
\]

\smallskip

Now write
\[
V_L(z)-V_L(z')=\sum_{n=1}^L g(U_{2n}), \qquad \text{with} \ g(u):=\one_{\{u=z\}}-\one_{\{u=z'\}}.
\]
Put also $f(u)=\one_{\{u=z\}}+\one_{\{u=z'\}}$ and let $\mathcal F_{2n}=\sigma(U_0,\ldots,U_{2n})$. For integers
$n,r\ge0$, the increment $U_{2(n+r)}-U_{2n}$ is independent of
$\mathcal F_{2n}$ and has mass function $k_r$. Thus
\begin{align*}
 E_{x_0}[g(U_{2(n+r)})\mid\mathcal F_{2n}]
 &=k_r(z-U_{2n})-k_r(z'-U_{2n}),\\
 E_{x_0}[f(U_{2(n+r)})\mid\mathcal F_{2n}]
 &=k_r(z-U_{2n})+k_r(z'-U_{2n}).
\end{align*}
Symmetry gives
$k_r(z-z')=k_r(z'-z)=k_r(\delta)$, and hence, for $u\in\{z,z'\}$, 
\[
 k_r(z-u)-k_r(z'-u)=b_r g(u),\qquad
 k_r(z-u)+k_r(z'-u)=a_r.
\]
Since $g^2=f$, multiplying the conditional expectations by
$g(U_{2n})$ therefore yields
\begin{align*}
 g(U_{2n})E_{x_0}[g(U_{2(n+r)})\mid\mathcal F_{2n}]
 &=b_r f(U_{2n}),\\
 g(U_{2n})E_{x_0}[f(U_{2(n+r)})\mid\mathcal F_{2n}]
 &=a_r g(U_{2n}).
\end{align*}
Both identities also hold when $U_{2n}\notin\{z,z'\}$, since both
sides then vanish.

Fix an integer $p\ge1$ and ordered times
$1\le n_1\le\cdots\le n_{2p}\le L$, and set
$r_j=n_j-n_{j-1}$ for $j\ge2$. When $p=1$, the first conditional
identity gives
\[
 E_{x_0}[g(U_{2n_1})g(U_{2n_2})]
 =b_{r_2}E_{x_0}[f(U_{2n_1})]=b_{r_2}c_{n_1}.
\]
For $p\ge2$, condition successively on $\mathcal F_{2n_{2p-1}}$
and $\mathcal F_{2n_{2p-2}}$. Using the two identities in that order
gives
\begin{align*}
& E_{x_0}\left[\prod_{j=1}^{2p}g(U_{2n_j})\right]=b_{r_{2p}}E_{x_0}\left[
       \left(\prod_{j=1}^{2p-2}g(U_{2n_j})\right)
             f(U_{2n_{2p-1}})\right]=b_{r_{2p}}a_{r_{2p-1}}
       E_{x_0}\left[\prod_{j=1}^{2p-2}g(U_{2n_j})\right].
\end{align*}
Iterating this reduction down to $p=1$ proves
\[
 E_{x_0}\left[\prod_{j=1}^{2p}g(U_{2n_j})\right]=c_{n_1}\prod_{j=1}^{p}b_{r_{2j}}\prod_{j=1}^{p-1}a_{r_{2j+1}}\ge0.
\]

This also holds when some times coincide, since $a_0=b_0=1$.
Expanding the $2p$th moment, ordering its time indices, and then
enlarging the sums over their gaps gives
\begin{align*}
 E_{x_0}\left[|V_L(z)-V_L(z')|^{2p}\right]
 &\le (2p)!\left(\sum_{r=1}^L c_r\right)
       \left(\sum_{r=0}^L a_r\right)^{p-1}
       \left(\sum_{r=0}^L b_r\right)^p\\
 &\le C_p L^{p/2}\delta^p.
\end{align*}
Taking the $2p$th root proves the estimate at even integer moments.
For any real $q\ge2$, choose $p=\lceil q/2\rceil$ and use
monotonicity of $L^q$ norms.
\end{proof}

For the following path estimates, consider the grid times $m/N$, with $m\ge0$, and the functionals $Y^N_{m/N}$ and $R^N_{m/N}$ of
\eqref{eq:YN}-\eqref{eq:RN}; both are zero when $m=0$.
For two independent walks define
\begin{equation}\label{eq:discrete-two-point-intersection}
 A^N_{m,n}(d)=\beta_N^2
   \sum_{i=1}^m\sum_{\ell=1}^n
       \gamma_{H_0}(m-n-i+\ell)
       \gamma_{H_1}(d+S_i-\widetilde S_\ell),\qquad d\in\mathbb Z.
\end{equation}

The covariance identity \eqref{eq:cell-covariance} extends to cells
anchored at different terminal times. Indeed, for any valid cell
indices and $t,s\ge0$, the changes of variables
$r=t-i/N+u/N$ and $r'=s-\ell/N+v/N$ give
\begin{align*}
 &N^{2H_0}\int_{\mathfrak I^N_{t,i}}
                \int_{\mathfrak I^N_{s,\ell}}
                   K_{H_0}(r-r')\dd r'\dd r\\
 &\qquad=\int_0^1\int_0^1
       K_{H_0}\bigl(N(t-s)-i+\ell+u-v\bigr)\dd v\dd u\\
 &\qquad=\gamma_{H_0}\bigl(N(t-s)-i+\ell\bigr),
\end{align*}
where the last equality is the tent representation
\eqref{fr:tent-representation}. Combining this with the spatial
cell integral computed in Proposition~\ref{prop:coupling} gives
\begin{equation}\label{eq:cross-time-cell-covariance}
 \ip{e^N_{t,i,a}}{e^N_{s,\ell,b}}_{\Hh}
 =\gamma_{H_0}\bigl(N(t-s)-i+\ell\bigr)\gamma_{H_1}(a-b).
\end{equation}
Taking $t=m/N$, $s=n/N$, $a=j+S_i$ and
$b=k+\widetilde S_\ell$, and summing according to
\eqref{eq:hNS-explicit}, therefore yields
\[
 A^N_{m,n}(j-k)
 =\ip{h^{S,N}_{m/N,j/\sqrt N}}
           {h^{\widetilde S,N}_{n/N,k/\sqrt N}}_{\Hh}.
\]

\

\begin{lemma}[Discrete path increments]
\label{lem:discrete-path-increments}
Let $\eta$ be the exponent defined in \eqref{eq:increment-exponent-used}.
For $|\theta|\le1$, set
\[
 A^N_{m;\theta}(d)=\beta_N^2
       \sum_{i,\ell=1}^m\gamma_{H_0}(\theta-i+\ell)
                         \gamma_{H_1}(d+S_i-\widetilde S_\ell).
\]
Fix $T>0$. For every $q\ge2$, the following bounds hold uniformly in
\begin{equation}\label{eq:discrete-path-parameter-range}
 N\in\N,\ N\ge1,\qquad
 m,m',n,n'\in\{0,\ldots,\lfloor NT\rfloor\},\qquad
 d,d'\in\mathbb Z,\qquad |\theta|\le1,
\end{equation}
with the additional restriction $d-d'\in2\mathbb Z$ in
\eqref{eq:discrete-mutual-space} when $H_1=\frac12$:
\begin{align}
 \|Y^N_{m/N}-Y^N_{n/N}\|_q+|R^N_{m/N}-R^N_{n/N}|
    &\le C_{q,h}\left|\frac{m-n}{N}\right|^{a_{h}},
                    \label{eq:discrete-self-time-increment}\\
 \sup_{d\in\mathbb Z}
       \bigl(\|A^N_{m,n}(d)\|_q+\|A^N_{m;\theta}(d)\|_q\bigr)
    &\le C_{h,T}\sqrt q, \label{eq:discrete-two-point-moment}\\
 \sup_{d\in\mathbb Z}\|A^N_{m,n}(d)-A^N_{m',n'}(d)\|_q
    &\le C_{q,h,T}\left(
          \left|\frac{m-m'}N\right|^{2H_0-1}
         +\left|\frac{n-n'}N\right|^{2H_0-1}\right),
                    \label{eq:discrete-mutual-time}\\
 \|A^N_{m,n}(d)-A^N_{m,n}(d')\|_q
    &\le C_{q,h,T}
            \left|\frac{d-d'}{\sqrt N}\right|^{2\eta},
                    \label{eq:discrete-mutual-space}\\
 \sup_{d\in\mathbb Z}\|A^N_{m;\theta}(d)-A^N_{m;0}(d)\|_q
    &\le C_{q,h,T}N^{-(2H_0-1)}|\theta|^{2H_0-1}.
                    \label{eq:discrete-anchor-increment}
\end{align}
Moreover, for every $\lambda>0$, uniformly under
\eqref{eq:discrete-path-parameter-range},
\begin{equation}\label{eq:discrete-path-exponential-moments}
\begin{aligned}
 E\bigl[e^{\lambda|Y^N_{m/N}|}\bigr]
    &\le C_{\lambda,h,T},\\
 \sup_{d\in\mathbb Z}E\bigl[e^{\lambda A^N_{m,n}(d)}\bigr]
 +\sup_{d\in\mathbb Z}E\bigl[e^{\lambda A^N_{m;\theta}(d)}\bigr]
    &\le C_{\lambda,h,T}.
\end{aligned}
\end{equation}
\end{lemma}

\begin{proof}

For the centered self-intersection, let $m=n+k$, $k\ge1$, and put
\[
 S^{n,-}_r=S_n-S_{n-r}\quad(0\le r\le n),\qquad
 S^{n,+}_v=S_n-S_{n+v}\quad(0\le v\le k).
\]

These are independent simple walks starting at zero. Splitting
\eqref{eq:GammaN} at $n$ and using \eqref{eq:discrete-scaling} gives
\begin{align*}
 Y^N_{m/N}-Y^N_{n/N}
   &=\left(\frac kN\right)^{a_h}Y^k_1(S^{n,+})
       +2N^{-a_h}\bigl(C_{n,k}-E[C_{n,k}]\bigr),\\
 C_{n,k}
   &=\sum_{r=0}^{n-1}\sum_{v=1}^k
       \gamma_{H_0}(r+v)\gamma_{H_1}(S^{n,-}_r-S^{n,+}_v).
\end{align*}
Decompose $C_{n,k}=C^{n,k}_{-1}+\sum_{j\ge0}C^{n,k}_j$, where, with $R_j=2^jk$,
\begin{align*}
 C^{n,k}_{-1}
   &=\sum_{r=0}^{\min(n,k)-1}\sum_{v=1}^k
       \gamma_{H_0}(r+v)\gamma_{H_1}(S^{n,-}_r-S^{n,+}_v),\\
 C^{n,k}_j
   &=\sum_{r=R_j}^{\min(n,2R_j)-1}\sum_{v=1}^k
       \gamma_{H_0}(r+v)\gamma_{H_1}(S^{n,-}_r-S^{n,+}_v).
\end{align*}
Empty sums are zero. Conditionally on $S^{n,+}_1$, and after extending $S^{n,-}$ beyond time $n$ by an independent walk, which by nonnegativity can only increase $C^{n,k}_{-1}$ to a full $k\times k$ block, Lemma~\ref{lem:weighted-block}, with $v=1+s$, yields
\[
 \|C^{n,k}_{-1}\|_q\le C_hq^{1-H_1}k^{a_h}.
\]
For each nonempty strip, \eqref{fr:lattice-envelope} and the
conditional occupation bound \eqref{eq:RW-occupation} give
\begin{align*}
 \|C^{n,k}_j\|_q
 &\le C_hR_j^{2H_0-2}\sum_{v=1}^k
       \left\|\sum_{r=R_j}^{\min(n,2R_j)-1}
          \gamma_{H_1}(S^{n,-}_r-S^{n,+}_v)\right\|_q\\
 &\le C_hq^{1-H_1}kR_j^{a_h-1}.
\end{align*}
Since $a_h<1$ and $q^{1-H_1}\le\sqrt q$,
\[
 \|C_{n,k}\|_q
 \le C_hq^{1-H_1}k^{a_h}
          \left(1+\sum_{j\ge0}2^{j(a_h-1)}\right)
 \le C_h\sqrt q\,k^{a_h}.
\]
Thus \eqref{eq:YN-subgaussian-Lp} implies
\[
 \|Y^N_{m/N}-Y^N_{n/N}\|_q
 \le \left(\frac kN\right)^{a_h}\|Y^k_1(S^{n,+})\|_q
       +4N^{-a_h}\|C_{n,k}\|_q
 \le C_h\sqrt q\,\left(\frac kN\right)^{a_h}.
\]
This includes $n=0$; the case $m=n$ is immediate, and $m<n$ follows
by interchanging the indices.

For the remainder put $w_r=\gamma_{H_0}(r)q_{H_1}(r)$ with
$q_{H_1}(r)=E_S[\gamma_{H_1}(S_r)]$, as in Proposition~\ref{prop:counterterm}.
Using \eqref{eq:Lambda-remainder} with $m$ in place of $N$ and the definition \eqref{eq:RN}, we obtain
\[
 R^N_{m/N}=-N^{-a_{h}}\sum_{r\ge1}\min(m,r)w_r.
\]
For $m\ge n$, the identity
\[
 \min(m,r)-\min(n,r)=\sum_{i=n}^{m-1}\one_{\{r>i\}}
\]
and the envelope $0\le w_r\le C_h(1+r)^{a_h-2}$ of
\eqref{eq:w-envelope} give
\begin{align*}
 |R^N_{m/N}-R^N_{n/N}|
 &=N^{-a_h}\sum_{r\ge1}
       \bigl(\min(m,r)-\min(n,r)\bigr)w_r\\
 &=N^{-a_h}\sum_{i=n}^{m-1}\sum_{r=i+1}^{\infty}w_r\le C_hN^{-a_h}\sum_{i=n}^{m-1}(1+i)^{a_h-1}\\
 &\hspace{2cm}\le C_hN^{-a_h}\sum_{j=0}^{m-n-1}(1+j)^{a_h-1}
 \le C_h\left(\frac{m-n}{N}\right)^{a_h}.
\end{align*}
This also covers $n=0$ and $m=n$; symmetry gives the case $m<n$.

For mutual intersections, put $M=\lfloor NT\rfloor$; if $M=0$,
all sums vanish. Suppose $M\ge1$. For $r\in\mathbb Z$, set
\[
 r_+=\max(r,0),\qquad r_-=\max(-r,0),\qquad
 \ell_{m,n}(r)=\max\{0,\min(m-r_+,n-r_-)\}.
\]
Grouping the pairs $(i,\ell)$ by their index difference $r=i-\ell$
gives
\begin{align*}
 D^r_j&=S_{r_++j}-\widetilde S_{r_-+j},\qquad j\ge0,\\
 L^{m,n}_r(d)
   &:=\sum_{j=1}^{\ell_{m,n}(r)}\gamma_{H_1}(d+D^r_j)
     =\sum_{\substack{1\le i\le m,\ 1\le\ell\le n\\i-\ell=r}}
                   \gamma_{H_1}(d+S_i-\widetilde S_\ell),\\
 A^N_{m,n}(d)
   &=N^{-a_h}\sum_{r\in\mathbb Z}
                   \gamma_{H_0}(m-n-r)L^{m,n}_r(d).
\end{align*}
Conditionally on $D^r_0$, the walk $D^r$ has the increment law of
Lemma~\ref{tp:difference-local-time}. Since
$0\le\ell_{m,n}(r)\le M$, \eqref{tp:difference-occupation-bound}
and Lemmas~\ref{lem:difference-kernel-increments}
and~\ref{lem:difference-local-time-increments} yield
\begin{align*}
 \|L^{m,n}_r(d)\|_q
   &\le C_hq^{1-H_1}M^{H_1},\\
 \|L^{m,n}_r(d)-L^{m,n}_r(d')\|_q
   &\le
   \begin{cases}
    C_h\sqrt q\,M^{1/2}|d-d'|^{2H_1-1},
       &H_1>\frac12,\\
    C_qM^{1/4}|d-d'|^{1/2},
       &H_1=\frac12,\ d-d'\in2\mathbb Z.
   \end{cases}
\end{align*}
Also $L^{m,n}_r=0$ for $|r|\ge M$. Minkowski's inequality and
\eqref{fr:lattice-envelope} therefore give
\begin{align*}
 \|A^N_{m,n}(d)\|_q
 &\le C_hq^{1-H_1}N^{-a_h}M^{H_1}
           \sum_{|r|<M}\gamma_{H_0}(m-n-r)\\
 &\le C_hq^{1-H_1}N^{-a_h}M^{H_1}(M+1)^{2H_0-1}
 \le C_{h,T}\sqrt q.
\end{align*}
The same summation gives
\[
 \|A^N_{m,n}(d)-A^N_{m,n}(d')\|_q
 \le C_{q,h}N^{-a_h}(M+1)^{2H_0-1}
 \begin{cases}
  M^{1/2}|d-d'|^{2H_1-1},&H_1>\frac12,\\
  M^{1/4}|d-d'|^{1/2},&H_1=\frac12,\ d-d'\in2\mathbb Z.
 \end{cases}
\]
With $\rho=|d-d'|/\sqrt N$, these estimates and
$a_h=2H_0+H_1-1$ imply
\[
 \|A^N_{m,n}(d)-A^N_{m,n}(d')\|_q
 \le C_{q,h,T}
 \begin{cases}
  1\wedge\rho^{2H_1-1},&H_1>\frac12,\\
  1\wedge\rho^{1/2},&H_1=\frac12,\ d-d'\in2\mathbb Z
 \end{cases}
 \le C_{q,h,T}\rho^{2\eta},
\]
which proves \eqref{eq:discrete-mutual-space}.

For the time increment, assume $m=m'+k$, $k\ge1$. Then
\begin{align*}
 0&\le\ell_{m,n}(r)-\ell_{m',n}(r)\le k,\\
 L^{m,n}_r(d)-L^{m',n}_r(d)
   &=\sum_{j=\ell_{m',n}(r)+1}^{\ell_{m,n}(r)}
                    \gamma_{H_1}(d+D^r_j),\\
 \|L^{m,n}_r(d)-L^{m',n}_r(d)\|_q
   &\le C_hq^{1-H_1}k^{H_1},
\end{align*}
the last bound following from \eqref{tp:difference-occupation-bound}
conditionally on $D^r_{\ell_{m',n}(r)}$. 

The exact decomposition
\begin{align*}
 A^N_{m,n}(d)-A^N_{m',n}(d)
 &=N^{-a_h}\sum_r
     \bigl[\gamma_{H_0}(m-n-r)-\gamma_{H_0}(m'-n-r)\bigr]
                  L^{m',n}_r(d)\\
 &\quad+N^{-a_h}\sum_r\gamma_{H_0}(m-n-r)
                     \bigl[L^{m,n}_r(d)-L^{m',n}_r(d)\bigr]
\end{align*}
and \eqref{fr:lattice-translation} and \eqref{fr:lattice-envelope} yield
\begin{align*}
 \|A^N_{m,n}(d)-A^N_{m',n}(d)\|_q
 &\le C_hq^{1-H_1}N^{-a_h}
   \left(M^{H_1}\sum_{r\in\mathbb Z}
            |\gamma_{H_0}(r+k)-\gamma_{H_0}(r)|
          +k^{H_1}(M+1)^{2H_0-1}\right)\\
 &\le C_{q,h,T}
       \left(\left(\frac kN\right)^{2H_0-1}
                      +\left(\frac kN\right)^{H_1}\right)
 \le C_{q,h,T}\left(\frac kN\right)^{2H_0-1},
\end{align*}
where $2H_0-1<1-H_1\le H_1$ by \eqref{fr:range}, and $k/N\le T$.
The joint symmetry
\[
 \bigl(A^N_{m,n}(d),A^N_{m,n'}(d)\bigr)
 \stackrel{\mathrm{law}}=
 \bigl(A^N_{n,m}(-d),A^N_{n',m}(-d)\bigr)
\]
follows by interchanging the walks. Together with the triangle
inequality, it proves \eqref{eq:discrete-mutual-time}.

For the shifted time cells,
\[
 A^N_{m;\theta}(d)
   =N^{-a_h}\sum_{r\in\mathbb Z}
                  \gamma_{H_0}(\theta-r)L^{m,m}_r(d).
\]
By \eqref{fr:lattice-translation} and \eqref{fr:lattice-envelope},
\begin{align*}
 \|A^N_{m;\theta}(d)-A^N_{m;0}(d)\|_q
 &\le C_hq^{1-H_1}M^{H_1}N^{-a_h}
       \sum_{r\in\mathbb Z}
           |\gamma_{H_0}(\theta-r)-\gamma_{H_0}(-r)|\\
 &\le C_{q,h,T}N^{-(2H_0-1)}|\theta|^{2H_0-1},\\
 \|A^N_{m;\theta}(d)\|_q
 &\le C_hq^{1-H_1}M^{H_1}N^{-a_h}
                    \sum_{|r|<M}\gamma_{H_0}(\theta-r)
 \le C_{h,T}\sqrt q.
\end{align*}
This proves \eqref{eq:discrete-anchor-increment} and completes
\eqref{eq:discrete-two-point-moment}.

Finally, Lemma~\ref{lem:discrete-scaling} and
\eqref{eq:YN-exp-unif} give, for every $\lambda>0$,
\[
 E\bigl[e^{\lambda|Y^N_{m/N}|}\bigr]
 \le E\bigl[e^{\lambda Y^N_{m/N}}\bigr]
      +E\bigl[e^{-\lambda Y^N_{m/N}}\bigr]
 \le 2e^{C_h\lambda^2T^{2a_h}}.
\]
All the mutual intersections above are nonnegative.
Their bounds by $C_{h,T}\sqrt q$, valid for every $q\ge2$,
and Corollary~\ref{tp:moment-to-exp} give the remaining bounds in \eqref{eq:discrete-path-exponential-moments}.
\end{proof}

\

We now turn the path estimates of Lemma~\ref{lem:discrete-path-increments}
into increment bounds for the partition field and its bilinear
interpolation $\operatorname{Interp}_N$, defined in
\eqref{eq:bilinear-interpolation}.

Recall from \eqref{eq:increment-exponent-used} that
\[
 \eta=
 \begin{cases}
  \eta_h=\min\{2H_0-1,H_1-\frac12\},&H_1>\frac12,\\[2pt]
  \eta_{H_0}=\min\{2H_0-1,\frac14\},&H_1=\frac12.
 \end{cases}
\]

\begin{proposition}[Discrete field increments]
\label{prop:discrete-increments}
For every finite $q\ge2$, $T,L>0$, and $(t,x),(s,y)$ in
$\mathcal K_{T,L}=[0,T]\times[-L,L]$, uniformly in $N$,
\begin{align}
 \|(\operatorname{Interp}_N\widehat{\mathsf Z}^N)(t,x)
       -(\operatorname{Interp}_N\widehat{\mathsf Z}^N)(s,y)\|_q
 &\le C_{q,h,T,L}
        \bigl(|t-s|^{\eta/2}+|x-y|^{\eta}\bigr),
                         \label{eq:discrete-field-increments}\\
 \|(\operatorname{Interp}_N\widehat{\mathsf Z}^N)(t,x)-\widehat{\mathsf Z}^N_t(x)\|_q
 &\le C_{q,h,T,L}N^{-\eta/2}.
                         \label{eq:interpolation-pointwise-error}
\end{align}
Both estimates also hold for $1\le q<2$ by monotonicity of norms.
\end{proposition}

\begin{proof}
For fixed $N$, write
$Q_m(j)=\widehat{\mathsf Z}^N_{m/N}\bigl(\frac j{\sqrt N}\bigr)$.
Fix $q\ge2$ and put $\tau_q=\frac12\log(q-1)$. As in the proof of
Proposition~\ref{prop:continuous-functional}, set
\[
 Q^{\tau_q}_m(j)=e^{R^N_{m/N}}
 E_S\left[e^{\frac12 Y^N_{m/N}}
          \cE(e^{\tau_q}h^{S,N}_{m/N,j/\sqrt N})\right].
\]
This is a finite sum of Gaussian exponentials, and
\eqref{eq:OU-exponential} gives $P_{\tau_q}Q^{\tau_q}_m(j)=Q_m(j)$. By
\eqref{eq:Nelson},
\begin{equation}\label{eq:discrete-U-by-V}
 \|Q_m(j)-Q_n(k)\|_q\le\|Q^{\tau_q}_m(j)-Q^{\tau_q}_n(k)\|_2.
\end{equation}
As in \eqref{eq:continuous-amplified-two-point},
\begin{equation}\label{eq:discrete-amplified-two-point}
 \E\bigl[Q^{\tau_q}_m(j)Q^{\tau_q}_n(k)\bigr]
 =e^{R^N_{m/N}+R^N_{n/N}}E\left[e^{\frac12 Y^N_{m/N}+\frac12\widetilde Y^N_{n/N}+(q-1)A^N_{m,n}(j-k)}\right].
\end{equation}

\

\noindent
\textit{First bounds on the grid.}
Put
\[
 F^{N,q}_{m,n}(d)
 :=R^N_{m/N}+R^N_{n/N}
   +\tfrac12\bigl(Y^N_{m/N}+\widetilde Y^N_{n/N}\bigr)
   +(q-1)A^N_{m,n}(d).
\]
Lemma~\ref{lem:discrete-path-increments} gives
\[
 \sup_{\substack{N\ge1,\ 0\le m,n\le\lfloor NT\rfloor\\d\in\mathbb Z}}
       \bigl\|e^{F^{N,q}_{m,n}(d)}\bigr\|_2
 \le C_{q,h,T}.
\]
For path variables $X,X'$, Cauchy-Schwarz yields
\[
 |E[e^X]-E[e^{X'}]|
 \le E\bigl[|X-X'|(e^X+e^{X'})\bigr]
 \le\|X-X'\|_2\bigl(\|e^X\|_2+\|e^{X'}\|_2\bigr).
\]
Hence \eqref{eq:discrete-amplified-two-point} and
\eqref{eq:discrete-self-time-increment}-\eqref{eq:discrete-mutual-time}
imply
\begin{align*}
 \|Q^{\tau_q}_m(j)-Q^{\tau_q}_n(j)\|_2^2
 &=E\bigl[e^{F^{N,q}_{m,m}(0)}
          +e^{F^{N,q}_{n,n}(0)}-2e^{F^{N,q}_{m,n}(0)}\bigr]\\
 &\le C_{q,h,T}\bigl(
       \|F^{N,q}_{m,m}(0)-F^{N,q}_{m,n}(0)\|_2
      +\|F^{N,q}_{n,n}(0)-F^{N,q}_{m,n}(0)\|_2\bigr)\\
 &\le C_{q,h,T}\left(
       \left|\frac{m-n}{N}\right|^{a_h}
      +\left|\frac{m-n}{N}\right|^{2H_0-1}\right)
 \le C_{q,h,T}\left|\frac{m-n}{N}\right|^{2H_0-1}.
\end{align*}
Similarly, \eqref{eq:discrete-mutual-space} gives
\begin{align*}
 \|Q^{\tau_q}_m(j)-Q^{\tau_q}_m(k)\|_2^2
 &=2E\bigl[e^{F^{N,q}_{m,m}(0)}
                  -e^{F^{N,q}_{m,m}(j-k)}\bigr]\\
 &\le C_{q,h,T}\|A^N_{m,m}(0)-A^N_{m,m}(j-k)\|_2\\
 &\le C_{q,h,T}\left|\frac{j-k}{\sqrt N}\right|^{2\eta},
\end{align*}
with $j-k\in2\mathbb Z$ if $H_1=\frac12$.
Taking square roots and using \eqref{eq:discrete-U-by-V} yields
\begin{align}
 \|Q_m(j)-Q_n(j)\|_q
       &\le C_{q,h,T}\left|\frac{m-n}N\right|^{H_0-\frac12},
                            \label{eq:discrete-grid-time}\\
 \|Q_m(j)-Q_m(k)\|_q
       &\le C_{q,h,T}\left|\frac{j-k}{\sqrt N}\right|^{\eta},
                            \label{eq:discrete-grid-even-space}
\end{align}
with the same parity restriction in the second estimate.

\

\noindent
\textit{Odd displacements when $H_1=\frac12$.}
For $m\ge1$, put
\[
 \omega_{N,m}(j)
 =W\left(N^{H_0+1/4}
       \one_{((m-1)/N,m/N]\times\mathfrak J_j^N}\right)
 \sim\mathcal N(0,1).
\]
Conditioning on the first walk step, with the noise fixed, gives
\begin{equation}\label{eq:discrete-one-step}
 Q_m(j)=\frac12\sum_{\sigma=\pm1}
       e^{\beta_N\omega_{N,m}(j+\sigma)-c_{\rm disc}\beta_N^2}
       Q_{m-1}(j+\sigma).
\end{equation}

For $\xi\sim\mathcal N(0,1)$ and $\beta_N\le1$,
\[
 \|e^{\beta_N\xi-c_{\rm disc}\beta_N^2}-1\|_{2q}
 \le\|\beta_N\xi-c_{\rm disc}\beta_N^2\|_{4q}
       \bigl(1+\|e^{\beta_N\xi-c_{\rm disc}\beta_N^2}\|_{4q}\bigr)
 \le C_{q,H_0}\beta_N.
\]
H\"older's inequality and Proposition~\ref{prop:discrete-uniform-moments}
therefore give, without any independence assumption on the factors,
\begin{align*}
 &\left\|Q_m(j)-\tfrac12
             \bigl(Q_{m-1}(j+1)+Q_{m-1}(j-1)\bigr)\right\|_q\\
 &\quad\le\frac12\sum_{\sigma=\pm1}
       \|e^{\beta_N\omega_{N,m}(j+\sigma)-c_{\rm disc}\beta_N^2}-1\|_{2q}
       \|Q_{m-1}(j+\sigma)\|_{2q}
 \le C_{q,H_0,T}\beta_N.
\end{align*}
Combining this with
\eqref{eq:discrete-grid-time}-\eqref{eq:discrete-grid-even-space},
\begin{equation}\label{eq:discrete-adjacent-space}
\begin{aligned}
 &\|Q_m(j+1)-Q_m(j)\|_q\\
 &\quad\le
   \left\|Q_m(j)-\tfrac12
       \bigl(Q_{m-1}(j+1)+Q_{m-1}(j-1)\bigr)\right\|_q\\
 &\qquad+\tfrac12\|Q_{m-1}(j+1)-Q_{m-1}(j-1)\|_q
          +\|Q_m(j+1)-Q_{m-1}(j+1)\|_q\\
 &\quad\le C_{q,H_0,T}
       \bigl(\beta_N+N^{-\eta/2}+N^{-(H_0-\frac12)}\bigr)
 \le C_{q,H_0,T}N^{-\eta/2},
\end{aligned}
\end{equation}
since $\beta_N=N^{-H_0+1/4}$ and $\eta\le2H_0-1$.
For odd $j-k$, choose $\sigma=\operatorname{sgn}(j-k)$; then
$j-k-\sigma\in2\mathbb Z$ and
\begin{align*}
 \|Q_m(j)-Q_m(k)\|_q
 &\le\|Q_m(j)-Q_m(k+\sigma)\|_q
       +\|Q_m(k+\sigma)-Q_m(k)\|_q\\
 &\le C_{q,H_0,T}
       \left(\left|\frac{j-k-\sigma}{\sqrt N}\right|^\eta
                         +N^{-\eta/2}\right)
 \le C_{q,H_0,T}\left|\frac{j-k}{\sqrt N}\right|^\eta.
\end{align*}
Together with \eqref{eq:discrete-grid-even-space} and $Q_0\equiv1$,
this proves, in both spatial regimes,
\begin{equation}\label{eq:discrete-grid-all-space}
 \|Q_m(j)-Q_m(k)\|_q
 \le C_{q,h,T}\left|\frac{j-k}{\sqrt N}\right|^\eta,
 \qquad j,k\in\mathbb Z.
\end{equation}

\

\noindent
\textit{Interpolation.} Put $F_N=\operatorname{Interp}_N\widehat{\mathsf Z}^N$ and use
\eqref{eq:discrete-grid-time} and \eqref{eq:discrete-grid-all-space}
with $T+1$ in place of $T$, to include the last interpolation cell.
For fixed $x$, put $j=\lfloor\sqrt N x\rfloor$ and
$\zeta=\sqrt N x-j$. If $s,t\in[m/N,(m+1)/N]$, the definition
\eqref{eq:bilinear-interpolation} gives
\begin{align*}
 \|F_N(t,x)-F_N(s,x)\|_q
 &=N|t-s|\left\|\sum_{b=0}^1\lambda_b(\zeta)
                    \bigl(Q_{m+1}(j+b)-Q_m(j+b)\bigr)\right\|_q\\
 &\le C_{q,h,T}N|t-s|N^{-(H_0-\frac12)}
 \le C_{q,h,T}|t-s|^{H_0-\frac12},
\end{align*}
since $N|t-s|\le1$. Similarly, for fixed $t$, put
$m=\lfloor Nt\rfloor$ and $\theta=Nt-m$. If
$x,y\in[j/\sqrt N,(j+1)/\sqrt N]$, then
\begin{align*}
 \|F_N(t,x)-F_N(t,y)\|_q
 &=\sqrt N|x-y|\left\|\sum_{a=0}^1\lambda_a(\theta)
                    \bigl(Q_{m+a}(j+1)-Q_{m+a}(j)\bigr)\right\|_q\\
 &\le C_{q,h,T}\sqrt N|x-y|N^{-\eta/2}
 \le C_{q,h,T}|x-y|^\eta,
\end{align*}
using $\sqrt N|x-y|\le1$ and $0<\eta\le1$.

At time grid points, convexity of the interpolation in space gives
\[
 \|F_N(m/N,x)-F_N(n/N,x)\|_q
 \le\sum_{b=0}^1\lambda_b(\zeta)
                   \|Q_m(j+b)-Q_n(j+b)\|_q
 \le C_{q,h,T}\left|\frac{m-n}{N}\right|^{H_0-\frac12};
\]
the analogous spatial bound at spatial grid points follows by
convexity in time. If $s\le t$ do not lie in one time grid interval,
let $u,v$ be the first and last time grid points in $[s,t]$. The
same-cell bounds and the preceding grid-point bound yield
\begin{align*}
 \|F_N(t,x)-F_N(s,x)\|_q
 &\le C_{q,h,T}\bigl(
       (u-s)^{H_0-\frac12}+(v-u)^{H_0-\frac12}
                              +(t-v)^{H_0-\frac12}\bigr)\\
 &\le C_{q,h,T}|t-s|^{H_0-\frac12}.
\end{align*}
Splitting $[x,y]$, for $x\le y$, at its first and last spatial grid
points gives the spatial bound in the same way. Thus, for arbitrary
$(t,x),(s,y)\in\mathcal K_{T,L}$, since $\eta\le2H_0-1$,
\begin{align*}
 \|F_N(t,x)-F_N(s,y)\|_q
 &\le\|F_N(t,x)-F_N(s,x)\|_q
       +\|F_N(s,x)-F_N(s,y)\|_q\\
 &\le C_{q,h,T,L}
          \bigl(|t-s|^{H_0-\frac12}+|x-y|^\eta\bigr)\\
 &\le C_{q,h,T,L}
          \bigl(|t-s|^{\eta/2}+|x-y|^\eta\bigr).
\end{align*}
This is \eqref{eq:discrete-field-increments}.

\

\noindent
\textit{Comparison with the pointwise approximation.}
For $m=N_t$, $j=\widetilde N_x$, set
$\theta=Nt-m\in[0,1)$ and $\zeta=\sqrt N x-j\in[0,1)$.
Translation of the time cells and
\eqref{eq:cross-time-cell-covariance} give
\begin{align*}
 Y^N_t&=Y^N_{m/N},\qquad R^N_t=R^N_{m/N},\\
 \ip{h^{S,N}_{t,x}}{h^{\widetilde S,N}_{m/N,j/\sqrt N}}_{\Hh}
     &=A^N_{m;\theta}(0),\\
 \ip{h^{S,N}_{t,x}}{h^{\widetilde S,N}_{t,x}}_{\Hh}
     &=\ip{h^{S,N}_{m/N,j/\sqrt N}}
            {h^{\widetilde S,N}_{m/N,j/\sqrt N}}_{\Hh}
       =A^N_{m;0}(0).
\end{align*}
Write
\[
 F^{N,q}_{m;\theta}
 :=2R^N_{m/N}+\tfrac12\bigl(Y^N_{m/N}+\widetilde Y^N_{m/N}\bigr)
                  +(q-1)A^N_{m;\theta}(0).
\]
Applying \eqref{eq:Nelson} after replacing each $\cE(h)$ by
$\cE(e^{\tau_q}h)$, and then the exponential comparison above, yields
\begin{align*}
 \|\widehat{\mathsf Z}^N_t(x)-Q_m(j)\|_q^2
 &\le2E\bigl[e^{F^{N,q}_{m;0}}-e^{F^{N,q}_{m;\theta}}\bigr]\\
 &\le C_{q,h,T}\|A^N_{m;0}(0)-A^N_{m;\theta}(0)\|_2\\
 &\le C_{q,h,T}N^{-(2H_0-1)}|\theta|^{2H_0-1},
\end{align*}
by \eqref{eq:discrete-anchor-increment} and
\eqref{eq:discrete-path-exponential-moments}. Thus
\begin{equation}\label{eq:raw-to-grid-time}
 \|\widehat{\mathsf Z}^N_t(x)-Q_m(j)\|_q
 \le C_{q,h,T}N^{-(H_0-\frac12)}|\theta|^{H_0-\frac12}.
\end{equation}
For $m=0$, both variables are $1$. Finally,
\eqref{eq:bilinear-interpolation},
\eqref{eq:discrete-grid-time} and \eqref{eq:discrete-grid-all-space}
give
\begin{align*}
 \|F_N(t,x)-Q_m(j)\|_q
 &\le\sum_{a,b\in\{0,1\}}\lambda_a(\theta)\lambda_b(\zeta)
                  \|Q_{m+a}(j+b)-Q_m(j)\|_q\\
 &\le C_{q,h,T,L}
              \bigl(N^{-(H_0-\frac12)}+N^{-\eta/2}\bigr).
\end{align*}
Consequently, by \eqref{eq:raw-to-grid-time} and since $\eta\le2H_0-1$,
\[
 \|F_N(t,x)-\widehat{\mathsf Z}^N_t(x)\|_q
 \le\|F_N(t,x)-Q_m(j)\|_q
       +\|Q_m(j)-\widehat{\mathsf Z}^N_t(x)\|_q
 \le C_{q,h,T,L}N^{-\eta/2},
\]
which is \eqref{eq:interpolation-pointwise-error}. 
\end{proof}

\

We now combine convergence of the $\mathcal S$-transforms and second
moments with the increment estimates to prove
Theorem~\ref{thm:coupled-L2}. We then transfer the result to the
original environment through \eqref{eq:interpolated-field-law},
completing the proof of Theorem~\ref{thm:fractional-main}.

\subsection{Proof of Theorem~\ref{thm:coupled-L2}}
\label{subsec:proof-coupled-L2}

\

\smallskip

\noindent
\textit{Pointwise convergence.}
Fix $t>0$ and $x\in\R$, the case $t=0$ being trivial.
Corollary~\ref{prop:second-moment-convergence} gives
\[
 M_{t,x}:=\sup_N\|\widehat{\mathsf Z}^N_t(x)\|_2+\|U_t(x)\|_2<\infty.
\]
Let $\mathscr V$ be the span in Lemma~\ref{lem:total}.
Corollary~\ref{prop:S-convergence} gives
$\E[(\widehat{\mathsf Z}^N_t(x)-U_t(x))\Phi]\to0$
for every $\Phi\in\mathscr V$.

Since $U_t(x)\in L^2(\Omega,\sigma(W),\Pp)$, for every
$\varepsilon>0$ we can choose $\Phi_\varepsilon\in\mathscr V$ with
$\|U_t(x)-\Phi_\varepsilon\|_2<\varepsilon$.
Cauchy-Schwarz gives
\[
 \left|\E\left[
       (\widehat{\mathsf Z}^N_t(x)-U_t(x))U_t(x)\right]\right|
 \le\left|\E\left[
       (\widehat{\mathsf Z}^N_t(x)-U_t(x))\Phi_\varepsilon\right]\right|
       +M_{t,x}\varepsilon.
\]
Taking first $N\to\infty$ and then $\varepsilon\downarrow0$ yields
\[
 \E[\widehat{\mathsf Z}^N_t(x)U_t(x)]
 \longrightarrow\E[U_t(x)^2].
\]
Together with Corollary~\ref{prop:second-moment-convergence},
this gives
\[
 \E\left[|\widehat{\mathsf Z}^N_t(x)-U_t(x)|^2\right]
 =\E\left[\widehat{\mathsf Z}^N_t(x)^2\right]
    +\E\left[U_t(x)^2\right]
    -2\E\left[\widehat{\mathsf Z}^N_t(x)U_t(x)\right]\longrightarrow0.
\]
This proves $L^2$ convergence at the fixed point $(t,x)$.
For $1\le p\le2$, it implies $L^p$ convergence directly. If $p>2$,
choose $q>p$ and $\nu\in(0,1)$ such that
$\frac1p=\frac\nu2+\frac{1-\nu}{q}$. Interpolation yields
\[
 \|\widehat{\mathsf Z}^N_t(x)-U_t(x)\|_p
 \le\|\widehat{\mathsf Z}^N_t(x)-U_t(x)\|_2^{\nu}
       \|\widehat{\mathsf Z}^N_t(x)-U_t(x)\|_q^{1-\nu}
 \longrightarrow0.
\]
The second factor is uniformly bounded by
Proposition~\ref{prop:discrete-uniform-moments} and the continuous
moment bounds in \eqref{eq:continuous-limit-moments}. This proves
\eqref{eq:main-L2} for every finite $p\ge1$, at every point
$(t,x)$.

\

\noindent
\textit{Functional convergence.}
Fix $T,L>0$ and set $\mathcal K_{T,L}=[0,T]\times[-L,L]$.
Let $\eta>0$ be the exponent in \eqref{eq:increment-exponent-used}.
Proposition~\ref{prop:discrete-increments} gives, for every finite
$p\ge1$ and every $(t,x)\in\mathcal K_{T,L}$,
\begin{align*}
 \|(\operatorname{Interp}_N\widehat{\mathsf Z}^N)(t,x)-U_t(x)\|_p
 &\le\|(\operatorname{Interp}_N\widehat{\mathsf Z}^N)(t,x)
                                      -\widehat{\mathsf Z}^N_t(x)\|_p
       +\|\widehat{\mathsf Z}^N_t(x)-U_t(x)\|_p\\
 &\le C_{p,h,T,L}N^{-\eta/2}
       +\|\widehat{\mathsf Z}^N_t(x)-U_t(x)\|_p
 \longrightarrow0.
\end{align*}
The same proposition supplies uniform increment bounds for the
continuous fields $\operatorname{Interp}_N\widehat{\mathsf Z}^N$,
whose values at $(0,0)$ are all $1$.
Lemma~\ref{lem:compact-upgrade}, applied to these fields and the
continuous field $U$, therefore yields
\[
 \|\operatorname{Interp}_N\widehat{\mathsf Z}^N-U
                  \|_{L^p(\Omega;C(\mathcal K_{T,L}))}
 \longrightarrow0,\qquad 1\le p<\infty,
\]
which is \eqref{eq:main-functional}.
Applying this conclusion on the rectangles $[0,k]\times[-k,k]$,
$k\in\N$, also proves convergence in probability in $C_{\rm loc}$.
\qed

\

\subsection{Proof of Theorem~\ref{thm:fractional-main}}

The continuous part of the theorem has already been established:
Theorem~\ref{thm:continuous-field} gives a nonnegative continuous
field $U$, measurable with respect to $W_h$, with $U_0(x)=1$,
and the convergence \eqref{fr:continuous-functional}.
Definition~\ref{def:continuous-solution} identifies $U$ as the
renormalized Stratonovich solution specified by
\eqref{fr:regularized-equation}; Proposition~\ref{prop:CDOT-comparison}
gives the identification with the solution of
\cite[Definition~3.4]{CDOT-KRP} when $2H_0+H_1>\frac53$.
It remains to prove \eqref{fr:original-functional}
and \eqref{fr:compact-moments} for the original discrete model.

By Theorem~\ref{thm:coupled-L2},
$\operatorname{Interp}_N\widehat{\mathsf Z}^N\to U$ in probability in
$C_{\rm loc}([0,\infty)\times\R)$, and hence in distribution in
that space. Equality of the field laws in
\eqref{eq:interpolated-field-law} transfers this convergence to
$\operatorname{Interp}_N\mathsf Z^N$, proving \eqref{fr:original-functional}.
Since restriction to a compact $\mathcal K\subset[0,\infty)\times\R$ is a continuous map into $C(\mathcal K)$, the continuous mapping theorem also gives convergence in law in $C(\mathcal K)$ for every such $\mathcal K$.

Fix such a $\mathcal K$, choose finite $T,L>0$ with
$\mathcal K\subset[0,T]\times[-L,L]$, and let $1\le p<\infty$.
The elementary norm inequality and \eqref{eq:main-functional} give
\[
 \left\|\|\operatorname{Interp}_N\widehat{\mathsf Z}^N\|_{C(\mathcal K)}
             -\|U\|_{C(\mathcal K)}\right\|_{L^p(\Omega)}
 \le\left\|\operatorname{Interp}_N\widehat{\mathsf Z}^N-U
                   \right\|_{L^p(\Omega;C([0,T]\times[-L,L]))}
 \longrightarrow0.
\]
Continuity of the $L^p$ norm yields convergence of the $p$th moments
of these supremum norms. Equality in law
\eqref{eq:interpolated-field-law} then proves
\eqref{fr:compact-moments}.
\qed

\subsection{Necessity of the lattice counterterm}
\label{subsec:discrete-necessity}

The convergence at a fixed point now gives the non-tightness assertion
of Proposition~\ref{prop:discrete-renormalization-needed}, by the
argument used for Proposition~\ref{prop:renormalization-needed}.

\begin{proof}[Proof of Proposition~\ref{prop:discrete-renormalization-needed}]
The covariances $\gamma_{H_0}$ and $\gamma_{H_1}$ are nonnegative, so
that $c_{\rm disc}(h)\ge\frac12$ by \eqref{fr:counterterm-discrete}.
For $t>0$ we have $N_t\sim Nt$, and $a_h<1$; hence
\[
 c_{\rm disc}(h)N_t\beta_{N,h}^2\ge\frac12N_tN^{-a_h}
 \longrightarrow+\infty.
\]
Fix $t>0$, $x\in\R$ and $K,\delta>0$. By
\eqref{eq:unrenormalized-partition}, as soon as
$Ke^{-c_{\rm disc}(h)N_t\beta_{N,h}^2}\le\delta$,
\[
 \Pp\bigl(\mathsf z^N_t(x)\le K\bigr)
 \le\Pp\bigl(\mathsf Z^N_t(x)\le\delta\bigr).
\]
By Lemma~\ref{lem:environment-law}, $\mathsf Z^N_t(x)$ has the law of
$\widehat{\mathsf Z}^N_t(x)$, which converges to $U_t(x)$ in
$L^1(\Omega)$ by Theorem~\ref{thm:coupled-L2}. Hence
$\mathsf Z^N_t(x)\to U_t(x)$ in law, and the Portmanteau theorem
applied to the closed set $(-\infty,\delta]$ gives
$\limsup_{N\to\infty}\Pp(\mathsf z^N_t(x)\le K)\le
\Pp(U_t(x)\le\delta)$. Letting $\delta\downarrow0$ and using
$U_t(x)\ge0$ yields the bound $\Pp(U_t(x)=0)$, which is smaller than
one because $\E[U_t(x)]>0$ (Theorem~\ref{thm:continuous-summary}).
Consequently,
$\liminf_{N\to\infty}\Pp(\mathsf z^N_t(x)>K)\ge1-\Pp(U_t(x)=0)>0$
for every $K$, so that no sequence $(\mathsf z^{N_k}_t(x))$ is tight,
and therefore none converges in law.
\end{proof}

\begin{remark}[Divergence of the mean]\label{rem:discrete-mean}
The divergence is already visible on the mean, by an elementary
computation which does not rely on Theorem~\ref{thm:coupled-L2}. Fix
$t>0$ and $x\in\R$. By \eqref{fr:discrete-covariance}, conditionally
on the walk, the exponent in \eqref{eq:unrenormalized-partition} is a
centered Gaussian variable, whose variance is the sum
$\Gamma^N_t(S)$ displayed in \eqref{eq:GammaN}. Fubini's theorem,
Jensen's inequality and \eqref{eq:RN} then give
\[
 \E\bigl[\mathsf z^N_t(x)\bigr]
 =E_S\Bigl[e^{\frac12\Gamma^N_t(S)}\Bigr]
 \ge e^{\frac12E_S[\Gamma^N_t(S)]}
 =e^{c_{\rm disc}(h)N_t\beta_{N,h}^2+R^N_t}.
\]
By Lemma~\ref{lem:discrete-scaling},
Proposition~\ref{prop:counterterm} and \eqref{eq:Yt-alphat-Rt},
\[
 R^N_t=\theta^N_tR^{N_t}_1\longrightarrow t^{a_h}R_1=R_t,
 \qquad
 c_{\rm disc}(h)N_t\beta_{N,h}^2\sim c_{\rm disc}(h)\,t\,N^{1-a_h}.
\]
The mean of the partition function without counterterm thus diverges,
and the counterterm removes exactly the divergent part of the exponent
$\frac12E_S[\Gamma^N_t(S)]$.
\end{remark}

\clearpage
\part{Appendices}\label{chap:appendices-simulations}

\appendix
\section{\texorpdfstring{Gaussian and lattice regularizations of $K_H$}{Gaussian and lattice regularizations of K\_H}}
\label{app:temporal}

For $\frac12<H<1$, recall the kernel $K_H(x)=c_H|x|^{2H-2}$ of
\eqref{fr:spatial-kernel} and its Gaussian regularizations
$K_H^{[v]}=K_H*\heat_v$ of \eqref{fr:kernels}, where $\heat_v$ is the
centered Gaussian density with variance $v$ of \eqref{eq:heat-kernel}.
Write $\kappa_H=E[K_H(Z)]$ for a standard
normal variable $Z$.

\begin{lemma}[Gaussian regularizations of $K_H$]\label{lem:kernel-regularization}
Let $\frac12<H<1$. There is a constant $C_H$ such that, for every
$v\ge0$, $x\ne0$, $R>0$ and $z\in\R$,
\begin{align}
 0\le K_H^{[v]}(x)&\le C_H(|x|+\sqrt v)^{2H-2}\le C_H|x|^{2H-2},
                                             \label{eq:KH-envelope}\\
 \int_{|x|\le R}K_H^{[v]}(x)\dd x&\le C_HR^{2H-1},
                                             \label{eq:KH-mass}\\
 \|K_H^{[v]}(\cdot+z)-K_H^{[v]}\|_{L^1(\R)}
 &\le C_H|z|^{2H-1}.                          \label{eq:KH-translation}
\end{align}
Moreover, $K_H^{[u]}(0)=\kappa_Hu^{H-1}$ for every $u>0$, and, for
every $\eta>0$,
\begin{equation}\label{eq:KH-convergence}
 \sup_{|x|\ge\eta}\bigl|K_H^{[v]}(x)-K_H(x)\bigr|
 \longrightarrow0\qquad(v\downarrow0).
\end{equation}
\end{lemma}

\begin{proof}
Put $\beta=2-2H\in(0,1)$, so that $K_H(x)=c_H|x|^{-\beta}$. For $v>0$,
the change of variables $y=\sqrt v\,w$ gives
\[
 K_H^{[v]}(x)=c_Hv^{-\beta/2}F(x/\sqrt v),
 \qquad F(u)=\int_\R |u-w|^{-\beta}\heat_1(w)\dd w.
\]
For $|u|\le2$, splitting at $|u-w|=1$ yields
\[
 F(u)\le\frac1{\sqrt{2\pi}}\int_{|y|\le1}|y|^{-\beta}\dd y+1
       =\frac{2}{\sqrt{2\pi}(1-\beta)}+1.
\]
For $|u|>2$, split instead at $|u-w|=|u|/2$. On $|u-w|\le|u|/2$ one
has $|w|\ge|u|/2$, so
\begin{align*}
 F(u)&\le \left(\frac2{|u|}\right)^\beta
       +\frac{e^{-u^2/8}}{\sqrt{2\pi}}
                     \int_{|y|\le|u|/2}|y|^{-\beta}\dd y\\
 &=\left(\frac2{|u|}\right)^\beta
       +\frac{2e^{-u^2/8}}{\sqrt{2\pi}(1-\beta)}
                     \left(\frac{|u|}{2}\right)^{1-\beta}
 \le C_H|u|^{-\beta}.
\end{align*}
Thus $F(u)\le C_H(1+|u|)^{-\beta}$, and scaling back gives
$K_H^{[v]}(x)\le C_H(|x|+\sqrt v)^{-\beta}$, which also holds for
$v=0$ and proves \eqref{eq:KH-envelope}. Integrating its last member
over $|x|\le R$ gives \eqref{eq:KH-mass}, the exponent $2H-1$ being
positive.

For \eqref{eq:KH-translation}, the case $z=0$ is immediate. For
$z\ne0$, homogeneity and symmetry give
\[
 \int_\R|K_H(x+z)-K_H(x)|\dd x
 =|z|^{2H-1}\int_\R|K_H(y+1)-K_H(y)|\dd y,
\]
and the last integral is finite: the singularities at $0$ and $-1$ are
integrable because $H>\frac12$, and the integrand is $O(|y|^{2H-3})$ at
infinity. Convolution with a probability density contracts this $L^1$
difference, which proves \eqref{eq:KH-translation} uniformly in $v$.
The value at the origin is
$K_H^{[u]}(0)=E[K_H(\sqrt u\,Z)]=u^{H-1}\kappa_H$.

Finally, fix $x$ with $|x|\ge\eta$ and put $d=|x|/2\ge\eta/2$. Since
$\int_\R\heat_v=1$,
\[
 |K_H^{[v]}(x)-K_H(x)|
 \le\int_\R|K_H(x-y)-K_H(x)|\heat_v(y)\dd y.
\]
On $|y|\le d$, every point between $x$ and $x-y$ has absolute value at
least $d$, and the mean value theorem bounds the integrand by
$c_H\beta d^{-\beta-1}|y|$. On $\{|y|>d,\ |x-y|\le d\}$, one has $|x-y|\le d<|x|$, hence $0\le K_H(x)\le K_H(x-y)$ and the integrand is at most $K_H(x-y)\heat_v(y)$; there the Gaussian factor is at most $(2\pi v)^{-1/2}e^{-d^2/(2v)}$, while
$\int_{|x-y|\le d}K_H(x-y)\dd y=2c_Hd^{1-\beta}/(1-\beta)$. On the
remaining set, both $K_H(x-y)$ and $K_H(x)$ are at most $c_Hd^{-\beta}$,
and $\int_{|y|>d}\heat_v\le d^{-1}\int_\R|y|\heat_v(y)\dd y$. With
$\int_\R|y|\heat_v(y)\dd y=\sqrt{2v/\pi}$, this gives
\[
 |K_H^{[v]}(x)-K_H(x)|
 \le C_H\left(d^{-\beta-1}\sqrt v
          +\frac{d^{1-\beta}}{\sqrt v}e^{-d^2/(2v)}\right).
\]
The function $\theta\mapsto\theta^{1-\beta}e^{-\theta^2/(2v)}$ is
nonincreasing on $[\sqrt v,\infty)$. For $v\le\eta^2/4$, the
supremum of the right-hand side over $d\ge\eta/2$ is therefore at most
\[
 C_H\left((\eta/2)^{-\beta-1}\sqrt v
   +\frac{(\eta/2)^{1-\beta}}{\sqrt v}e^{-\eta^2/(8v)}\right)
 \longrightarrow0\qquad(v\downarrow0),
\]
which proves \eqref{eq:KH-convergence}.
\end{proof}

The lattice kernel is obtained by averaging $K_H$ against the tent
function. The following estimates apply to both the temporal index
$H_0$ and the spatial index $H_1$.

\begin{lemma}[Bounds and translations of the lattice kernel]
\label{fr:kernel-lemma}
Let $\frac12<H<1$ and extend the kernel of \eqref{fr:fgn} to real
arguments by
\[
 \gamma_H(x)=\frac12\bigl(|x+1|^{2H}+|x-1|^{2H}-2|x|^{2H}\bigr),
 \qquad x\in\R.
\]
This function is nonnegative and continuously differentiable. There
is a constant $C_H$ such that, for every $x,z\in\R$ and $R\ge1$,
\begin{equation}\label{fr:lattice-envelope}
\begin{split}
 0\le\gamma_H(x)&\le C_H(1+|x|)^{2H-2},\\
 \sup_{v\in\R}\sum_{|j|\le R}\gamma_H(j+v)&\le C_HR^{2H-1}.
\end{split}
\end{equation}
Moreover,
\begin{equation}\label{fr:lattice-translation}
 \sum_{j\in\mathbb Z}|\gamma_H(j+z)-\gamma_H(j)|
 \le C_H|z|^{2H-1},
\end{equation}
and
\begin{equation}\label{fr:lattice-asymptotic}
\bigl|\gamma_H(x)-c_H|x|^{2H-2}\bigr|\le C_H|x|^{2H-4} \quad \text{for} \ |x|\ge2.
\end{equation}
At $H=\frac12$, one has $\gamma_{1/2}(x)=(1-|x|)_+$. The two bounds
in \eqref{fr:lattice-envelope} remain valid, with $2H-1=0$ in the
mass bound, and the translation estimate takes the explicit form
\[
 \sum_{j\in\mathbb Z}|\gamma_{1/2}(j+z)-\gamma_{1/2}(j)|
 =2\min\{|z|,1\}.
\]
In particular, on the lattice $\gamma_{1/2}(j)=\one_{\{j=0\}}$.
\end{lemma}

\begin{proof}
For $H>\frac12$, the second-difference identity gives
\begin{equation}\label{fr:tent-representation}
 \gamma_H(x)=\int_{-1}^1(1-|u|)K_H(x+u)\dd u,
 \qquad x\in\R.
\end{equation}
This follows by integrating twice the locally integrable second
derivative of $|x|^{2H}$. It proves nonnegativity. The explicit
second-difference formula also shows that $\gamma_H$ is continuously
differentiable, since $2H>1$.
For $|x|\ge2$, one has $|x+u|\ge|x|/2$ on the integration interval,
which gives the envelope in \eqref{fr:lattice-envelope}; on $[-2,2]$
it follows from continuity.

For the mass bound, put $\tau(x)=(1-|x|)_+$. Its integer
translates form a partition of unity, and the sum over $|j|\le R$
vanishes outside $[-R-1,R+1]$. Hence, for every $x\in\R$,
\begin{equation}\label{property-tau}
 \sum_{j\in\mathbb Z}\tau(x-j)=1,\qquad
 0\le\sum_{|j|\le R}\tau(x-j)
       \le\one_{\{|x|\le R+1\}}.
\end{equation}
Using \eqref{fr:tent-representation} and changing variables in
each integral, we obtain, for every $v\in\R$,
\begin{align*}
 \sum_{|j|\le R}\gamma_H(j+v)
 &=\sum_{|j|\le R}\int_\R\tau(u)K_H(j+v+u)\dd u\\
 &=\int_\R\left(\sum_{|j|\le R}\tau(u-j)\right)K_H(v+u)\dd u\\
 &\le\int_{-(R+1)}^{R+1}K_H(v+u)\dd u.
\end{align*}
Put $a=R+1$ and $\alpha=2H-1\in(0,1)$. Direct integration gives
\[
 \int_{-a}^{a}K_H(v+u)\dd u
 =\frac{c_H}{\alpha}
 \begin{cases}
  (a+v)^\alpha+(a-v)^\alpha,& |v|\le a,\\[1ex]
  (|v|+a)^\alpha-(|v|-a)^\alpha,& |v|>a.
 \end{cases}
\]
For $|v|\le a$, concavity gives $(a+v)^\alpha+(a-v)^\alpha\le2a^\alpha$. For $|v|>a$, subadditivity of $x\mapsto x^\alpha$ on $[0,\infty)$ gives
\[
 (|v|+a)^\alpha-(|v|-a)^\alpha
 \le  2 a^\alpha.
\]
Consequently,
\[
 \sup_{v\in\R}\int_{-(R+1)}^{R+1}K_H(v+u)\dd u
 \le\frac{2c_H}{2H-1}(R+1)^{2H-1}
 \le C_HR^{2H-1},
\]
where the last inequality uses $R\ge1$. Together with the preceding
bound on the lattice sum, this proves \eqref{fr:lattice-envelope}.

For the asymptotic expansion, Taylor's formula, uniformly for
$|u|\le1$ and $|x|\ge2$, gives
\[
 K_H(x+u)=K_H(x)+uK_H'(x)+O_H(u^2|x|^{2H-4}).
\]
The tent function has integral one and vanishing first moment.
Integrating this identity proves \eqref{fr:lattice-asymptotic}.

Finally, changing variables in
\eqref{fr:tent-representation}, then applying Tonelli's theorem, \eqref{property-tau} and
\eqref{eq:KH-translation} with $v=0$, gives
\begin{align*}
 \sum_j|\gamma_H(j+z)-\gamma_H(j)|
 &\le\int_\R|K_H(y+z)-K_H(y)|\sum_j\tau(y-j)\dd y\\
 &=\|K_H(\cdot+z)-K_H\|_{L^1(\R)}
 \le C_H|z|^{2H-1}.
\end{align*}
This proves \eqref{fr:lattice-translation} for every real $z$.

At $H=\frac12$, the second difference is exactly $\tau$. The
envelope is immediate, and the mass bound follows from
$\sum_j\tau(j+v)=1$. Since $\tau(j)=\one_{\{j=0\}}$, the translation
sum equals
$1-\tau(z)+\sum_{j\ne0}\tau(j+z)=2(1-\tau(z))
=2\min\{|z|,1\}$.
\end{proof}

\section{Sub-Gaussian moment and sum estimates}
\label{app:subgaussian}

For completeness, we prove the elementary estimate used to sum
independent centered blocks in the Brownian and discrete constructions.

\begin{proof}[Proof of Lemma~\ref{lem:subgaussian}]
Let $Z'$ be an independent copy of $Z$. Conditional Jensen gives
$e^{\lambda Z}\le E[e^{\lambda(Z-Z')}\mid Z]$.
The difference $Z-Z'$ is symmetric, and monotone convergence for
the series of $\cosh$ yields
\[
 E[e^{\lambda Z}]
 \le E[\cosh(\lambda(Z-Z'))]
 =1+\sum_{m\ge1}
       \frac{\lambda^{2m}E[(Z-Z')^{2m}]}{(2m)!}.
\]
These statements hold for extended nonnegative expectations;
the following estimate proves finiteness.
By Minkowski's inequality,
\[
 E[(Z-Z')^{2m}]
 \le(2c\sqrt{2m}\,s)^{2m}=(8c^2s^2m)^m.
\]
Since $\frac{(2m)!}{m!}=(m+1)\cdots(2m)\ge m^m$, the series
is bounded by
$\sum_{m\ge0}\frac{(8c^2s^2\lambda^2)^m}{m!}
=\exp(8c^2s^2\lambda^2)$.
This proves \eqref{eq:subgaussian-mgf}.

For the sum, independence gives
\[
 E[e^{\lambda\Sigma_n}]
 =\prod_{i=1}^n E[e^{\lambda Z_i}]
 \le\prod_{i=1}^n e^{8c^2\lambda^2s_i^2}
 =e^{8c^2\lambda^2V_n^2},
\]
which is \eqref{eq:subgaussian-sum-mgf}.

Assume $V_n>0$. For $u,\lambda>0$, exponential Markov gives
\[
 \Pp(\Sigma_n\ge u)
 \le\exp(-\lambda u+8c^2\lambda^2V_n^2).
\]
Taking $\lambda=\frac{u}{16c^2V_n^2}$ bounds this probability
by $\exp(-\frac{u^2}{32c^2V_n^2})$.
Apply the same argument to $-\Sigma_n$ and add the two bounds:
\begin{equation}\label{eq:subgaussian-sum-tail}
 \Pp(|\Sigma_n|\ge u)
 \le2\exp\left(-\frac{u^2}{32c^2V_n^2}\right),\qquad u>0.
\end{equation}
For an integer $m\ge1$, integration of this tail gives
\begin{align*}
 E[|\Sigma_n|^{2m}]
 &\le4m\int_0^\infty u^{2m-1}
                    e^{-\frac{u^2}{32c^2V_n^2}}\dd u\\
 &=2(32c^2V_n^2)^m m!.
\end{align*}
Given $p\ge2$, choose $m=\lceil\frac p2\rceil$ and use
$m!\le m^m$:
\[
 \|\Sigma_n\|_p\le\|\Sigma_n\|_{2m}
 \le 2^{\frac1{2m}}\sqrt{32}\,c\,V_n(m!)^{\frac1{2m}}
\le8c\sqrt m\,V_n\le8c\sqrt p\,V_n,
\]
since $m\le p$. This proves \eqref{eq:subgaussian-sum-Lp}.
\end{proof}

\section{Feynman-Kac formula for an unbounded potential}
\label{app:FK-unbounded}

We prove Lemma~\ref{lem:FK-unbounded} by combining Friedman's
uniqueness theorem with a direct verification of the Brownian
representation. The derivative bounds in
\eqref{eq:FK-linear-upper-bound} are used to justify differentiation
under the expectation.

\begin{proof}[Proof of Lemma~\ref{lem:FK-unbounded}]
Write $A=A_T$ for the constant in
\eqref{eq:FK-linear-upper-bound}.

\smallskip
\textit{Step 1: uniqueness.} With the growth assumption \eqref{eq:FK-linear-upper-bound} on $q$, uniqueness is an immediate consequence of \cite[Chapter~2, Section~4, Theorem~10]{Friedman-parabolic}.

\smallskip
\textit{Step 2: construction and classical regularity.}
Put $M_T=\sup_{0\le r\le T}|B_r|$. The reflection principle gives
\[
 P_B(M_T>z)\le4e^{-z^2/(2T)},\qquad z>0,
\]
and hence $E_B[e^{\lambda M_T}]<\infty$ for every $\lambda\ge0$.
By changing variables in \eqref{eq:bounded-potential-FK}, write
\[
 I_t(x)=\int_0^t q(t-r,x+B_r)\dd r,
 \qquad u(t,x)=E_B[e^{I_t(x)}].
\]
The bound $|I_t(x)|\le AT(1+|x|+M_T)$ yields
\begin{equation}\label{eq:FK-exponential-bound}
 0<u(t,x)\le e^{AT(1+|x|)}E_B[e^{AT M_T}]
 =C_Te^{AT|x|},\qquad 0\le t\le T.
\end{equation}
For each continuous Brownian path, $I_t(x)$ is continuous in
$(t,x)$, and its exponential is dominated on
$[0,T]\times[-L,L]$ by $e^{AT(1+L+M_T)}$. Thus
$u\in C([0,T]\times\R)$ and $u(0,x)=1$.

For a fixed continuous path, ordinary differentiation gives
\begin{align*}
 \partial_x I_t(x)&=\int_0^t q_x(t-r,x+B_r)\dd r,\\
 \partial_{xx} I_t(x)&=\int_0^t q_{xx}(t-r,x+B_r)\dd r,\\
 \partial_t I_t(x)&=q(0,x+B_t)
                 +\int_0^t q_t(t-r,x+B_r)\dd r.
\end{align*}
These derivatives are jointly continuous (with one-sided time
derivatives at the endpoints). By
\eqref{eq:FK-linear-upper-bound}, the exponential and its first
time derivative and first two spatial derivatives are dominated
on $[0,T]\times[-L,L]$ by an integrable random variable of the form
\[
 C_{T,L}(1+M_T)^2e^{C_{T,L}(1+M_T)}.
\]
Differentiation under the expectation and dominated convergence
therefore give $u\in C^{1,2}((0,T]\times\R)$ and
\begin{align*}
 \partial_xu(t,x)&=E_B[e^{I_t(x)}\partial_xI_t(x)],\\
 \partial_{xx}u(t,x)&=E_B\left[e^{I_t(x)}
       \big((\partial_xI_t(x))^2+\partial_{xx}I_t(x)\big)\right],\\
 \partial_tu(t,x)&=E_B[e^{I_t(x)}\partial_tI_t(x)].
\end{align*}
The same bounds also show that
\begin{equation}\label{eq:FK-second-derivative-growth}
 |\partial_{xx}u(t,x)|\le C_T(1+|x|)^2e^{AT|x|},
 \qquad 0\le t\le T.
\end{equation}

\smallskip
\textit{Step 3: verification of the equation.}
Fix $0<t<T$. For $0<\delta\le T-t$, the Markov property at time
$\delta$ gives
\[
 u(t+\delta,x)=E_B\left[e^{J_\delta(t,x)}u(t,x+B_\delta)\right],
 \qquad
 J_\delta(t,x)=\int_0^\delta q(t+\delta-r,x+B_r)\dd r.
\]
Writing $P_\delta f(x)=E_B[f(x+B_\delta)]$, we obtain
\begin{align*}
 \frac{u(t+\delta,x)-u(t,x)}{\delta}
 ={}&\frac{P_\delta u(t,\cdot)(x)-u(t,x)}{\delta}\\
 &+E_B\left[\frac{e^{J_\delta(t,x)}-1}{\delta}
                         u(t,x+B_\delta)\right].
\end{align*}
For a standard normal variable $Z$, Taylor's formula and
$E[Z]=0$ yield
\begin{align*}
 \frac{P_\delta u(t,\cdot)(x)-u(t,x)}{\delta}
 &=E\left[Z^2\int_0^1(1-\theta)
       \partial_{xx}u(t,x+\theta\sqrt\delta\,Z)\dd\theta\right]\\
 &\longrightarrow\tfrac12\partial_{xx}u(t,x).
\end{align*}
Dominated convergence here follows from
\eqref{eq:FK-second-derivative-growth} and Gaussian exponential
moments. For the other term, continuity gives
$J_\delta(t,x)/\delta\to q(t,x)$ almost surely. Moreover,
\[
 |J_\delta(t,x)|/\delta\le A(1+|x|+M_T).
\]
Together with \eqref{eq:FK-exponential-bound}, this implies
\[
 \left|\frac{e^{J_\delta(t,x)}-1}{\delta}
                        u(t,x+B_\delta)\right|
 \le C_{T,x}(1+M_T)e^{C_{T,x}(1+M_T)}.
\]
Another application of dominated convergence shows that this term
tends to $q(t,x)u(t,x)$. Hence
$\partial_tu=\frac12\partial_{xx}u+qu$ on $(0,T)\times\R$,
and also at $t=T$ by continuity. Together with the initial
condition and \eqref{eq:FK-exponential-bound}, this completes the
existence proof.
\end{proof}

\section{From pointwise to compact-uniform convergence}
\label{app:compact-upgrade}

We prove Lemma~\ref{lem:compact-upgrade}, used in the continuous
construction and in the functional convergence of the discrete model.
The H\"older bounds follow from the Garsia-Rodemich-Rumsey lemma, and
the convergence from a Cauchy argument in $L^p(\Omega;C(\mathcal K))$.

\begin{proof}[Proof of Lemma~\ref{lem:compact-upgrade}]
\textit{Step 0: reduction.} For $z=(t,x)$ and $z'=(s,y)$ in $\mathcal K$, one has $|x-y|\le2L$, so that
$|x-y|^\eta\le(2L)^{\eta/2}|x-y|^{\eta/2}$. With $|\cdot|$ the Euclidean
norm on $\R^2$, \eqref{eq:compact-increment-criterion} therefore gives,
uniformly in $\iota$,
\begin{equation}\label{eq:chaining-increment}
 \|X_\iota(z)-X_\iota(z')\|_q\le C_{q,\mathcal K}|z-z'|^{\eta/2},
 \qquad z,z'\in\mathcal K.
\end{equation}
Fix $p$, choose $q>\max\{p,\frac4\eta\}$ and
$0<\gamma<\frac\eta2-\frac2q$. All constants below depend on $q$,
$\gamma$, $\eta$ and $\mathcal K$, but not on $\iota$. We write
$[f]_\gamma=\sup_{z\ne z'}|f(z)-f(z')|/|z-z'|^\gamma$ for the
$\gamma$-H\"older seminorm on $\mathcal K$, and
$\|f\|_{C^\gamma(\mathcal K)}=\|f\|_{C(\mathcal K)}+[f]_\gamma$.

\textit{Step 1: H\"older bounds.}
The Garsia-Rodemich-Rumsey lemma \cite{GRR} gives, for
every continuous function $f$ on $\mathcal K$,
\begin{equation}\label{eq:GRR}
 [f]_\gamma
 \le C\left(\int_{\mathcal K}\int_{\mathcal K}
       \frac{|f(z)-f(z')|^q}{|z-z'|^{\gamma q+4}}\dd z\dd z'\right)^{\frac1q}.
\end{equation}
Apply \eqref{eq:GRR} to the continuous function $f=X_\iota$, take the
$q$th moment, and use Tonelli's theorem and
\eqref{eq:chaining-increment}:
\[
 \E\bigl[[X_\iota]_\gamma^q\bigr]
 \le C\int_{\mathcal K}\int_{\mathcal K}
     \frac{\E\bigl[|X_\iota(z)-X_\iota(z')|^q\bigr]}
          {|z-z'|^{\gamma q+4}}\dd z\dd z'
 \le C\int_{\mathcal K}\int_{\mathcal K}
     |z-z'|^{\frac{q\eta}2-\gamma q-4}\dd z\dd z'.
\]
The last integral is finite because the exponent exceeds $-2$, that is,
because $\gamma<\frac\eta2-\frac2q$. Since
$\|X_\iota\|_{C(\mathcal K)}\le|X_\iota(0,0)|+C[X_\iota]_\gamma$ and
$q\ge p$, this proves
$\sup_\iota\E\bigl[\|X_\iota\|_{C^\gamma(\mathcal K)}^p\bigr]<\infty$,
the first part of \eqref{eq:compact-holder-control}.

\textit{Step 2: Cauchy property and continuous modification.} Fix
$\delta>0$ and a finite set $D_\delta\subset\mathcal K$ such that every
point of $\mathcal K$ lies within distance $\delta$ of $D_\delta$. For
all $\iota,\iota'$,
\[
 \|X_\iota-X_{\iota'}\|_{C(\mathcal K)}
 \le\max_{z\in D_\delta}|X_\iota(z)-X_{\iota'}(z)|
   +\bigl([X_\iota]_\gamma+[X_{\iota'}]_\gamma\bigr)\delta^\gamma,
\]
so that, by Step~1,
\[
 \|X_\iota-X_{\iota'}\|_{L^p(\Omega;C(\mathcal K))}
 \le\Bigl\|\max_{z\in D_\delta}|X_\iota(z)-X_{\iota'}(z)|\Bigr\|_p
   +C\delta^\gamma.
\]
The first term tends to zero as $\iota$ and $\iota'$ tend to their
limit, since it involves finitely many points at which $X_\iota$
converges in $L^p$. Letting then $\delta\downarrow0$ shows that
$(X_\iota)$ is a Cauchy family in the Banach space
$L^p(\Omega;C(\mathcal K))$. It therefore converges in this space to a
random continuous field $\widetilde X$. For each $z\in\mathcal K$, both
$X_\iota(z)\to\widetilde X(z)$ and $X_\iota(z)\to X(z)$ in $L^p$, so
that $\widetilde X(z)=X(z)$ almost surely. Thus $\widetilde X$ is a
continuous modification of $X$, and
$\|X_\iota-\widetilde X\|_{L^p(\Omega;C(\mathcal K))}\to0$.

\textit{Step 3: H\"older bound for the limit.} Along a sequence
$\iota_n$, $X_{\iota_n}\to\widetilde X$ uniformly on $\mathcal K$ almost
surely. Hence
$\|\widetilde X\|_{C^\gamma(\mathcal K)}
\le\liminf_n\|X_{\iota_n}\|_{C^\gamma(\mathcal K)}$ almost surely, and
Fatou's lemma together with Step~1 gives
$\E\bigl[\|\widetilde X\|_{C^\gamma(\mathcal K)}^p\bigr]
\le\sup_\iota\E\bigl[\|X_\iota\|_{C^\gamma(\mathcal K)}^p\bigr]<\infty$,
which completes \eqref{eq:compact-holder-control}.
\end{proof}

\section*{Acknowledgments}
I would like to thank Jian Song for stimulating discussions on the topic.

\section*{Disclosure of AI assistance}
I developed the main mathematical ideas and the proof strategy for this
paper. The Brownian Feynman-Kac representation of the regularized
solution, which is central to the argument, already appears in my joint
work with Chen, Ouyang and Tindel \cite[Section~4.2, Eq.~(4.10)]{CDOT-moments}.
I used ChatGPT6 (OpenAI) and Claude Fable (Anthropic) to help work out
the details of the proofs.
I also used AI assistants (ChatGPT, OpenAI; Claude, Anthropic) to draft and revise parts of the text. I checked and revised all of it, and I take full responsibility for the final manuscript.

\section*{Data availability}
This article has no associated data set. The figures of
Section~\ref{sec:numerics} display simulated realizations, obtained by
the procedure described in Section~\ref{subsec:numerical-method}.

\section*{Conflict of interest}
The author has no competing interests to declare that are relevant to
the content of this article.


\end{document}